\documentclass[sigconf]{acmart}
\AtBeginDocument{%
  \providecommand\BibTeX{{%
    \normalfont B\kern-0.5em{\scshape i\kern-0.25em b}\kern-0.8em\TeX}}}

\copyrightyear{2022}
\acmYear{2022}
\setcopyright{acmcopyright}\acmConference[e-Energy '22]{The Thirteenth ACM International Conference on Future Energy Systems}{June 28-July 1, 2022}{Virtual Event, USA}
\acmBooktitle{The Thirteenth ACM International Conference on Future Energy Systems (e-Energy '22), June 28-July 1, 2022, Virtual Event, USA}
\acmPrice{15.00}

\acmConference[ACM e-Energy '22]{ACM e-Energy 2022}{June 28 - July 1, 2022}{virtual event}

\usepackage{algorithm}
\usepackage{algorithmic}
\DeclareMathOperator{\diag}{diag}

\begin{document}

%
%
%
%



\newcommand{\hermitiantranspose}{\text{H}}
\newcommand{\transpose}{\text{T}}
\newcommand{\complexconj}{*}

\newcommand{\Arrow}[1]{%
\parbox{#1}{\tikz{\draw[->](0,0)--(#1,0);}}
}
\newcommand{\Arrowright}[1]{%
\parbox{#1}{\tikz{\draw[<-](0,0)--(#1,0);}}
}

\newcommand{\indexScenario}{s}
\newcommand{\indexZone}{z}
\newcommand{\indexStakeholder}{b}
\newcommand{\indexUnitTech}{{x}_{\unitss}}
\newcommand{\indexPNTech}{{x}_{\pnss}}
\newcommand{\indexGITech}{{x}_{\giss}}
\newcommand{\indexGETech}{{x}_{\gess}}
\newcommand{\indexUnitRating}{{g}_{\unitss}}
\newcommand{\indexPNRating}{{g}_{\pnss}}
\newcommand{\indexGIRating}{{g}_{\giss}}
\newcommand{\indexGERating}{{g}_{\gess}}
\newcommand{\indexNodeRating}{{g}_{\nodess}}
\newcommand{\indexGridNode}{{i}}
\newcommand{\indexGridNodeTwo}{{j}}
\newcommand{\indexGridNodeThree}{{k}}
\newcommand{\indexGridLines}{{l}}
\newcommand{\indexGridLinesTwo}{{m}}
\newcommand{\indexGridLoad}{{h}}
\newcommand{\indexPhases}{p}
\newcommand{\indexPhasesTwo}{q}
\newcommand{\indexTimestep}{{k}}
\newcommand{\indexDay}{{d}}
\newcommand{\indexWeek}{{w}}
\newcommand{\indexUnit}{{u}}
\newcommand{\indexShunt}{{h}}
\newcommand{\indexPowerNode}{{e}}
\newcommand{\indexAdd}{{a}}
\newcommand{\indexRemove}{{r}}
\newcommand{\indexCirculation}{{\indexAdd\indexRemove}}
\newcommand{\indexFlows}{{f}}
\newcommand{\indexConductor}{n}
\newcommand{\indexdim}{n}
\newcommand{\indexdimtwo}{m}
\newcommand{\indexIteration}{k}
\newcommand{\indexContingency}{c}


\newcommand{\symPower}{P}
\newcommand{\symPowerToEnergy}{\gamma}
\newcommand{\symReactivePower}{Q}
\newcommand{\symApparentPower}{S}
\newcommand{\symVoltage}{U}
\newcommand{\symVoltageSOCP}{W}
\newcommand{\symVoltageAngle}{\theta}
\newcommand{\symPhaseDifference}{\varphi}
\newcommand{\symAnnualCost}{K}
\newcommand{\symPenalty}{Y}
\newcommand{\symPenaltyWeight}{N}
\newcommand{\symMass}{m}
\newcommand{\symSpecificMass}{\dot{\symMass}}
\newcommand{\symVolume}{V}
\newcommand{\symSpecificVolume}{v}
\newcommand{\symEnergy}{E}
\newcommand{\symEnergyFlow}{\dot{E}}
\newcommand{\symCurrent}{I}
\newcommand{\symCurrentSOCP}{L}
\newcommand{\symCurrentSquared}{L}
\newcommand{\symLineCurrent}{J}
\newcommand{\symObjectiveWeight}{w}
\newcommand{\symProbability}{\lambda}
\newcommand{\symBinary}{\alpha}
\newcommand{\symDecision}{{d}}
\newcommand{\symPhaseA}{{a}}
\newcommand{\symPhaseB}{{b}}
\newcommand{\symPhaseC}{{c}}
\newcommand{\symPhaseP}{{p}}
\newcommand{\symPhaseN}{{n}}
\newcommand{\symPhaseG}{{g}}
\newcommand{\symHarmonic}{{h}}
\newcommand{\symLossFactor}{{\rho}}
\newcommand{\symAvailability}{{b}}
\newcommand{\symDirectionality}{{d}}
\newcommand{\symPrice}{p}
\newcommand{\symInvestment}{C}
\newcommand{\symSpecificInvestment}{c}
\newcommand{\symLifetime}{\tau}
\newcommand{\symResistance}{r}
\newcommand{\symReactance}{x}
\newcommand{\symImpedance}{z}
\newcommand{\symAdmittance}{y}
\newcommand{\symCapacitance}{C}
\newcommand{\symSusceptance}{b}
\newcommand{\symConductance}{g}
\newcommand{\symLength}{l}
\newcommand{\symTime}{t}
\newcommand{\symRatio}{T}
\newcommand{\symLocation}{l}
\newcommand{\symEnthalpy}{H}
\newcommand{\symConnectivity}{C}

\newcommand{\symVolumeFlow}{\bar{\symVolume}}
\newcommand{\symMassFlow}{\bar{\symMass}}
\newcommand{\symDensity}{{\rho}}
\newcommand{\symGravitationalConstant}{\textcolor{\paramcolor}{g}}
\newcommand{\symElevation}{{h}}
\newcommand{\symElevationDifference}{\Delta{\symElevation}}
\newcommand{\symSetting}{\sigma}

\newcommand{\symRampRate}{\dot{\symPower}}
\newcommand{\symSpecRampRate}{{\zeta}}

\newcommand{\symPressure}[0]{p}
\newcommand{\symPressureSquared}[0]{\beta}
\newcommand{\symTemperature}[0]{   T }
\newcommand{\symHeatFlow}[0]{   Q }
\newcommand{\symViscosity}[0]{   \mu }
\newcommand{\symHeatResistance}[0]{   R }
\newcommand{\symFluidPipeDiameter}[0]{  \textcolor{\paramcolor}{D }    }
\newcommand{\symFluidPipeArea}[0]{  A     }
\newcommand{\symFluidPipeFriction}[0]{  \textcolor{\paramcolor}{f }    }
\newcommand{\symHeatCapacity}[0]{ c   }

\newcommand{\symRadiality}[0]{ \beta   }
\newcommand{\symSlack}{\epsilon}



\newcommand{\basess}{\text{{base}}}
\newcommand{\seriesss}{\text{{s}}}
\newcommand{\shuntss}{\text{{sh}}}
\newcommand{\invss}{\text{{inv}}}
\newcommand{\unbss}{\text{{unb}}}
\newcommand{\totss}{\text{{tot}}}
\newcommand{\batss}{\text{{bat}}}
\newcommand{\gridss}{\text{{grid}}}
\newcommand{\usabless}{\text{{usable}}}
\newcommand{\giss}{\text{{gi}}}
\newcommand{\pnss}{\text{{pn}}}
\newcommand{\lossss}{\text{{loss}}}
\newcommand{\cyclesss}{\text{{cycles}}}
\newcommand{\standbyss}{\text{{standby}}}
\newcommand{\operationss}{\text{{op}}}
\newcommand{\Pss}{\text{\sc{\symPower}}}
\newcommand{\Vss}{\text{\sc{\symVoltage}}}
\newcommand{\Iss}{\text{\sc{\symCurrent}}}
\newcommand{\Sss}{\text{\sc{\symApparentPower}}}
\newcommand{\Qss}{\text{\sc{\symReactivePower}}}
\newcommand{\Zss}{\text{\sc{\symImpedance}}}
\newcommand{\Ess}{\text{\sc{\symEnergy}}}
\newcommand{\nomss}{\text{nom}}
\newcommand{\ratedss}{\text{rated}}
\newcommand{\absss}{\text{abs}}
\newcommand{\maxss}{\text{max}}
\newcommand{\minss}{\text{min}}
\newcommand{\effss}{\text{eff}}
\newcommand{\chargess}{\text{c}}
\newcommand{\dischargess}{\text{d}}
\newcommand{\exoaddss}{\text{+}}
\newcommand{\exoremss}{\text{-}}
\newcommand{\orthss}{\perp}
\newcommand{\deprss}{\text{depr}}
\newcommand{\PFconvexss}{\text{PFconvex}}
\newcommand{\GIconvexss}{\text{GIconvex}}
\newcommand{\iterationss}{{(q)}}
\newcommand{\refss}{\text{ref}}
\newcommand{\addss}{\indexAdd,\text{add}}
\newcommand{\addtextss}{\text{add}}
\newcommand{\remss}{\indexRemove,\text{rem}}
\newcommand{\remtextss}{\text{rem}}
\newcommand{\circss}{\indexAdd\indexRemove,\text{circ}}
\newcommand{\flowss}{\indexFlows,\text{flow}}
\newcommand{\switchss}{\text{switch}}
\newcommand{\rotatedss}{\text{rot}}
\newcommand{\circumss}{\text{circum}}
\newcommand{\sbss}{\text{SB}}
\newcommand{\nonsbss}{\text{nonSB}}
\newcommand{\gess}{\text{ge}}
\newcommand{\unitss}{\text{unit}}
\newcommand{\nodess}{\text{node}}
\newcommand{\linss}{\text{lin}}
\newcommand{\realss}{\text{re}}
\newcommand{\imagss}{\text{im}}
\newcommand{\hagenposeuilless}[0]{  \text{ HP }}
\newcommand{\darcyweisbachss}[0]{ \text{  DC }}
\newcommand{\bentallinearization}[0]{ \text{  BenTal }}

\newcommand{\fortescuess}{\text{012}}

\newcommand{\plusss}{\text{+}}
\newcommand{\subsss}{\text{-}}
\newcommand{\combustionss}{\text{comb}}
\newcommand{\formationss}{\text{form}}



\newcommand{\PowerBase}{\symApparentPower_{\basess}}
\newcommand{\CurrentBase}{\symCurrent_{\basess}}
\newcommand{\VoltageBase}{\symVoltage_{\basess}^{\text{line}}}
\newcommand{\VoltageBasePhase}{\symVoltage_{\basess}^{\text{phase}}}
\newcommand{\TimeBase}{\symTime_{\basess}}
\newcommand{\EnergyBase}{\symEnergy_{\basess}}
\newcommand{\ImpedanceBase}{\symImpedance_{\basess}}
\newcommand{\AdmittanceBase}{\symAdmittance_{\basess}}
\newcommand{\CurrencyBase}{\euro_{\basess}}

\newcommand{\Knul}{\textcolor{\paramcolor}{\symAnnualCost_{0}}}
\newcommand{\Ynul}{\textcolor{\paramcolor}{\symPenalty_{0}}}



\newcommand{\yearunit}{a}
\newcommand{\VA}{VA}
\newcommand{\var}{var}
\newcommand{\kWh}{kWh}
\newcommand{\MWh}{MWh}
\newcommand{\si}[1]{#1}
\newcommand{\Si}[2]{\unit{#1}{#2}}
\newcommand{\percent}{\%}
\newcommand{\num}{}
\newcommand{\invtan}{\mathrm{atan2}}
\newcommand{\onen}[1]{\textcolor{\paramcolor}{\mathbf{1}_{#1}}}
\newcommand{\onenon}{\textcolor{\paramcolor}{\mathbf{1}}}

\newcommand{\nk}{n_{\setTimesteps}}
\newcommand{\nunits}{n_\setUnits}
\newcommand{\nn}{n_{\omegaGridNodes}}
\newcommand{\nc}{n_{\omegaGridLines}}
\newcommand{\nslackbuses}{n_{\omegaSlackbuses}}
\newcommand{\nnonslackbuses}{n_{\omegaNonSlackbuses}}

\newcommand{\setTimesteps}{\mathcal{K}}
\newcommand{\setTimestepsprev}{\mathcal{K}_0}
\newcommand{\omegaPNTech}{\mathcal{X}_{\pnss}}
\newcommand{\omegaGITech}{\mathcal{X}_{\giss}}
\newcommand{\omegaGETech}{\mathcal{X}_{\gess}}
\newcommand{\omegaUnitTech}{\mathcal{X}_{\unitss}}
\newcommand{\setPNRating}{\mathcal{G}_{\pnss}}
\newcommand{\setGIRating}{\mathcal{G}_{\giss}}
\newcommand{\setGERating}{\mathcal{G}_{\gess}}
\newcommand{\setUnitRating}{\mathcal{G}_{\unitss}}  
\newcommand{\setNodeRating}{\mathcal{G}_{\nodess}}  %
\newcommand{\omegaTransformers}{\textcolor{\paramcolor}{\mathcal{F}}}
\newcommand{\omegaGridLines}{\textcolor{\paramcolor}{\mathcal{J}}}
\newcommand{\omegaGridNodes}{\textcolor{\paramcolor}{\mathcal{I}}}
\newcommand{\omegaGridNodesDomain}[1]{\omegaGridNodes_{#1}}
\newcommand{\omegaSlackbuses}{\omegaGridNodes_{\indexDomains, \sbss}}
\newcommand{\omegaNonSlackbuses}{\omegaGridNodes_{\indexDomains,\nonsbss}}
\newcommand{\setReferenceNodes}{\omegaGridNodes_{\text{ref}}}
\newcommand{\setReferenceNodesPrime}{\omegaGridNodes'_{\text{ref}}}
\newcommand{\setReferenceNodesRankOne}{\omegaGridNodes^{\text{R-1}}_{\text{ref}}}
\newcommand{\setReferenceNodesRankOnePrime}{\omegaGridNodes'^{\text{R-1}}_{\text{ref}}}
\newcommand{\setReferenceNodesBalanced}{\omegaGridNodes^{\text{bal}}_{\text{ref}}}
\newcommand{\setReferenceNodesFixed}{\omegaGridNodes^{\text{fix}}_{\text{ref}}}

\newcommand{\setHarmonics}{\mathcal{H}}
\newcommand{\setHarmonicsTwo}{\mathcal{H}_1}

\newcommand{\setPowerNodes}{\mathcal{E}}
\newcommand{\setAdd}{\mathcal{A}}
\newcommand{\setRemove}{\mathcal{R}}
\newcommand{\setCirculation}{\mathcal{C}}
\newcommand{\setFlows}{\mathcal{F}}
\newcommand{\setContingency}{\mathcal{C}}
\newcommand{\setTopology}{\textcolor{\paramcolor}{\mathcal{T}^{\Arrow{.1cm}}}}
\newcommand{\setTopologyBoth}{\textcolor{\paramcolor}{\mathcal{T}}}
\newcommand{\setTopologyReverse}{\textcolor{\paramcolor}{\mathcal{T}^{\Arrowright{.1cm}}}}
\newcommand{\setUnits}{\textcolor{\paramcolor}{\mathcal{U}}}
\newcommand{\setUnitTopology}{\textcolor{\paramcolor}{\mathcal{T^{\text{units}}}}}
\newcommand{\setShuntTopology}{\textcolor{\paramcolor}{\mathcal{T^{\text{shunts}}}}}

\newcommand{\setBinary}[0]{\textcolor{\setscolor}{ \left\{ 0,1 \right\} }}
\newcommand{\setComplex}[0]{\textcolor{\setscolor}{ \mathbb{C}   }}
\newcommand{\setReal}[0]{\textcolor{\setscolor}{ \mathbb{R}   }}
\newcommand{\setIntegers}[0]{\textcolor{\setscolor}{ \mathbb{Z}   }}
\newcommand{\setNaturalNumbers}[0]{\textcolor{\setscolor}{ \mathbb{N}   }}
\newcommand{\setNaturalNumbersPos}[0]{\textcolor{\setscolor}{ \mathbb{N}^{+}   }}
\newcommand{\setRealMatrix}[2]{\textcolor{\setscolor}{ \mathbb{R}^{#1 \times #2}   }}
\newcommand{\setComplexMatrix}[2]{\textcolor{\setscolor}{ \mathbb{C}^{#1\times #2}   }}
\newcommand{\setHermitianMatrix}[1]{\textcolor{\setscolor}{ \mathbb{H}^{#1}   }}
\newcommand{\setRealPos}[0]{\textcolor{\setscolor}{  \setReal_{+}    }}

\newcommand{\setSymmetricRealMatrix}[1]{\textcolor{\setscolor}{ \mathbb{P}^{#1}   }}

\newcommand{\setPhases}{\textcolor{\paramcolor}{\mathcal{P}}}
\newcommand{\setPhasesTwo}{{\Phi}}
\newcommand{\setFortescue}{\textcolor{\paramcolor}{\mathcal{F}}}
\newcommand{\setCircuit}{\textcolor{\paramcolor}{\mathcal{C}^{4W}}}
\newcommand{\setConductors}{\textcolor{\paramcolor}{\mathcal{N}}}
\newcommand{\setBuspairs}{\textcolor{\paramcolor}{\mathcal{B}}}
\newcommand{\setZIP}{\mathcal{Z}}

\newcommand{\setCircuitTW}{\textcolor{\paramcolor}{\mathcal{C}^{3W}}}

\newcommand{\setBalanced}[0]{\textcolor{\setscolor}{  \mathbb{B}^{\indexdim}    }}

\newcommand{\setPSDConen}[1]{\textcolor{\setscolor}{  \mathbb{S}_{+}^{#1}    }}
\newcommand{\setNNOrthant}[0]{\textcolor{\setscolor}{  \setReal_{+}^{\indexdim}    }}
\newcommand{\setSOCone}[0]{\textcolor{\setscolor}{  \mathbb{Q}_{}^{\indexdim}    }}
\newcommand{\setPSDCone}[0]{\textcolor{\setscolor}{  \mathbb{S}_{+}^{\indexdim}    }}
\newcommand{\setPDCone}[0]{\textcolor{\setscolor}{  \mathbb{S}_{++}^{\indexdim}    }}
\newcommand{\setScenario}{\mathcal{S}}
\newcommand{\setZone}{\mathcal{Z}}
\newcommand{\setStakeholder}{\mathcal{B}}

\newcommand{\ineuro}[1]{\EUR{\ensuremath{#1}}}
\newcommand{\eurounit}{\euro}
\newcommand{\ineuroperyear}[1]{\ensuremath{#1}\,\eurounit \si{\per\yearunit}}

\newcommand{\variables}{\textcolor{black}{real-valued variables}}
\newcommand{\parameters}{\textcolor{\paramcolor}{real-valued parameters}}
\newcommand{\binaryparameters}{\textcolor{\decisioncolor}{binary parameters}}
\newcommand{\bounds}{\textcolor{\boundscolor}{generic bounds}}
\newcommand{\sizingvariables}{\textcolor{\sizingcolor}{sizing parameters}}
\newcommand{\sizingbounds}{\textcolor{\ratingboundscolor}{sizing bounds}}
\newcommand{\sizingparameters}{\textcolor{\investmentcolor}{sizing parameters}}
\newcommand{\binaryvariables}{\textcolor{\binarycolor}{binary variables}}
\newcommand{\timeparameters}{\textcolor{\timecolor}{time parameters}}
\newcommand{\priceparameters}{\textcolor{\pricecolor}{price parameters}}
\newcommand{\complexvariables}{\textcolor{\complexcolor}{complex variables}}
\newcommand{\complexparameters}{\textcolor{\complexparamcolor}{complex parameters}}

\newcommand{\paramcolor}{red}
\newcommand{\sizingcolor}{red}
\newcommand{\investmentcolor}{teal}
\newcommand{\decisioncolor}{brown}
\newcommand{\boundscolor}{violet}
\newcommand{\binarycolor}{darkgray}
\newcommand{\timecolor}{cyan}
\newcommand{\pricecolor}{brown}
\newcommand{\ratingboundscolor}{olive}
\newcommand{\setscolor}{darkgray}
\newcommand{\complexcolor}{blue}  
\newcommand{\complexparamcolor}{brown}  


\newcommand{\costheta}{{c}_{\indexGridNode \indexGridNodeTwo}   }
\newcommand{\sintheta}{{s}_{\indexGridNode \indexGridNodeTwo}   }
\newcommand{\usquared}{{u}_{\indexGridNode \indexGridNode}   }
\newcommand{\umult}{{u}_{\indexGridNode \indexGridNodeTwo}   }
\newcommand{\uucos}{{r}_{\indexGridNode \indexGridNodeTwo}   }
\newcommand{\uusin}{{t}_{\indexGridNode \indexGridNodeTwo}   }

\newcommand{\phasesi}{\textcolor{\paramcolor}{{\setPhasesTwo}_{\indexGridNode }   }}
\newcommand{\phasesj}{\textcolor{\paramcolor}{{\setPhasesTwo}_{\indexGridNodeTwo }   }}
\newcommand{\phasesij}{\textcolor{\paramcolor}{{\setPhasesTwo}_{\indexGridNode \indexGridNodeTwo}   }}

\newcommand{\fortescue}{\textcolor{\complexparamcolor}{\mathbf{F}}}
\newcommand{\alpharot}{\textcolor{\complexparamcolor}{\alpha}}

\newcommand{\Ts}{\textcolor{\timecolor}{T_{\text{s}} }}
\newcommand{\thor}{\textcolor{\timecolor}{t_{\text{h}}}}
\newcommand{\etac}{\textcolor{\paramcolor}{\eta_{\chargess} }}
\newcommand{\etad}{\textcolor{\paramcolor}{\eta_{\dischargess} }}
\newcommand{\etap}{\textcolor{\paramcolor}{\eta_{\exoaddss} }}
\newcommand{\etam}{\textcolor{\paramcolor}{\eta_{\exoremss} }}

\newcommand{\etaadd}{\textcolor{\paramcolor}{\eta_{\addss} }}
\newcommand{\etarem}{\textcolor{\paramcolor}{\eta_{\remss} }}


\newcommand{\Pck}[0]{\symPower_{\chargess} }
\newcommand{\Pdk}[0]{\symPower_{\dischargess} }
\newcommand{\Ppk}[0]{\symPower_{\exoaddss} }
\newcommand{\Pmk}[0]{\symPower_{\exoremss} }

\newcommand{\Paddk}[0]{\symPower_{\addss} }
\newcommand{\Premk}[0]{\symPower_{\remss} }
\newcommand{\Pflowk}[0]{\symPower_{\flowss} }
\newcommand{\Pflowkprev}[0]{\symPower_{\flowss} prev}
\newcommand{\Pflowrampmax}[0]{\textcolor{\paramcolor}{\symRampRate_{\flowss}^{\maxss}}}
\newcommand{\Pflowrampspecmax}[0]{\textcolor{\paramcolor}{\symSpecRampRate_{\flowss}}}

\newcommand{\Pclossk}[0]{\symPower_{\chargess, \lossss} }
\newcommand{\Pdlossk}[0]{\symPower_{\dischargess, \lossss} }
\newcommand{\Pplossk}[0]{\symPower_{\exoaddss, \lossss} }
\newcommand{\Pmlossk}[0]{\symPower_{\exoremss, \lossss} }
\newcommand{\Pgilossk}[0]{\symPower_{\exoremss, \lossss} }
\newcommand{\PPNlossk}[0]{\symPower_{\indexPowerNode, \lossss} }
\newcommand{\PPNexolossk}[0]{\symPower_{\indexPowerNode, \text{exo},\lossss} }
\newcommand{\PPNeleclossk}[0]{\symPower_{\indexPowerNode, \text{elec}, \lossss} }
\newcommand{\Punitlossk}[0]{\symPower_{\indexUnit, \lossss} }

\newcommand{\Paddlossk}[0]{\symPower_{\addss, \lossss} }
\newcommand{\Premlossk}[0]{\symPower_{\remss, \lossss} }

\newcommand{\Kswitch}[0]{\symAnnualCost^{\switchss}_{}}
\newcommand{\Kswitchl}[0]{\symAnnualCost^{\switchss}_{\indexGridLines}}

\newcommand{\Kgridloss}[0]{\symAnnualCost^{\gridss\lossss}_{}}
\newcommand{\Kgridlossl}[0]{\symAnnualCost^{\gridss\lossss}_{\indexGridLines}}

\newcommand{\Kc}[0]{\symAnnualCost^{\chargess}_{}}
\newcommand{\Kd}[0]{\symAnnualCost^{\dischargess}_{}}
\newcommand{\Kp}[0]{\symAnnualCost^{\exoaddss}_{}}
\newcommand{\Km}[0]{\symAnnualCost^{\exoremss}_{}}
\newcommand{\KE}[0]{\symAnnualCost^{\Ess}_{}}
\newcommand{\Kgi}[0]{\symAnnualCost^{\giss}_{}}
\newcommand{\KPN}[0]{\symAnnualCost_{\indexPowerNode}}
\newcommand{\Kunit}[0]{\symAnnualCost_{\indexUnit}}
\newcommand{\Kunits}[0]{\symAnnualCost_{\setUnits}}

\newcommand{\Kflow}[0]{\symAnnualCost^{\flowss}_{}}

\newcommand{\Kopex}[0]{\symAnnualCost^{\operationss}_{}}

\newcommand{\Kcoper}[0]{\symAnnualCost^{\chargess,\operationss}}
\newcommand{\Kdoper}[0]{\symAnnualCost^{\dischargess,\operationss}}
\newcommand{\Kctotoper}[0]{\symAnnualCost^{\chargess,\totss,\operationss}}
\newcommand{\Kdtotoper}[0]{\symAnnualCost^{\dischargess,\totss,\operationss}}
\newcommand{\Kpoper}[0]{\symAnnualCost^{\exoaddss,\operationss}}
\newcommand{\Kmoper}[0]{\symAnnualCost^{\exoremss,\operationss}}
\newcommand{\KPNoper}[0]{\symAnnualCost_{\indexPowerNode}^{\operationss}}
\newcommand{\Kunitoper}[0]{\symAnnualCost_{\indexUnit}^{\operationss}}
\newcommand{\Kunitsoper}[0]{\symAnnualCost_{\setUnits}^{\operationss}}

\newcommand{\Kflowoper}[0]{\symAnnualCost^{\flowss,\operationss}}

\newcommand{\Ec}[0]{\symEnergy_{\chargess}}
\newcommand{\Ed}[0]{\symEnergy_{\dischargess}}
\newcommand{\Ep}[0]{\symEnergy_{\exoaddss}}
\newcommand{\Em}[0]{\symEnergy_{\exoremss}}
\newcommand{\Egi}[0]{\symEnergy_{\giss}}
\newcommand{\EE}[0]{n_{\cyclesss}}

\newcommand{\Eflow}[0]{\symEnergy_{\flowss}}

\newcommand{\Ecflow}[0]{\symEnergyFlow_{\chargess}}
\newcommand{\Edflow}[0]{\symEnergyFlow_{\dischargess}}
\newcommand{\Epflow}[0]{\symEnergyFlow_{\exoaddss}}
\newcommand{\Emflow}[0]{\symEnergyFlow_{\exoremss}}
\newcommand{\Egiflow}[0]{\symEnergyFlow_{\giss}}
\newcommand{\EEflow}[0]{\dot{n}_{\cyclesss}}

\newcommand{\Eflowflow}[0]{\symEnergyFlow_{\flowss}}

\newcommand{\Ecmax}[0]{\textcolor{\paramcolor}{\symEnergy_{\chargess}^{\maxss}}}
\newcommand{\Edmax}[0]{\textcolor{\paramcolor}{\symEnergy_{\dischargess}^{\maxss}}}
\newcommand{\Epmax}[0]{\textcolor{\paramcolor}{\symEnergy_{\exoaddss}^{\maxss}}}
\newcommand{\Emmax}[0]{\textcolor{\paramcolor}{\symEnergy_{\exoremss}^{\maxss}}}
\newcommand{\Egimax}[0]{\textcolor{\paramcolor}{\symEnergy_{\giss}^{\maxss}}}
\newcommand{\EEmax}[0]{\textcolor{\paramcolor}{n_{\cyclesss}^{\maxss}}}

\newcommand{\Eflowmax}[0]{\textcolor{\paramcolor}{\symEnergy_{\flowss}^{\maxss}}}

\newcommand{\Ecmaxflow}[0]{\textcolor{\paramcolor}{\dot{\symEnergy}_{\chargess}^{\maxss}}}
\newcommand{\Edmaxflow}[0]{\textcolor{\paramcolor}{\dot{\symEnergy}_{\dischargess}^{\maxss}}}
\newcommand{\Epmaxflow}[0]{\textcolor{\paramcolor}{\dot{\symEnergy}_{\exoaddss}^{\maxss}}}
\newcommand{\Emmaxflow}[0]{\textcolor{\paramcolor}{\dot{\symEnergy}_{\exoremss}^{\maxss}}}
\newcommand{\Egimaxflow}[0]{\textcolor{\paramcolor}{\dot{\symEnergy}_{\giss}^{\maxss}}}
\newcommand{\EEmaxflow}[0]{\textcolor{\paramcolor}{\dot{n}_{\cyclesss}^{\maxss}}}

\newcommand{\Eflowmaxflow}[0]{\textcolor{\paramcolor}{\dot{\symEnergy}_{\flowss}^{\maxss}}}

\newcommand{\Kcdepr}[0]{\symAnnualCost^{\chargess,\deprss}}
\newcommand{\Kddepr}[0]{\symAnnualCost^{\dischargess,\deprss}}
\newcommand{\Kpdepr}[0]{\symAnnualCost^{\exoaddss,\deprss}}
\newcommand{\Kmdepr}[0]{\symAnnualCost^{\exoremss,\deprss}}
\newcommand{\KEdepr}[0]{\symAnnualCost^{\Ess,\deprss}}
\newcommand{\Kgidepr}[0]{\symAnnualCost^{\giss,\deprss}}
\newcommand{\KPNdepr}[0]{\symAnnualCost_{\indexPowerNode}^{\deprss}}
\newcommand{\Kunitdepr}[0]{\symAnnualCost_{\indexUnit}^{\deprss}}
\newcommand{\Kunitsdepr}[0]{\symAnnualCost_{\setUnits}^{\deprss}}

\newcommand{\Kflowdepr}[0]{\symAnnualCost^{\flowss,\deprss}}

\newcommand{\Cc}[0]{\textcolor{\investmentcolor}{\symInvestment_{\chargess}}}
\newcommand{\Cd}[0]{\textcolor{\investmentcolor}{\symInvestment_{\dischargess}}}
\newcommand{\Cp}[0]{\textcolor{\investmentcolor}{\symInvestment_{\exoaddss}}}
\newcommand{\Cm}[0]{\textcolor{\investmentcolor}{\symInvestment_{\exoremss}}}
\newcommand{\Cgi}[0]{\textcolor{\investmentcolor}{\symInvestment_{\Sss}}}
\newcommand{\CE}[0]{\textcolor{\investmentcolor}{\symInvestment_{\Ess}}}
\newcommand{\Cunit}[0]{\textcolor{\investmentcolor}{\symInvestment_{\indexUnit}}}
\newcommand{\CPN}[0]{\textcolor{\investmentcolor}{\symInvestment_{\indexPowerNode}}}
\newcommand{\Cunitmax}[0]{\textcolor{\boundscolor}{\symInvestment^{\maxss}_{\indexUnit}}}

\newcommand{\Cflow}[0]{\textcolor{\investmentcolor}{\symInvestment_{\flowss}}}

\newcommand{\cpc}[0]{\textcolor{\pricecolor}{\symSpecificInvestment_{\chargess}}}
\newcommand{\cpd}[0]{\textcolor{\pricecolor}{\symSpecificInvestment_{\dischargess}}}
\newcommand{\cpp}[0]{\textcolor{\pricecolor}{\symSpecificInvestment_{\exoaddss}}}
\newcommand{\cpm}[0]{\textcolor{\pricecolor}{\symSpecificInvestment_{\exoremss}}}
\newcommand{\cpgi}[0]{\textcolor{\pricecolor}{\symSpecificInvestment_{\Sss}}}
\newcommand{\cpE}[0]{\textcolor{\pricecolor}{\symSpecificInvestment_{\Ess}}}

\newcommand{\cpflow}[0]{\textcolor{\pricecolor}{\symSpecificInvestment_{\flowss}}}

\newcommand{\Vc}[0]{\textcolor{\investmentcolor}{\symVolume_{\chargess}}}
\newcommand{\Vd}[0]{\textcolor{\investmentcolor}{\symVolume_{\dischargess}}}
\newcommand{\Vp}[0]{\textcolor{\investmentcolor}{\symVolume_{\exoaddss}}}
\newcommand{\Vm}[0]{\textcolor{\investmentcolor}{\symVolume_{\exoremss}}}
\newcommand{\Vgi}[0]{\textcolor{\investmentcolor}{\symVolume_{\Sss}}}
\newcommand{\VE}[0]{\textcolor{\investmentcolor}{\symVolume_{\Ess}}}
\newcommand{\Vunit}[0]{\textcolor{\investmentcolor}{\symVolume_{\indexUnit}}}
\newcommand{\Vunitmax}[0]{\textcolor{\boundscolor}{\symVolume^{\maxss}_{\indexUnit}}}
\newcommand{\VPN}[0]{\textcolor{\investmentcolor}{\symVolume_{\indexPowerNode}}}

\newcommand{\Vflow}[0]{\textcolor{\investmentcolor}{\symVolume_{\flowss}}}

\newcommand{\vpc}[0]{\textcolor{\pricecolor}{\symSpecificVolume_{\chargess}}}
\newcommand{\vpd}[0]{\textcolor{\pricecolor}{\symSpecificVolume_{\dischargess}}}
\newcommand{\vpp}[0]{\textcolor{\pricecolor}{\symSpecificVolume_{\exoaddss}}}
\newcommand{\vpm}[0]{\textcolor{\pricecolor}{\symSpecificVolume_{\exoremss}}}
\newcommand{\vpgi}[0]{\textcolor{\pricecolor}{\symSpecificVolume_{\Sss}}}
\newcommand{\vpE}[0]{\textcolor{\pricecolor}{\symSpecificVolume_{\Ess}}}

\newcommand{\vpflow}[0]{\textcolor{\pricecolor}{\symSpecificVolume_{\flowss}}}

\newcommand{\Md}[0]{\textcolor{\investmentcolor}{\symMass_{\dischargess}}}
\newcommand{\Mp}[0]{\textcolor{\investmentcolor}{\symMass_{\exoaddss}}}
\newcommand{\Mm}[0]{\textcolor{\investmentcolor}{\symMass_{\exoremss}}}
\newcommand{\Mgi}[0]{\textcolor{\investmentcolor}{\symMass_{\Sss}}}
\newcommand{\ME}[0]{\textcolor{\investmentcolor}{\symMass_{\Ess}}}
\newcommand{\Munit}[0]{\textcolor{\investmentcolor}{\symMass_{\indexUnit}}}
\newcommand{\Munitmax}[0]{\textcolor{\boundscolor}{\symMass^{\maxss}_{\indexUnit}}}
\newcommand{\MPN}[0]{\textcolor{\investmentcolor}{\symMass_{\indexPowerNode}}}

\newcommand{\Mflow}[0]{\textcolor{\investmentcolor}{\symMass_{\flowss}}}

\newcommand{\mpc}[0]{\textcolor{\pricecolor}{\symSpecificMass_{\chargess}}}
\newcommand{\mpd}[0]{\textcolor{\pricecolor}{\symSpecificMass_{\dischargess}}}
\newcommand{\mpp}[0]{\textcolor{\pricecolor}{\symSpecificMass_{\exoaddss}}}
\newcommand{\mpm}[0]{\textcolor{\pricecolor}{\symSpecificMass_{\exoremss}}}
\newcommand{\mpgi}[0]{\textcolor{\pricecolor}{\symSpecificMass_{\Sss}}}
\newcommand{\mpE}[0]{\textcolor{\pricecolor}{\symSpecificMass_{\Ess}}}

\newcommand{\mpflow}[0]{\textcolor{\pricecolor}{\symSpecificMass_{\flowss}}}

\newcommand{\tauc}[0]{\textcolor{black}{\symLifetime_{\chargess}}}
\newcommand{\taud}[0]{\textcolor{black}{\symLifetime_{\dischargess}}}
\newcommand{\taup}[0]{\textcolor{black}{\symLifetime_{\exoaddss}}}
\newcommand{\taum}[0]{\textcolor{black}{\symLifetime_{\exoremss}}}
\newcommand{\taugi}[0]{\textcolor{black}{\symLifetime_{\Sss}}}
\newcommand{\tauE}[0]{\textcolor{black}{\symLifetime{_\Ess}}}

\newcommand{\tauflow}[0]{\textcolor{black}{\symLifetime_{\flowss}}}

\newcommand{\taucmax}[0]{\textcolor{\paramcolor}{\symLifetime^{\maxss}_{\chargess}}}
\newcommand{\taudmax}[0]{\textcolor{\paramcolor}{\symLifetime{^{\maxss}_\dischargess}}}
\newcommand{\taupmax}[0]{\textcolor{\paramcolor}{\symLifetime^{\maxss}_{\exoaddss}}}
\newcommand{\taummax}[0]{\textcolor{\paramcolor}{\symLifetime^{\maxss}_{\exoremss}}}
\newcommand{\taugimax}[0]{\textcolor{\paramcolor}{\symLifetime^{\maxss}_{\Sss}}}
\newcommand{\tauEmax}[0]{\textcolor{\paramcolor}{\symLifetime^{\maxss}_{\Ess}}}

\newcommand{\tauflowmax}[0]{\textcolor{\paramcolor}{\symLifetime^{\maxss}_{\flowss}}}

\newcommand{\ppc}[0]{\textcolor{\pricecolor}{\symPrice_{\chargess} }}
\newcommand{\ppd}[0]{\textcolor{\pricecolor}{\symPrice_{\dischargess} }}
\newcommand{\ppp}[0]{\textcolor{\pricecolor}{\symPrice_{\exoaddss} }}
\newcommand{\ppm}[0]{\textcolor{\pricecolor}{\symPrice_{\exoremss} }}

\newcommand{\ppflow}[0]{\textcolor{\pricecolor}{\symPrice_{\flowss} }}

\newcommand{\ppswitch}[0]{\textcolor{\pricecolor}{\symPrice_{\switchss} }}
\newcommand{\ppswitchl}[0]{\textcolor{\pricecolor}{\symPrice_{\indexGridLines,\switchss} }}

\newcommand{\ppPgridlossl}[0]{\textcolor{\pricecolor}{\symPrice_{\indexGridLines\gridss\lossss}^{\Pss}    }}
\newcommand{\ppQgridlossl}[0]{\textcolor{\pricecolor}{\symPrice_{\indexGridLines\gridss\lossss}^{\Qss} }}
\newcommand{\ppQabsgridlossl}[0]{\textcolor{\pricecolor}{\symPrice_{\indexGridLines\gridss\lossss}^{|\Qss|} }}

\newcommand{\Pckref}[0]{\textcolor{\paramcolor}{\symPower^{\refss}_{\chargess} }}
\newcommand{\Pdkref}[0]{\textcolor{\paramcolor}{\symPower^{\refss}_{\dischargess} }}
\newcommand{\Ppkref}[0]{\textcolor{\paramcolor}{\symPower^{\refss}_{\exoaddss} }}
\newcommand{\Pmkref}[0]{\textcolor{\paramcolor}{\symPower^{\refss}_{\exoremss} }}

\newcommand{\Pckmax}[0]{\textcolor{\boundscolor}{\symPower^{\maxss}_{\chargess} }}
\newcommand{\Pdkmax}[0]{\textcolor{\boundscolor}{\symPower^{\maxss}_{\dischargess} }}
\newcommand{\Ppkmax}[0]{\textcolor{\boundscolor}{\symPower^{\maxss}_{\exoaddss} }}
\newcommand{\Pmkmax}[0]{\textcolor{\boundscolor}{\symPower^{\maxss}_{\exoremss} }}

\newcommand{\Pckmin}[0]{\textcolor{\boundscolor}{\symPower^{\minss}_{\chargess} }}
\newcommand{\Pdkmin}[0]{\textcolor{\boundscolor}{\symPower^{\minss}_{\dischargess} }}
\newcommand{\Ppkmin}[0]{\textcolor{\boundscolor}{\symPower^{\minss}_{\exoaddss} }}
\newcommand{\Pmkmin}[0]{\textcolor{\boundscolor}{\symPower^{\minss}_{\exoremss} }}
\newcommand{\Pflowkref}[0]{\textcolor{\paramcolor}{\symPower^{\refss}_{\flowss} }}

\newcommand{\Paddkmin}[0]{\textcolor{\boundscolor}{\symPower^{\minss}_{\addss} }}
\newcommand{\Paddkmax}[0]{\textcolor{\boundscolor}{\symPower^{\maxss}_{\addss} }}
\newcommand{\Premkmin}[0]{\textcolor{\boundscolor}{\symPower^{\minss}_{\remss} }}
\newcommand{\Premkmax}[0]{\textcolor{\boundscolor}{\symPower^{\maxss}_{\remss} }}
\newcommand{\Pflowkmin}[0]{\textcolor{\boundscolor}{\symPower^{\minss}_{\flowss} }}
\newcommand{\Pflowkmax}[0]{\textcolor{\boundscolor}{\symPower^{\maxss}_{\flowss} }}

\newcommand{\Pcratedmax}[0]{\textcolor{\ratingboundscolor}{\symPower'^{\maxss}_{\chargess}}}
\newcommand{\Pdratedmax}[0]{\textcolor{\ratingboundscolor}{\symPower'^{\maxss}_{\dischargess}}}
\newcommand{\Ppratedmax}[0]{\textcolor{\ratingboundscolor}{\symPower'^{\maxss}_{\exoaddss}}}
\newcommand{\Pmratedmax}[0]{\textcolor{\ratingboundscolor}{\symPower'^{\maxss}_{\exoremss}}}

\newcommand{\Pcratedmin}[0]{\textcolor{\ratingboundscolor}{\symPower'^{\minss}_{\chargess}}}
\newcommand{\Pdratedmin}[0]{\textcolor{\ratingboundscolor}{\symPower'^{\minss}_{\dischargess}}}
\newcommand{\Ppratedmin}[0]{\textcolor{\ratingboundscolor}{\symPower'^{\minss}_{\exoaddss}}}
\newcommand{\Pmratedmin}[0]{\textcolor{\ratingboundscolor}{\symPower'^{\minss}_{\exoremss}}}

\newcommand{\Paddratedmin}[0]{\textcolor{\ratingboundscolor}{\symPower'^{\minss}_{\addss}}}
\newcommand{\Paddratedmax}[0]{\textcolor{\ratingboundscolor}{\symPower'^{\maxss}_{\addss}}}
\newcommand{\Premratedmin}[0]{\textcolor{\ratingboundscolor}{\symPower'^{\minss}_{\remss}}}
\newcommand{\Premratedmax}[0]{\textcolor{\ratingboundscolor}{\symPower'^{\maxss}_{\remss}}}

\newcommand{\Pflowrated}[0]{\textcolor{\sizingcolor}{\symPower'^{}_{\flowss}}}
\newcommand{\Pflowratedmin}[0]{\textcolor{\ratingboundscolor}{\symPower'^{\minss}_{\flowss}}}
\newcommand{\Pflowratedmax}[0]{\textcolor{\ratingboundscolor}{\symPower'^{\maxss}_{\flowss}}}

\newcommand{\Eratedmaxmax}[0]{\textcolor{\ratingboundscolor}{\symEnergy'^{\maxss\maxss}}}
\newcommand{\Eratedmaxmin}[0]{\textcolor{\ratingboundscolor}{\symEnergy'^{\maxss\minss}}}
\newcommand{\Eratedmax}[0]{\textcolor{\sizingcolor}{\symEnergy'^{\maxss}}}
\newcommand{\Eratedminmax}[0]{\textcolor{\ratingboundscolor}{\symEnergy'^{\minss\maxss}}}
\newcommand{\Eratedminmin}[0]{\textcolor{\ratingboundscolor}{\symEnergy'^{\minss\minss}}}
\newcommand{\Eratedmin}[0]{\textcolor{\sizingcolor}{\symEnergy'^{\minss}}}

\newcommand{\Pcrated}[0]{\textcolor{\sizingcolor}{\symPower'^{}_{\chargess}}}
\newcommand{\Pdrated}[0]{\textcolor{\sizingcolor}{\symPower'^{}_{\dischargess}}}
\newcommand{\Pprated}[0]{\textcolor{\sizingcolor}{\symPower'^{}_{\exoaddss}}}
\newcommand{\Pmrated}[0]{\textcolor{\sizingcolor}{\symPower'^{}_{\exoremss}}}

\newcommand{\Paddrated}[0]{\textcolor{\sizingcolor}{\symPower'^{}_{\addss}}}
\newcommand{\Premrated}[0]{\textcolor{\sizingcolor}{\symPower'^{}_{\remss}}}

\newcommand{\Pcorthk}[0]{\symPower^{\orthss}_{\chargess} }
\newcommand{\Pdorthk}[0]{\symPower^{\orthss}_{\dischargess} }
\newcommand{\Pporthk}[0]{\symPower^{\orthss}_{\exoaddss} }
\newcommand{\Pmorthk}[0]{\symPower^{\orthss}_{\exoremss} }
\newcommand{\Porthk}[0]{\symPower^{\orthss}_{} }

\newcommand{\Paddorthk}[0]{\symPower^{\orthss}_{\addss} }
\newcommand{\Premorthk}[0]{\symPower^{\orthss}_{\remss} }
\newcommand{\Pfloworthk}[0]{\symPower^{\orthss}_{\flowss} }
\newcommand{\Paddorthtotk}[0]{\symPower^{\orthss}_{\addtextss} }
\newcommand{\Premorthtotk}[0]{\symPower^{\orthss}_{\remtextss} }

\newcommand{\Pcorthkmax}[0]{\textcolor{\boundscolor}{\symPower^{\orthss, \maxss}_{\chargess} }}
\newcommand{\Pdorthkmax}[0]{\textcolor{\boundscolor}{\symPower^{\orthss, \maxss}_{\dischargess} }}
\newcommand{\Pporthkmax}[0]{\textcolor{\boundscolor}{\symPower^{\orthss, \maxss}_{\exoaddss} }}
\newcommand{\Pmorthkmax}[0]{\textcolor{\boundscolor}{\symPower^{\orthss, \maxss}_{\exoremss} }}

\newcommand{\Paddorthkmax}[0]{\textcolor{\boundscolor}{\symPower^{\orthss, \maxss}_{\addss} }}
\newcommand{\Premorthkmax}[0]{\textcolor{\boundscolor}{\symPower^{\orthss, \maxss}_{\remss} }}
\newcommand{\Pfloworthkmax}[0]{\textcolor{\boundscolor}{\symPower^{\orthss, \maxss}_{\flowss} }}

\newcommand{\Pcdk}[0]{\symPower_{\chargess\dischargess} }
\newcommand{\Pcmk}[0]{\symPower_{\chargess\exoremss} }
\newcommand{\Ppmk}[0]{\symPower_{\exoaddss\exoremss} }
\newcommand{\Ppdk}[0]{\symPower_{\exoaddss\dischargess} }

\newcommand{\Pcirck}[0]{\symPower_{\circss} }
\newcommand{\Pcircmax}[0]{\symPower^{\maxss}_{\circss}}

\newcommand{\Pcdmax}[0]{\symPower^{\maxss}_{\chargess\dischargess}}
\newcommand{\Pcmmax}[0]{\symPower^{\maxss}_{\chargess\exoremss}}
\newcommand{\Ppmmax}[0]{\symPower^{\maxss}_{\exoaddss\exoremss}}
\newcommand{\Ppdmax}[0]{\symPower^{\maxss}_{\exoaddss\dischargess}}

\newcommand{\Pcpk}[0]{\symPower_{\chargess\exoaddss} }
\newcommand{\Pdmk}[0]{\symPower_{\dischargess\exoremss} }
\newcommand{\gammacp}[0]{\textcolor{\paramcolor}{\symPowerToEnergy_{\chargess\exoaddss} }}
\newcommand{\gammadm}[0]{\textcolor{\paramcolor}{\symPowerToEnergy_{\dischargess\exoremss} }}

\newcommand{\gammaadd}[0]{\textcolor{\paramcolor}{\symPowerToEnergy_{\addss} }}
\newcommand{\gammarem}[0]{\textcolor{\paramcolor}{\symPowerToEnergy_{\remss} }}

\newcommand{\dcd}[0]{\textcolor{\decisioncolor}{\symDecision_{\chargess\dischargess} }}
\newcommand{\dcm}[0]{\textcolor{\decisioncolor}{\symDecision_{\chargess\exoremss} }}
\newcommand{\dpm}[0]{\textcolor{\decisioncolor}{\symDecision_{\exoaddss\exoremss} }}
\newcommand{\dpd}[0]{\textcolor{\decisioncolor}{\symDecision_{\exoaddss\dischargess} }}

\newcommand{\dcirc}[0]{\textcolor{\decisioncolor}{\symDecision_{\circss} }}

\newcommand{\dcorth}[0]{\textcolor{\decisioncolor}{\symDecision^{\orthss}_{\chargess} }}
\newcommand{\ddorth}[0]{\textcolor{\decisioncolor}{\symDecision^{\orthss}_{\dischargess} }}
\newcommand{\dporth}[0]{\textcolor{\decisioncolor}{\symDecision^{\orthss}_{\exoaddss} }}
\newcommand{\dmorth}[0]{\textcolor{\decisioncolor}{\symDecision^{\orthss}_{\exoremss} }}

\newcommand{\daddorth}[0]{\textcolor{\decisioncolor}{\symDecision^{\orthss}_{\addss} }}
\newcommand{\dremorth}[0]{\textcolor{\decisioncolor}{\symDecision^{\orthss}_{\remss} }}
\newcommand{\dfloworth}[0]{\textcolor{\decisioncolor}{\symDecision^{\orthss}_{\flowss} }}

\newcommand{\dclevels}[0]{\textcolor{\decisioncolor}{\symDecision^{}_{\chargess}  }}
\newcommand{\ddlevels}[0]{\textcolor{\decisioncolor}{\symDecision^{}_{\dischargess}  }}
\newcommand{\dplevels}[0]{\textcolor{\decisioncolor}{\symDecision^{}_{\exoaddss}  }}
\newcommand{\dmlevels}[0]{\textcolor{\decisioncolor}{\symDecision^{}_{\exoremss}  }}
\newcommand{\dflowlevels}[0]{\textcolor{\decisioncolor}{\symDecision^{}_{\flowss}  }}

\newcommand{\dcmaxlevels}[0]{\textcolor{\paramcolor}{\symDecision^{\maxss}_{\chargess}  }}
\newcommand{\ddmaxlevels}[0]{\textcolor{\paramcolor}{\symDecision^{\maxss}_{\dischargess}  }}
\newcommand{\dpmaxlevels}[0]{\textcolor{\paramcolor}{\symDecision^{\maxss}_{\exoaddss}  }}
\newcommand{\dmmaxlevels}[0]{\textcolor{\paramcolor}{\symDecision^{\maxss}_{\exoremss}  }}
\newcommand{\dflowmaxlevels}[0]{\textcolor{\paramcolor}{\symDecision^{\maxss}_{\flowss}  }}


\newcommand{\dE}[0]{\textcolor{\decisioncolor}{\symDecision_{\Ess}}}

\newcommand{\Emax}[0]{\textcolor{\sizingcolor}{\symEnergy'^{\maxss}}}
\newcommand{\Emin}[0]{\textcolor{\sizingcolor}{\symEnergy'^{\minss}}}
\newcommand{\Eeff}[0]{\textcolor{\sizingcolor}{\symEnergy'^{\usabless}}}
\newcommand{\Ekmax}[0]{\textcolor{\boundscolor}{\symEnergy^{\maxss} }}
\newcommand{\Ekmin}[0]{\textcolor{\boundscolor}{\symEnergy^{\minss} }}
\newcommand{\Ekeff}[0]{\textcolor{\boundscolor}{\symEnergy^{\usabless} }}
\newcommand{\Ek}[0]{\symEnergy }
\newcommand{\Ekprev}[0]{\symEnergy prev}
\newcommand{\iEk}[0]{\textcolor{\binarycolor}{\symBinary_{\Ess} }}
\newcommand{\iEkprev}[0]{\textcolor{\binarycolor}{\symBinary_{\Ess} prev}}

\newcommand{\iijk}[0]{\textcolor{\binarycolor}{\symBinary_{\indexGridNode \indexGridNodeTwo} }}
\newcommand{\ijik}[0]{\textcolor{\binarycolor}{\symBinary_{\indexGridNodeTwo \indexGridNode} }}
\newcommand{\ilk}[0]{\textcolor{\binarycolor}{\symBinary_{\indexGridLines} }}
\newcommand{\ilkprev}[0]{\textcolor{\binarycolor}{\symBinary_{\indexGridLines} prev}}
\newcommand{\ilkmin}[0]{\textcolor{\boundscolor}{\symBinary^{\minss}_{\indexGridLines} }}
\newcommand{\ilkmax}[0]{\textcolor{\boundscolor}{\symBinary^{\maxss}_{\indexGridLines} }}

\newcommand{\betaijk}[0]{\textcolor{\binarycolor}{\symRadiality_{\indexGridNode \indexGridNodeTwo} }}
\newcommand{\betajik}[0]{\textcolor{\binarycolor}{\symRadiality_{\indexGridNodeTwo \indexGridNode} }}


\newcommand{\Sgik}[0]{\textcolor{\complexcolor}{\symApparentPower_{\giss} }}
\newcommand{\Sauxk}[0]{\symApparentPower_{\text{aux}} }
\newcommand{\SauxkA}[0]{\symApparentPower_{\text{aux},\symPhaseA} }
\newcommand{\SauxkB}[0]{\symApparentPower_{\text{aux},\symPhaseB} }
\newcommand{\SauxkC}[0]{\symApparentPower_{\text{aux},\symPhaseC} }
\newcommand{\Sgirated}[0]{\textcolor{\sizingcolor}{\symApparentPower_{\giss}^{'}}}
\newcommand{\lgik}[0]{\textcolor{\paramcolor}{\symLocation_{\giss} }}
\newcommand{\Pgik}[0]{\symPower_{\giss} }
\newcommand{\Qgik}[0]{\symReactivePower_{\giss} }
\newcommand{\Qgikref}[0]{\textcolor{\boundscolor}{\symReactivePower^{\refss}_{\giss} }}
\newcommand{\Qgikmax}[0]{\textcolor{\boundscolor}{\symReactivePower^{\maxss}_{\giss} }}
\newcommand{\Qgikmin}[0]{\textcolor{\boundscolor}{\symReactivePower^{\minss}_{\giss} }}
\newcommand{\Pctotk}[0]{\symPower_{\chargess,\totss} }
\newcommand{\Pdtotk}[0]{\symPower_{\dischargess,\totss} }
\newcommand{\Pstandbyk}[0]{\symPower_{\standbyss} }
\newcommand{\PSlossk}[0]{\symPower_{\Sss,\lossss} }
\newcommand{\Punblossk}[0]{\symPower_{\unbss,\lossss} }
\newcommand{\etagi}{\textcolor{\paramcolor}{\eta_\giss  }}

\newcommand{\Pstandbykmin}[0]{\textcolor{\boundscolor}{\symPower_{\standbyss}^{\minss}  }}
\newcommand{\Pstandbykmax}[0]{\textcolor{\boundscolor}{\symPower_{\standbyss}^{\maxss}  }}

\newcommand{\aPZk}[0]{\textcolor{\paramcolor}{a_{  }^{ \Zss\Pss }  }}
\newcommand{\aPIk}[0]{\textcolor{\paramcolor}{a_{  }^{ \Iss\Pss }  }}
\newcommand{\aPPk}[0]{\textcolor{\paramcolor}{a_{  }^{ \Pss\Pss }  }}
\newcommand{\aQZk}[0]{\textcolor{\paramcolor}{a_{  }^{ \Zss\Qss }  }}
\newcommand{\aQIk}[0]{\textcolor{\paramcolor}{a_{  }^{ \Iss\Qss }  }}
\newcommand{\aQPk}[0]{\textcolor{\paramcolor}{a_{  }^{ \Pss\Qss }  }}

\newcommand{\aPZ}[0]{\textcolor{\paramcolor}{a_{ \indexZIP }^{ \Zss \Pss} }}
\newcommand{\aPI}[0]{\textcolor{\paramcolor}{a_{ \indexZIP }^{ \Iss \Pss} }}
\newcommand{\aPP}[0]{\textcolor{\paramcolor}{a_{\indexZIP  }^{ \Pss\Pss } }}
\newcommand{\aQZ}[0]{\textcolor{\paramcolor}{a_{  \indexZIP}^{ \Zss\Qss } }}
\newcommand{\aQI}[0]{\textcolor{\paramcolor}{a_{ \indexZIP }^{ \Iss \Qss} }}
\newcommand{\aQP}[0]{\textcolor{\paramcolor}{a_{ \indexZIP }^{ \Pss\Qss } }}

\newcommand{\PgikA}[0]{\symPower_{\giss,\symPhaseA} }
\newcommand{\QgikA}[0]{\symReactivePower_{\giss,\symPhaseA} }
\newcommand{\SgikA}[0]{\textcolor{\complexcolor}{\symApparentPower_{\giss,\symPhaseA} }}
\newcommand{\PgikB}[0]{\symPower_{\giss,\symPhaseB} }
\newcommand{\QgikB}[0]{\symReactivePower_{\giss,\symPhaseB} }
\newcommand{\SgikB}[0]{\textcolor{\complexcolor}{\symApparentPower_{\giss,\symPhaseB} }}
\newcommand{\PgikC}[0]{\symPower_{\giss,\symPhaseC} }
\newcommand{\QgikC}[0]{\symReactivePower_{\giss,\symPhaseC} }
\newcommand{\SgikC}[0]{\textcolor{\complexcolor}{\symApparentPower_{\giss,\symPhaseC} }}
\newcommand{\SgiratedA}[0]{\textcolor{\sizingcolor}{\symApparentPower_{\giss,\symPhaseA}^{'}}}
\newcommand{\SgiratedB}[0]{\textcolor{\sizingcolor}{\symApparentPower_{\giss,\symPhaseB}^{'}}}
\newcommand{\SgiratedC}[0]{\textcolor{\sizingcolor}{\symApparentPower_{\giss,\symPhaseC}^{'}}}

\newcommand{\Sgiratedmin}[0]{\textcolor{\ratingboundscolor}{\symApparentPower_{\giss}^{\minss}}}
\newcommand{\SgiratedminA}[0]{\textcolor{\ratingboundscolor}{\symApparentPower_{\giss,\symPhaseA}^{\minss}}}
\newcommand{\SgiratedminB}[0]{\textcolor{\ratingboundscolor}{\symApparentPower_{\giss,\symPhaseB}^{\minss}}}
\newcommand{\SgiratedminC}[0]{\textcolor{\ratingboundscolor}{\symApparentPower_{\giss,\symPhaseC}^{\minss}}}
\newcommand{\Sgiratedmax}[0]{\textcolor{\ratingboundscolor}{\symApparentPower_{\giss}^{\maxss}}}
\newcommand{\SgiratedmaxA}[0]{\textcolor{\ratingboundscolor}{\symApparentPower_{\giss,\symPhaseA}^{\maxss}}}
\newcommand{\SgiratedmaxB}[0]{\textcolor{\ratingboundscolor}{\symApparentPower_{\giss,\symPhaseB}^{\maxss}}}
\newcommand{\SgiratedmaxC}[0]{\textcolor{\ratingboundscolor}{\symApparentPower_{\giss,\symPhaseC}^{\maxss}}}

\newcommand{\Pgikunb}[0]{\symPower_{\giss,\unbss} }
\newcommand{\Qgikunb}[0]{\symReactivePower_{\giss,\unbss} }

\newcommand{\rhoPgikunb}[0]{\textcolor{\paramcolor}{\symLossFactor_{\giss,\Pss, \unbss} }}
\newcommand{\rhoQgikunb}[0]{\textcolor{\paramcolor}{\symLossFactor_{\giss,\Qss,\unbss} }}

\newcommand{\bgik}[0]{\textcolor{\paramcolor}{\symAvailability_{\giss} }}
\newcommand{\bgikA}[0]{\textcolor{\paramcolor}{\symAvailability_{\giss, \symPhaseA} }}
\newcommand{\bgikB}[0]{\textcolor{\paramcolor}{\symAvailability_{\giss, \symPhaseB} }}
\newcommand{\bgikC}[0]{\textcolor{\paramcolor}{\symAvailability_{\giss, \symPhaseC} }}

\newcommand{\dPtogik}[0]{\textcolor{\paramcolor}{\symDirectionality_{\giss}^{\Pss\plusss}   }}
\newcommand{\dPtogikA}[0]{\textcolor{\paramcolor}{\symDirectionality_{\giss, \symPhaseA}^{\Pss\plusss} }}
\newcommand{\dPtogikB}[0]{\textcolor{\paramcolor}{\symDirectionality_{\giss, \symPhaseB}^{\Pss\plusss} }}
\newcommand{\dPtogikC}[0]{\textcolor{\paramcolor}{\symDirectionality_{\giss, \symPhaseC}^{\Pss\plusss} }}

\newcommand{\dQtogik}[0]{\textcolor{\paramcolor}{\symDirectionality_{\giss}^{\Qss\plusss}  }}
\newcommand{\dQtogikA}[0]{\textcolor{\paramcolor}{\symDirectionality_{\giss, \symPhaseA}^{\Qss\plusss}  }}
\newcommand{\dQtogikB}[0]{\textcolor{\paramcolor}{\symDirectionality_{\giss, \symPhaseB}^{\Qss\plusss}  }}
\newcommand{\dQtogikC}[0]{\textcolor{\paramcolor}{\symDirectionality_{\giss, \symPhaseC}^{\Qss\plusss}  }}

\newcommand{\dPfromgik}[0]{\textcolor{\paramcolor}{\symDirectionality_{\giss}^{\Pss\subsss} }}
\newcommand{\dPfromgikA}[0]{\textcolor{\paramcolor}{\symDirectionality_{\giss, \symPhaseA}^{\Pss\subsss} }}
\newcommand{\dPfromgikB}[0]{\textcolor{\paramcolor}{\symDirectionality_{\giss, \symPhaseB}^{\Pss\subsss} }}
\newcommand{\dPfromgikC}[0]{\textcolor{\paramcolor}{\symDirectionality_{\giss, \symPhaseC}^{\Pss\subsss} }}

\newcommand{\dQfromgik}[0]{\textcolor{\paramcolor}{\symDirectionality_{\giss}^{\Qss\subsss}  }}
\newcommand{\dQfromgikA}[0]{\textcolor{\paramcolor}{\symDirectionality_{\giss, \symPhaseA}^{\Qss\subsss}  }}
\newcommand{\dQfromgikB}[0]{\textcolor{\paramcolor}{\symDirectionality_{\giss, \symPhaseB}^{\Qss\subsss}  }}
\newcommand{\dQfromgikC}[0]{\textcolor{\paramcolor}{\symDirectionality_{\giss, \symPhaseC}^{\Qss\subsss}  }}

\newcommand{\rhoDODmax}[0]{\textcolor{\paramcolor}{\symLossFactor_{\Ess}^{\maxss}}}
\newcommand{\rhoDODmin}[0]{\textcolor{\paramcolor}{\symLossFactor_{\Ess}^{\minss}}}
\newcommand{\rhoDODeff}[0]{\textcolor{\paramcolor}{\symLossFactor_{\Ess}^{\usabless}}}

\newcommand{\Eijk}[0]{\textcolor{\complexcolor}{\symEnergy_{\indexGridNode \indexGridNodeTwo}^{}  }}
\newcommand{\Eijkre}[0]{\textcolor{black}{\symEnergy_{\indexGridNode \indexGridNodeTwo}^{\realss}  }}

\newcommand{\Sij}[0]{\textcolor{\complexcolor}{\symApparentPower_{\indexGridNode \indexGridNodeTwo}^{} }}

\newcommand{\Pzref}[0]{\textcolor{\paramcolor}{\symPower_{\indexZIP }^{\refss} }}
\newcommand{\Pinode}[0]{\textcolor{black}{\symPower_{\indexGridNode}^{} }}
\newcommand{\Pz}[0]{\textcolor{black}{\symPower_{\indexZIP}^{} }}

\newcommand{\Pij}[0]{\textcolor{black}{\symPower_{\indexGridNode \indexGridNodeTwo}^{} }}
\newcommand{\Pji}[0]{\textcolor{black}{\symPower_{\indexGridNodeTwo \indexGridNode }^{} }}
\newcommand{\Pijloss}[0]{\textcolor{black}{\symPower_{ \indexGridLines}^{\lossss} }}
\newcommand{\Pijmax}[0]{\textcolor{\paramcolor}{\symPower_{\indexGridNode \indexGridNodeTwo}^{\maxss} }}
\newcommand{\Pijacc}[0]{\textcolor{black}{\symPower_{\indexGridNode \indexGridNodeTwo}^{\star} }}
\newcommand{\Pijrated}[0]{\textcolor{\sizingcolor}{\symPower_{\indexGridNode \indexGridNodeTwo}^{\ratedss} }}
\newcommand{\Pjirated}[0]{\textcolor{\sizingcolor}{\symPower_{\indexGridNodeTwo\indexGridNode }^{\ratedss} }}

\newcommand{\Qij}[0]{\textcolor{black}{\symReactivePower_{\indexGridNode \indexGridNodeTwo}^{} }}
\newcommand{\Qji}[0]{\textcolor{black}{\symReactivePower_{\indexGridNodeTwo \indexGridNode }^{} }}
\newcommand{\Qijrated}[0]{\textcolor{\sizingcolor}{\symReactivePower_{\indexGridNode \indexGridNodeTwo}^{\ratedss} }}
\newcommand{\Qjirated}[0]{\textcolor{\sizingcolor}{\symReactivePower_{\indexGridNodeTwo\indexGridNode }^{\ratedss} }}
\newcommand{\Qijmax}[0]{\textcolor{\paramcolor}{\symReactivePower_{\indexGridNode \indexGridNodeTwo}^{\maxss} }}
\newcommand{\Qijmin}[0]{\textcolor{\paramcolor}{\symReactivePower_{\indexGridNode \indexGridNodeTwo}^{\minss} }}
\newcommand{\Qjimax}[0]{\textcolor{\paramcolor}{\symReactivePower_{\indexGridNodeTwo \indexGridNode }^{\maxss} }}
\newcommand{\Qjimin}[0]{\textcolor{\paramcolor}{\symReactivePower_{\indexGridNodeTwo\indexGridNode }^{\minss} }}
\newcommand{\Qijloss}[0]{\textcolor{black}{\symReactivePower_{ \indexGridLines}^{\lossss} }}

\newcommand{\Sjiks}[0]{\textcolor{\complexcolor}{\symApparentPower_{\indexGridNodeTwo \indexGridNode , \seriesss}^{}  }}
\newcommand{\Sijks}[0]{\textcolor{\complexcolor}{\symApparentPower_{\indexGridNode \indexGridNodeTwo, \seriesss}^{}  }}
\newcommand{\Pijks}[0]{\textcolor{black}{\symPower_{\indexGridNode \indexGridNodeTwo, \seriesss}^{}  }}
\newcommand{\Qijks}[0]{\textcolor{black}{\symReactivePower_{\indexGridNode \indexGridNodeTwo, \seriesss}^{}  }}
\newcommand{\Pjiks}[0]{\textcolor{black}{\symPower_{\indexGridNodeTwo \indexGridNode , \seriesss}^{}  }}
\newcommand{\Qjiks}[0]{\textcolor{black}{\symReactivePower_{\indexGridNodeTwo \indexGridNode , \seriesss}^{}  }}

\newcommand{\Sijlossk}[0]{\textcolor{\complexcolor}{\symApparentPower_{ \indexGridLines}^{\lossss}  }}
\newcommand{\Slossk}[0]{\textcolor{black}{\symApparentPower_{}^{\lossss}  }}

\newcommand{\Sijlosssk}[0]{\textcolor{\complexcolor}{\symApparentPower_{ \indexGridLines, \seriesss}^{\lossss}  }}

\newcommand{\Plossk}[0]{\textcolor{black}{\symPower_{}^{\lossss}  }}
\newcommand{\Pijk}[0]{\textcolor{black}{\symPower_{\indexGridNode \indexGridNodeTwo}^{}  }}
\newcommand{\Pijlossk}[0]{\textcolor{black}{\symPower_{ \indexGridLines}^{\lossss}  }}
\newcommand{\Pijlosssk}[0]{\textcolor{black}{\symPower_{ \indexGridLines, \seriesss}^{\lossss}  }}
\newcommand{\Pijlossshk}[0]{\textcolor{black}{\symPower_{ \indexGridLines, \shuntss}^{\lossss}  }}
\newcommand{\Pilossshk}[0]{\textcolor{black}{\symPower_{\indexGridNode \indexGridNodeTwo, \shuntss }^{\lossss}  }}
\newcommand{\Pjlossshk}[0]{\textcolor{black}{\symPower_{ \indexGridNodeTwo\indexGridNode, \shuntss}^{\lossss}  }}
\newcommand{\Pjik}[0]{\textcolor{black}{\symPower_{\indexGridNodeTwo \indexGridNode}^{}    }}
\newcommand{\Pinodek}[0]{\textcolor{black}{\symPower_{ \indexGridNode}^{}    }}
\newcommand{\Qijk}[0]{\textcolor{black}{\symReactivePower_{\indexGridNode \indexGridNodeTwo}^{}  }}
\newcommand{\Qijkdelta}[0]{\textcolor{black}{\symReactivePower_{\indexGridNode \indexGridNodeTwo}^{\Delta}  }}
\newcommand{\Qijlossk}[0]{\textcolor{black}{\symReactivePower_{ \indexGridLines}^{\lossss}  }}
\newcommand{\Qijlosssk}[0]{\textcolor{black}{\symReactivePower_{ \indexGridLines,\seriesss}^{\lossss}  }}
\newcommand{\Qijlossshk}[0]{\textcolor{black}{\symReactivePower_{ \indexGridLines,\shuntss}^{\lossss}  }}
\newcommand{\Qilossshk}[0]{\textcolor{black}{\symReactivePower_{\indexGridNode\indexGridNodeTwo,\shuntss }^{\lossss}  }}
\newcommand{\Qjlossshk}[0]{\textcolor{black}{\symReactivePower_{ \indexGridNodeTwo\indexGridNode,\shuntss}^{\lossss}  }}
\newcommand{\Qjik}[0]{\textcolor{black}{\symReactivePower_{\indexGridNodeTwo \indexGridNode}^{}    }}
\newcommand{\Qinodek}[0]{\textcolor{black}{\symReactivePower_{ \indexGridNode}^{}    }}
\newcommand{\Sinodek}[0]{\textcolor{\complexcolor}{\symApparentPower_{ \indexGridNode}^{}    }}

\newcommand{\Wijk}[0]{\textcolor{\complexcolor}{W_{\indexGridNode \indexGridNodeTwo}^{}  }} 
\newcommand{\WijkSDP}[0]{\textcolor{\complexcolor}{\mathbf{W}_{}^{}  }} 
\newcommand{\Wiik}[0]{\textcolor{black}{W_{\indexGridNode \indexGridNode}^{}  }} 
\newcommand{\Wjjk}[0]{\textcolor{black}{W_{\indexGridNodeTwo \indexGridNodeTwo}^{}  }} 
\newcommand{\Wjik}[0]{\textcolor{\complexcolor}{W_{\indexGridNodeTwo \indexGridNode}^{}  }} 

\newcommand{\Rijk}[0]{\textcolor{black}{R_{\indexGridNode \indexGridNodeTwo}^{}  }} 
\newcommand{\Tijk}[0]{\textcolor{black}{T_{\indexGridNode \indexGridNodeTwo}^{}  }}  
\newcommand{\Rjik}[0]{\textcolor{black}{R_{\indexGridNodeTwo\indexGridNode }^{}  }} 
\newcommand{\Tjik}[0]{\textcolor{black}{T_{\indexGridNodeTwo\indexGridNode }^{}  }} 
\newcommand{\Rlk}[0]{\textcolor{black}{R_{\indexGridLines }^{}  }}    
\newcommand{\rijk}[0]{\textcolor{black}{r_{\indexGridNode \indexGridNodeTwo}^{}  }}   

\newcommand{\SikA}[0]{\textcolor{\complexcolor}{\symApparentPower_{\indexGridNode , \symPhaseA}^{}  }}
\newcommand{\SikB}[0]{\textcolor{\complexcolor}{\symApparentPower_{\indexGridNode , \symPhaseB}^{}  }}
\newcommand{\SikC}[0]{\textcolor{\complexcolor}{\symApparentPower_{\indexGridNode , \symPhaseC}^{}  }}

\newcommand{\PijkA}[0]{\textcolor{black}{\symPower_{\indexGridNode \indexGridNodeTwo, \symPhaseA}^{}  }}
\newcommand{\PijkB}[0]{\textcolor{black}{\symPower_{\indexGridNode \indexGridNodeTwo, \symPhaseB}^{}  }}
\newcommand{\PijkC}[0]{\textcolor{black}{\symPower_{\indexGridNode \indexGridNodeTwo, \symPhaseC}^{}  }}
\newcommand{\QijkA}[0]{\textcolor{black}{\symReactivePower_{\indexGridNode \indexGridNodeTwo, \symPhaseA}^{}  }}
\newcommand{\QijkB}[0]{\textcolor{black}{\symReactivePower_{\indexGridNode \indexGridNodeTwo, \symPhaseB}^{}  }}
\newcommand{\QijkC}[0]{\textcolor{black}{\symReactivePower_{\indexGridNode \indexGridNodeTwo, \symPhaseC}^{}  }}

\newcommand{\PjikA}[0]{\textcolor{black}{\symPower_{ \indexGridNodeTwo\indexGridNode, \symPhaseA}^{}  }}
\newcommand{\PjikB}[0]{\textcolor{black}{\symPower_{ \indexGridNodeTwo\indexGridNode, \symPhaseB}^{}  }}
\newcommand{\PjikC}[0]{\textcolor{black}{\symPower_{ \indexGridNodeTwo\indexGridNode, \symPhaseC}^{}  }}
\newcommand{\QjikA}[0]{\textcolor{black}{\symReactivePower_{ \indexGridNodeTwo\indexGridNode, \symPhaseA}^{}  }}
\newcommand{\QjikB}[0]{\textcolor{black}{\symReactivePower_{ \indexGridNodeTwo\indexGridNode, \symPhaseB}^{}  }}
\newcommand{\QjikC}[0]{\textcolor{black}{\symReactivePower_{ \indexGridNodeTwo\indexGridNode, \symPhaseC}^{}  }}

\newcommand{\PijlosskA}[0]{\textcolor{black}{\symPower_{\indexGridNode \indexGridNodeTwo,\symPhaseA}^{\lossss}  }}
\newcommand{\PijlosskB}[0]{\textcolor{black}{\symPower_{\indexGridNode \indexGridNodeTwo,\symPhaseB}^{\lossss}  }}
\newcommand{\PijlosskC}[0]{\textcolor{black}{\symPower_{\indexGridNode \indexGridNodeTwo,\symPhaseC}^{\lossss}  }}
\newcommand{\QijlosskA}[0]{\textcolor{black}{\symReactivePower_{\indexGridNode \indexGridNodeTwo,\symPhaseA}^{\lossss}  }}
\newcommand{\QijlosskB}[0]{\textcolor{black}{\symReactivePower_{\indexGridNode \indexGridNodeTwo,\symPhaseB}^{\lossss}  }}
\newcommand{\QijlosskC}[0]{\textcolor{black}{\symReactivePower_{\indexGridNode \indexGridNodeTwo,\symPhaseC}^{\lossss}  }}

\newcommand{\PijlosskminA}[0]{\textcolor{\boundscolor}{\symPower_{\indexGridNode \indexGridNodeTwo,\symPhaseA}^{\lossss,\minss}  }}
\newcommand{\PijlosskminB}[0]{\textcolor{\boundscolor}{\symPower_{\indexGridNode \indexGridNodeTwo,\symPhaseB}^{\lossss,\minss}  }}
\newcommand{\PijlosskminC}[0]{\textcolor{\boundscolor}{\symPower_{\indexGridNode \indexGridNodeTwo,\symPhaseC}^{\lossss,\minss}  }}
\newcommand{\QijlosskminA}[0]{\textcolor{\boundscolor}{\symReactivePower_{\indexGridNode \indexGridNodeTwo,\symPhaseA}^{\lossss,\minss}  }}
\newcommand{\QijlosskminB}[0]{\textcolor{\boundscolor}{\symReactivePower_{\indexGridNode \indexGridNodeTwo,\symPhaseB}^{\lossss,\minss}  }}
\newcommand{\QijlosskminC}[0]{\textcolor{\boundscolor}{\symReactivePower_{\indexGridNode \indexGridNodeTwo,\symPhaseC}^{\lossss,\minss}  }}

\newcommand{\PijlosskmaxA}[0]{\textcolor{\boundscolor}{\symPower_{\indexGridNode \indexGridNodeTwo,\symPhaseA}^{\lossss,\maxss}  }}
\newcommand{\PijlosskmaxB}[0]{\textcolor{\boundscolor}{\symPower_{\indexGridNode \indexGridNodeTwo,\symPhaseB}^{\lossss,\maxss}  }}
\newcommand{\PijlosskmaxC}[0]{\textcolor{\boundscolor}{\symPower_{\indexGridNode \indexGridNodeTwo,\symPhaseC}^{\lossss,\maxss}  }}
\newcommand{\QijlosskmaxA}[0]{\textcolor{\boundscolor}{\symReactivePower_{\indexGridNode \indexGridNodeTwo,\symPhaseA}^{\lossss,\maxss}  }}
\newcommand{\QijlosskmaxB}[0]{\textcolor{\boundscolor}{\symReactivePower_{\indexGridNode \indexGridNodeTwo,\symPhaseB}^{\lossss,\maxss}  }}
\newcommand{\QijlosskmaxC}[0]{\textcolor{\boundscolor}{\symReactivePower_{\indexGridNode \indexGridNodeTwo,\symPhaseC}^{\lossss,\maxss}  }}

\newcommand{\Iijkabs}[0]{\textcolor{black}{|\symCurrent_{\indexGridLines\indexGridNode \indexGridNodeTwo}^{}  |}}
\newcommand{\Sijkabs}[0]{\textcolor{black}{|\symApparentPower_{\indexGridLines\indexGridNode \indexGridNodeTwo}^{}  |}}
\newcommand{\Sijk}[0]{\textcolor{\complexcolor}{\symApparentPower_{\indexGridLines\indexGridNode \indexGridNodeTwo}^{}  }}

\newcommand{\SijkSDPrated}[0]{\textcolor{\paramcolor}{\mathbf{\symApparentPower}_{\indexGridLines\indexGridNode \indexGridNodeTwo}^{\ratedss}  }}
\newcommand{\IijkSDPrated}[0]{\textcolor{\paramcolor}{\mathbf{\symCurrent}_{\indexGridLines\indexGridNode \indexGridNodeTwo}^{\ratedss}  }}
\newcommand{\IijksSDPrated}[0]{\textcolor{\paramcolor}{\mathbf{\symCurrent}_{\indexGridLines }^{\seriesss, \ratedss}  }}

\newcommand{\SijkSDPseq}[0]{\textcolor{\complexcolor}{\mathbf{\symApparentPower}_{\indexGridLines\indexGridNode \indexGridNodeTwo}^{\fortescuess}  }}
\newcommand{\SjikSDPseq}[0]{\textcolor{\complexcolor}{\mathbf{\symApparentPower}_{\indexGridLines\indexGridNodeTwo \indexGridNode }^{\fortescuess} }}
\newcommand{\SijkSDPH}[0]{\textcolor{\complexcolor}{(\mathbf{\symApparentPower}_{\indexGridLines\indexGridNode \indexGridNodeTwo})^{\hermitiantranspose}  }}
\newcommand{\SijkSDP}[0]{\textcolor{\complexcolor}{\mathbf{\symApparentPower}_{\indexGridLines\indexGridNode \indexGridNodeTwo}^{}  }}
\newcommand{\SjikSDP}[0]{\textcolor{\complexcolor}{\mathbf{\symApparentPower}_{\indexGridLines\indexGridNodeTwo \indexGridNode }^{}  }}
\newcommand{\SjkkSDP}[0]{\textcolor{\complexcolor}{\mathbf{\symApparentPower}_{ \indexGridLines\indexGridNodeTwo \indexGridNodeThree}^{}  }}
\newcommand{\PjkkSDP}[0]{\textcolor{black}{\mathbf{\symPower}_{ \indexGridLines\indexGridNodeTwo \indexGridNodeThree}^{}  }}
\newcommand{\QjkkSDP}[0]{\textcolor{black}{\mathbf{\symReactivePower}_{ \indexGridLines\indexGridNodeTwo \indexGridNodeThree}^{}  }}

\newcommand{\SijkSDPkron}[0]{\textcolor{\complexcolor}{\mathbf{\symApparentPower}_{\indexGridLines\indexGridNode \indexGridNodeTwo}^{\text{Kr}}  }}
\newcommand{\SjikSDPkron}[0]{\textcolor{\complexcolor}{\mathbf{\symApparentPower}_{\indexGridLines\indexGridNodeTwo \indexGridNode }^{\text{Kr}}  }}

\newcommand{\SijkSDPAdot}[0]{\textcolor{\complexcolor}{\mathbf{\symApparentPower}_{\indexGridLines\indexGridNode \indexGridNodeTwo}^{\symPhaseA \! \cdot}  }}
\newcommand{\SijkSDPBdot}[0]{\textcolor{\complexcolor}{\mathbf{\symApparentPower}_{\indexGridLines\indexGridNode \indexGridNodeTwo}^{\symPhaseB \! \cdot}  }}
\newcommand{\SijkSDPCdot}[0]{\textcolor{\complexcolor}{\mathbf{\symApparentPower}_{\indexGridLines\indexGridNode \indexGridNodeTwo}^{\symPhaseC \! \cdot}  }}
\newcommand{\SijkSDPNdot}[0]{\textcolor{\complexcolor}{\mathbf{\symApparentPower}_{\indexGridLines\indexGridNode \indexGridNodeTwo}^{\symPhaseN \! \cdot}  }}
\newcommand{\SijkSDPGdot}[0]{\textcolor{\complexcolor}{\mathbf{\symApparentPower}_{\indexGridLines\indexGridNode \indexGridNodeTwo}^{\symPhaseG \! \cdot}  }}

\newcommand{\SjikSDPAdot}[0]{\textcolor{\complexcolor}{\mathbf{\symApparentPower}_{\indexGridLines\indexGridNodeTwo\indexGridNode }^{\symPhaseA \! \cdot}  }}
\newcommand{\SjikSDPBdot}[0]{\textcolor{\complexcolor}{\mathbf{\symApparentPower}_{\indexGridLines\indexGridNodeTwo\indexGridNode }^{\symPhaseB \! \cdot}  }}
\newcommand{\SjikSDPCdot}[0]{\textcolor{\complexcolor}{\mathbf{\symApparentPower}_{\indexGridLines\indexGridNodeTwo\indexGridNode }^{\symPhaseC \! \cdot}  }}
\newcommand{\SjikSDPNdot}[0]{\textcolor{\complexcolor}{\mathbf{\symApparentPower}_{\indexGridLines\indexGridNodeTwo\indexGridNode }^{\symPhaseN \! \cdot}  }}

\newcommand{\PijkSDPAdot}[0]{\textcolor{black}{\mathbf{\symPower}_{\indexGridLines\indexGridNode \indexGridNodeTwo}^{\symPhaseA \! \cdot}  }}
\newcommand{\PijkSDPBdot}[0]{\textcolor{black}{\mathbf{\symPower}_{\indexGridLines\indexGridNode \indexGridNodeTwo}^{\symPhaseB \! \cdot}  }}
\newcommand{\PijkSDPCdot}[0]{\textcolor{black}{\mathbf{\symPower}_{\indexGridLines\indexGridNode \indexGridNodeTwo}^{\symPhaseC \! \cdot}  }}
\newcommand{\PijkSDPNdot}[0]{\textcolor{black}{\mathbf{\symPower}_{\indexGridLines\indexGridNode \indexGridNodeTwo}^{\symPhaseN \! \cdot}  }}

\newcommand{\QijkSDPAdot}[0]{\textcolor{black}{\mathbf{\symReactivePower}_{\indexGridLines\indexGridNode \indexGridNodeTwo}^{\symPhaseA \! \cdot}  }}
\newcommand{\QijkSDPBdot}[0]{\textcolor{black}{\mathbf{\symReactivePower}_{\indexGridLines\indexGridNode \indexGridNodeTwo}^{\symPhaseB \! \cdot}  }}
\newcommand{\QijkSDPCdot}[0]{\textcolor{black}{\mathbf{\symReactivePower}_{\indexGridLines\indexGridNode \indexGridNodeTwo}^{\symPhaseC \! \cdot}  }}
\newcommand{\QijkSDPNdot}[0]{\textcolor{black}{\mathbf{\symReactivePower}_{\indexGridLines\indexGridNode \indexGridNodeTwo}^{\symPhaseN \! \cdot}  }}

\newcommand{\PijksSDPAdot}[0]{\textcolor{black}{\mathbf{\symPower}_{\indexGridLines\indexGridNode \indexGridNodeTwo}^{\seriesss,\symPhaseA \! \cdot}  }}
\newcommand{\PijksSDPBdot}[0]{\textcolor{black}{\mathbf{\symPower}_{\indexGridLines\indexGridNode \indexGridNodeTwo}^{\seriesss,\symPhaseB \! \cdot}  }}
\newcommand{\PijksSDPCdot}[0]{\textcolor{black}{\mathbf{\symPower}_{\indexGridLines\indexGridNode \indexGridNodeTwo}^{\seriesss,\symPhaseC \! \cdot}  }}
\newcommand{\PijksSDPNdot}[0]{\textcolor{black}{\mathbf{\symPower}_{\indexGridLines\indexGridNode \indexGridNodeTwo}^{\seriesss,\symPhaseN \! \cdot}  }}

\newcommand{\QijksSDPAdot}[0]{\textcolor{black}{\mathbf{\symReactivePower}_{\indexGridLines\indexGridNode \indexGridNodeTwo}^{\seriesss,\symPhaseA \! \cdot}  }}
\newcommand{\QijksSDPBdot}[0]{\textcolor{black}{\mathbf{\symReactivePower}_{\indexGridLines\indexGridNode \indexGridNodeTwo}^{\seriesss,\symPhaseB \! \cdot}  }}
\newcommand{\QijksSDPCdot}[0]{\textcolor{black}{\mathbf{\symReactivePower}_{\indexGridLines\indexGridNode \indexGridNodeTwo}^{\seriesss,\symPhaseC \! \cdot}  }}
\newcommand{\QijksSDPNdot}[0]{\textcolor{black}{\mathbf{\symReactivePower}_{\indexGridLines\indexGridNode \indexGridNodeTwo}^{\seriesss,\symPhaseN \! \cdot}  }}

\newcommand{\SijkSDPdotA}[0]{\textcolor{\complexcolor}{\mathbf{\symApparentPower}_{\indexGridLines\indexGridNode \indexGridNodeTwo}^{\cdot\!\symPhaseA}  }}
\newcommand{\SijkSDPdotB}[0]{\textcolor{\complexcolor}{\mathbf{\symApparentPower}_{\indexGridLines\indexGridNode \indexGridNodeTwo}^{\cdot\!\symPhaseB }  }}
\newcommand{\SijkSDPdotC}[0]{\textcolor{\complexcolor}{\mathbf{\symApparentPower}_{\indexGridLines\indexGridNode \indexGridNodeTwo}^{\cdot\!\symPhaseC }  }}
\newcommand{\SijkSDPdotN}[0]{\textcolor{\complexcolor}{\mathbf{\symApparentPower}_{\indexGridLines\indexGridNode \indexGridNodeTwo}^{\cdot \!\symPhaseN }  }}

\newcommand{\PijkSDPdotA}[0]{\textcolor{black}{\mathbf{\symPower}_{\indexGridLines\indexGridNode \indexGridNodeTwo}^{\cdot\!\symPhaseA}  }}
\newcommand{\PijkSDPdotB}[0]{\textcolor{black}{\mathbf{\symPower}_{\indexGridLines\indexGridNode \indexGridNodeTwo}^{\cdot\!\symPhaseB }  }}
\newcommand{\PijkSDPdotC}[0]{\textcolor{black}{\mathbf{\symPower}_{\indexGridLines\indexGridNode \indexGridNodeTwo}^{\cdot\!\symPhaseC }  }}
\newcommand{\PijkSDPdotN}[0]{\textcolor{black}{\mathbf{\symPower}_{\indexGridLines\indexGridNode \indexGridNodeTwo}^{\cdot \!\symPhaseN }  }}

\newcommand{\QijkSDPdotA}[0]{\textcolor{black}{\mathbf{\symReactivePower}_{\indexGridLines\indexGridNode \indexGridNodeTwo}^{\cdot\!\symPhaseA}  }}
\newcommand{\QijkSDPdotB}[0]{\textcolor{black}{\mathbf{\symReactivePower}_{\indexGridLines\indexGridNode \indexGridNodeTwo}^{\cdot\!\symPhaseB }  }}
\newcommand{\QijkSDPdotC}[0]{\textcolor{black}{\mathbf{\symReactivePower}_{\indexGridLines\indexGridNode \indexGridNodeTwo}^{\cdot\!\symPhaseC }  }}
\newcommand{\QijkSDPdotN}[0]{\textcolor{black}{\mathbf{\symReactivePower}_{\indexGridLines\indexGridNode \indexGridNodeTwo}^{\cdot \!\symPhaseN }  }}

\newcommand{\PijksSDPdotA}[0]{\textcolor{black}{\mathbf{\symPower}_{\indexGridLines\indexGridNode \indexGridNodeTwo}^{\seriesss,\cdot\!\symPhaseA}  }}
\newcommand{\PijksSDPdotB}[0]{\textcolor{black}{\mathbf{\symPower}_{\indexGridLines\indexGridNode \indexGridNodeTwo}^{\seriesss,\cdot\!\symPhaseB }  }}
\newcommand{\PijksSDPdotC}[0]{\textcolor{black}{\mathbf{\symPower}_{\indexGridLines\indexGridNode \indexGridNodeTwo}^{\seriesss,\cdot\!\symPhaseC }  }}
\newcommand{\PijksSDPdotN}[0]{\textcolor{black}{\mathbf{\symPower}_{\indexGridLines\indexGridNode \indexGridNodeTwo}^{\seriesss,\cdot \!\symPhaseN }  }}

\newcommand{\QijksSDPdotA}[0]{\textcolor{black}{\mathbf{\symReactivePower}_{\indexGridLines\indexGridNode \indexGridNodeTwo}^{\seriesss,\cdot\!\symPhaseA}  }}
\newcommand{\QijksSDPdotB}[0]{\textcolor{black}{\mathbf{\symReactivePower}_{\indexGridLines\indexGridNode \indexGridNodeTwo}^{\seriesss,\cdot\!\symPhaseB }  }}
\newcommand{\QijksSDPdotC}[0]{\textcolor{black}{\mathbf{\symReactivePower}_{\indexGridLines\indexGridNode \indexGridNodeTwo}^{\seriesss,\cdot\!\symPhaseC }  }}
\newcommand{\QijksSDPdotN}[0]{\textcolor{black}{\mathbf{\symReactivePower}_{\indexGridLines\indexGridNode \indexGridNodeTwo}^{\seriesss,\cdot \!\symPhaseN }  }}

\newcommand{\SijskSDPAdot}[0]{\textcolor{\complexcolor}{\mathbf{\symApparentPower}_{\indexGridLines\indexGridNode \indexGridNodeTwo}^{\seriesss, \symPhaseA \! \cdot}  }}
\newcommand{\SijskSDPBdot}[0]{\textcolor{\complexcolor}{\mathbf{\symApparentPower}_{\indexGridLines\indexGridNode \indexGridNodeTwo}^{\seriesss, \symPhaseB \! \cdot}  }}
\newcommand{\SijskSDPCdot}[0]{\textcolor{\complexcolor}{\mathbf{\symApparentPower}_{\indexGridLines\indexGridNode \indexGridNodeTwo}^{\seriesss, \symPhaseC \! \cdot}  }}
\newcommand{\SijskSDPNdot}[0]{\textcolor{\complexcolor}{\mathbf{\symApparentPower}_{\indexGridLines\indexGridNode \indexGridNodeTwo}^{\seriesss, \symPhaseN \! \cdot}  }}
\newcommand{\SijskSDPGdot}[0]{\textcolor{\complexcolor}{\mathbf{\symApparentPower}_{\indexGridLines\indexGridNode \indexGridNodeTwo}^{\seriesss, \symPhaseG \! \cdot}  }}

\newcommand{\SijskSDPdotA}[0]{\textcolor{\complexcolor}{\mathbf{\symApparentPower}_{\indexGridLines\indexGridNode \indexGridNodeTwo}^{\seriesss, \cdot\!\symPhaseA}  }}
\newcommand{\SijskSDPdotB}[0]{\textcolor{\complexcolor}{\mathbf{\symApparentPower}_{\indexGridLines\indexGridNode \indexGridNodeTwo}^{\seriesss, \cdot\!\symPhaseB }  }}
\newcommand{\SijskSDPdotC}[0]{\textcolor{\complexcolor}{\mathbf{\symApparentPower}_{\indexGridLines\indexGridNode \indexGridNodeTwo}^{\seriesss, \cdot\!\symPhaseC }  }}
\newcommand{\SijskSDPdotN}[0]{\textcolor{\complexcolor}{\mathbf{\symApparentPower}_{\indexGridLines\indexGridNode \indexGridNodeTwo}^{\seriesss, \cdot \!\symPhaseN }  }}

\newcommand{\SijksSDPmax}[0]{\textcolor{\paramcolor}{\mathbf{\symApparentPower}_{\indexGridLines \indexGridNode \indexGridNodeTwo}^{\seriesss, \maxss}  }}

\newcommand{\SijksSDP}[0]{\textcolor{\complexcolor}{\mathbf{\symApparentPower}_{\indexGridLines \indexGridNode \indexGridNodeTwo}^{\seriesss}  }}
\newcommand{\SjiksSDP}[0]{\textcolor{\complexcolor}{\mathbf{\symApparentPower}_{\indexGridLines \indexGridNodeTwo\indexGridNode }^{\seriesss}  }}

\newcommand{\SijksSDPH}[0]{\textcolor{\complexcolor}{\mathbf{(\symApparentPower}_{\indexGridLines \indexGridNode \indexGridNodeTwo}^{\seriesss})^{\hermitiantranspose}  }}

\newcommand{\SijksSDPdiag}[0]{\textcolor{\complexcolor}{\mathbf{\symApparentPower}_{\indexGridLines\indexGridNode \indexGridNodeTwo}^{\seriesss, \text{diag}}  }}
\newcommand{\PijksSDPdiag}[0]{\textcolor{black}{\mathbf{\symPower}_{\indexGridLines\indexGridNode \indexGridNodeTwo}^{\seriesss, \text{diag}}  }}
\newcommand{\QijksSDPdiag}[0]{\textcolor{black}{\mathbf{\symReactivePower}_{\indexGridLines\indexGridNode \indexGridNodeTwo}^{\seriesss, \text{diag}}  }}

\newcommand{\PijkSDP}[0]{\textcolor{black}{\mathbf{\symPower}_{\indexGridLines\indexGridNode \indexGridNodeTwo}^{}  }}
\newcommand{\QijkSDP}[0]{\textcolor{black}{\mathbf{\symReactivePower}_{\indexGridLines\indexGridNode \indexGridNodeTwo}^{}  }}
\newcommand{\PijkSDPH}[0]{\textcolor{black}{\mathbf{\symPower}_{\indexGridLines\indexGridNode \indexGridNodeTwo}^{\hermitiantranspose}  }}
\newcommand{\QijkSDPH}[0]{\textcolor{black}{\mathbf{\symReactivePower}_{\indexGridLines\indexGridNode \indexGridNodeTwo}^{\hermitiantranspose}  }}
\newcommand{\PjikSDP}[0]{\textcolor{black}{\mathbf{\symPower}_{\indexGridLines\indexGridNodeTwo\indexGridNode }^{}  }}
\newcommand{\QjikSDP}[0]{\textcolor{black}{\mathbf{\symReactivePower}_{\indexGridLines\indexGridNodeTwo\indexGridNode }^{}  }}

\newcommand{\PijksSDPH}[0]{\textcolor{black}{\mathbf{\symPower}_{\indexGridLines\indexGridNode \indexGridNodeTwo, \seriesss}^{\hermitiantranspose}  }}
\newcommand{\QijksSDPH}[0]{\textcolor{black}{\mathbf{\symReactivePower}_{\indexGridLines\indexGridNode \indexGridNodeTwo, \seriesss}^{\hermitiantranspose}  }}

\newcommand{\PijksSDP}[0]{\textcolor{black}{\mathbf{\symPower}_{\indexGridLines\indexGridNode \indexGridNodeTwo}^{\seriesss}  }}
\newcommand{\QijksSDP}[0]{\textcolor{black}{\mathbf{\symReactivePower}_{\indexGridLines\indexGridNode \indexGridNodeTwo}^{\seriesss}  }}
\newcommand{\PjiksSDP}[0]{\textcolor{black}{\mathbf{\symPower}_{\indexGridLines\indexGridNodeTwo\indexGridNode }^{\seriesss}  }}
\newcommand{\QjiksSDP}[0]{\textcolor{black}{\mathbf{\symReactivePower}_{\indexGridLines\indexGridNodeTwo\indexGridNode }^{\seriesss}  }}

\newcommand{\SijklossSDP}[0]{\textcolor{\complexcolor}{\mathbf{\symApparentPower}_{\indexGridNode \indexGridNodeTwo}^{\lossss}  }}
\newcommand{\SjiklossSDP}[0]{\textcolor{\complexcolor}{\mathbf{\symApparentPower}_{\indexGridNodeTwo\indexGridNode }^{\lossss}  }}
\newcommand{\SlklossSDP}[0]{\textcolor{\complexcolor}{\mathbf{\symApparentPower}_{\indexGridLines}^{\lossss}  }}
\newcommand{\PijklossSDP}[0]{\textcolor{black}{\mathbf{\symPower}_{\indexGridNode \indexGridNodeTwo}^{\lossss}  }}
\newcommand{\QijklossSDP}[0]{\textcolor{black}{\mathbf{\symReactivePower}_{\indexGridNode \indexGridNodeTwo}^{\lossss}  }}

\newcommand{\PijkslossSDP}[0]{\textcolor{black}{\mathbf{\symPower}_{\indexGridNode \indexGridNodeTwo, \seriesss}^{\lossss}  }}
\newcommand{\QijkslossSDP}[0]{\textcolor{black}{\mathbf{\symReactivePower}_{\indexGridNode \indexGridNodeTwo, \seriesss}^{\lossss}  }}

\newcommand{\SlklossSDPkron}[0]{\textcolor{\complexcolor}{\mathbf{\symApparentPower}_{\indexGridLines}^{\lossss, \text{Kr}}  }}

\newcommand{\SijkshlossSDP}[0]{\textcolor{\complexcolor}{\mathbf{\symApparentPower}_{\indexGridLines\indexGridNode \indexGridNodeTwo}^{\lossss, \shuntss}  }}
\newcommand{\SjikshlossSDP}[0]{\textcolor{\complexcolor}{\mathbf{\symApparentPower}_{\indexGridLines\indexGridNodeTwo\indexGridNode}^{\lossss, \shuntss}  }}
\newcommand{\SlkshlossSDP}[0]{\textcolor{\complexcolor}{\mathbf{\symApparentPower}_{\indexGridLines}^{\lossss, \shuntss}  }}

\newcommand{\SijkslossSDP}[0]{\textcolor{\complexcolor}{\mathbf{\symApparentPower}_{\indexGridLines }^{\lossss,\seriesss}  }}
\newcommand{\SjikslossSDP}[0]{\textcolor{\complexcolor}{\mathbf{\symApparentPower}_{\indexGridLines}^{\lossss,\seriesss}  }}
\newcommand{\SlkslossSDP}[0]{\textcolor{\complexcolor}{\mathbf{\symApparentPower}_{\indexGridLines}^{\lossss,\seriesss}  }}

\newcommand{\Sinode}[0]{\textcolor{black}{\symApparentPower_{\indexGridNode }^{}  }}
\newcommand{\Sjnode}[0]{\textcolor{black}{\symApparentPower_{ \indexGridNodeTwo}^{}  }}
\newcommand{\SjnodeSDP}[0]{\textcolor{\complexcolor}{\mathbf{\symApparentPower}_{ \indexGridNodeTwo}^{}  }}
\newcommand{\QjnodeSDP}[0]{\textcolor{black}{\mathbf{\symReactivePower}_{ \indexGridNodeTwo}^{}  }}
\newcommand{\PjnodeSDP}[0]{\textcolor{black}{\mathbf{\symPower}_{ \indexGridNodeTwo}^{}  }}

\newcommand{\SinodeSDP}[0]{\textcolor{\complexcolor}{\mathbf{\symApparentPower}_{ \indexGridNode}^{}  }}
\newcommand{\QinodeSDP}[0]{\textcolor{black}{\mathbf{\symReactivePower}_{ \indexGridNode}^{}  }}
\newcommand{\PinodeSDP}[0]{\textcolor{black}{\mathbf{\symPower}_{ \indexGridNode}^{}  }}

\newcommand{\IinodeSDP}[0]{\textcolor{\complexcolor}{\mathbf{\symCurrent}_{ \indexGridNode}^{}  }}

\newcommand{\IuunitSDPrated}[0]{\textcolor{\paramcolor}{\mathbf{\symCurrent}_{ \indexUnit}^{\ratedss}  }}
\newcommand{\IuunitSDPdeltarated}[0]{\textcolor{\paramcolor}{\mathbf{\symCurrent}_{ \indexUnit}^{\Delta\ratedss}  }}

\newcommand{\IuunitSDPref}[0]{\textcolor{\complexparamcolor}{\mathbf{\symCurrent}_{ \indexUnit}^{\refss}  }}

\newcommand{\IgenSDP}[0]{\textcolor{\complexcolor}{\mathbf{\symCurrent}_{ g}  }}
\newcommand{\IloadSDP}[0]{\textcolor{\complexcolor}{\mathbf{\symCurrent}_{ d} }}
\newcommand{\IgroundSDP}[0]{\textcolor{\complexcolor}{\mathbf{\symCurrent}_{i}^{\text{g}} }}

\newcommand{\Iloada}[0]{\textcolor{\complexcolor}{{\symCurrent}_{ d,a} }}
\newcommand{\Iloadb}[0]{\textcolor{\complexcolor}{{\symCurrent}_{ d,b} }}
\newcommand{\Iloadc}[0]{\textcolor{\complexcolor}{{\symCurrent}_{ d,c} }}

\newcommand{\Sloada}[0]{\textcolor{\complexcolor}{{\symApparentPower}_{ d,a} }}
\newcommand{\Sloadb}[0]{\textcolor{\complexcolor}{{\symApparentPower}_{ d,b} }}
\newcommand{\Sloadc}[0]{\textcolor{\complexcolor}{{\symApparentPower}_{ d,c} }}

\newcommand{\Sloadaref}[0]{\textcolor{\complexparamcolor}{{\symApparentPower}_{ d,a}^{\refss} }}
\newcommand{\Sloadbref}[0]{\textcolor{\complexparamcolor}{{\symApparentPower}_{ d,b}^{\refss} }}
\newcommand{\Sloadcref}[0]{\textcolor{\complexparamcolor}{{\symApparentPower}_{ d,c}^{\refss} }}

\newcommand{\SloadarefMAT}[0]{\textcolor{\complexparamcolor}{\mathbf{\symApparentPower}_{ d}^{\refss}  }}

\newcommand{\IuunitSDP}[0]{\textcolor{\complexcolor}{\mathbf{\symCurrent}_{ \indexUnit}^{}  }}
\newcommand{\IuunitSDPreal}[0]{\textcolor{black}{\mathbf{\symCurrent}_{ \indexUnit}^{\text{re}}  }}
\newcommand{\IuunitSDPimag}[0]{\textcolor{black}{\mathbf{\symCurrent}_{ \indexUnit}^{\text{im}}  }}

\newcommand{\IuunitSDPdelta}[0]{\textcolor{\complexcolor}{\mathbf{\symCurrent}_{ \indexUnit}^{\Delta}  }}
\newcommand{\IuunitSDPdeltareal}[0]{\textcolor{black}{\mathbf{\symCurrent}_{ \indexUnit}^{\Delta, \realss}  }}
\newcommand{\IuunitSDPdeltaimag}[0]{\textcolor{black}{\mathbf{\symCurrent}_{ \indexUnit}^{\Delta, \imagss}  }}

\newcommand{\UuunitSDP}[0]{\textcolor{\complexcolor}{\mathbf{\symVoltage}_{ \indexUnit}^{}  }}

\newcommand{\IbSDP}[0]{\textcolor{\complexcolor}{\mathbf{\symCurrent}_{ \indexShunt}^{}  }}
\newcommand{\IbSDPreal}[0]{\textcolor{black}{\mathbf{\symCurrent}_{ \indexShunt}^{\text{re}}  }}
\newcommand{\IbSDPimag}[0]{\textcolor{black}{\mathbf{\symCurrent}_{ \indexShunt}^{\text{im}}  }}

\newcommand{\SbSDP}[0]{\textcolor{\complexcolor}{\mathbf{\symApparentPower}_{ \indexShunt} }}
\newcommand{\PbSDP}[0]{\textcolor{black}{\mathbf{\symPower}_{ \indexShunt} }}
\newcommand{\QbSDP}[0]{\textcolor{black}{\mathbf{\symReactivePower}_{ \indexShunt} }}

\newcommand{\SuunitSDPint}[0]{\textcolor{\complexcolor}{\mathbf{\symApparentPower}_{ \indexUnit}^{'}  }}

\newcommand{\SuunitSDPref}[0]{\textcolor{\complexparamcolor}{\mathbf{\symApparentPower}_{ \indexUnit}^{\refss}  }}
\newcommand{\SuunitSDPrefdelta}[0]{\textcolor{\complexparamcolor}{\mathbf{\symApparentPower}_{ \indexUnit}^{\Delta\refss}  }}
\newcommand{\PuunitSDPref}[0]{\textcolor{\paramcolor}{\mathbf{\symPower}_{ \indexUnit}^{\refss}  }}
\newcommand{\QuunitSDPref}[0]{\textcolor{\paramcolor}{\mathbf{\symReactivePower}_{ \indexUnit}^{\refss}  }}

\newcommand{\PuunitSDPmin}[0]{\textcolor{\paramcolor}{\mathbf{\symPower}_{ \indexUnit}^{\minss}  }}
\newcommand{\QuunitSDPmin}[0]{\textcolor{\paramcolor}{\mathbf{\symReactivePower}_{ \indexUnit}^{\minss}  }}
\newcommand{\PuunitSDPmax}[0]{\textcolor{\paramcolor}{\mathbf{\symPower}_{ \indexUnit}^{\maxss}  }}
\newcommand{\QuunitSDPmax}[0]{\textcolor{\paramcolor}{\mathbf{\symReactivePower}_{ \indexUnit}^{\maxss}  }}
\newcommand{\SuunitSDPmax}[0]{\textcolor{\paramcolor}{\mathbf{\symApparentPower}_{ \indexUnit}^{\maxss}  }}

\newcommand{\PuunitSDPmindelta}[0]{\textcolor{\paramcolor}{\mathbf{\symPower}_{ \indexUnit}^{\Delta\minss}  }}
\newcommand{\QuunitSDPmindelta}[0]{\textcolor{\paramcolor}{\mathbf{\symReactivePower}_{ \indexUnit}^{\Delta\minss}  }}
\newcommand{\PuunitSDPmaxdelta}[0]{\textcolor{\paramcolor}{\mathbf{\symPower}_{ \indexUnit}^{\Delta\maxss}  }}
\newcommand{\QuunitSDPmaxdelta}[0]{\textcolor{\paramcolor}{\mathbf{\symReactivePower}_{ \indexUnit}^{\Delta\maxss}  }}

\newcommand{\SuunitSDP}[0]{\textcolor{\complexcolor}{\mathbf{\symApparentPower}_{ \indexUnit}^{}  }}
\newcommand{\PuunitSDP}[0]{\textcolor{black}{\mathbf{\symPower}_{ \indexUnit}^{}  }}
\newcommand{\QuunitSDP}[0]{\textcolor{black}{\mathbf{\symReactivePower}_{ \indexUnit}^{}  }}

\newcommand{\SuunitSDPdelta}[0]{\textcolor{\complexcolor}{\mathbf{\symApparentPower}_{ \indexUnit}^{\Delta}  }}
\newcommand{\PuunitSDPdelta}[0]{\textcolor{black}{\mathbf{\symPower}_{ \indexUnit}^{\Delta}  }}
\newcommand{\QuunitSDPdelta}[0]{\textcolor{black}{\mathbf{\symReactivePower}_{ \indexUnit}^{\Delta}  }}

\newcommand{\Tuunitdelta}[0]{\textcolor{\paramcolor}{\mathbf{T}^{\Delta}  }}
\newcommand{\Nwye}[0]{\textcolor{\paramcolor}{\mathbf{N}^{\text{pn}}  }}

\newcommand{\XljiSDPdelta}[0]{\textcolor{\complexcolor}{\mathbf{X}_{ lji}  }}
\newcommand{\XlijSDPdelta}[0]{\textcolor{\complexcolor}{\mathbf{X}_{ lij}  }}

\newcommand{\XuunitSDPdelta}[0]{\textcolor{\complexcolor}{\mathbf{X}_{ \indexUnit}^{\Delta}  }}
\newcommand{\XuunitSDPdeltareal}[0]{\textcolor{black}{\mathbf{X}_{ \indexUnit}^{\Delta\realss}  }}
\newcommand{\XuunitSDPdeltaimag}[0]{\textcolor{black}{\mathbf{X}_{ \indexUnit}^{\Delta\imagss}  }}

\newcommand{\IuunitAN}[0]{\textcolor{\complexcolor}{{\symCurrent}_{ \indexUnit, \symPhaseA}^{}  }}
\newcommand{\IuunitBN}[0]{\textcolor{\complexcolor}{{\symCurrent}_{ \indexUnit, \symPhaseB}^{}  }}
\newcommand{\IuunitCN}[0]{\textcolor{\complexcolor}{{\symCurrent}_{ \indexUnit, \symPhaseC}^{}  }}
\newcommand{\IuunitN}[0]{\textcolor{\complexcolor}{{\symCurrent}_{ \indexUnit, \symPhaseN}^{}  }}
\newcommand{\IuunitG}[0]{\textcolor{\complexcolor}{{\symCurrent}_{ \indexUnit, \symPhaseG}^{}  }}

\newcommand{\IuunitArated}[0]{\textcolor{\paramcolor}{{\symCurrent}_{ \indexUnit, \symPhaseA}^{\ratedss}  }}
\newcommand{\IuunitBrated}[0]{\textcolor{\paramcolor}{{\symCurrent}_{ \indexUnit, \symPhaseB}^{\ratedss}  }}
\newcommand{\IuunitCrated}[0]{\textcolor{\paramcolor}{{\symCurrent}_{ \indexUnit, \symPhaseC}^{\ratedss}  }}
\newcommand{\IuunitNrated}[0]{\textcolor{\paramcolor}{{\symCurrent}_{ \indexUnit, \symPhaseN}^{\ratedss}  }}

\newcommand{\IuunitABrated}[0]{\textcolor{\paramcolor}{{\symCurrent}_{ \indexUnit, \symPhaseA\symPhaseB}^{\ratedss}  }}
\newcommand{\IuunitBCrated}[0]{\textcolor{\paramcolor}{{\symCurrent}_{ \indexUnit, \symPhaseB\symPhaseC}^{\ratedss}  }}
\newcommand{\IuunitCArated}[0]{\textcolor{\paramcolor}{{\symCurrent}_{ \indexUnit, \symPhaseC\symPhaseA}^{\ratedss}  }}

\newcommand{\IbA}[0]{\textcolor{\complexcolor}{{\symCurrent}_{ \indexShunt, \symPhaseA}^{}  }}
\newcommand{\IbB}[0]{\textcolor{\complexcolor}{{\symCurrent}_{ \indexShunt, \symPhaseB}^{}  }}
\newcommand{\IbC}[0]{\textcolor{\complexcolor}{{\symCurrent}_{ \indexShunt, \symPhaseC}^{}  }}
\newcommand{\IbN}[0]{\textcolor{\complexcolor}{{\symCurrent}_{ \indexShunt, \symPhaseN}^{}  }}

\newcommand{\IuunitNreal}[0]{\textcolor{black}{{\symCurrent}_{ \indexUnit, \symPhaseN}^{\text{re}}  }}
\newcommand{\IuunitNimag}[0]{\textcolor{black}{{\symCurrent}_{ \indexUnit, \symPhaseN}^{\text{im}}  }}

\newcommand{\IuunitAB}[0]{\textcolor{\complexcolor}{{\symCurrent}_{ \indexUnit, \symPhaseA\symPhaseB}^{\Delta}  }}
\newcommand{\IuunitBC}[0]{\textcolor{\complexcolor}{{\symCurrent}_{ \indexUnit, \symPhaseB\symPhaseC}^{\Delta}  }}
\newcommand{\IuunitCA}[0]{\textcolor{\complexcolor}{{\symCurrent}_{ \indexUnit, \symPhaseC\symPhaseA}^{\Delta}  }}

\newcommand{\IuunitBA}[0]{\textcolor{\complexcolor}{{\symCurrent}_{ \indexUnit,\symPhaseB \symPhaseA}^{\Delta}  }}
\newcommand{\IuunitCB}[0]{\textcolor{\complexcolor}{{\symCurrent}_{ \indexUnit, \symPhaseC\symPhaseB}^{\Delta}  }}
\newcommand{\IuunitAC}[0]{\textcolor{\complexcolor}{{\symCurrent}_{ \indexUnit, \symPhaseA\symPhaseC}^{\Delta}  }}

\newcommand{\SuunitAAdelta}[0]{\textcolor{\complexcolor}{{\symApparentPower}_{ \indexUnit, \symPhaseA\symPhaseA}^{\Delta}  }}
\newcommand{\SuunitBBdelta}[0]{\textcolor{\complexcolor}{{\symApparentPower}_{ \indexUnit, \symPhaseB\symPhaseB}^{\Delta}  }}
\newcommand{\SuunitCCdelta}[0]{\textcolor{\complexcolor}{{\symApparentPower}_{ \indexUnit, \symPhaseC\symPhaseC}^{\Delta}  }}

\newcommand{\PuunitAAdelta}[0]{\textcolor{black}{{\symPower}_{ \indexUnit, \symPhaseA\symPhaseA}^{\Delta}  }}
\newcommand{\PuunitBBdelta}[0]{\textcolor{black}{{\symPower}_{ \indexUnit, \symPhaseB\symPhaseB}^{\Delta}  }}
\newcommand{\PuunitCCdelta}[0]{\textcolor{black}{{\symPower}_{ \indexUnit, \symPhaseC\symPhaseC}^{\Delta}  }}

\newcommand{\QuunitAAdelta}[0]{\textcolor{black}{{\symReactivePower}_{ \indexUnit, \symPhaseA\symPhaseA}^{\Delta}  }}
\newcommand{\QuunitBBdelta}[0]{\textcolor{black}{{\symReactivePower}_{ \indexUnit, \symPhaseB\symPhaseB}^{\Delta}  }}
\newcommand{\QuunitCCdelta}[0]{\textcolor{black}{{\symReactivePower}_{ \indexUnit, \symPhaseC\symPhaseC}^{\Delta}  }}

\newcommand{\SuunitABdelta}[0]{\textcolor{\complexcolor}{{\symApparentPower}_{ \indexUnit, \symPhaseA\symPhaseB}^{\Delta}  }}
\newcommand{\SuunitBCdelta}[0]{\textcolor{\complexcolor}{{\symApparentPower}_{ \indexUnit, \symPhaseB\symPhaseC}^{\Delta}  }}
\newcommand{\SuunitCAdelta}[0]{\textcolor{\complexcolor}{{\symApparentPower}_{ \indexUnit, \symPhaseC\symPhaseA}^{\Delta}  }}

\newcommand{\SuunitBAdelta}[0]{\textcolor{\complexcolor}{{\symApparentPower}_{ \indexUnit,\symPhaseB \symPhaseA}^{\Delta}  }}
\newcommand{\SuunitCBdelta}[0]{\textcolor{\complexcolor}{{\symApparentPower}_{ \indexUnit, \symPhaseC\symPhaseB}^{\Delta}  }}
\newcommand{\SuunitACdelta}[0]{\textcolor{\complexcolor}{{\symApparentPower}_{ \indexUnit, \symPhaseA\symPhaseC}^{\Delta}  }}

\newcommand{\SuunitABCdelta}[0]{\textcolor{\complexcolor}{{\symApparentPower}_{ \indexUnit, \symPhaseA\symPhaseB\symPhaseC}^{\Delta}  }}
\newcommand{\SuunitBCAdelta}[0]{\textcolor{\complexcolor}{{\symApparentPower}_{ \indexUnit, \symPhaseB\symPhaseC\symPhaseA}^{\Delta}  }}
\newcommand{\SuunitCABdelta}[0]{\textcolor{\complexcolor}{{\symApparentPower}_{ \indexUnit, \symPhaseC\symPhaseA\symPhaseB}^{\Delta}  }}

\newcommand{\PuunitABdeltamax}[0]{\textcolor{\paramcolor}{{\symPower}_{ \indexUnit, \symPhaseA\symPhaseB}^{\Delta, \maxss}  }}
\newcommand{\PuunitBCdeltamax}[0]{\textcolor{\paramcolor}{{\symPower}_{ \indexUnit, \symPhaseB\symPhaseC}^{\Delta, \maxss}  }}
\newcommand{\PuunitCAdeltamax}[0]{\textcolor{\paramcolor}{{\symPower}_{ \indexUnit, \symPhaseC\symPhaseA}^{\Delta, \maxss}  }}

\newcommand{\QuunitABdeltamax}[0]{\textcolor{\paramcolor}{{\symReactivePower}_{ \indexUnit, \symPhaseA\symPhaseB}^{\Delta, \maxss}  }}
\newcommand{\QuunitBCdeltamax}[0]{\textcolor{\paramcolor}{{\symReactivePower}_{ \indexUnit, \symPhaseB\symPhaseC}^{\Delta, \maxss}  }}
\newcommand{\QuunitCAdeltamax}[0]{\textcolor{\paramcolor}{{\symReactivePower}_{ \indexUnit, \symPhaseC\symPhaseA}^{\Delta, \maxss}  }}

\newcommand{\PuunitABdeltamin}[0]{\textcolor{\paramcolor}{{\symPower}_{ \indexUnit, \symPhaseA\symPhaseB}^{\Delta, \minss}  }}
\newcommand{\PuunitBCdeltamin}[0]{\textcolor{\paramcolor}{{\symPower}_{ \indexUnit, \symPhaseB\symPhaseC}^{\Delta, \minss}  }}
\newcommand{\PuunitCAdeltamin}[0]{\textcolor{\paramcolor}{{\symPower}_{ \indexUnit, \symPhaseC\symPhaseA}^{\Delta, \minss}  }}

\newcommand{\QuunitABdeltamin}[0]{\textcolor{\paramcolor}{{\symReactivePower}_{ \indexUnit, \symPhaseA\symPhaseB}^{\Delta, \minss}  }}
\newcommand{\QuunitBCdeltamin}[0]{\textcolor{\paramcolor}{{\symReactivePower}_{ \indexUnit, \symPhaseB\symPhaseC}^{\Delta, \minss}  }}
\newcommand{\QuunitCAdeltamin}[0]{\textcolor{\paramcolor}{{\symReactivePower}_{ \indexUnit, \symPhaseC\symPhaseA}^{\Delta, \minss}  }}

\newcommand{\PuunitAAmax}[0]{\textcolor{\paramcolor}{{\symPower}_{ \indexUnit, \symPhaseA\symPhaseA}^{\maxss}  }}
\newcommand{\PuunitBBmax}[0]{\textcolor{\paramcolor}{{\symPower}_{ \indexUnit, \symPhaseB\symPhaseB}^{ \maxss}  }}
\newcommand{\PuunitCCmax}[0]{\textcolor{\paramcolor}{{\symPower}_{ \indexUnit, \symPhaseC\symPhaseC}^{ \maxss}  }}

\newcommand{\QuunitAAmax}[0]{\textcolor{\paramcolor}{{\symReactivePower}_{ \indexUnit, \symPhaseA\symPhaseA}^{\maxss}  }}
\newcommand{\QuunitBBmax}[0]{\textcolor{\paramcolor}{{\symReactivePower}_{ \indexUnit, \symPhaseB\symPhaseB}^{ \maxss}  }}
\newcommand{\QuunitCCmax}[0]{\textcolor{\paramcolor}{{\symReactivePower}_{ \indexUnit, \symPhaseC\symPhaseC}^{ \maxss}  }}

\newcommand{\PuunitAAmin}[0]{\textcolor{\paramcolor}{{\symPower}_{ \indexUnit, \symPhaseA\symPhaseA}^{\minss}  }}
\newcommand{\PuunitBBmin}[0]{\textcolor{\paramcolor}{{\symPower}_{ \indexUnit, \symPhaseB\symPhaseB}^{ \minss}  }}
\newcommand{\PuunitCCmin}[0]{\textcolor{\paramcolor}{{\symPower}_{ \indexUnit, \symPhaseC\symPhaseC}^{ \minss}  }}

\newcommand{\QuunitAAmin}[0]{\textcolor{\paramcolor}{{\symReactivePower}_{ \indexUnit, \symPhaseA\symPhaseA}^{\minss}  }}
\newcommand{\QuunitBBmin}[0]{\textcolor{\paramcolor}{{\symReactivePower}_{ \indexUnit, \symPhaseB\symPhaseB}^{ \minss}  }}
\newcommand{\QuunitCCmin}[0]{\textcolor{\paramcolor}{{\symReactivePower}_{ \indexUnit, \symPhaseC\symPhaseC}^{ \minss}  }}

\newcommand{\PuunitABdelta}[0]{\textcolor{black}{{\symPower}_{ \indexUnit, \symPhaseA\symPhaseB}^{\Delta}  }}
\newcommand{\PuunitBCdelta}[0]{\textcolor{black}{{\symPower}_{ \indexUnit, \symPhaseB\symPhaseC}^{\Delta}  }}
\newcommand{\PuunitCAdelta}[0]{\textcolor{black}{{\symPower}_{ \indexUnit, \symPhaseC\symPhaseA}^{\Delta}  }}

\newcommand{\PuunitBAdelta}[0]{\textcolor{black}{{\symPower}_{ \indexUnit,\symPhaseB \symPhaseA}^{\Delta}  }}
\newcommand{\PuunitCBdelta}[0]{\textcolor{black}{{\symPower}_{ \indexUnit, \symPhaseC\symPhaseB}^{\Delta}  }}
\newcommand{\PuunitACdelta}[0]{\textcolor{black}{{\symPower}_{ \indexUnit, \symPhaseA\symPhaseC}^{\Delta}  }}

\newcommand{\QuunitABdelta}[0]{\textcolor{black}{{\symReactivePower}_{ \indexUnit, \symPhaseA\symPhaseB}^{\Delta}  }}
\newcommand{\QuunitBCdelta}[0]{\textcolor{black}{{\symReactivePower}_{ \indexUnit, \symPhaseB\symPhaseC}^{\Delta}  }}
\newcommand{\QuunitCAdelta}[0]{\textcolor{black}{{\symReactivePower}_{ \indexUnit, \symPhaseC\symPhaseA}^{\Delta}  }}

\newcommand{\QuunitBAdelta}[0]{\textcolor{black}{{\symReactivePower}_{ \indexUnit,\symPhaseB \symPhaseA}^{\Delta}  }}
\newcommand{\QuunitCBdelta}[0]{\textcolor{black}{{\symReactivePower}_{ \indexUnit, \symPhaseC\symPhaseB}^{\Delta}  }}
\newcommand{\QuunitACdelta}[0]{\textcolor{black}{{\symReactivePower}_{ \indexUnit, \symPhaseA\symPhaseC}^{\Delta}  }}

\newcommand{\PuunitAN}[0]{\textcolor{black}{{\symPower}_{ \indexUnit, \symPhaseA}^{}  }}
\newcommand{\PuunitBN}[0]{\textcolor{black}{{\symPower}_{ \indexUnit, \symPhaseB}^{}  }}
\newcommand{\PuunitCN}[0]{\textcolor{black}{{\symPower}_{ \indexUnit, \symPhaseC}^{}  }}

\newcommand{\SuunitrefAA}[0]{\textcolor{\complexparamcolor}{{\symApparentPower}_{ \indexUnit, \symPhaseA\symPhaseA}^{\refss}  }}
\newcommand{\SuunitrefBB}[0]{\textcolor{\complexparamcolor}{{\symApparentPower}_{ \indexUnit, \symPhaseB\symPhaseB}^{\refss}  }}
\newcommand{\SuunitrefCC}[0]{\textcolor{\complexparamcolor}{{\symApparentPower}_{ \indexUnit, \symPhaseC\symPhaseC}^{\refss}  }}

\newcommand{\SuunitrefAB}[0]{\textcolor{\complexparamcolor}{{\symApparentPower}_{ \indexUnit, \symPhaseA\symPhaseB}^{\refss}  }}
\newcommand{\SuunitrefBC}[0]{\textcolor{\complexparamcolor}{{\symApparentPower}_{ \indexUnit, \symPhaseB\symPhaseC}^{\refss}  }}
\newcommand{\SuunitrefCA}[0]{\textcolor{\complexparamcolor}{{\symApparentPower}_{ \indexUnit, \symPhaseC\symPhaseA}^{\refss}  }}

\newcommand{\PuunitrefAB}[0]{\textcolor{\paramcolor}{{\symPower}_{ \indexUnit, \symPhaseA\symPhaseB}^{\refss}  }}
\newcommand{\PuunitrefBC}[0]{\textcolor{\paramcolor}{{\symPower}_{ \indexUnit, \symPhaseB\symPhaseC}^{\refss}  }}
\newcommand{\PuunitrefCA}[0]{\textcolor{\paramcolor}{{\symPower}_{ \indexUnit, \symPhaseC\symPhaseA}^{\refss}  }}

\newcommand{\QuunitrefAB}[0]{\textcolor{\paramcolor}{{\symReactivePower}_{ \indexUnit, \symPhaseA\symPhaseB}^{\refss}  }}
\newcommand{\QuunitrefBC}[0]{\textcolor{\paramcolor}{{\symReactivePower}_{ \indexUnit, \symPhaseB\symPhaseC}^{\refss}  }}
\newcommand{\QuunitrefCA}[0]{\textcolor{\paramcolor}{{\symReactivePower}_{ \indexUnit, \symPhaseC\symPhaseA}^{\refss}  }}

\newcommand{\PuunitrefAA}[0]{\textcolor{\paramcolor}{{\symPower}_{ \indexUnit, \symPhaseA\symPhaseA}^{\refss}  }}
\newcommand{\PuunitrefBB}[0]{\textcolor{\paramcolor}{{\symPower}_{ \indexUnit, \symPhaseB\symPhaseB}^{\refss}  }}
\newcommand{\PuunitrefCC}[0]{\textcolor{\paramcolor}{{\symPower}_{ \indexUnit, \symPhaseC\symPhaseC}^{\refss}  }}

\newcommand{\QuunitrefAA}[0]{\textcolor{\paramcolor}{{\symReactivePower}_{ \indexUnit, \symPhaseA\symPhaseA}^{\refss}  }}
\newcommand{\QuunitrefBB}[0]{\textcolor{\paramcolor}{{\symReactivePower}_{ \indexUnit, \symPhaseB\symPhaseB}^{\refss}  }}
\newcommand{\QuunitrefCC}[0]{\textcolor{\paramcolor}{{\symReactivePower}_{ \indexUnit, \symPhaseC\symPhaseC}^{\refss}  }}

\newcommand{\SuunitAdot}[0]{\textcolor{\complexcolor}{\mathbf{\symApparentPower}_{ \indexUnit}^{\symPhaseA\!\cdot}  }}
\newcommand{\SuunitBdot}[0]{\textcolor{\complexcolor}{\mathbf{\symApparentPower}_{ \indexUnit}^{\symPhaseB\!\cdot}  }}
\newcommand{\SuunitCdot}[0]{\textcolor{\complexcolor}{\mathbf{\symApparentPower}_{ \indexUnit}^{\symPhaseC\!\cdot}  }}
\newcommand{\SuunitNdot}[0]{\textcolor{\complexcolor}{\mathbf{\symApparentPower}_{ \indexUnit}^{\symPhaseN\!\cdot}  }}
\newcommand{\SuunitGdot}[0]{\textcolor{\complexcolor}{\mathbf{\symApparentPower}_{ \indexUnit}^{\symPhaseG\!\cdot}  }}

\newcommand{\SuunitdotN}[0]{\textcolor{\complexcolor}{\mathbf{\symApparentPower}_{ \indexUnit}^{\cdot\! \symPhaseN}  }}

\newcommand{\PuunitNdot}[0]{\textcolor{black}{\mathbf{\symPower}_{ \indexUnit}^{\symPhaseN\!\cdot}  }}
\newcommand{\QuunitNdot}[0]{\textcolor{black}{\mathbf{\symReactivePower}_{ \indexUnit}^{\symPhaseN\!\cdot}  }}

\newcommand{\XuunitAA}[0]{\textcolor{\complexcolor}{{X}_{ \indexUnit, \symPhaseA\symPhaseA}^{}  }}
\newcommand{\XuunitAB}[0]{\textcolor{\complexcolor}{{X}_{ \indexUnit, \symPhaseA\symPhaseB}^{}  }}
\newcommand{\XuunitAC}[0]{\textcolor{\complexcolor}{{X}_{ \indexUnit, \symPhaseA\symPhaseC}^{}  }}

\newcommand{\XuunitBA}[0]{\textcolor{\complexcolor}{{X}_{ \indexUnit, \symPhaseB\symPhaseA}^{}  }}
\newcommand{\XuunitBB}[0]{\textcolor{\complexcolor}{{X}_{ \indexUnit, \symPhaseB\symPhaseB}^{}  }}
\newcommand{\XuunitBC}[0]{\textcolor{\complexcolor}{{X}_{ \indexUnit, \symPhaseB\symPhaseC}^{}  }}

\newcommand{\XuunitCA}[0]{\textcolor{\complexcolor}{{X}_{ \indexUnit, \symPhaseC\symPhaseA}^{}  }}
\newcommand{\XuunitCB}[0]{\textcolor{\complexcolor}{{X}_{ \indexUnit, \symPhaseC\symPhaseB}^{}  }}
\newcommand{\XuunitCC}[0]{\textcolor{\complexcolor}{{X}_{ \indexUnit, \symPhaseC\symPhaseC}^{}  }}

\newcommand{\XuunitAAreal}[0]{\textcolor{black}{{X}_{ \indexUnit, \symPhaseA\symPhaseA}^{\realss}  }}
\newcommand{\XuunitABreal}[0]{\textcolor{black}{{X}_{ \indexUnit, \symPhaseA\symPhaseB}^{\realss}  }}
\newcommand{\XuunitACreal}[0]{\textcolor{black}{{X}_{ \indexUnit, \symPhaseA\symPhaseC}^{\realss}  }}

\newcommand{\XuunitBAreal}[0]{\textcolor{black}{{X}_{ \indexUnit, \symPhaseB\symPhaseA}^{\realss}  }}
\newcommand{\XuunitBBreal}[0]{\textcolor{black}{{X}_{ \indexUnit, \symPhaseB\symPhaseB}^{\realss}  }}
\newcommand{\XuunitBCreal}[0]{\textcolor{black}{{X}_{ \indexUnit, \symPhaseB\symPhaseC}^{\realss}  }}

\newcommand{\XuunitCAreal}[0]{\textcolor{black}{{X}_{ \indexUnit, \symPhaseC\symPhaseA}^{\realss}  }}
\newcommand{\XuunitCBreal}[0]{\textcolor{black}{{X}_{ \indexUnit, \symPhaseC\symPhaseB}^{\realss}  }}
\newcommand{\XuunitCCreal}[0]{\textcolor{black}{{X}_{ \indexUnit, \symPhaseC\symPhaseC}^{\realss}  }}

\newcommand{\XuunitAAimag}[0]{\textcolor{black}{{X}_{ \indexUnit, \symPhaseA\symPhaseA}^{\imagss}  }}
\newcommand{\XuunitABimag}[0]{\textcolor{black}{{X}_{ \indexUnit, \symPhaseA\symPhaseB}^{\imagss}  }}
\newcommand{\XuunitACimag}[0]{\textcolor{black}{{X}_{ \indexUnit, \symPhaseA\symPhaseC}^{\imagss}  }}

\newcommand{\XuunitBAimag}[0]{\textcolor{black}{{X}_{ \indexUnit, \symPhaseB\symPhaseA}^{\imagss}  }}
\newcommand{\XuunitBBimag}[0]{\textcolor{black}{{X}_{ \indexUnit, \symPhaseB\symPhaseB}^{\imagss}  }}
\newcommand{\XuunitBCimag}[0]{\textcolor{black}{{X}_{ \indexUnit, \symPhaseB\symPhaseC}^{\imagss}  }}

\newcommand{\XuunitCAimag}[0]{\textcolor{black}{{X}_{ \indexUnit, \symPhaseC\symPhaseA}^{\imagss}  }}
\newcommand{\XuunitCBimag}[0]{\textcolor{black}{{X}_{ \indexUnit, \symPhaseC\symPhaseB}^{\imagss}  }}
\newcommand{\XuunitCCimag}[0]{\textcolor{black}{{X}_{ \indexUnit, \symPhaseC\symPhaseC}^{\imagss}  }}

\newcommand{\SuunitAA}[0]{\textcolor{\complexcolor}{{\symApparentPower}_{ \indexUnit, \symPhaseA\symPhaseA}^{}  }}
\newcommand{\SuunitAB}[0]{\textcolor{\complexcolor}{{\symApparentPower}_{ \indexUnit, \symPhaseA\symPhaseB}^{}  }}
\newcommand{\SuunitAC}[0]{\textcolor{\complexcolor}{{\symApparentPower}_{ \indexUnit, \symPhaseA\symPhaseC}^{}  }}
\newcommand{\SuunitAN}[0]{\textcolor{\complexcolor}{{\symApparentPower}_{ \indexUnit, \symPhaseA\symPhaseN}^{}  }}

\newcommand{\SuunitBA}[0]{\textcolor{\complexcolor}{{\symApparentPower}_{ \indexUnit, \symPhaseB\symPhaseA}^{}  }}
\newcommand{\SuunitBB}[0]{\textcolor{\complexcolor}{{\symApparentPower}_{ \indexUnit, \symPhaseB\symPhaseB}^{}  }}
\newcommand{\SuunitBC}[0]{\textcolor{\complexcolor}{{\symApparentPower}_{ \indexUnit, \symPhaseB\symPhaseC}^{}  }}
\newcommand{\SuunitBN}[0]{\textcolor{\complexcolor}{{\symApparentPower}_{ \indexUnit, \symPhaseB\symPhaseN}^{}  }}

\newcommand{\SuunitCA}[0]{\textcolor{\complexcolor}{{\symApparentPower}_{ \indexUnit, \symPhaseC\symPhaseA}^{}  }}
\newcommand{\SuunitCB}[0]{\textcolor{\complexcolor}{{\symApparentPower}_{ \indexUnit, \symPhaseC\symPhaseB}^{}  }}
\newcommand{\SuunitCC}[0]{\textcolor{\complexcolor}{{\symApparentPower}_{ \indexUnit, \symPhaseC\symPhaseC}^{}  }}
\newcommand{\SuunitCN}[0]{\textcolor{\complexcolor}{{\symApparentPower}_{ \indexUnit, \symPhaseC\symPhaseN}^{}  }}
\newcommand{\SuunitNN}[0]{\textcolor{\complexcolor}{{\symApparentPower}_{ \indexUnit, \symPhaseN\symPhaseN}^{}  }}
\newcommand{\SuunitGG}[0]{\textcolor{\complexcolor}{{\symApparentPower}_{ \indexUnit, \symPhaseG\symPhaseG}^{}  }}

\newcommand{\SuunitNA}[0]{\textcolor{\complexcolor}{{\symApparentPower}_{ \indexUnit, \symPhaseN\symPhaseA}^{}  }}
\newcommand{\SuunitNB}[0]{\textcolor{\complexcolor}{{\symApparentPower}_{ \indexUnit, \symPhaseN\symPhaseB}^{}  }}
\newcommand{\SuunitNC}[0]{\textcolor{\complexcolor}{{\symApparentPower}_{ \indexUnit, \symPhaseN\symPhaseC}^{}  }}

\newcommand{\PuunitNN}[0]{\textcolor{black}{{\symPower}_{ \indexUnit, \symPhaseN\symPhaseN}^{}  }}
\newcommand{\PuunitNA}[0]{\textcolor{black}{{\symPower}_{ \indexUnit, \symPhaseN\symPhaseA}^{}  }}
\newcommand{\PuunitNB}[0]{\textcolor{black}{{\symPower}_{ \indexUnit, \symPhaseN\symPhaseB}^{}  }}
\newcommand{\PuunitNC}[0]{\textcolor{black}{{\symPower}_{ \indexUnit, \symPhaseN\symPhaseC}^{}  }}

\newcommand{\PuunitAA}[0]{\textcolor{black}{{\symPower}_{ \indexUnit, \symPhaseA\symPhaseA}^{}  }}
\newcommand{\PuunitAB}[0]{\textcolor{black}{{\symPower}_{ \indexUnit, \symPhaseA\symPhaseB}^{}  }}
\newcommand{\PuunitAC}[0]{\textcolor{black}{{\symPower}_{ \indexUnit, \symPhaseA\symPhaseC}^{}  }}

\newcommand{\PuunitBA}[0]{\textcolor{black}{{\symPower}_{ \indexUnit, \symPhaseB\symPhaseA}^{}  }}
\newcommand{\PuunitBB}[0]{\textcolor{black}{{\symPower}_{ \indexUnit, \symPhaseB\symPhaseB}^{}  }}
\newcommand{\PuunitBC}[0]{\textcolor{black}{{\symPower}_{ \indexUnit, \symPhaseB\symPhaseC}^{}  }}

\newcommand{\PuunitCA}[0]{\textcolor{black}{{\symPower}_{ \indexUnit, \symPhaseC\symPhaseA}^{}  }}
\newcommand{\PuunitCB}[0]{\textcolor{black}{{\symPower}_{ \indexUnit, \symPhaseC\symPhaseB}^{}  }}
\newcommand{\PuunitCC}[0]{\textcolor{black}{{\symPower}_{ \indexUnit, \symPhaseC\symPhaseC}^{}  }}

\newcommand{\PuunitAAref}[0]{\textcolor{\paramcolor}{{\symPower}_{ \indexUnit, \symPhaseA\symPhaseA}^{\refss}  }}
\newcommand{\PuunitBBref}[0]{\textcolor{\paramcolor}{{\symPower}_{ \indexUnit, \symPhaseB\symPhaseB}^{\refss}  }}
\newcommand{\PuunitCCref}[0]{\textcolor{\paramcolor}{{\symPower}_{ \indexUnit, \symPhaseC\symPhaseC}^{\refss}  }}

\newcommand{\QuunitAAref}[0]{\textcolor{\paramcolor}{{\symReactivePower}_{ \indexUnit, \symPhaseA\symPhaseA}^{\refss}  }}
\newcommand{\QuunitBBref}[0]{\textcolor{\paramcolor}{{\symReactivePower}_{ \indexUnit, \symPhaseB\symPhaseB}^{\refss}  }}
\newcommand{\QuunitCCref}[0]{\textcolor{\paramcolor}{{\symReactivePower}_{ \indexUnit, \symPhaseC\symPhaseC}^{\refss}  }}

\newcommand{\QuunitNN}[0]{\textcolor{black}{{\symReactivePower}_{ \indexUnit, \symPhaseN\symPhaseN}^{}  }}
\newcommand{\QuunitNA}[0]{\textcolor{black}{{\symReactivePower}_{ \indexUnit, \symPhaseN\symPhaseA}^{}  }}
\newcommand{\QuunitNB}[0]{\textcolor{black}{{\symReactivePower}_{ \indexUnit, \symPhaseN\symPhaseB}^{}  }}
\newcommand{\QuunitNC}[0]{\textcolor{black}{{\symReactivePower}_{ \indexUnit, \symPhaseN\symPhaseC}^{}  }}

\newcommand{\QuunitAA}[0]{\textcolor{black}{{\symReactivePower}_{ \indexUnit, \symPhaseA\symPhaseA}^{}  }}
\newcommand{\QuunitAB}[0]{\textcolor{black}{{\symReactivePower}_{ \indexUnit, \symPhaseA\symPhaseB}^{}  }}
\newcommand{\QuunitAC}[0]{\textcolor{black}{{\symReactivePower}_{ \indexUnit, \symPhaseA\symPhaseC}^{}  }}

\newcommand{\QuunitBA}[0]{\textcolor{black}{{\symReactivePower}_{ \indexUnit, \symPhaseB\symPhaseA}^{}  }}
\newcommand{\QuunitBB}[0]{\textcolor{black}{{\symReactivePower}_{ \indexUnit, \symPhaseB\symPhaseB}^{}  }}
\newcommand{\QuunitBC}[0]{\textcolor{black}{{\symReactivePower}_{ \indexUnit, \symPhaseB\symPhaseC}^{}  }}

\newcommand{\QuunitCA}[0]{\textcolor{black}{{\symReactivePower}_{ \indexUnit, \symPhaseC\symPhaseA}^{}  }}
\newcommand{\QuunitCB}[0]{\textcolor{black}{{\symReactivePower}_{ \indexUnit, \symPhaseC\symPhaseB}^{}  }}
\newcommand{\QuunitCC}[0]{\textcolor{black}{{\symReactivePower}_{ \indexUnit, \symPhaseC\symPhaseC}^{}  }}

\newcommand{\QuunitAN}[0]{\textcolor{black}{{\symReactivePower}_{ \indexUnit, \symPhaseA}^{}  }}
\newcommand{\QuunitBN}[0]{\textcolor{black}{{\symReactivePower}_{ \indexUnit, \symPhaseB}^{}  }}
\newcommand{\QuunitCN}[0]{\textcolor{black}{{\symReactivePower}_{ \indexUnit, \symPhaseC}^{}  }}

\newcommand{\SuunitrefSDP}[0]{\textcolor{\complexparamcolor}{\mathbf{\symApparentPower}_{ \indexUnit}^{\refss}  }}
\newcommand{\PuunitrefSDP}[0]{\textcolor{\paramcolor}{\mathbf{\symPower}_{ \indexUnit}^{\refss}  }}
\newcommand{\QuunitrefSDP}[0]{\textcolor{\paramcolor}{\mathbf{\symReactivePower}_{ \indexUnit}^{\refss}  }}

\newcommand{\PuunitrefAN}[0]{\textcolor{\paramcolor}{{\symPower}_{ \indexUnit, \symPhaseA}^{\refss}  }}
\newcommand{\PuunitrefBN}[0]{\textcolor{\paramcolor}{{\symPower}_{ \indexUnit, \symPhaseB}^{\refss}  }}
\newcommand{\PuunitrefCN}[0]{\textcolor{\paramcolor}{{\symPower}_{ \indexUnit, \symPhaseC}^{\refss}  }}

\newcommand{\QuunitrefAN}[0]{\textcolor{\paramcolor}{{\symReactivePower}_{ \indexUnit, \symPhaseA}^{\refss}  }}
\newcommand{\QuunitrefBN}[0]{\textcolor{\paramcolor}{{\symReactivePower}_{ \indexUnit, \symPhaseB}^{\refss}  }}
\newcommand{\QuunitrefCN}[0]{\textcolor{\paramcolor}{{\symReactivePower}_{ \indexUnit, \symPhaseC}^{\refss}  }}

\newcommand{\SbNN}[0]{\textcolor{\complexcolor}{{\symApparentPower}_{ \indexShunt, \symPhaseN\symPhaseN}^{}  }}

\newcommand{\PijkrefAA}[0]{\textcolor{\paramcolor}{\symPower_{\indexGridLines\indexGridNode \indexGridNodeTwo, \symPhaseA\symPhaseA}^{\refss}  }}
\newcommand{\PijkrefBB}[0]{\textcolor{\paramcolor}{\symPower_{\indexGridLines\indexGridNode \indexGridNodeTwo, \symPhaseB\symPhaseB}^{\refss}  }}
\newcommand{\PijkrefCC}[0]{\textcolor{\paramcolor}{\symPower_{\indexGridLines\indexGridNode \indexGridNodeTwo, \symPhaseC\symPhaseC}^{\refss}  }}

\newcommand{\Pijkpp}[0]{\textcolor{black}{\symPower_{\indexGridLines\indexGridNode \indexGridNodeTwo, \indexPhases\indexPhases}^{}  }}
\newcommand{\Qijkpp}[0]{\textcolor{black}{\symReactivePower_{\indexGridLines\indexGridNode \indexGridNodeTwo, \indexPhases\indexPhases}^{}  }}

\newcommand{\PjikAA}[0]{\textcolor{black}{\symPower_{\indexGridLines\indexGridNodeTwo\indexGridNode , \symPhaseA\symPhaseA}^{}  }}
\newcommand{\PijkAAloss}[0]{\textcolor{black}{\symPower_{\indexGridLines\indexGridNode \indexGridNodeTwo, \symPhaseA\symPhaseA}^{\lossss}  }}
\newcommand{\PjikAAloss}[0]{\textcolor{black}{\symPower_{\indexGridLines\indexGridNodeTwo\indexGridNode , \symPhaseA\symPhaseA}^{\lossss}  }}

\newcommand{\SlklossAA}[0]{\textcolor{\complexcolor}{\symApparentPower_{\indexGridLines , \symPhaseA\symPhaseA}^{\lossss}  }}
\newcommand{\SlklossBB}[0]{\textcolor{\complexcolor}{\symApparentPower_{\indexGridLines , \symPhaseB\symPhaseB}^{\lossss}  }}
\newcommand{\SlklossCC}[0]{\textcolor{\complexcolor}{\symApparentPower_{\indexGridLines , \symPhaseC\symPhaseC}^{\lossss}  }}
\newcommand{\SlklossNN}[0]{\textcolor{\complexcolor}{\symApparentPower_{\indexGridLines , \symPhaseN\symPhaseN}^{\lossss}  }}
\newcommand{\SlklossGG}[0]{\textcolor{\complexcolor}{\symApparentPower_{\indexGridLines , \symPhaseG\symPhaseG}^{\lossss}  }}
 \newcommand{\Slklosstot}[0]{\textcolor{\complexcolor}{\symApparentPower_{\indexGridLines}^{\text{tot}}  }}

\newcommand{\SlklossAAkron}[0]{\textcolor{\complexcolor}{\symApparentPower_{\indexGridLines , \symPhaseA\symPhaseA}^{\lossss, \text{Kr}}  }}
\newcommand{\SlklossBBkron}[0]{\textcolor{\complexcolor}{\symApparentPower_{\indexGridLines , \symPhaseB\symPhaseB}^{\lossss, \text{Kr}}  }}
\newcommand{\SlklossCCkron}[0]{\textcolor{\complexcolor}{\symApparentPower_{\indexGridLines , \symPhaseC\symPhaseC}^{\lossss, \text{Kr}}  }}
 \newcommand{\Slklosstotkron}[0]{\textcolor{\complexcolor}{\symApparentPower_{\indexGridLines}^{\text{tot, Kr}}  }}

\newcommand{\PlklossAA}[0]{\textcolor{black}{\symPower_{\indexGridLines , \symPhaseA\symPhaseA}^{\lossss}  }}
\newcommand{\PlklossBB}[0]{\textcolor{black}{\symPower_{\indexGridLines , \symPhaseB\symPhaseB}^{\lossss}  }}
\newcommand{\PlklossCC}[0]{\textcolor{black}{\symPower_{\indexGridLines , \symPhaseC\symPhaseC}^{\lossss}  }}
\newcommand{\PlklossNN}[0]{\textcolor{black}{\symPower_{\indexGridLines , \symPhaseN\symPhaseN}^{\lossss}  }}
\newcommand{\PlklossGG}[0]{\textcolor{black}{\symPower_{\indexGridLines , \symPhaseG\symPhaseG}^{\lossss}  }}
 \newcommand{\Plklosstot}[0]{\textcolor{black}{\symPower_{\indexGridLines}^{\text{tot}}  }}

\newcommand{\QlklossAA}[0]{\textcolor{black}{\symReactivePower_{\indexGridLines , \symPhaseA\symPhaseA}^{\lossss}  }}
\newcommand{\QlklossBB}[0]{\textcolor{black}{\symReactivePower_{\indexGridLines , \symPhaseB\symPhaseB}^{\lossss}  }}
\newcommand{\QlklossCC}[0]{\textcolor{black}{\symReactivePower_{\indexGridLines , \symPhaseC\symPhaseC}^{\lossss}  }}
\newcommand{\QlklossNN}[0]{\textcolor{black}{\symReactivePower_{\indexGridLines , \symPhaseN\symPhaseN}^{\lossss}  }}
\newcommand{\QlklossGG}[0]{\textcolor{black}{\symReactivePower_{\indexGridLines , \symPhaseG\symPhaseG}^{\lossss}  }}
 \newcommand{\Qlklosstot}[0]{\textcolor{black}{\symReactivePower_{\indexGridLines}^{\text{tot}}  }}

\newcommand{\SijkAA}[0]{\textcolor{\complexcolor}{\symApparentPower_{\indexGridLines\indexGridNode \indexGridNodeTwo, \symPhaseA\symPhaseA}^{}  }}
\newcommand{\SijkAB}[0]{\textcolor{\complexcolor}{\symApparentPower_{\indexGridLines\indexGridNode \indexGridNodeTwo, \symPhaseA\symPhaseB}^{}  }}
\newcommand{\SijkAC}[0]{\textcolor{\complexcolor}{\symApparentPower_{\indexGridLines\indexGridNode \indexGridNodeTwo, \symPhaseA\symPhaseC}^{}  }}
\newcommand{\SijkAN}[0]{\textcolor{\complexcolor}{\symApparentPower_{\indexGridLines\indexGridNode \indexGridNodeTwo, \symPhaseA\symPhaseN}^{}  }}
\newcommand{\SijkAG}[0]{\textcolor{\complexcolor}{\symApparentPower_{\indexGridLines\indexGridNode \indexGridNodeTwo, \symPhaseA\symPhaseG}^{}  }}
\newcommand{\SijkBA}[0]{\textcolor{\complexcolor}{\symApparentPower_{\indexGridLines\indexGridNode \indexGridNodeTwo, \symPhaseB\symPhaseA}^{}  }}
\newcommand{\SijkBB}[0]{\textcolor{\complexcolor}{\symApparentPower_{\indexGridLines\indexGridNode \indexGridNodeTwo, \symPhaseB\symPhaseB}^{}  }}
\newcommand{\SijkBC}[0]{\textcolor{\complexcolor}{\symApparentPower_{\indexGridLines\indexGridNode \indexGridNodeTwo, \symPhaseB\symPhaseC}^{}  }}
\newcommand{\SijkBN}[0]{\textcolor{\complexcolor}{\symApparentPower_{\indexGridLines\indexGridNode \indexGridNodeTwo, \symPhaseB\symPhaseN}^{}  }}
\newcommand{\SijkBG}[0]{\textcolor{\complexcolor}{\symApparentPower_{\indexGridLines\indexGridNode \indexGridNodeTwo, \symPhaseB\symPhaseG}^{}  }}
\newcommand{\SijkCA}[0]{\textcolor{\complexcolor}{\symApparentPower_{\indexGridLines\indexGridNode \indexGridNodeTwo, \symPhaseC\symPhaseA}^{}  }}
\newcommand{\SijkCB}[0]{\textcolor{\complexcolor}{\symApparentPower_{\indexGridLines\indexGridNode \indexGridNodeTwo, \symPhaseC\symPhaseB}^{}  }}
\newcommand{\SijkCC}[0]{\textcolor{\complexcolor}{\symApparentPower_{\indexGridLines\indexGridNode \indexGridNodeTwo, \symPhaseC\symPhaseC}^{}  }}
\newcommand{\SijkCN}[0]{\textcolor{\complexcolor}{\symApparentPower_{\indexGridLines\indexGridNode \indexGridNodeTwo, \symPhaseC\symPhaseN}^{}  }}
\newcommand{\SijkCG}[0]{\textcolor{\complexcolor}{\symApparentPower_{\indexGridLines\indexGridNode \indexGridNodeTwo, \symPhaseC\symPhaseG}^{}  }}

\newcommand{\SijkNA}[0]{\textcolor{\complexcolor}{\symApparentPower_{\indexGridLines\indexGridNode \indexGridNodeTwo, \symPhaseN\symPhaseA}^{}  }}
\newcommand{\SijkNB}[0]{\textcolor{\complexcolor}{\symApparentPower_{\indexGridLines\indexGridNode \indexGridNodeTwo, \symPhaseN\symPhaseB}^{}  }}
\newcommand{\SijkNC}[0]{\textcolor{\complexcolor}{\symApparentPower_{\indexGridLines\indexGridNode \indexGridNodeTwo, \symPhaseN\symPhaseC}^{}  }}
\newcommand{\SijkNN}[0]{\textcolor{\complexcolor}{\symApparentPower_{\indexGridLines\indexGridNode \indexGridNodeTwo, \symPhaseN\symPhaseN}^{}  }}
\newcommand{\SijkNG}[0]{\textcolor{\complexcolor}{\symApparentPower_{\indexGridLines\indexGridNode \indexGridNodeTwo, \symPhaseN\symPhaseG}^{}  }}
\newcommand{\PijkNG}[0]{\textcolor{black}{\symPower_{\indexGridLines\indexGridNode \indexGridNodeTwo, \symPhaseN\symPhaseG}^{}  }}
\newcommand{\QijkNG}[0]{\textcolor{black}{\symReactivePower_{\indexGridLines\indexGridNode \indexGridNodeTwo, \symPhaseN\symPhaseG}^{}  }}

\newcommand{\SijkGA}[0]{\textcolor{\complexcolor}{\symApparentPower_{\indexGridLines\indexGridNode \indexGridNodeTwo, \symPhaseG\symPhaseA}^{}  }}
\newcommand{\SijkGB}[0]{\textcolor{\complexcolor}{\symApparentPower_{\indexGridLines\indexGridNode \indexGridNodeTwo, \symPhaseG\symPhaseB}^{}  }}
\newcommand{\SijkGC}[0]{\textcolor{\complexcolor}{\symApparentPower_{\indexGridLines\indexGridNode \indexGridNodeTwo, \symPhaseG\symPhaseC}^{}  }}
\newcommand{\SijkGN}[0]{\textcolor{\complexcolor}{\symApparentPower_{\indexGridLines\indexGridNode \indexGridNodeTwo, \symPhaseG\symPhaseN}^{}  }}
\newcommand{\SijkGG}[0]{\textcolor{\complexcolor}{\symApparentPower_{\indexGridLines\indexGridNode \indexGridNodeTwo, \symPhaseG\symPhaseG}^{}  }}

\newcommand{\SjikAA}[0]{\textcolor{\complexcolor}{\symApparentPower_{\indexGridLines\indexGridNodeTwo\indexGridNode , \symPhaseA\symPhaseA}^{}  }}
\newcommand{\SjikAB}[0]{\textcolor{\complexcolor}{\symApparentPower_{\indexGridLines\indexGridNodeTwo\indexGridNode , \symPhaseA\symPhaseB}^{}  }}
\newcommand{\SjikAC}[0]{\textcolor{\complexcolor}{\symApparentPower_{\indexGridLines\indexGridNodeTwo\indexGridNode , \symPhaseA\symPhaseC}^{}  }}
\newcommand{\SjikAN}[0]{\textcolor{\complexcolor}{\symApparentPower_{\indexGridLines\indexGridNodeTwo\indexGridNode , \symPhaseA\symPhaseN}^{}  }}
\newcommand{\SjikBB}[0]{\textcolor{\complexcolor}{\symApparentPower_{\indexGridLines\indexGridNodeTwo\indexGridNode , \symPhaseB\symPhaseB}^{}  }}
\newcommand{\SjikCC}[0]{\textcolor{\complexcolor}{\symApparentPower_{\indexGridLines\indexGridNodeTwo\indexGridNode , \symPhaseC\symPhaseC}^{}  }}
\newcommand{\SjikNA}[0]{\textcolor{\complexcolor}{\symApparentPower_{\indexGridLines\indexGridNodeTwo\indexGridNode , \symPhaseN\symPhaseA}^{}  }}
\newcommand{\SjikNB}[0]{\textcolor{\complexcolor}{\symApparentPower_{\indexGridLines\indexGridNodeTwo\indexGridNode , \symPhaseN\symPhaseB}^{}  }}
\newcommand{\SjikNC}[0]{\textcolor{\complexcolor}{\symApparentPower_{\indexGridLines\indexGridNodeTwo\indexGridNode , \symPhaseN\symPhaseC}^{}  }}
\newcommand{\SjikNN}[0]{\textcolor{\complexcolor}{\symApparentPower_{\indexGridLines\indexGridNodeTwo\indexGridNode , \symPhaseN\symPhaseN}^{}  }}

\newcommand{\SikNG}[0]{\textcolor{\complexcolor}{\symApparentPower_{\indexGridNode, \symPhaseN\symPhaseG}^{}  }}
\newcommand{\PikNG}[0]{\textcolor{black}{\symPower_{\indexGridNode, \symPhaseN\symPhaseG}^{}  }}
\newcommand{\QikNG}[0]{\textcolor{black}{\symReactivePower_{\indexGridNode, \symPhaseN\symPhaseG}^{}  }}

\newcommand{\PijkAA}[0]{\textcolor{black}{\symPower_{\indexGridLines\indexGridNode \indexGridNodeTwo, \symPhaseA\symPhaseA}^{}  }}
\newcommand{\PijkAB}[0]{\textcolor{black}{\symPower_{\indexGridLines\indexGridNode \indexGridNodeTwo, \symPhaseA\symPhaseB}^{}  }}
\newcommand{\PijkAC}[0]{\textcolor{black}{\symPower_{\indexGridLines\indexGridNode \indexGridNodeTwo, \symPhaseA\symPhaseC}^{}  }}
\newcommand{\PijkBA}[0]{\textcolor{black}{\symPower_{\indexGridLines\indexGridNode \indexGridNodeTwo, \symPhaseB\symPhaseA}^{}  }}
\newcommand{\PijkBB}[0]{\textcolor{black}{\symPower_{\indexGridLines\indexGridNode \indexGridNodeTwo, \symPhaseB\symPhaseB}^{}  }}
\newcommand{\PijkBC}[0]{\textcolor{black}{\symPower_{\indexGridLines\indexGridNode \indexGridNodeTwo, \symPhaseB\symPhaseC}^{}  }}
\newcommand{\PijkCA}[0]{\textcolor{black}{\symPower_{\indexGridLines\indexGridNode \indexGridNodeTwo, \symPhaseC\symPhaseA}^{}  }}
\newcommand{\PijkCB}[0]{\textcolor{black}{\symPower_{\indexGridLines\indexGridNode \indexGridNodeTwo, \symPhaseC\symPhaseB}^{}  }}
\newcommand{\PijkCC}[0]{\textcolor{black}{\symPower_{\indexGridLines\indexGridNode \indexGridNodeTwo, \symPhaseC\symPhaseC}^{}  }}
\newcommand{\PijkNN}[0]{\textcolor{black}{\symPower_{\indexGridLines\indexGridNode \indexGridNodeTwo, \symPhaseN\symPhaseN}^{}  }}
\newcommand{\QijkAA}[0]{\textcolor{black}{\symReactivePower_{\indexGridLines\indexGridNode \indexGridNodeTwo, \symPhaseA\symPhaseA}^{}  }}
\newcommand{\QijkAB}[0]{\textcolor{black}{\symReactivePower_{\indexGridLines\indexGridNode \indexGridNodeTwo, \symPhaseA\symPhaseB}^{}  }}
\newcommand{\QijkAC}[0]{\textcolor{black}{\symReactivePower_{\indexGridLines\indexGridNode \indexGridNodeTwo, \symPhaseA\symPhaseC}^{}  }}
\newcommand{\QijkBA}[0]{\textcolor{black}{\symReactivePower_{\indexGridLines\indexGridNode \indexGridNodeTwo, \symPhaseB\symPhaseA}^{}  }}
\newcommand{\QijkBB}[0]{\textcolor{black}{\symReactivePower_{\indexGridLines\indexGridNode \indexGridNodeTwo, \symPhaseB\symPhaseB}^{}  }}
\newcommand{\QijkBC}[0]{\textcolor{black}{\symReactivePower_{\indexGridLines\indexGridNode \indexGridNodeTwo, \symPhaseB\symPhaseC}^{}  }}
\newcommand{\QijkCA}[0]{\textcolor{black}{\symReactivePower_{\indexGridLines\indexGridNode \indexGridNodeTwo, \symPhaseC\symPhaseA}^{}  }}
\newcommand{\QijkCB}[0]{\textcolor{black}{\symReactivePower_{\indexGridLines\indexGridNode \indexGridNodeTwo, \symPhaseC\symPhaseB}^{}  }}
\newcommand{\QijkCC}[0]{\textcolor{black}{\symReactivePower_{\indexGridLines\indexGridNode \indexGridNodeTwo, \symPhaseC\symPhaseC}^{}  }}
\newcommand{\QijkNN}[0]{\textcolor{black}{\symReactivePower_{\indexGridLines\indexGridNode \indexGridNodeTwo, \symPhaseN\symPhaseN}^{}  }}

\newcommand{\PjikAB}[0]{\textcolor{black}{\symPower_{\indexGridLines \indexGridNodeTwo \indexGridNode, \symPhaseA\symPhaseB}^{}  }}
\newcommand{\QjikAB}[0]{\textcolor{black}{\symReactivePower_{\indexGridLines \indexGridNodeTwo \indexGridNode, \symPhaseA\symPhaseB}^{}  }}

\newcommand{\PijksAA}[0]{\textcolor{black}{\symPower_{\indexGridLines\indexGridNode \indexGridNodeTwo, \symPhaseA\symPhaseA }^{\seriesss}  }}
\newcommand{\QijksAA}[0]{\textcolor{black}{\symReactivePower_{\indexGridLines\indexGridNode \indexGridNodeTwo, \symPhaseA\symPhaseA }^{\seriesss}  }}
\newcommand{\SijksAA}[0]{\textcolor{\complexcolor}{\symApparentPower_{\indexGridLines\indexGridNode \indexGridNodeTwo, \symPhaseA\symPhaseA }^{\seriesss}  }}
\newcommand{\PijksAB}[0]{\textcolor{black}{\symPower_{\indexGridLines\indexGridNode \indexGridNodeTwo, \symPhaseA\symPhaseB }^{\seriesss}  }}
\newcommand{\QijksAB}[0]{\textcolor{black}{\symReactivePower_{\indexGridLines\indexGridNode \indexGridNodeTwo, \symPhaseA\symPhaseB }^{\seriesss}  }}
\newcommand{\SijksAB}[0]{\textcolor{\complexcolor}{\symApparentPower_{\indexGridLines\indexGridNode \indexGridNodeTwo, \symPhaseA\symPhaseB }^{\seriesss}  }}
\newcommand{\PijksAC}[0]{\textcolor{black}{\symPower_{\indexGridLines\indexGridNode \indexGridNodeTwo, \symPhaseA\symPhaseC }^{\seriesss}  }}
\newcommand{\QijksAC}[0]{\textcolor{black}{\symReactivePower_{\indexGridLines\indexGridNode \indexGridNodeTwo, \symPhaseA\symPhaseC }^{\seriesss}  }}
\newcommand{\SijksAC}[0]{\textcolor{\complexcolor}{\symApparentPower_{\indexGridLines\indexGridNode \indexGridNodeTwo, \symPhaseA\symPhaseC }^{\seriesss}  }}
\newcommand{\PijksBA}[0]{\textcolor{black}{\symPower_{\indexGridLines\indexGridNode \indexGridNodeTwo, \symPhaseB\symPhaseA}^{\seriesss}  }}
\newcommand{\QijksBA}[0]{\textcolor{black}{\symReactivePower_{\indexGridLines\indexGridNode \indexGridNodeTwo, \symPhaseB\symPhaseA}^{\seriesss}  }}
\newcommand{\SijksBA}[0]{\textcolor{\complexcolor}{\symApparentPower_{\indexGridLines\indexGridNode \indexGridNodeTwo, \symPhaseB\symPhaseA}^{\seriesss}  }}
\newcommand{\PijksBB}[0]{\textcolor{black}{\symPower_{\indexGridLines\indexGridNode \indexGridNodeTwo, \symPhaseB\symPhaseB}^{\seriesss}  }}
\newcommand{\QijksBB}[0]{\textcolor{black}{\symReactivePower_{\indexGridLines\indexGridNode \indexGridNodeTwo, \symPhaseB\symPhaseB}^{\seriesss}  }}
\newcommand{\SijksBB}[0]{\textcolor{\complexcolor}{\symApparentPower_{\indexGridLines\indexGridNode \indexGridNodeTwo, \symPhaseB\symPhaseB}^{\seriesss}  }}
\newcommand{\PijksBC}[0]{\textcolor{black}{\symPower_{\indexGridLines\indexGridNode \indexGridNodeTwo, \symPhaseB\symPhaseC}^{\seriesss}  }}
\newcommand{\QijksBC}[0]{\textcolor{black}{\symReactivePower_{\indexGridLines\indexGridNode \indexGridNodeTwo, \symPhaseB\symPhaseC}^{\seriesss}  }}
\newcommand{\SijksBC}[0]{\textcolor{\complexcolor}{\symApparentPower_{\indexGridLines\indexGridNode \indexGridNodeTwo, \symPhaseB\symPhaseC}^{\seriesss}  }}
\newcommand{\PijksCA}[0]{\textcolor{black}{\symPower_{\indexGridLines\indexGridNode \indexGridNodeTwo, \symPhaseC\symPhaseA}^{\seriesss}  }}
\newcommand{\QijksCA}[0]{\textcolor{black}{\symReactivePower_{\indexGridLines\indexGridNode \indexGridNodeTwo, \symPhaseC\symPhaseA}^{\seriesss}  }}
\newcommand{\SijksCA}[0]{\textcolor{\complexcolor}{\symApparentPower_{\indexGridLines\indexGridNode \indexGridNodeTwo, \symPhaseC\symPhaseA}^{\seriesss}  }}
\newcommand{\PijksCB}[0]{\textcolor{black}{\symPower_{\indexGridLines\indexGridNode \indexGridNodeTwo, \symPhaseC\symPhaseB}^{\seriesss}  }}
\newcommand{\QijksCB}[0]{\textcolor{black}{\symReactivePower_{\indexGridLines\indexGridNode \indexGridNodeTwo, \symPhaseC\symPhaseB}^{\seriesss}  }}
\newcommand{\SijksCB}[0]{\textcolor{\complexcolor}{\symApparentPower_{\indexGridLines\indexGridNode \indexGridNodeTwo, \symPhaseC\symPhaseB}^{\seriesss}  }}
\newcommand{\PijksCC}[0]{\textcolor{black}{\symPower_{\indexGridLines\indexGridNode \indexGridNodeTwo, \symPhaseC\symPhaseC}^{\seriesss}  }}
\newcommand{\QijksCC}[0]{\textcolor{black}{\symReactivePower_{\indexGridLines\indexGridNode \indexGridNodeTwo, \symPhaseC\symPhaseC}^{\seriesss}  }}
\newcommand{\SijksCC}[0]{\textcolor{\complexcolor}{\symApparentPower_{\indexGridLines\indexGridNode \indexGridNodeTwo, \symPhaseC\symPhaseC}^{\seriesss}  }}

\newcommand{\PijksAN}[0]{\textcolor{black}{\symPower_{\indexGridLines\indexGridNode \indexGridNodeTwo, \symPhaseA\symPhaseN }^{\seriesss}  }}
\newcommand{\QijksAN}[0]{\textcolor{black}{\symReactivePower_{\indexGridLines\indexGridNode \indexGridNodeTwo, \symPhaseA\symPhaseN }^{\seriesss}  }}
\newcommand{\SijksAN}[0]{\textcolor{\complexcolor}{\symApparentPower_{\indexGridLines\indexGridNode \indexGridNodeTwo, \symPhaseA\symPhaseN }^{\seriesss}  }}
\newcommand{\PijksBN}[0]{\textcolor{black}{\symPower_{\indexGridLines\indexGridNode \indexGridNodeTwo, \symPhaseB\symPhaseN }^{\seriesss}  }}
\newcommand{\QijksBN}[0]{\textcolor{black}{\symReactivePower_{\indexGridLines\indexGridNode \indexGridNodeTwo, \symPhaseB\symPhaseN }^{\seriesss}  }}
\newcommand{\SijksBN}[0]{\textcolor{\complexcolor}{\symApparentPower_{\indexGridLines\indexGridNode \indexGridNodeTwo, \symPhaseB\symPhaseN }^{\seriesss}  }}
\newcommand{\PijksCN}[0]{\textcolor{black}{\symPower_{\indexGridLines\indexGridNode \indexGridNodeTwo, \symPhaseC\symPhaseN }^{\seriesss}  }}
\newcommand{\QijksCN}[0]{\textcolor{black}{\symReactivePower_{\indexGridLines\indexGridNode \indexGridNodeTwo, \symPhaseC\symPhaseN }^{\seriesss}  }}
\newcommand{\SijksCN}[0]{\textcolor{\complexcolor}{\symApparentPower_{\indexGridLines\indexGridNode \indexGridNodeTwo, \symPhaseC\symPhaseN }^{\seriesss}  }}
\newcommand{\PijksNN}[0]{\textcolor{black}{\symPower_{\indexGridLines\indexGridNode \indexGridNodeTwo, \symPhaseN\symPhaseN }^{\seriesss}  }}
\newcommand{\QijksNN}[0]{\textcolor{black}{\symReactivePower_{\indexGridLines\indexGridNode \indexGridNodeTwo, \symPhaseN\symPhaseN }^{\seriesss}  }}
\newcommand{\SijksNN}[0]{\textcolor{\complexcolor}{\symApparentPower_{\indexGridLines\indexGridNode \indexGridNodeTwo, \symPhaseN\symPhaseN }^{\seriesss}  }}
\newcommand{\PijksNA}[0]{\textcolor{black}{\symPower_{\indexGridLines\indexGridNode \indexGridNodeTwo, \symPhaseN\symPhaseA }^{\seriesss}  }}
\newcommand{\QijksNA}[0]{\textcolor{black}{\symReactivePower_{\indexGridLines\indexGridNode \indexGridNodeTwo, \symPhaseN\symPhaseA }^{\seriesss}  }}
\newcommand{\SijksNA}[0]{\textcolor{\complexcolor}{\symApparentPower_{\indexGridLines\indexGridNode \indexGridNodeTwo, \symPhaseN\symPhaseA }^{\seriesss}  }}
\newcommand{\PijksNB}[0]{\textcolor{black}{\symPower_{\indexGridLines\indexGridNode \indexGridNodeTwo, \symPhaseN\symPhaseB }^{\seriesss}  }}
\newcommand{\QijksNB}[0]{\textcolor{black}{\symReactivePower_{\indexGridLines\indexGridNode \indexGridNodeTwo, \symPhaseN\symPhaseB }^{\seriesss}  }}
\newcommand{\SijksNB}[0]{\textcolor{\complexcolor}{\symApparentPower_{\indexGridLines\indexGridNode \indexGridNodeTwo, \symPhaseN\symPhaseB }^{\seriesss}  }}
\newcommand{\PijksNC}[0]{\textcolor{black}{\symPower_{\indexGridLines\indexGridNode \indexGridNodeTwo, \symPhaseN\symPhaseC }^{\seriesss}  }}
\newcommand{\QijksNC}[0]{\textcolor{black}{\symReactivePower_{\indexGridLines\indexGridNode \indexGridNodeTwo, \symPhaseN\symPhaseC }^{\seriesss}  }}
\newcommand{\SijksNC}[0]{\textcolor{\complexcolor}{\symApparentPower_{\indexGridLines\indexGridNode \indexGridNodeTwo, \symPhaseN\symPhaseC }^{\seriesss}  }}

\newcommand{\PijklossAA}[0]{\textcolor{black}{\symPower_{\indexGridNode \indexGridNodeTwo, \symPhaseA\symPhaseA}^{\lossss}  }}
\newcommand{\PijklossAB}[0]{\textcolor{black}{\symPower_{\indexGridNode \indexGridNodeTwo, \symPhaseA\symPhaseB}^{\lossss}  }}
\newcommand{\PijklossAC}[0]{\textcolor{black}{\symPower_{\indexGridNode \indexGridNodeTwo, \symPhaseA\symPhaseC}^{\lossss}  }}
\newcommand{\PijklossBA}[0]{\textcolor{black}{\symPower_{\indexGridNode \indexGridNodeTwo, \symPhaseB\symPhaseA}^{\lossss}  }}
\newcommand{\PijklossBB}[0]{\textcolor{black}{\symPower_{\indexGridNode \indexGridNodeTwo, \symPhaseB\symPhaseB}^{\lossss}  }}
\newcommand{\PijklossBC}[0]{\textcolor{black}{\symPower_{\indexGridNode \indexGridNodeTwo, \symPhaseB\symPhaseC}^{\lossss}  }}
\newcommand{\PijklossCA}[0]{\textcolor{black}{\symPower_{\indexGridNode \indexGridNodeTwo, \symPhaseC\symPhaseA}^{\lossss}  }}
\newcommand{\PijklossCB}[0]{\textcolor{black}{\symPower_{\indexGridNode \indexGridNodeTwo, \symPhaseC\symPhaseB}^{\lossss}  }}
\newcommand{\PijklossCC}[0]{\textcolor{black}{\symPower_{\indexGridNode \indexGridNodeTwo, \symPhaseC\symPhaseC}^{\lossss}  }}
\newcommand{\QijklossAA}[0]{\textcolor{black}{\symReactivePower_{\indexGridNode \indexGridNodeTwo, \symPhaseA\symPhaseA}^{\lossss}  }}
\newcommand{\QijklossAB}[0]{\textcolor{black}{\symReactivePower_{\indexGridNode \indexGridNodeTwo, \symPhaseA\symPhaseB}^{\lossss}  }}
\newcommand{\QijklossAC}[0]{\textcolor{black}{\symReactivePower_{\indexGridNode \indexGridNodeTwo, \symPhaseA\symPhaseC}^{\lossss}  }}
\newcommand{\QijklossBA}[0]{\textcolor{black}{\symReactivePower_{\indexGridNode \indexGridNodeTwo, \symPhaseB\symPhaseA}^{\lossss}  }}
\newcommand{\QijklossBB}[0]{\textcolor{black}{\symReactivePower_{\indexGridNode \indexGridNodeTwo, \symPhaseB\symPhaseB}^{\lossss}  }}
\newcommand{\QijklossBC}[0]{\textcolor{black}{\symReactivePower_{\indexGridNode \indexGridNodeTwo, \symPhaseB\symPhaseC}^{\lossss}  }}
\newcommand{\QijklossCA}[0]{\textcolor{black}{\symReactivePower_{\indexGridNode \indexGridNodeTwo, \symPhaseC\symPhaseA}^{\lossss}  }}
\newcommand{\QijklossCB}[0]{\textcolor{black}{\symReactivePower_{\indexGridNode \indexGridNodeTwo, \symPhaseC\symPhaseB}^{\lossss}  }}
\newcommand{\QijklossCC}[0]{\textcolor{black}{\symReactivePower_{\indexGridNode \indexGridNodeTwo, \symPhaseC\symPhaseC}^{\lossss}  }}


\newcommand{\Iik}[0]{\textcolor{\complexcolor}{\symCurrent_{\indexGridNode }^{}  }}

\newcommand{\Ijik}[0]{\textcolor{\complexcolor}{\symCurrent_{\indexGridLines\indexGridNodeTwo  \indexGridNode }^{}  }}

\newcommand{\Iijkre}[0]{\textcolor{black}{\symCurrent_{\indexGridLines\indexGridNode \indexGridNodeTwo}^{\realss}  }}
\newcommand{\Iijkim}[0]{\textcolor{black}{\symCurrent_{\indexGridLines\indexGridNode \indexGridNodeTwo}^{\imagss}  }}
\newcommand{\Iijkangle}[0]{\textcolor{black}{\symCurrent_{\indexGridLines\indexGridNode \indexGridNodeTwo}^{\imagangle}  }}
\newcommand{\Iijkanglemin}[0]{\textcolor{\paramcolor}{\symCurrent_{\indexGridLines\indexGridNode \indexGridNodeTwo}^{\imagangle\minss}  }}
\newcommand{\Iijkanglemax}[0]{\textcolor{\paramcolor}{\symCurrent_{\indexGridLines\indexGridNode \indexGridNodeTwo}^{\imagangle\maxss}  }}

\newcommand{\IijkSDPseq}[0]{\textcolor{\complexcolor}{\mathbf{\symCurrent}_{\indexGridLines\indexGridNode \indexGridNodeTwo   }^{\fortescuess}  }}

\newcommand{\IijkSDPreal}[0]{\textcolor{black}{\mathbf{\symCurrent}_{\indexGridLines\indexGridNode \indexGridNodeTwo   }^{\realss}  }}
\newcommand{\IijkSDPimag}[0]{\textcolor{black}{\mathbf{\symCurrent}_{\indexGridLines\indexGridNode \indexGridNodeTwo   }^{\imagss}  }}

\newcommand{\IijkSDP}[0]{\textcolor{\complexcolor}{\mathbf{\symCurrent}_{\indexGridLines\indexGridNode \indexGridNodeTwo   }^{}  }}
\newcommand{\IijkSDPH}[0]{\textcolor{\complexcolor}{(\mathbf{\symCurrent}_{\indexGridLines\indexGridNode \indexGridNodeTwo   })^{\hermitiantranspose}  }}

\newcommand{\IijkSDPref}[0]{\textcolor{\complexparamcolor}{\mathbf{\symCurrent}_{\indexGridLines\indexGridNode \indexGridNodeTwo   }^{\refss}  }}
\newcommand{\IijkSDPrefre}[0]{\textcolor{\paramcolor}{\mathbf{\symCurrent}_{\indexGridLines\indexGridNode \indexGridNodeTwo   }^{\refss, \realss}  }}
\newcommand{\IijkSDPrefim}[0]{\textcolor{\paramcolor}{\mathbf{\symCurrent}_{\indexGridLines\indexGridNode \indexGridNodeTwo   }^{\refss, \imagss}  }}

\newcommand{\IijksSDPref}[0]{\textcolor{\complexparamcolor}{\mathbf{\symCurrent}_{\indexGridLines\indexGridNode \indexGridNodeTwo   }^{\seriesss, \refss}  }}
\newcommand{\IijksSDPrefre}[0]{\textcolor{\paramcolor}{\mathbf{\symCurrent}_{\indexGridLines\indexGridNode \indexGridNodeTwo   }^{\seriesss, \refss, \realss}  }}
\newcommand{\IijksSDPrefim}[0]{\textcolor{\paramcolor}{\mathbf{\symCurrent}_{\indexGridLines\indexGridNode \indexGridNodeTwo   }^{\seriesss, \refss, \imagss}  }}

\newcommand{\IijkshSDPref}[0]{\textcolor{\complexparamcolor}{\mathbf{\symCurrent}_{\indexGridLines\indexGridNode \indexGridNodeTwo   }^{\shuntss, \refss}  }}
\newcommand{\IijkshSDPrefre}[0]{\textcolor{\paramcolor}{\mathbf{\symCurrent}_{\indexGridLines\indexGridNode \indexGridNodeTwo   }^{\shuntss, \refss, \realss}  }}
\newcommand{\IijkshSDPrefim}[0]{\textcolor{\paramcolor}{\mathbf{\symCurrent}_{\indexGridLines\indexGridNode \indexGridNodeTwo   }^{\shuntss, \refss, \imagss}  }}

\newcommand{\IijkSDPdelta}[0]{\textcolor{\complexcolor}{\mathbf{\symCurrent}_{\indexGridLines\indexGridNode \indexGridNodeTwo   }^{\Delta}  }}
\newcommand{\IijkSDPdeltare}[0]{\textcolor{black}{\mathbf{\symCurrent}_{\indexGridLines\indexGridNode \indexGridNodeTwo   }^{\Delta, \realss}  }}
\newcommand{\IijkSDPdeltaim}[0]{\textcolor{black}{\mathbf{\symCurrent}_{\indexGridLines\indexGridNode \indexGridNodeTwo   }^{\Delta, \imagss}  }}

\newcommand{\IijksSDPdelta}[0]{\textcolor{\complexcolor}{\mathbf{\symCurrent}_{\indexGridLines\indexGridNode \indexGridNodeTwo   }^{\seriesss, \Delta}  }}
\newcommand{\IijksSDPdeltare}[0]{\textcolor{black}{\mathbf{\symCurrent}_{\indexGridLines\indexGridNode \indexGridNodeTwo   }^{\seriesss, \Delta, \realss}  }}
\newcommand{\IijksSDPdeltaim}[0]{\textcolor{black}{\mathbf{\symCurrent}_{\indexGridLines\indexGridNode \indexGridNodeTwo   }^{\seriesss, \Delta, \imagss}  }}

\newcommand{\IsqijkSDPdelta}[0]{\textcolor{\complexcolor}{\mathbf{\symCurrentSquared}_{\indexGridLines\indexGridNode \indexGridNodeTwo   }^{\Delta}  }}
\newcommand{\IsqijkSDPdeltare}[0]{\textcolor{black}{\mathbf{\symCurrentSquared}_{\indexGridLines\indexGridNode \indexGridNodeTwo   }^{\Delta, \realss}  }}
\newcommand{\IsqijkSDPdeltaim}[0]{\textcolor{black}{\mathbf{\symCurrentSquared}_{\indexGridLines\indexGridNode \indexGridNodeTwo   }^{\Delta, \imagss}  }}

\newcommand{\IsqijksSDPdelta}[0]{\textcolor{\complexcolor}{\mathbf{\symCurrentSquared}_{\indexGridLines\indexGridNode \indexGridNodeTwo   }^{\seriesss, \Delta}  }}
\newcommand{\IsqijksSDPdeltare}[0]{\textcolor{black}{\mathbf{\symCurrentSquared}_{\indexGridLines\indexGridNode \indexGridNodeTwo   }^{\seriesss, \Delta, \realss}  }}
\newcommand{\IsqijksSDPdeltaim}[0]{\textcolor{black}{\mathbf{\symCurrentSquared}_{\indexGridLines\indexGridNode \indexGridNodeTwo   }^{\seriesss, \Delta, \imagss}  }}

\newcommand{\IjikSDP}[0]{\textcolor{\complexcolor}{\mathbf{\symCurrent}_{\indexGridLines\indexGridNodeTwo\indexGridNode    }^{}  }}
\newcommand{\IjikSDPreal}[0]{\textcolor{black}{\mathbf{\symCurrent}_{\indexGridLines\indexGridNodeTwo\indexGridNode    }^{\realss}  }}
\newcommand{\IjikSDPimag}[0]{\textcolor{black}{\mathbf{\symCurrent}_{\indexGridLines\indexGridNodeTwo\indexGridNode    }^{\imagss}  }}

\newcommand{\IijkrefSDP}[0]{\textcolor{\complexparamcolor}{\mathbf{\symCurrent}_{\indexGridLines\indexGridNode \indexGridNodeTwo   }^{\refss}  }}

\newcommand{\IijksrefSDP}[0]{\textcolor{\complexparamcolor}{\mathbf{\symCurrent}_{\indexGridLines\indexGridNode \indexGridNodeTwo, \seriesss   }^{\refss}  }}

\newcommand{\IijskSDPseq}[0]{\textcolor{\complexcolor}{\mathbf{\symCurrent}_{\indexGridLines\indexGridNode \indexGridNodeTwo   }^{\seriesss, \fortescuess}  }}

\newcommand{\IijskSDPH}[0]{\textcolor{\complexcolor}{(\mathbf{\symCurrent}_{\indexGridLines\indexGridNode \indexGridNodeTwo   }^{\seriesss})^{\hermitiantranspose}  }}
\newcommand{\IijskSDP}[0]{\textcolor{\complexcolor}{\mathbf{\symCurrent}_{\indexGridLines\indexGridNode \indexGridNodeTwo   }^{\seriesss}  }}
\newcommand{\IijskSDPreal}[0]{\textcolor{black}{\mathbf{\symCurrent}_{\indexGridLines\indexGridNode \indexGridNodeTwo   }^{\seriesss, \realss}  }}
\newcommand{\IijskSDPimag}[0]{\textcolor{black}{\mathbf{\symCurrent}_{\indexGridLines\indexGridNode \indexGridNodeTwo   }^{\seriesss, \imagss}  }}

\newcommand{\IjiskSDP}[0]{\textcolor{\complexcolor}{\mathbf{\symCurrent}_{\indexGridLines\indexGridNodeTwo\indexGridNode    }^{\seriesss}  }}
\newcommand{\IjiskSDPreal}[0]{\textcolor{black}{\mathbf{\symCurrent}_{\indexGridLines\indexGridNodeTwo\indexGridNode    }^{\seriesss, \realss}  }}
\newcommand{\IjiskSDPimag}[0]{\textcolor{black}{\mathbf{\symCurrent}_{\indexGridLines\indexGridNodeTwo\indexGridNode    }^{\seriesss, \imagss}  }}

\newcommand{\IijshkSDP}[0]{\textcolor{\complexcolor}{\mathbf{\symCurrent}_{\indexGridLines\indexGridNode \indexGridNodeTwo   }^{\shuntss}  }}
\newcommand{\IjishkSDP}[0]{\textcolor{\complexcolor}{\mathbf{\symCurrent}_{ \indexGridLines\indexGridNodeTwo \indexGridNode   }^{\shuntss}  }}

\newcommand{\IijshkSDPreal}[0]{\textcolor{black}{\mathbf{\symCurrent}_{\indexGridLines\indexGridNode \indexGridNodeTwo   }^{\shuntss, \realss}  }}

\newcommand{\IijshkSDPimag}[0]{\textcolor{black}{\mathbf{\symCurrent}_{\indexGridLines\indexGridNode \indexGridNodeTwo   }^{\shuntss, \imagss}  }}

\newcommand{\Iijk}[0]{\textcolor{\complexcolor}{\symCurrent_{\indexGridLines\indexGridNode \indexGridNodeTwo}^{}  }}
\newcommand{\Iijsk}[0]{\textcolor{\complexcolor}{\symCurrent_{\indexGridLines\indexGridNode \indexGridNodeTwo,\seriesss}^{}  }}
\newcommand{\Iijshk}[0]{\textcolor{\complexcolor}{\symCurrent_{\indexGridLines\indexGridNode \indexGridNodeTwo,\shuntss}^{}  }}
\newcommand{\Iijsqk}[0]{\textcolor{black}{\symCurrent_{\indexGridLines\indexGridNode \indexGridNodeTwo}^{2}  }}
\newcommand{\Iijsqsk}[0]{\textcolor{black}{\symCurrent_{\indexGridLines\indexGridLines,\seriesss}^{2}  }}
\newcommand{\Iijsqshk}[0]{\textcolor{black}{\symCurrent_{\indexGridLines\indexGridNode \indexGridNodeTwo,\shuntss}^{2}  }}

\newcommand{\IijkH}[0]{\textcolor{\complexcolor}{\symCurrent_{\indexGridLines\indexGridNode \indexGridNodeTwo}^{\hermitiantranspose}  }}
\newcommand{\Ijisk}[0]{\textcolor{\complexcolor}{\symCurrent_{\indexGridLines\indexGridNodeTwo\indexGridNode ,\seriesss}^{}  }}
\newcommand{\Ijishk}[0]{\textcolor{\complexcolor}{\symCurrent_{\indexGridLines\indexGridNodeTwo\indexGridNode ,\shuntss}^{}  }}

\newcommand{\Iijkmax}[0]{\textcolor{\boundscolor}{\symCurrent_{\indexGridLines\indexGridNode \indexGridNodeTwo}^{\maxss}  }}
\newcommand{\Iijkmin}[0]{\textcolor{\boundscolor}{\symCurrent_{\indexGridLines\indexGridNode \indexGridNodeTwo}^{\minss}  }}
\newcommand{\Iijrated}[0]{\textcolor{\sizingcolor}{\symCurrent_{\indexGridLines\indexGridNode \indexGridNodeTwo}^{\ratedss} }}
\newcommand{\Ijirated}[0]{\textcolor{\sizingcolor}{\symCurrent_{\indexGridLines\indexGridNodeTwo  \indexGridNode }^{\ratedss} }}
\newcommand{\IijratedA}[0]{\textcolor{\sizingcolor}{\symCurrent_{\indexGridLines\indexGridNode \indexGridNodeTwo,\symPhaseA}^{\ratedss} }}
\newcommand{\IjiratedA}[0]{\textcolor{\sizingcolor}{\symCurrent_{\indexGridLines\indexGridNodeTwo  \indexGridNode,\symPhaseA}^{\ratedss} }}
\newcommand{\IijratedB}[0]{\textcolor{\sizingcolor}{\symCurrent_{\indexGridLines\indexGridNode \indexGridNodeTwo,\symPhaseB}^{\ratedss} }}
\newcommand{\IjiratedB}[0]{\textcolor{\sizingcolor}{\symCurrent_{\indexGridLines\indexGridNodeTwo  \indexGridNode,\symPhaseB}^{\ratedss} }}
\newcommand{\IijratedC}[0]{\textcolor{\sizingcolor}{\symCurrent_{\indexGridLines\indexGridNode \indexGridNodeTwo,\symPhaseC}^{\ratedss} }}
\newcommand{\IjiratedC}[0]{\textcolor{\sizingcolor}{\symCurrent_{\indexGridLines\indexGridNodeTwo  \indexGridNode,\symPhaseC}^{\ratedss} }}

\newcommand{\IijratedN}[0]{\textcolor{\sizingcolor}{\symCurrent_{\indexGridLines\indexGridNode \indexGridNodeTwo,\symPhaseN}^{\ratedss} }}

\newcommand{\IijsqskSDPreal}[0]{\textcolor{black}{\mathbf{\symCurrentSOCP}_{\indexGridLines  }^{\seriesss,\realss}  }}
\newcommand{\IijsqskSDPimag}[0]{\textcolor{black}{\mathbf{\symCurrentSOCP}_{\indexGridLines  }^{\seriesss,\imagss}  }}

\newcommand{\IijsqkSDPreal}[0]{\textcolor{black}{\mathbf{\symCurrentSOCP}_{\indexGridLines\indexGridNode \indexGridNodeTwo }^{\realss}  }}
\newcommand{\IijsqkSDPimag}[0]{\textcolor{black}{\mathbf{\symCurrentSOCP}_{\indexGridLines\indexGridNode \indexGridNodeTwo }^{\imagss}  }}

\newcommand{\IijsqkSDPAdot}[0]{\textcolor{\complexcolor}{\mathbf{\symCurrentSOCP}_{\indexGridLines\indexGridNode \indexGridNodeTwo }^{\symPhaseA\!\cdot}  }}
\newcommand{\IijsqkSDPBdot}[0]{\textcolor{\complexcolor}{\mathbf{\symCurrentSOCP}_{\indexGridLines\indexGridNode \indexGridNodeTwo }^{\symPhaseB\!\cdot}  }}
\newcommand{\IijsqkSDPCdot}[0]{\textcolor{\complexcolor}{\mathbf{\symCurrentSOCP}_{\indexGridLines\indexGridNode \indexGridNodeTwo }^{\symPhaseC\!\cdot}  }}

\newcommand{\IjisqkSDP}[0]{\textcolor{\complexcolor}{\mathbf{\symCurrentSOCP}_{\indexGridLines \indexGridNodeTwo \indexGridNode  }^{}  }}
\newcommand{\IijsqkSDP}[0]{\textcolor{\complexcolor}{\mathbf{\symCurrentSOCP}_{\indexGridLines\indexGridNode \indexGridNodeTwo }^{}  }}
\newcommand{\IijsqkmaxSDP}[0]{\textcolor{\complexparamcolor}{\mathbf{\symCurrentSOCP}_{\indexGridLines\indexGridNode \indexGridNodeTwo }^{\maxss}  }}

\newcommand{\IijkmaxSDP}[0]{\textcolor{\paramcolor}{\mathbf{\symCurrent}_{\indexGridLines\indexGridNode \indexGridNodeTwo }^{\maxss}  }}

\newcommand{\IjisqkSDPreal}[0]{\textcolor{black}{\mathbf{\symCurrentSOCP}_{\indexGridLines \indexGridNodeTwo \indexGridNode  }^{\realss}  }}

\newcommand{\IijsqshkSDP}[0]{\textcolor{\complexcolor}{\mathbf{\symCurrentSOCP}_{\indexGridLines\indexGridNode \indexGridNodeTwo }^{\shuntss}  }}
\newcommand{\IjisqshkSDP}[0]{\textcolor{\complexcolor}{\mathbf{\symCurrentSOCP}_{\indexGridLines\indexGridNodeTwo \indexGridNode  }^{\shuntss}  }}

\newcommand{\IusqkSDPAdot}[0]{\textcolor{\complexcolor}{\mathbf{\symCurrentSOCP}_{\indexUnit  }^{\symPhaseA\!\cdot}  }}
\newcommand{\IusqkSDPBdot}[0]{\textcolor{\complexcolor}{\mathbf{\symCurrentSOCP}_{\indexUnit  }^{\symPhaseB\!\cdot}  }}
\newcommand{\IusqkSDPCdot}[0]{\textcolor{\complexcolor}{\mathbf{\symCurrentSOCP}_{\indexUnit  }^{\symPhaseC\!\cdot}  }}
\newcommand{\IusqkSDPNdot}[0]{\textcolor{\complexcolor}{\mathbf{\symCurrentSOCP}_{\indexUnit  }^{\symPhaseN\!\cdot}  }}

\newcommand{\IusqkSDPdotA}[0]{\textcolor{\complexcolor}{\mathbf{\symCurrentSOCP}_{\indexUnit  }^{\!\cdot\!\symPhaseA}  }}
\newcommand{\IusqkSDPdotB}[0]{\textcolor{\complexcolor}{\mathbf{\symCurrentSOCP}_{\indexUnit  }^{\!\cdot\!\symPhaseB}  }}
\newcommand{\IusqkSDPdotC}[0]{\textcolor{\complexcolor}{\mathbf{\symCurrentSOCP}_{\indexUnit  }^{\!\cdot\!\symPhaseC}  }}
\newcommand{\IusqkSDPdotN}[0]{\textcolor{\complexcolor}{\mathbf{\symCurrentSOCP}_{\indexUnit  }^{\!\cdot\!\symPhaseN}  }}

\newcommand{\IijsqskSDPseq}[0]{\textcolor{\complexcolor}{\mathbf{\symCurrentSOCP}_{\indexGridLines  }^{\seriesss, \fortescuess}  }}
\newcommand{\IusqkSDP}[0]{\textcolor{\complexcolor}{\mathbf{\symCurrentSOCP}_{\indexUnit  }  }}
\newcommand{\IusqkSDPreal}[0]{\textcolor{black}{\mathbf{\symCurrentSOCP}_{\indexUnit  }^{\realss}  }}
\newcommand{\IusqkSDPimag}[0]{\textcolor{black}{\mathbf{\symCurrentSOCP}_{\indexUnit  }^{\imagss}  }}

\newcommand{\IusqkSDPdelta}[0]{\textcolor{\complexcolor}{\mathbf{\symCurrentSOCP}_{\indexUnit  }^{\Delta}  }}
\newcommand{\IusqkSDPdeltareal}[0]{\textcolor{black}{\mathbf{\symCurrentSOCP}_{\indexUnit  }^{\Delta\realss}  }}
\newcommand{\IusqkSDPdeltaimag}[0]{\textcolor{black}{\mathbf{\symCurrentSOCP}_{\indexUnit  }^{\Delta\imagss}  }}

\newcommand{\IijsqskSDP}[0]{\textcolor{\complexcolor}{\mathbf{\symCurrentSOCP}_{\indexGridLines  }^{\seriesss}  }}
\newcommand{\IjisqskSDP}[0]{\textcolor{\complexcolor}{\mathbf{\symCurrentSOCP}_{\indexGridLines  }^{\seriesss}  }}
\newcommand{\Iijsqsksocp}[0]{\textcolor{black}{\symCurrentSOCP_{\indexGridLines\indexGridNode \indexGridNodeTwo ,\seriesss}  }}

\newcommand{\Iijsqksocp}[0]{\textcolor{black}{\symCurrentSOCP_{\indexGridLines\indexGridNode \indexGridNodeTwo}  }}

\newcommand{\IijsqsksocpAA}[0]{\textcolor{black}{\symCurrentSOCP_{\indexGridLines\indexGridNode \indexGridNodeTwo,\seriesss, \symPhaseA\symPhaseA}  }}
\newcommand{\IijsqsksocpAB}[0]{\textcolor{\complexcolor}{\symCurrentSOCP_{\indexGridLines\indexGridNode \indexGridNodeTwo,\seriesss, \symPhaseA\symPhaseB}  }}
\newcommand{\IijsqsksocpAC}[0]{\textcolor{\complexcolor}{\symCurrentSOCP_{\indexGridLines\indexGridNode \indexGridNodeTwo,\seriesss, \symPhaseA\symPhaseC}  }}
\newcommand{\IijsqsksocpBA}[0]{\textcolor{\complexcolor}{\symCurrentSOCP_{\indexGridLines\indexGridNode \indexGridNodeTwo,\seriesss, \symPhaseB\symPhaseA}  }}
\newcommand{\IijsqsksocpBB}[0]{\textcolor{black}{\symCurrentSOCP_{\indexGridLines\indexGridNode \indexGridNodeTwo,\seriesss, \symPhaseB\symPhaseB}  }}
\newcommand{\IijsqsksocpBC}[0]{\textcolor{\complexcolor}{\symCurrentSOCP_{\indexGridLines\indexGridNode \indexGridNodeTwo,\seriesss, \symPhaseB\symPhaseC}  }}
\newcommand{\IijsqsksocpCA}[0]{\textcolor{\complexcolor}{\symCurrentSOCP_{\indexGridLines\indexGridNode \indexGridNodeTwo,\seriesss, \symPhaseC\symPhaseA}  }}
\newcommand{\IijsqsksocpCB}[0]{\textcolor{\complexcolor}{\symCurrentSOCP_{\indexGridLines\indexGridNode \indexGridNodeTwo,\seriesss, \symPhaseC\symPhaseB}  }}
\newcommand{\IijsqsksocpCC}[0]{\textcolor{black}{\symCurrentSOCP_{\indexGridLines\indexGridNode \indexGridNodeTwo,\seriesss, \symPhaseC\symPhaseC}  }}

\newcommand{\IijsqksocpAAreal}[0]{\textcolor{black}{\symCurrentSOCP_{\indexGridLines\indexGridNode \indexGridNodeTwo, \symPhaseA\symPhaseA}^{\realss}  }}
\newcommand{\IijsqksocpBBreal}[0]{\textcolor{black}{\symCurrentSOCP_{\indexGridLines\indexGridNode \indexGridNodeTwo, \symPhaseB\symPhaseB}^{\realss}  }}
\newcommand{\IijsqksocpCCreal}[0]{\textcolor{black}{\symCurrentSOCP_{\indexGridLines\indexGridNode \indexGridNodeTwo, \symPhaseC\symPhaseC}^{\realss}  }}

\newcommand{\IijsqksocpNNreal}[0]{\textcolor{black}{\symCurrentSOCP_{\indexGridLines\indexGridNode \indexGridNodeTwo, \symPhaseN\symPhaseN}^{\realss}  }}
\newcommand{\IijsqksocpGGreal}[0]{\textcolor{black}{\symCurrentSOCP_{\indexGridLines\indexGridNode \indexGridNodeTwo, \symPhaseG\symPhaseG}^{\realss}  }}
\newcommand{\kifsqksocpGGreal}[0]{\textcolor{black}{\symCurrentSOCP_{k\indexGridNode f, \symPhaseG\symPhaseG}^{\realss}  }}

\newcommand{\IijingsqksocpGG}[0]{\textcolor{\complexcolor}{\symCurrentSOCP_{\indexGridLines\indexGridNode \indexGridNodeTwo, \symPhaseN\symPhaseG}  }}

\newcommand{\kifingsqksocpGG}[0]{\textcolor{\complexcolor}{\symCurrentSOCP_{k\indexGridNode f, \symPhaseN\symPhaseG}  }}
    \newcommand{\parm}{\mathord{\color{black!33}\bullet}}%

\newcommand{\IijsqsksocpAAreal}[0]{\textcolor{black}{\symCurrentSOCP_{\indexGridLines\indexGridNode \indexGridNodeTwo, \symPhaseA\symPhaseA}^{\seriesss,\realss}  }}
\newcommand{\IijsqsksocpABreal}[0]{\textcolor{black}{\symCurrentSOCP_{\indexGridLines\indexGridNode \indexGridNodeTwo, \symPhaseA\symPhaseB}^{\seriesss,\realss}  }}
\newcommand{\IijsqsksocpACreal}[0]{\textcolor{black}{\symCurrentSOCP_{\indexGridLines\indexGridNode \indexGridNodeTwo, \symPhaseA\symPhaseC}^{\seriesss,\realss}  }}
\newcommand{\IijsqsksocpBBreal}[0]{\textcolor{black}{\symCurrentSOCP_{\indexGridLines\indexGridNode \indexGridNodeTwo, \symPhaseB\symPhaseB}^{\seriesss,\realss}  }}
\newcommand{\IijsqsksocpBCreal}[0]{\textcolor{black}{\symCurrentSOCP_{\indexGridLines\indexGridNode \indexGridNodeTwo, \symPhaseB\symPhaseC}^{\seriesss,\realss}  }}
\newcommand{\IijsqsksocpCCreal}[0]{\textcolor{black}{\symCurrentSOCP_{\indexGridLines\indexGridNode \indexGridNodeTwo, \symPhaseC\symPhaseC}^{\seriesss,\realss}  }}
\newcommand{\IijsqsksocpNNreal}[0]{\textcolor{black}{\symCurrentSOCP_{\indexGridLines\indexGridNode \indexGridNodeTwo, \symPhaseN\symPhaseN}^{\seriesss,\realss}  }}

\newcommand{\IijsqsksocpANreal}[0]{\textcolor{black}{\symCurrentSOCP_{\indexGridLines\indexGridNode \indexGridNodeTwo, \symPhaseA\symPhaseN}^{\seriesss,\realss}  }}
\newcommand{\IijsqsksocpBNreal}[0]{\textcolor{black}{\symCurrentSOCP_{\indexGridLines\indexGridNode \indexGridNodeTwo, \symPhaseB\symPhaseN}^{\seriesss,\realss}  }}
\newcommand{\IijsqsksocpCNreal}[0]{\textcolor{black}{\symCurrentSOCP_{\indexGridLines\indexGridNode \indexGridNodeTwo, \symPhaseC\symPhaseN}^{\seriesss,\realss}  }}

\newcommand{\IijsqsksocpANimag}[0]{\textcolor{black}{\symCurrentSOCP_{\indexGridLines\indexGridNode \indexGridNodeTwo, \symPhaseA\symPhaseN}^{\seriesss,\imagss}  }}
\newcommand{\IijsqsksocpBNimag}[0]{\textcolor{black}{\symCurrentSOCP_{\indexGridLines\indexGridNode \indexGridNodeTwo, \symPhaseB\symPhaseN}^{\seriesss,\imagss}  }}
\newcommand{\IijsqsksocpCNimag}[0]{\textcolor{black}{\symCurrentSOCP_{\indexGridLines\indexGridNode \indexGridNodeTwo, \symPhaseC\symPhaseN}^{\seriesss,\imagss}  }}

\newcommand{\IijsqsksocpGGreal}[0]{\textcolor{black}{\symCurrentSOCP_{\indexGridLines\indexGridNode \indexGridNodeTwo, \symPhaseG\symPhaseG}^{\seriesss,\realss}  }}
\newcommand{\IisqsksocpNGreal}[0]{\textcolor{black}{\symCurrentSOCP_{\indexGridNode , \symPhaseN\symPhaseG}  }}

\newcommand{\IjisqsksocpBBreal}[0]{\textcolor{black}{\symCurrentSOCP_{\indexGridLines\indexGridNodeTwo\indexGridNode , \symPhaseB\symPhaseB}^{\seriesss,\realss}  }}

\newcommand{\IijsqsksocpABimag}[0]{\textcolor{black}{\symCurrentSOCP_{\indexGridLines\indexGridNode \indexGridNodeTwo, \symPhaseA\symPhaseB}^{\seriesss,\imagss}  }}
\newcommand{\IijsqsksocpACimag}[0]{\textcolor{black}{\symCurrentSOCP_{\indexGridLines\indexGridNode \indexGridNodeTwo, \symPhaseA\symPhaseC}^{\seriesss,\imagss}  }}
\newcommand{\IijsqsksocpBCimag}[0]{\textcolor{black}{\symCurrentSOCP_{\indexGridLines\indexGridNode \indexGridNodeTwo, \symPhaseB\symPhaseC}^{\seriesss,\imagss}  }}

\newcommand{\VjsqksocpAAreal}[0]{\textcolor{black}{\symVoltageSOCP_{\indexGridNodeTwo, \symPhaseA\symPhaseA}^{\realss}  }}
\newcommand{\VjsqksocpBBreal}[0]{\textcolor{black}{\symVoltageSOCP_{\indexGridNodeTwo, \symPhaseB\symPhaseB}^{\realss}  }}
\newcommand{\VjsqksocpCCreal}[0]{\textcolor{black}{\symVoltageSOCP_{\indexGridNodeTwo, \symPhaseC\symPhaseC}^{\realss}  }}

\newcommand{\VisqksocpNNreal}[0]{\textcolor{black}{\symVoltageSOCP_{\indexGridNode, \symPhaseN\symPhaseN}^{\realss}  }}
\newcommand{\VisqksocpGGreal}[0]{\textcolor{black}{\symVoltageSOCP_{\indexGridNode, \symPhaseG\symPhaseG}^{\realss}  }}

\newcommand{\VisqksocpAAreal}[0]{\textcolor{black}{\symVoltageSOCP_{\indexGridNode, \symPhaseA\symPhaseA}^{\realss}  }}
\newcommand{\VisqksocpABreal}[0]{\textcolor{black}{\symVoltageSOCP_{\indexGridNode, \symPhaseA\symPhaseB}^{\realss}  }}
\newcommand{\VisqksocpACreal}[0]{\textcolor{black}{\symVoltageSOCP_{\indexGridNode, \symPhaseA\symPhaseC}^{\realss}  }}

\newcommand{\VisqksocpANreal}[0]{\textcolor{black}{\symVoltageSOCP_{\indexGridNode, \symPhaseA\symPhaseN}^{\realss}  }}
\newcommand{\VisqksocpBNreal}[0]{\textcolor{black}{\symVoltageSOCP_{\indexGridNode, \symPhaseB\symPhaseN}^{\realss}  }}
\newcommand{\VisqksocpCNreal}[0]{\textcolor{black}{\symVoltageSOCP_{\indexGridNode, \symPhaseC\symPhaseN}^{\realss}  }}
\newcommand{\VisqksocpANimag}[0]{\textcolor{black}{\symVoltageSOCP_{\indexGridNode, \symPhaseA\symPhaseN}^{\imagss}  }}
\newcommand{\VisqksocpBNimag}[0]{\textcolor{black}{\symVoltageSOCP_{\indexGridNode, \symPhaseB\symPhaseN}^{\imagss}  }}
\newcommand{\VisqksocpCNimag}[0]{\textcolor{black}{\symVoltageSOCP_{\indexGridNode, \symPhaseC\symPhaseN}^{\imagss}  }}

\newcommand{\VisqksocpBBreal}[0]{\textcolor{black}{\symVoltageSOCP_{\indexGridNode, \symPhaseB\symPhaseB}^{\realss}  }}
\newcommand{\VisqksocpBCreal}[0]{\textcolor{black}{\symVoltageSOCP_{\indexGridNode, \symPhaseB\symPhaseC}^{\realss}  }}
\newcommand{\VisqksocpCCreal}[0]{\textcolor{black}{\symVoltageSOCP_{\indexGridNode, \symPhaseC\symPhaseC}^{\realss}  }}

\newcommand{\VisqksocpABimag}[0]{\textcolor{black}{\symVoltageSOCP_{\indexGridNode, \symPhaseA\symPhaseB}^{\imagss}  }}
\newcommand{\VisqksocpACimag}[0]{\textcolor{black}{\symVoltageSOCP_{\indexGridNode, \symPhaseA\symPhaseC}^{\imagss}  }}
\newcommand{\VisqksocpBCimag}[0]{\textcolor{black}{\symVoltageSOCP_{\indexGridNode, \symPhaseB\symPhaseC}^{\imagss}  }}

\newcommand{\VisqksocpAA}[0]{\textcolor{black}{\symVoltageSOCP_{\indexGridNode, \symPhaseA\symPhaseA} }}
\newcommand{\VisqksocpAB}[0]{\textcolor{\complexcolor}{\symVoltageSOCP_{\indexGridNode, \symPhaseA\symPhaseB}}}
\newcommand{\VisqksocpAC}[0]{\textcolor{\complexcolor}{\symVoltageSOCP_{\indexGridNode, \symPhaseA\symPhaseC}  }}
\newcommand{\VisqksocpBA}[0]{\textcolor{\complexcolor}{\symVoltageSOCP_{\indexGridNode, \symPhaseB\symPhaseA} }}
\newcommand{\VisqksocpBB}[0]{\textcolor{black}{\symVoltageSOCP_{\indexGridNode, \symPhaseB\symPhaseB} }}
\newcommand{\VisqksocpBC}[0]{\textcolor{\complexcolor}{\symVoltageSOCP_{\indexGridNode, \symPhaseB\symPhaseC}}}
\newcommand{\VisqksocpCA}[0]{\textcolor{\complexcolor}{\symVoltageSOCP_{\indexGridNode, \symPhaseC\symPhaseA}  }}
\newcommand{\VisqksocpCB}[0]{\textcolor{\complexcolor}{\symVoltageSOCP_{\indexGridNode, \symPhaseC\symPhaseB}  }}
\newcommand{\VisqksocpCC}[0]{\textcolor{black}{\symVoltageSOCP_{\indexGridNode, \symPhaseC\symPhaseC}  }}

\newcommand{\VijsqksocpAAreal}[0]{\textcolor{black}{\symVoltageSOCP_{\indexGridNode\indexGridNodeTwo, \symPhaseA\symPhaseA}^{\realss}  }}
\newcommand{\VijsqksocpABreal}[0]{\textcolor{black}{\symVoltageSOCP_{\indexGridNode\indexGridNodeTwo, \symPhaseA\symPhaseB}^{\realss}  }}
\newcommand{\VijsqksocpACreal}[0]{\textcolor{black}{\symVoltageSOCP_{\indexGridNode\indexGridNodeTwo, \symPhaseA\symPhaseC}^{\realss}  }}

\newcommand{\VijsqksocpBAreal}[0]{\textcolor{black}{\symVoltageSOCP_{\indexGridNode\indexGridNodeTwo, \symPhaseB\symPhaseA}^{\realss}  }}
\newcommand{\VijsqksocpBBreal}[0]{\textcolor{black}{\symVoltageSOCP_{\indexGridNode\indexGridNodeTwo, \symPhaseB\symPhaseB}^{\realss}  }}
\newcommand{\VijsqksocpBCreal}[0]{\textcolor{black}{\symVoltageSOCP_{\indexGridNode\indexGridNodeTwo, \symPhaseB\symPhaseC}^{\realss}  }}

\newcommand{\VijsqksocpCAreal}[0]{\textcolor{black}{\symVoltageSOCP_{\indexGridNode\indexGridNodeTwo, \symPhaseC\symPhaseA}^{\realss}  }}
\newcommand{\VijsqksocpCBreal}[0]{\textcolor{black}{\symVoltageSOCP_{\indexGridNode\indexGridNodeTwo, \symPhaseC\symPhaseB}^{\realss}  }}
\newcommand{\VijsqksocpCCreal}[0]{\textcolor{black}{\symVoltageSOCP_{\indexGridNode\indexGridNodeTwo, \symPhaseC\symPhaseC}^{\realss}  }}

\newcommand{\VjsqksocpABimag}[0]{\textcolor{black}{\symVoltageSOCP_{\indexGridNodeTwo, \symPhaseA\symPhaseB}^{\imagss}  }}
\newcommand{\VjsqksocpACimag}[0]{\textcolor{black}{\symVoltageSOCP_{\indexGridNodeTwo, \symPhaseA\symPhaseC}^{\imagss}  }}
\newcommand{\VjsqksocpBCimag}[0]{\textcolor{black}{\symVoltageSOCP_{\indexGridNodeTwo, \symPhaseB\symPhaseC}^{\imagss}  }}

\newcommand{\VjsqksocpABreal}[0]{\textcolor{black}{\symVoltageSOCP_{\indexGridNodeTwo, \symPhaseA\symPhaseB}^{\realss}  }}
\newcommand{\VjsqksocpACreal}[0]{\textcolor{black}{\symVoltageSOCP_{\indexGridNodeTwo, \symPhaseA\symPhaseC}^{\realss}  }}
\newcommand{\VjsqksocpBCreal}[0]{\textcolor{black}{\symVoltageSOCP_{\indexGridNodeTwo, \symPhaseB\symPhaseC}^{\realss}  }}

\newcommand{\VijsqksocpAAimag}[0]{\textcolor{black}{\symVoltageSOCP_{\indexGridNode\indexGridNodeTwo, \symPhaseA\symPhaseA}^{\imagss}  }}
\newcommand{\VijsqksocpABimag}[0]{\textcolor{black}{\symVoltageSOCP_{\indexGridNode\indexGridNodeTwo, \symPhaseA\symPhaseB}^{\imagss}  }}
\newcommand{\VijsqksocpACimag}[0]{\textcolor{black}{\symVoltageSOCP_{\indexGridNode\indexGridNodeTwo, \symPhaseA\symPhaseC}^{\imagss}  }}

\newcommand{\VijsqksocpBAimag}[0]{\textcolor{black}{\symVoltageSOCP_{\indexGridNode\indexGridNodeTwo, \symPhaseB\symPhaseA}^{\imagss}  }}
\newcommand{\VijsqksocpBBimag}[0]{\textcolor{black}{\symVoltageSOCP_{\indexGridNode\indexGridNodeTwo, \symPhaseB\symPhaseB}^{\imagss}  }}
\newcommand{\VijsqksocpBCimag}[0]{\textcolor{black}{\symVoltageSOCP_{\indexGridNode\indexGridNodeTwo, \symPhaseB\symPhaseC}^{\imagss}  }}

\newcommand{\VijsqksocpCAimag}[0]{\textcolor{black}{\symVoltageSOCP_{\indexGridNode\indexGridNodeTwo, \symPhaseC\symPhaseA}^{\imagss}  }}
\newcommand{\VijsqksocpCBimag}[0]{\textcolor{black}{\symVoltageSOCP_{\indexGridNode\indexGridNodeTwo, \symPhaseC\symPhaseB}^{\imagss}  }}
\newcommand{\VijsqksocpCCimag}[0]{\textcolor{black}{\symVoltageSOCP_{\indexGridNode\indexGridNodeTwo, \symPhaseC\symPhaseC}^{\imagss}  }}

\newcommand{\IijsqskA}[0]{\textcolor{black}{\symCurrent_{\indexGridLines\indexGridNode \indexGridNodeTwo,\seriesss, \symPhaseA}^{2}  }}
\newcommand{\IijsqskB}[0]{\textcolor{black}{\symCurrent_{\indexGridLines\indexGridNode \indexGridNodeTwo,\seriesss, \symPhaseB}^{2}  }}
\newcommand{\IijsqskC}[0]{\textcolor{black}{\symCurrent_{\indexGridLines\indexGridNode \indexGridNodeTwo,\seriesss, \symPhaseC}^{2}  }}

\newcommand{\Sijrated}[0]{\textcolor{\sizingcolor}{\symApparentPower_{\indexGridLines\indexGridNode \indexGridNodeTwo}^{\ratedss} }}
\newcommand{\Sjirated}[0]{\textcolor{\sizingcolor}{\symApparentPower_{\indexGridLines\indexGridNodeTwo \indexGridNode }^{\ratedss} }}
\newcommand{\SijratedA}[0]{\textcolor{\sizingcolor}{\symApparentPower_{\indexGridLines\indexGridNode \indexGridNodeTwo,\symPhaseA}^{\ratedss} }}
\newcommand{\SjiratedA}[0]{\textcolor{\sizingcolor}{\symApparentPower_{\indexGridLines\indexGridNodeTwo \indexGridNode ,\symPhaseA}^{\ratedss} }}
\newcommand{\SijratedB}[0]{\textcolor{\sizingcolor}{\symApparentPower_{\indexGridLines\indexGridNode \indexGridNodeTwo,\symPhaseB}^{\ratedss} }}
\newcommand{\SjiratedB}[0]{\textcolor{\sizingcolor}{\symApparentPower_{\indexGridLines\indexGridNodeTwo \indexGridNode ,\symPhaseB}^{\ratedss} }}
\newcommand{\SijratedC}[0]{\textcolor{\sizingcolor}{\symApparentPower_{\indexGridLines\indexGridNode \indexGridNodeTwo,\symPhaseC}^{\ratedss} }}
\newcommand{\SijratedN}[0]{\textcolor{\sizingcolor}{\symApparentPower_{\indexGridLines\indexGridNode \indexGridNodeTwo,\symPhaseN}^{\ratedss} }}
\newcommand{\SjiratedC}[0]{\textcolor{\sizingcolor}{\symApparentPower_{\indexGridLines\indexGridNodeTwo \indexGridNode ,\symPhaseC}^{\ratedss} }}

\newcommand{\VikSDPmin}[0]{\textcolor{\paramcolor}{\mathbf{\symVoltage}_{\indexGridNode }^{\minss}  }}
\newcommand{\VikSDPmax}[0]{\textcolor{\paramcolor}{\mathbf{\symVoltage}_{\indexGridNode }^{\maxss}  }}
\newcommand{\VikSDPref}[0]{\textcolor{\paramcolor}{\mathbf{\symVoltage}_{\indexGridNode }^{\refss}  }}

\newcommand{\VikSDPwye}[0]{\textcolor{\complexcolor}{\mathbf{\symVoltage}_{\indexGridNode }^{\text{pn}}  }}
\newcommand{\VjkSDPwye}[0]{\textcolor{\complexcolor}{\mathbf{\symVoltage}_{\indexGridNodeTwo }^{\text{pn}}  }}

\newcommand{\VikSDPminwye}[0]{\textcolor{\paramcolor}{\mathbf{\symVoltage}_{\indexGridNode }^{\text{pn},\minss}  }}
\newcommand{\VikSDPmaxwye}[0]{\textcolor{\paramcolor}{\mathbf{\symVoltage}_{\indexGridNode }^{\text{pn},\maxss}  }}

\newcommand{\VikSDPmindelta}[0]{\textcolor{\paramcolor}{\mathbf{\symVoltage}_{\indexGridNode }^{\Delta\minss}  }}
\newcommand{\VikSDPmaxdelta}[0]{\textcolor{\paramcolor}{\mathbf{\symVoltage}_{\indexGridNode }^{\Delta\maxss}  }}

\newcommand{\VjkSDPmin}[0]{\textcolor{\paramcolor}{\mathbf{\symVoltage}_{\indexGridNodeTwo }^{\minss}  }}
\newcommand{\VjkSDPmax}[0]{\textcolor{\paramcolor}{\mathbf{\symVoltage}_{\indexGridNodeTwo }^{\maxss}  }}

\newcommand{\VjkSDPseq}[0]{\textcolor{\complexcolor}{\mathbf{\symVoltage}_{\indexGridNodeTwo }^{\fortescuess}  }}

\newcommand{\VikSDPacc}[0]{\textcolor{\complexcolor}{\mathbf{\bar{\symVoltage}}_{\indexGridNode }  }}
\newcommand{\VikSDPaccreal}[0]{\textcolor{black}{\mathbf{\bar{\symVoltage}}_{\indexGridNode }^{\realss}  }}
\newcommand{\VikSDPaccimag}[0]{\textcolor{black}{\mathbf{\bar{\symVoltage}}_{\indexGridNode }^{\imagss}  }}

\newcommand{\VikSDPdelta}[0]{\textcolor{\complexcolor}{\mathbf{\symVoltage}_{\indexGridNode }^{\Delta}  }}
\newcommand{\VikSDPdeltareal}[0]{\textcolor{black}{\mathbf{\symVoltage}_{\indexGridNode }^{\Delta, \realss}  }}
\newcommand{\VikSDPdeltaimag}[0]{\textcolor{black}{\mathbf{\symVoltage}_{\indexGridNode }^{\Delta, \imagss}  }}

\newcommand{\VikSDPseq}[0]{\textcolor{\complexcolor}{\mathbf{\symVoltage}_{\indexGridNode }^{\fortescuess}  }}
\newcommand{\VikSDP}[0]{\textcolor{\complexcolor}{\mathbf{\symVoltage}_{\indexGridNode }^{}  }}
\newcommand{\VikSDPreal}[0]{\textcolor{black}{\mathbf{\symVoltage}_{\indexGridNode }^{\realss}  }}
\newcommand{\VikSDPimag}[0]{\textcolor{black}{\mathbf{\symVoltage}_{\indexGridNode }^{\imagss}  }}
\newcommand{\VjkSDPreal}[0]{\textcolor{black}{\mathbf{\symVoltage}_{\indexGridNodeTwo }^{\realss}  }}
\newcommand{\VjkSDPimag}[0]{\textcolor{black}{\mathbf{\symVoltage}_{\indexGridNodeTwo }^{\imagss}  }}
\newcommand{\VikSDPH}[0]{\textcolor{\complexcolor}{\mathbf{\symVoltage}_{\indexGridNode }^{\hermitiantranspose}  }}
\newcommand{\VjkSDP}[0]{\textcolor{\complexcolor}{\mathbf{\symVoltage}_{\indexGridNodeTwo }^{}  }}
\newcommand{\VzkSDP}[0]{\textcolor{\complexcolor}{\mathbf{\symVoltage}_{z }^{}  }}
\newcommand{\VjkSDPH}[0]{\textcolor{\complexcolor}{\mathbf{\symVoltage}_{\indexGridNodeTwo }^{\hermitiantranspose}  }}
\newcommand{\Vik}[0]{\textcolor{black}{\symVoltage_{\indexGridNode }^{}  }}
\newcommand{\Vjk}[0]{\textcolor{black}{\symVoltage_{\indexGridNodeTwo }^{}  }}
\newcommand{\Vikdelta}[0]{\textcolor{black}{\Delta\symVoltage_{\indexGridNode }^{}  }}
\newcommand{\Vjkdelta}[0]{\textcolor{black}{\Delta\symVoltage_{\indexGridNodeTwo }^{}  }}
\newcommand{\Viacck}[0]{\textcolor{black}{\symVoltage_{\indexGridNode }^{'}  }}
\newcommand{\Vjacck}[0]{\textcolor{black}{\symVoltage_{\indexGridNodeTwo }^{'}  }}
\newcommand{\Viacckl}[0]{\textcolor{black}{\symVoltage_{\indexGridNode }^{\star}  }}
\newcommand{\Vjacckl}[0]{\textcolor{black}{\symVoltage_{\indexGridNodeTwo }^{\star}  }}

\newcommand{\VikSDPmag}[0]{\textcolor{black}{\mathbf{\symVoltage}_{\indexGridNode }^{\text{mag}}  }}
\newcommand{\VikSDPang}[0]{\textcolor{black}{\mathbf{\Theta}_{\indexGridNode }  }}

\newcommand{\VjkSDPmag}[0]{\textcolor{black}{\mathbf{\symVoltage}_{\indexGridNodeTwo }^{\text{mag}}  }}
\newcommand{\VjkSDPang}[0]{\textcolor{black}{\mathbf{\Theta}_{\indexGridNodeTwo }  }}

\newcommand{\VukrefSDP}[0]{\textcolor{\complexparamcolor}{\mathbf{\symVoltage}_{\indexUnit }^{\refss}  }}

\newcommand{\IurefSDP}[0]{\textcolor{\complexparamcolor}{\mathbf{\symCurrent}_{\indexGridNode }^{\refss}  }}

\newcommand{\VikrefSDP}[0]{\textcolor{\complexparamcolor}{\mathbf{\symVoltage}_{\indexGridNode }^{\refss}  }}
\newcommand{\VikrefSDPreal}[0]{\textcolor{\paramcolor}{\mathbf{\symVoltage}_{\indexGridNode }^{\refss,\realss}  }}
\newcommand{\VikrefSDPimag}[0]{\textcolor{\paramcolor}{\mathbf{\symVoltage}_{\indexGridNode }^{\refss,\imagss}  }}
\newcommand{\VikdeltaSDP}[0]{\textcolor{\complexcolor}{\mathbf{\symVoltage}_{\indexGridNode }^{\Delta}  }}
\newcommand{\VikdeltaSDPreal}[0]{\textcolor{black}{\mathbf{\symVoltage}_{\indexGridNode }^{\Delta,\realss}  }}
\newcommand{\VikdeltaSDPimag}[0]{\textcolor{black}{\mathbf{\symVoltage}_{\indexGridNode }^{\Delta,\imagss}  }}

\newcommand{\VsqikdeltaSDP}[0]{\textcolor{\complexcolor}{\mathbf{\symVoltageSOCP}_{\indexGridNode }^{\Delta}  }}

\newcommand{\VikaccSDP}[0]{\textcolor{\complexcolor}{\mathbf{\symVoltage}_{\indexGridNode }^{'}  }}
\newcommand{\VjkaccSDP}[0]{\textcolor{\complexcolor}{\mathbf{\symVoltage}_{\indexGridNodeTwo }^{'}  }}
\newcommand{\Cdelta}[0]{\textcolor{\complexcolor}{\mathbf{\symConnectivity}_{\Delta }^{}  }}
\newcommand{\Czigzag}[0]{\textcolor{\complexcolor}{\mathbf{\symConnectivity}_{Z }^{}  }}
\newcommand{\Cwye}[0]{\textcolor{\complexcolor}{\mathbf{\symConnectivity}_{Y }^{}  }}

\newcommand{\CijkSDP}[0]{\textcolor{\complexcolor}{\mathbf{\symConnectivity}_{\indexGridNode\indexGridNodeTwo }^{}  }}
\newcommand{\CjikSDP}[0]{\textcolor{\complexcolor}{\mathbf{\symConnectivity}_{\indexGridNodeTwo\indexGridNode }^{}  }}

\newcommand{\Visqk}[0]{\textcolor{black}{\symVoltage_{\indexGridNode }^{2}  }}
\newcommand{\Vjsqk}[0]{\textcolor{black}{\symVoltage_{\indexGridNodeTwo }^{2}  }}

\newcommand{\Vjh}[0]{\textcolor{\complexcolor}{\mathbf{\symVoltage}_{\indexGridNodeTwo ,\symHarmonic}   }}

\newcommand{\VjAh}[0]{\textcolor{\complexcolor}{\symVoltage_{\indexGridNodeTwo,\symPhaseA ,\symHarmonic}  }}
\newcommand{\VjBh}[0]{\textcolor{\complexcolor}{\symVoltage_{\indexGridNodeTwo,\symPhaseB,\symHarmonic }  }}
\newcommand{\VjCh}[0]{\textcolor{\complexcolor}{\symVoltage_{\indexGridNodeTwo,\symPhaseC,\symHarmonic }   }}

\newcommand{\Vih}[0]{\textcolor{\complexcolor}{\mathbf{\symVoltage}_{\indexGridNode ,\symHarmonic}   }}

\newcommand{\ViAh}[0]{\textcolor{\complexcolor}{\symVoltage_{\indexGridNode,\symPhaseA ,\symHarmonic}  }}
\newcommand{\ViBh}[0]{\textcolor{\complexcolor}{\symVoltage_{\indexGridNode,\symPhaseB,\symHarmonic }  }}
\newcommand{\ViCh}[0]{\textcolor{\complexcolor}{\symVoltage_{\indexGridNode,\symPhaseC,\symHarmonic }   }}

\newcommand{\VmagiAh}[0]{\symVoltage_{\indexGridNode,\symPhaseA ,\symHarmonic}^{\text{mag}}  }
\newcommand{\VmagiBh}[0]{\symVoltage_{\indexGridNode,\symPhaseB,\symHarmonic }^{\text{mag}}  }
\newcommand{\VmagiCh}[0]{\symVoltage_{\indexGridNode,\symPhaseC,\symHarmonic }^{\text{mag}}  }

\newcommand{\VangiAh}[0]{\symVoltageAngle_{\indexGridNode,\symPhaseA ,\symHarmonic}  }
\newcommand{\VangiBh}[0]{\symVoltageAngle_{\indexGridNode,\symPhaseB,\symHarmonic }  }
\newcommand{\VangiCh}[0]{\symVoltageAngle_{\indexGridNode,\symPhaseC,\symHarmonic }  }

\newcommand{\VmagiAhmax}[0]{\textcolor{\paramcolor}{\symVoltage_{\indexGridNode,\symPhaseA ,\symHarmonic}^{\text{max}}  }}
\newcommand{\VmagiBhmax}[0]{\textcolor{\paramcolor}{\symVoltage_{\indexGridNode,\symPhaseB,\symHarmonic }^{\text{max}}  }}
\newcommand{\VmagiChmax}[0]{\textcolor{\paramcolor}{\symVoltage_{\indexGridNode,\symPhaseC,\symHarmonic }^{\text{max}}  }}

\newcommand{\VmagiAthd}[0]{\textcolor{\paramcolor}{\symVoltage_{\indexGridNode,\symPhaseA }^{\text{THD}}  }}
\newcommand{\VmagiBthd}[0]{\textcolor{\paramcolor}{\symVoltage_{\indexGridNode,\symPhaseB }^{\text{THD}}  }}
\newcommand{\VmagiCthd}[0]{\textcolor{\paramcolor}{\symVoltage_{\indexGridNode,\symPhaseC }^{\text{THD}}  }}

\newcommand{\VikAN}[0]{\textcolor{\complexcolor}{\symVoltage_{\indexGridNode,\symPhaseA\symPhaseN }^{}  }}
\newcommand{\VikBN}[0]{\textcolor{\complexcolor}{\symVoltage_{\indexGridNode,\symPhaseB \symPhaseN}^{}  }}
\newcommand{\VikCN}[0]{\textcolor{\complexcolor}{\symVoltage_{\indexGridNode,\symPhaseC\symPhaseN }^{}  }}

\newcommand{\VikAB}[0]{\textcolor{\complexcolor}{\symVoltage_{\indexGridNode,\symPhaseA\symPhaseB }^{\Delta}  }}
\newcommand{\VikBC}[0]{\textcolor{\complexcolor}{\symVoltage_{\indexGridNode,\symPhaseB\symPhaseC }^{\Delta}  }}
\newcommand{\VikCA}[0]{\textcolor{\complexcolor}{\symVoltage_{\indexGridNode,\symPhaseC\symPhaseA }^{\Delta}  }}

\newcommand{\Vikzero}[0]{\textcolor{\complexcolor}{\symVoltage_{\indexGridNode,0 }^{}  }}
\newcommand{\Vikone}[0]{\textcolor{\complexcolor}{\symVoltage_{\indexGridNode,1 }^{}  }}
\newcommand{\Viktwo}[0]{\textcolor{\complexcolor}{\symVoltage_{\indexGridNode,2 }^{}  }}

\newcommand{\VikAreal}[0]{\textcolor{black}{\symVoltage_{\indexGridNode,\symPhaseA }^{\realss}  }}
\newcommand{\VikBreal}[0]{\textcolor{black}{\symVoltage_{\indexGridNode,\symPhaseB }^{\realss}  }}
\newcommand{\VikCreal}[0]{\textcolor{black}{\symVoltage_{\indexGridNode,\symPhaseC }^{\realss}  }}
\newcommand{\VikAimag}[0]{\textcolor{black}{\symVoltage_{\indexGridNode,\symPhaseA }^{\imagss}  }}
\newcommand{\VikBimag}[0]{\textcolor{black}{\symVoltage_{\indexGridNode,\symPhaseB }^{\imagss}  }}
\newcommand{\VikCimag}[0]{\textcolor{black}{\symVoltage_{\indexGridNode,\symPhaseC }^{\imagss}  }}

\newcommand{\VikANreal}[0]{\textcolor{black}{\symVoltage_{\indexGridNode,\symPhaseA\symPhaseN }^{\realss}  }}
\newcommand{\VikBNreal}[0]{\textcolor{black}{\symVoltage_{\indexGridNode,\symPhaseB\symPhaseN }^{\realss}  }}
\newcommand{\VikCNreal}[0]{\textcolor{black}{\symVoltage_{\indexGridNode,\symPhaseC\symPhaseN }^{\realss}  }}
\newcommand{\VikNreal}[0]{\textcolor{black}{\symVoltage_{\indexGridNode,\symPhaseN }^{\realss}  }}

\newcommand{\VikANimag}[0]{\textcolor{black}{\symVoltage_{\indexGridNode,\symPhaseA }^{\imagss}  }}
\newcommand{\VikBNimag}[0]{\textcolor{black}{\symVoltage_{\indexGridNode,\symPhaseB }^{\imagss}  }}
\newcommand{\VikCNimag}[0]{\textcolor{black}{\symVoltage_{\indexGridNode,\symPhaseC }^{\imagss}  }}
\newcommand{\VikNimag}[0]{\textcolor{black}{\symVoltage_{\indexGridNode,\symPhaseN }^{\imagss}  }}

\newcommand{\VikpNreal}[0]{\textcolor{black}{\symVoltage_{\indexGridNode,\indexPhases }^{\realss}  }}
\newcommand{\VikpNimag}[0]{\textcolor{black}{\symVoltage_{\indexGridNode,\indexPhases }^{\imagss}  }}
\newcommand{\VikppNreal}[0]{\textcolor{black}{\symVoltage_{\indexGridNode,\indexPhasesTwo }^{\realss}  }}
\newcommand{\VikppNimag}[0]{\textcolor{black}{\symVoltage_{\indexGridNode,\indexPhasesTwo }^{\imagss}  }}
\newcommand{\VjkppNreal}[0]{\textcolor{black}{\symVoltage_{\indexGridNodeTwo,\indexPhasesTwo }^{\realss}  }}
\newcommand{\VjkppNimag}[0]{\textcolor{black}{\symVoltage_{\indexGridNodeTwo,\indexPhasesTwo }^{\imagss}  }}

\newcommand{\VikpNmag}[0]{\textcolor{black}{\symVoltage_{\indexGridNode,\indexPhases }^{\text{mag}}  }}
\newcommand{\VikppNmag}[0]{\textcolor{black}{\symVoltage_{\indexGridNode,\indexPhasesTwo }^{\text{mag}}  }}
\newcommand{\VjkppNmag}[0]{\textcolor{black}{\symVoltage_{\indexGridNodeTwo,\indexPhasesTwo }^{\text{mag}}  }}

\newcommand{\VikANmag}[0]{\textcolor{black}{\symVoltage_{\indexGridNode,\symPhaseA }^{\text{mag}}  }}
\newcommand{\VikBNmag}[0]{\textcolor{black}{\symVoltage_{\indexGridNode,\symPhaseB }^{\text{mag}}  }}
\newcommand{\VikCNmag}[0]{\textcolor{black}{\symVoltage_{\indexGridNode,\symPhaseC }^{\text{mag}}  }}
\newcommand{\Vikmag}[0]{\textcolor{black}{\symVoltage_{\indexGridNode }^{\text{mag}}  }}

\newcommand{\VjkANmag}[0]{\textcolor{black}{\symVoltage_{\indexGridNodeTwo,\symPhaseA }^{\text{mag}}  }}
\newcommand{\VjkBNmag}[0]{\textcolor{black}{\symVoltage_{\indexGridNodeTwo,\symPhaseB }^{\text{mag}}  }}
\newcommand{\VjkCNmag}[0]{\textcolor{black}{\symVoltage_{\indexGridNodeTwo,\symPhaseC }^{\text{mag}}  }}

\newcommand{\VikANangle}[0]{\textcolor{black}{\symVoltageAngle_{\indexGridNode,\symPhaseA }^{}  }}
\newcommand{\VikBNangle}[0]{\textcolor{black}{\symVoltageAngle_{\indexGridNode,\symPhaseB }^{}  }}
\newcommand{\VikCNangle}[0]{\textcolor{black}{\symVoltageAngle_{\indexGridNode,\symPhaseC }^{}  }}
\newcommand{\VikNangle}[0]{\textcolor{black}{\symVoltageAngle_{\indexGridNode,\symPhaseN }^{}  }}

\newcommand{\VjkANangle}[0]{\textcolor{black}{\symVoltageAngle_{\indexGridNodeTwo,\symPhaseA }^{}  }}
\newcommand{\VjkBNangle}[0]{\textcolor{black}{\symVoltageAngle_{\indexGridNodeTwo,\symPhaseB }^{}  }}
\newcommand{\VjkCNangle}[0]{\textcolor{black}{\symVoltageAngle_{\indexGridNodeTwo,\symPhaseC }^{}  }}
\newcommand{\VjkNangle}[0]{\textcolor{black}{\symVoltageAngle_{\indexGridNodeTwo,\symPhaseN }^{}  }}

\newcommand{\VikpNangle}[0]{\textcolor{black}{\symVoltageAngle_{\indexGridNode,\indexPhases }^{}  }}
\newcommand{\VikhNangle}[0]{\textcolor{black}{\symVoltageAngle_{\indexGridNode,\indexPhasesTwo }^{}  }}
\newcommand{\VjkpNangle}[0]{\textcolor{black}{\symVoltageAngle_{\indexGridNodeTwo,\indexPhases }^{}  }}
\newcommand{\VjkhNangle}[0]{\textcolor{black}{\symVoltageAngle_{\indexGridNodeTwo,\indexPhasesTwo }^{}  }}

\newcommand{\VikrefAN}[0]{\textcolor{\complexparamcolor}{\symVoltage_{\indexGridNode,\symPhaseA }^{\refss}  }}
\newcommand{\VikrefBN}[0]{\textcolor{\complexparamcolor}{\symVoltage_{\indexGridNode,\symPhaseB }^{\refss}  }}
\newcommand{\VikrefCN}[0]{\textcolor{\complexparamcolor}{\symVoltage_{\indexGridNode,\symPhaseC  }^{\refss}  }}

\newcommand{\VikANref}[0]{\textcolor{\complexparamcolor}{\symVoltage_{\indexGridNode,\symPhaseA }^{\refss}  }}
\newcommand{\VikBNref}[0]{\textcolor{\complexparamcolor}{\symVoltage_{\indexGridNode,\symPhaseB }^{\refss}  }}
\newcommand{\VikCNref}[0]{\textcolor{\complexparamcolor}{\symVoltage_{\indexGridNode,\symPhaseC }^{\refss}  }}
\newcommand{\VikNref}[0]{\textcolor{\complexparamcolor}{\symVoltage_{\indexGridNode,n }^{\refss}  }}

\newcommand{\thetaikANref}[0]{\textcolor{\paramcolor}{\symVoltageAngle_{\indexGridNode,\symPhaseA }^{\refss}  }}
\newcommand{\thetaikBNref}[0]{\textcolor{\paramcolor}{\symVoltageAngle_{\indexGridNode,\symPhaseB }^{\refss}  }}
\newcommand{\thetaikCNref}[0]{\textcolor{\paramcolor}{\symVoltageAngle_{\indexGridNode,\symPhaseC }^{\refss}  }}
\newcommand{\thetaikNref}[0]{\textcolor{\paramcolor}{\symVoltageAngle_{\indexGridNode,n }^{\refss}  }}

\newcommand{\VikAmin}[0]{\textcolor{\paramcolor}{\symVoltage_{\indexGridNode,\symPhaseA }^{\minss}  }}
\newcommand{\VikBmin}[0]{\textcolor{\paramcolor}{\symVoltage_{\indexGridNode,\symPhaseB }^{\minss}  }}
\newcommand{\VikCmin}[0]{\textcolor{\paramcolor}{\symVoltage_{\indexGridNode,\symPhaseC }^{\minss}  }}
\newcommand{\VikAmax}[0]{\textcolor{\paramcolor}{\symVoltage_{\indexGridNode,\symPhaseA }^{\maxss}  }}
\newcommand{\VikBmax}[0]{\textcolor{\paramcolor}{\symVoltage_{\indexGridNode,\symPhaseB }^{\maxss}  }}
\newcommand{\VikCmax}[0]{\textcolor{\paramcolor}{\symVoltage_{\indexGridNode,\symPhaseC }^{\maxss}  }}

\newcommand{\VikANmin}[0]{\textcolor{\paramcolor}{\symVoltage_{\indexGridNode,\symPhaseA\symPhaseN }^{\minss}  }}
\newcommand{\VikBNmin}[0]{\textcolor{\paramcolor}{\symVoltage_{\indexGridNode,\symPhaseB\symPhaseN }^{\minss}  }}
\newcommand{\VikCNmin}[0]{\textcolor{\paramcolor}{\symVoltage_{\indexGridNode,\symPhaseC\symPhaseN }^{\minss}  }}
\newcommand{\VikNmin}[0]{\textcolor{\paramcolor}{\symVoltage_{\indexGridNode,\symPhaseN }^{\minss}  }}
\newcommand{\VikANmax}[0]{\textcolor{\paramcolor}{\symVoltage_{\indexGridNode,\symPhaseA\symPhaseN }^{\maxss}  }}
\newcommand{\VikBNmax}[0]{\textcolor{\paramcolor}{\symVoltage_{\indexGridNode,\symPhaseB\symPhaseN }^{\maxss}  }}
\newcommand{\VikCNmax}[0]{\textcolor{\paramcolor}{\symVoltage_{\indexGridNode,\symPhaseC\symPhaseN }^{\maxss}  }}
\newcommand{\VikNmax}[0]{\textcolor{\paramcolor}{\symVoltage_{\indexGridNode,\symPhaseN }^{\maxss}  }}

\newcommand{\VikABmax}[0]{\textcolor{\paramcolor}{\symVoltage_{\indexGridNode,\symPhaseA\symPhaseB }^{\maxss}  }}
\newcommand{\VikBCmax}[0]{\textcolor{\paramcolor}{\symVoltage_{\indexGridNode,\symPhaseB\symPhaseC }^{\maxss}  }}
\newcommand{\VikCAmax}[0]{\textcolor{\paramcolor}{\symVoltage_{\indexGridNode,\symPhaseC\symPhaseA }^{\maxss}  }}
\newcommand{\VikABmin}[0]{\textcolor{\paramcolor}{\symVoltage_{\indexGridNode,\symPhaseA\symPhaseB }^{\minss}  }}
\newcommand{\VikBCmin}[0]{\textcolor{\paramcolor}{\symVoltage_{\indexGridNode,\symPhaseB\symPhaseC }^{\minss}  }}
\newcommand{\VikCAmin}[0]{\textcolor{\paramcolor}{\symVoltage_{\indexGridNode,\symPhaseC\symPhaseA }^{\minss}  }}

\newcommand{\VjkAN}[0]{\textcolor{\complexcolor}{\symVoltage_{\indexGridNodeTwo,\symPhaseA\symPhaseN }^{}  }}
\newcommand{\VjkBN}[0]{\textcolor{\complexcolor}{\symVoltage_{\indexGridNodeTwo,\symPhaseB\symPhaseN }^{}  }}
\newcommand{\VjkCN}[0]{\textcolor{\complexcolor}{\symVoltage_{\indexGridNodeTwo,\symPhaseC\symPhaseN }^{}  }}
\newcommand{\VjkANmax}[0]{\textcolor{\paramcolor}{\symVoltage_{\indexGridNodeTwo,\symPhaseA }^{\maxss}  }}

\newcommand{\VjkA}[0]{\textcolor{\complexcolor}{\symVoltage_{\indexGridNodeTwo,\symPhaseA }^{}  }}
\newcommand{\VjkB}[0]{\textcolor{\complexcolor}{\symVoltage_{\indexGridNodeTwo,\symPhaseB }^{}  }}
\newcommand{\VjkC}[0]{\textcolor{\complexcolor}{\symVoltage_{\indexGridNodeTwo,\symPhaseC }^{}  }}
\newcommand{\VjkN}[0]{\textcolor{\complexcolor}{\symVoltage_{\indexGridNodeTwo,\symPhaseN }^{}  }}
\newcommand{\VjkG}[0]{\textcolor{\complexcolor}{\symVoltage_{\indexGridNodeTwo,\symPhaseG }^{}  }}

\newcommand{\VikA}[0]{\textcolor{\complexcolor}{\symVoltage_{\indexGridNode,\symPhaseA }^{}  }}
\newcommand{\VikB}[0]{\textcolor{\complexcolor}{\symVoltage_{\indexGridNode,\symPhaseB }^{}  }}
\newcommand{\VikC}[0]{\textcolor{\complexcolor}{\symVoltage_{\indexGridNode,\symPhaseC }^{}  }}
\newcommand{\VikN}[0]{\textcolor{\complexcolor}{\symVoltage_{\indexGridNode,\symPhaseN }^{}  }}
\newcommand{\VikG}[0]{\textcolor{\complexcolor}{\symVoltage_{\indexGridNode,\symPhaseG }^{}  }}

\newcommand{\Vikmin}[0]{\textcolor{\boundscolor}{\symVoltage^{\minss}_{\indexGridNode }  }}
\newcommand{\Vjkmin}[0]{\textcolor{\boundscolor}{\symVoltage^{\minss}_{\indexGridNodeTwo }  }}
\newcommand{\Vikmax}[0]{\textcolor{\boundscolor}{\symVoltage^{\maxss}_{\indexGridNode }  }}
\newcommand{\Vjkmax}[0]{\textcolor{\boundscolor}{\symVoltage^{\maxss}_{\indexGridNodeTwo }  }}
\newcommand{\bigMikmax}[0]{\textcolor{\boundscolor}{M^{}_{\indexGridNode }  }}
\newcommand{\bigMjkmax}[0]{\textcolor{\boundscolor}{M^{}_{\indexGridNodeTwo }  }}
\newcommand{\Virated}[0]{\textcolor{\sizingcolor}{\symVoltage^{\ratedss}_{\indexGridNode } }}
\newcommand{\Vjrated}[0]{\textcolor{\sizingcolor}{\symVoltage^{\ratedss}_{\indexGridNodeTwo } }}
\newcommand{\Vijrated}[0]{\textcolor{\sizingcolor}{\symVoltage^{\ratedss}_{\indexGridNode \indexGridNodeTwo} }}
\newcommand{\Vjirated}[0]{\textcolor{\sizingcolor}{\symVoltage^{\ratedss}_{\indexGridNodeTwo \indexGridNode } }}

\newcommand{\ViratedAN}[0]{\textcolor{\sizingcolor}{\symVoltage^{\ratedss}_{\indexGridNode,\symPhaseA\symPhaseN } }}
\newcommand{\ViratedBN}[0]{\textcolor{\sizingcolor}{\symVoltage^{\ratedss}_{\indexGridNode,\symPhaseB\symPhaseN } }}
\newcommand{\ViratedCN}[0]{\textcolor{\sizingcolor}{\symVoltage^{\ratedss}_{\indexGridNode,\symPhaseC\symPhaseN } }}

\newcommand{\VijratedAN}[0]{\textcolor{\sizingcolor}{\symVoltage^{\ratedss}_{\indexGridNode \indexGridNodeTwo,\symPhaseA\symPhaseN} }}
\newcommand{\VjiratedAN}[0]{\textcolor{\sizingcolor}{\symVoltage^{\ratedss}_{\indexGridNodeTwo \indexGridNode,\symPhaseA\symPhaseN } }}
\newcommand{\VijratedBN}[0]{\textcolor{\sizingcolor}{\symVoltage^{\ratedss}_{\indexGridNode \indexGridNodeTwo,\symPhaseB\symPhaseN} }}
\newcommand{\VjiratedBN}[0]{\textcolor{\sizingcolor}{\symVoltage^{\ratedss}_{\indexGridNodeTwo \indexGridNode ,\symPhaseB\symPhaseN} }}
\newcommand{\VijratedCN}[0]{\textcolor{\sizingcolor}{\symVoltage^{\ratedss}_{\indexGridNode \indexGridNodeTwo,\symPhaseC\symPhaseN} }}
\newcommand{\VjiratedCN}[0]{\textcolor{\sizingcolor}{\symVoltage^{\ratedss}_{\indexGridNodeTwo \indexGridNode,\symPhaseC\symPhaseN } }}

\newcommand{\Visqksocplin}[0]{\textcolor{black}{\symVoltage_{\indexGridNode }^{\linss}  }}

\newcommand{\VisqkSDPseq}[0]{\textcolor{\complexcolor}{\mathbf{\symVoltageSOCP}_{\indexGridNode }^{\fortescuess} }}

\newcommand{\VusqkSDP}[0]{\textcolor{\complexcolor}{\mathbf{\symVoltageSOCP}_{\indexUnit }  }}

\newcommand{\VisqkSDPdelta}[0]{\textcolor{\complexcolor}{\mathbf{\symVoltageSOCP}_{\indexGridNode }^{\Delta}  }}
\newcommand{\VisqkSDPdeltareal}[0]{\textcolor{black}{\mathbf{\symVoltageSOCP}_{\indexGridNode }^{\Delta \realss}  }}
\newcommand{\VisqkSDPdeltaimag}[0]{\textcolor{black}{\mathbf{\symVoltageSOCP}_{\indexGridNode }^{\Delta \imagss}  }}

\newcommand{\VisqkSDP}[0]{\textcolor{\complexcolor}{\mathbf{\symVoltageSOCP}_{\indexGridNode }  }}
\newcommand{\VjsqkSDP}[0]{\textcolor{\complexcolor}{\mathbf{\symVoltageSOCP}_{\indexGridNodeTwo }  }}
\newcommand{\VisqkSDPreal}[0]{\textcolor{black}{\mathbf{\symVoltageSOCP}_{\indexGridNode }^{\realss}  }}
\newcommand{\VisqkSDPimag}[0]{\textcolor{black}{\mathbf{\symVoltageSOCP}_{\indexGridNode }^{\imagss}  }}
\newcommand{\VjsqkSDPreal}[0]{\textcolor{black}{\mathbf{\symVoltageSOCP}_{\indexGridNodeTwo }^{\realss}  }}
\newcommand{\VjsqkSDPimag}[0]{\textcolor{black}{\mathbf{\symVoltageSOCP}_{\indexGridNodeTwo }^{\imagss}  }}
\newcommand{\VisqkaccSDP}[0]{\textcolor{\complexcolor}{\mathbf{\symVoltageSOCP}_{\indexGridNode }^{'}  }}

\newcommand{\VzsqkSDP}[0]{\textcolor{\complexcolor}{\mathbf{\symVoltageSOCP}_{z }  }}
\newcommand{\VsqkSDP}[0]{\textcolor{\complexcolor}{\mathbf{M}  }}
\newcommand{\VsqkSDPtworeal}[0]{\textcolor{black}{\mathbf{M}^{2\realss}  }}
\newcommand{\VsqkSDPreal}[0]{\textcolor{black}{\mathbf{M}^{\realss}  }}
\newcommand{\VsqkSDPimag}[0]{\textcolor{black}{\mathbf{M}^{\imagss}  }}
\newcommand{\MisqkSDP}[0]{\textcolor{\complexcolor}{\mathbf{M}_i  }}
\newcommand{\MisqkSDPtworeal}[0]{\textcolor{black}{\mathbf{M}^{2\realss}_{\indexGridNode}  }}

\newcommand{\VbpsqkSDP}[0]{\textcolor{\complexcolor}{\mathbf{M}_{\indexGridNode\indexGridNodeTwo}  }}
\newcommand{\VbpsqkSDPtworeal}[0]{\textcolor{black}{\mathbf{M}^{2\realss}_{\indexGridNode\indexGridNodeTwo}  }}
\newcommand{\VbpsqkSDPreal}[0]{\textcolor{black}{\mathbf{M}_{\indexGridNode\indexGridNodeTwo} ^{\realss}  }}
\newcommand{\VbpsqkSDPimag}[0]{\textcolor{black}{\mathbf{M}_{\indexGridNode\indexGridNodeTwo} ^{\imagss}  }}

\newcommand{\VlinesqkSDP}[0]{\textcolor{\complexcolor}{\mathbf{M}_{\indexGridLines\indexGridNode\indexGridNodeTwo}  }}
\newcommand{\VlinesqkSDPreal}[0]{\textcolor{black}{\mathbf{M}_{\indexGridLines\indexGridNode\indexGridNodeTwo}^{\realss}  }}
\newcommand{\VlinesqkSDPrealtwo}[0]{\textcolor{black}{\mathbf{M}_{\indexGridLines\indexGridNode\indexGridNodeTwo}^{2\realss}  }}

\newcommand{\VlinetosqkSDP}[0]{\textcolor{\complexcolor}{\mathbf{M}_{\indexGridLines\indexGridNodeTwo\indexGridNode}  }}

\newcommand{\VunitsqkSDP}[0]{\textcolor{\complexcolor}{\mathbf{M}_{\indexGridNode\indexUnit}  }}
\newcommand{\VunitsqkSDPvoltage}[0]{\textcolor{\complexcolor}{\mathbf{M}^{\symVoltage}_{\indexGridNode\indexUnit}  }}
\newcommand{\VunitsqkSDPcurrent}[0]{\textcolor{\complexcolor}{\mathbf{M}_{\indexGridNode\indexUnit}  }}
\newcommand{\VunitsqkSDPcurrenttworeal}[0]{\textcolor{black}{\mathbf{M}^{2\realss}_{\indexGridNode\indexUnit}  }}

\newcommand{\VunitsqkSDPcurrentdelta}[0]{\textcolor{\complexcolor}{\mathbf{M}^{\Delta}_{\indexGridNode\indexUnit}  }}
\newcommand{\VunitsqkSDPcurrentdeltap}[0]{\textcolor{\complexcolor}{\mathbf{M}^{\Delta'}_{\indexGridNode\indexUnit}  }}

\newcommand{\VjzsqkSDP}[0]{\textcolor{\complexcolor}{\mathbf{\symVoltageSOCP}_{jz }  }}
\newcommand{\VizsqkSDP}[0]{\textcolor{\complexcolor}{\mathbf{\symVoltageSOCP}_{iz }  }}
\newcommand{\VzjsqkSDP}[0]{\textcolor{\complexcolor}{\mathbf{\symVoltageSOCP}_{zj }  }}
\newcommand{\VzisqkSDP}[0]{\textcolor{\complexcolor}{\mathbf{\symVoltageSOCP}_{zi }  }}

\newcommand{\ViusqkSDP}[0]{\textcolor{\complexcolor}{\mathbf{\symVoltageSOCP}_{\indexGridNode \indexUnit}  }}
\newcommand{\VuisqkSDP}[0]{\textcolor{\complexcolor}{\mathbf{\symVoltageSOCP}_{ \indexUnit \indexGridNode}  }}
\newcommand{\ViusqkSDPreal}[0]{\textcolor{black}{\mathbf{\symVoltageSOCP}_{\indexGridNode \indexUnit}^{\realss}  }}
\newcommand{\VuisqkSDPreal}[0]{\textcolor{black}{\mathbf{\symVoltageSOCP}_{ \indexUnit \indexGridNode}^{\realss}  }}
\newcommand{\ViusqkSDPimag}[0]{\textcolor{black}{\mathbf{\symVoltageSOCP}_{\indexGridNode \indexUnit}^{\imagss}   }}
\newcommand{\VuisqkSDPimag}[0]{\textcolor{black}{\mathbf{\symVoltageSOCP}_{ \indexUnit \indexGridNode}^{\imagss}   }}

\newcommand{\VijsqkSDP}[0]{\textcolor{\complexcolor}{\mathbf{\symVoltageSOCP}_{\indexGridNode \indexGridNodeTwo}  }}
\newcommand{\VjisqkSDP}[0]{\textcolor{\complexcolor}{\mathbf{\symVoltageSOCP}_{\indexGridNodeTwo \indexGridNode}  }}
\newcommand{\VijsqkSDPreal}[0]{\textcolor{black}{\mathbf{\symVoltageSOCP}_{\indexGridNode \indexGridNodeTwo}^{\realss}  }}
\newcommand{\VijsqkSDPimag}[0]{\textcolor{black}{\mathbf{\symVoltageSOCP}_{\indexGridNode \indexGridNodeTwo}^{\imagss}  }}

\newcommand{\ratioij}[0]{\textcolor{\complexcolor}{\mathbf{\symRatio}_{\indexGridLines }  }}
\newcommand{\ratiosqij}[0]{\textcolor{black}{\mathbf{H}_{\indexGridLines }  }}
\newcommand{\ratioijA}[0]{\textcolor{\complexcolor}{{\symRatio}_{\indexGridLines , \symPhaseA}  }}
\newcommand{\ratioijB}[0]{\textcolor{\complexcolor}{{\symRatio}_{\indexGridLines , \symPhaseB}  }}
\newcommand{\ratioijC}[0]{\textcolor{\complexcolor}{{\symRatio}_{\indexGridLines , \symPhaseC}  }}

\newcommand{\ratioijreal}[0]{\textcolor{black}{\mathbf{\symRatio}_{\indexGridLines }^{\realss}  }}
\newcommand{\ratioijimag}[0]{\textcolor{black}{\mathbf{\symRatio}_{\indexGridLines }^{\imagss}  }}
\newcommand{\ratioijmag}[0]{\textcolor{black}{\mathbf{\symRatio}_{\indexGridLines }^{\text{mag}}  }}

\newcommand{\ratiomijA}[0]{\textcolor{black}{{\symRatio}^{\text{mag}}_{\indexGridLines , \symPhaseA}  }}
\newcommand{\ratiomijB}[0]{\textcolor{black}{{\symRatio}^{\text{mag}}_{\indexGridLines , \symPhaseB}  }}
\newcommand{\ratiomijC}[0]{\textcolor{black}{{\symRatio}^{\text{mag}}_{\indexGridLines , \symPhaseC}  }}

\newcommand{\ratiomminijA}[0]{\textcolor{\paramcolor}{{\symRatio}^{\text{min}}_{\indexGridLines , \symPhaseA}  }}
\newcommand{\ratiomminijB}[0]{\textcolor{\paramcolor}{{\symRatio}^{\text{min}}_{\indexGridLines , \symPhaseB}  }}
\newcommand{\ratiomminijC}[0]{\textcolor{\paramcolor}{{\symRatio}^{\text{min}}_{\indexGridLines , \symPhaseC}  }}

\newcommand{\ratiommaxijA}[0]{\textcolor{\paramcolor}{{\symRatio}^{\text{max}}_{\indexGridLines , \symPhaseA}  }}
\newcommand{\ratiommaxijB}[0]{\textcolor{\paramcolor}{{\symRatio}^{\text{max}}_{\indexGridLines , \symPhaseB}  }}
\newcommand{\ratiommaxijC}[0]{\textcolor{\paramcolor}{{\symRatio}^{\text{max}}_{\indexGridLines , \symPhaseC}  }}

\newcommand{\ratioaijA}[0]{\textcolor{black}{{\symRatio}^{\angle}_{\indexGridLines , \symPhaseA}  }}
\newcommand{\ratioaijB}[0]{\textcolor{black}{{\symRatio}^{\angle}_{\indexGridLines , \symPhaseB}  }}
\newcommand{\ratioaijC}[0]{\textcolor{black}{{\symRatio}^{\angle}_{\indexGridLines , \symPhaseC}  }}

\newcommand{\ratioaijAmin}[0]{\textcolor{\paramcolor}{{\symRatio}^{\angle}_{\indexGridLines , \symPhaseA}  }}
\newcommand{\ratioaijBmin}[0]{\textcolor{\paramcolor}{{\symRatio}^{\angle}_{\indexGridLines , \symPhaseB}  }}
\newcommand{\ratioaijCmin}[0]{\textcolor{\paramcolor}{{\symRatio}^{\angle}_{\indexGridLines , \symPhaseC}  }}

\newcommand{\ratioaijAmax}[0]{\textcolor{\paramcolor}{{\symRatio}^{\angle}_{\indexGridLines , \symPhaseA}  }}
\newcommand{\ratioaijBmax}[0]{\textcolor{\paramcolor}{{\symRatio}^{\angle}_{\indexGridLines , \symPhaseB}  }}
\newcommand{\ratioaijCmax}[0]{\textcolor{\paramcolor}{{\symRatio}^{\angle}_{\indexGridLines , \symPhaseC}  }}

\newcommand{\ratioaijref}[0]{\textcolor{\paramcolor}{{\symRatio}^{\angle \text{ref}}_{\indexGridLines\indexGridNode \indexGridNodeTwo}  }}

\newcommand{\ratiomaxij}[0]{\textcolor{\paramcolor}{\mathbf{\symRatio}^{\text{max}}_{\indexGridLines }  }}
\newcommand{\ratiominij}[0]{\textcolor{\paramcolor}{\mathbf{\symRatio}^{\text{min}}_{\indexGridLines }  }}

\newcommand{\ratioamaxij}[0]{\textcolor{\paramcolor}{\mathbf{\symRatio}^{\angle\text{max}}_{\indexGridLines }  }}
\newcommand{\ratioaminij}[0]{\textcolor{\paramcolor}{\mathbf{\symRatio}^{\angle\text{min}}_{\indexGridLines }  }}

\newcommand{\Visqksocp}[0]{\textcolor{black}{\symVoltageSOCP_{\indexGridNode }  }}
\newcommand{\Vjsqksocp}[0]{\textcolor{black}{\symVoltageSOCP_{\indexGridNodeTwo } }}
\newcommand{\Visqksocpacc}[0]{\textcolor{black}{\symVoltageSOCP_{\indexGridNode }^{'}  }}
\newcommand{\Vjsqksocpacc}[0]{\textcolor{black}{\symVoltageSOCP_{\indexGridNodeTwo }^{'} }}
\newcommand{\Visqksocpaccl}[0]{\textcolor{black}{\symVoltageSOCP_{\indexGridNode }^{\star}  }}
\newcommand{\Vjsqksocpaccl}[0]{\textcolor{black}{\symVoltageSOCP_{\indexGridNodeTwo }^{\star} }}
\newcommand{\Visqksocpjabr}[0]{\textcolor{black}{\symVoltageSOCP^{\star}_{\indexGridNode }  }}
\newcommand{\Vjsqksocpjabr}[0]{\textcolor{black}{\symVoltageSOCP^{\star}_{\indexGridNodeTwo } }}
\newcommand{\Visqksocpjabrl}[0]{\textcolor{black}{\symVoltageSOCP^{\star\indexGridLines}_{\indexGridNode }  }}
\newcommand{\Vjsqksocpjabrl}[0]{\textcolor{black}{\symVoltageSOCP^{\star\indexGridLines}_{\indexGridNodeTwo } }}

\newcommand{\Ajabr}[0]{\textcolor{\paramcolor}{A^{'}_{\indexGridNode \indexGridNodeTwo } }}
\newcommand{\Bjabr}[0]{\textcolor{\paramcolor}{B^{'}_{\indexGridNode \indexGridNodeTwo } }}
\newcommand{\Cjabr}[0]{\textcolor{\paramcolor}{C^{'}_{\indexGridNode \indexGridNodeTwo } }}
\newcommand{\Djabr}[0]{\textcolor{\paramcolor}{D^{'}_{\indexGridNode \indexGridNodeTwo } }}

\newcommand{\bigM}[0]{\textcolor{\paramcolor}{M^{} }}

\newcommand{\thetaikref}[0]{\textcolor{\paramcolor}{\symVoltageAngle_{\indexGridNode }^{\refss}  }}
\newcommand{\thetaikrefone}[0]{\textcolor{\paramcolor}{\symVoltageAngle_{\indexGridNode,1 }^{\refss}  }}
\newcommand{\Vikref}[0]{\textcolor{\paramcolor}{\symVoltage_{\indexGridNode }^{\refss}  }}
\newcommand{\Viref}[0]{\textcolor{\paramcolor}{\symVoltage_{\indexGridNode }^{\refss} }}

\newcommand{\Vjkref}[0]{\textcolor{\paramcolor}{\symVoltage_{\indexGridNodeTwo }^{\refss}  }}
\newcommand{\Vjref}[0]{\textcolor{\paramcolor}{\symVoltage_{\indexGridNodeTwo }^{\refss} }}

\newcommand{\thetaik}[0]{\textcolor{black}{\symVoltageAngle_{\indexGridNode }^{}  }}
\newcommand{\thetajk}[0]{\textcolor{black}{\symVoltageAngle_{\indexGridNodeTwo }^{}  }}
\newcommand{\thetaijk}[0]{\textcolor{black}{\symVoltageAngle_{\indexGridNode\indexGridNodeTwo }^{}  }}
\newcommand{\thetaijacckmin}[0]{\textcolor{\paramcolor}{\symVoltageAngle_{\indexGridNode\indexGridNodeTwo }^{'\minss}  }}
\newcommand{\thetaijacckmax}[0]{\textcolor{\paramcolor}{\symVoltageAngle_{\indexGridNode\indexGridNodeTwo }^{'\maxss}  }}
\newcommand{\thetaikmin}[0]{\textcolor{\boundscolor}{\symVoltageAngle_{\indexGridNode }^{\minss}  }}
\newcommand{\thetaikmax}[0]{\textcolor{\boundscolor}{\symVoltageAngle_{\indexGridNode }^{\maxss}  }}

\newcommand{\thetaijkmax}[0]{\textcolor{\boundscolor}{\symVoltageAngle_{\indexGridNode \indexGridNodeTwo}^{\maxss}  }}
\newcommand{\thetaijkmin}[0]{\textcolor{\boundscolor}{\symVoltageAngle_{\indexGridNode \indexGridNodeTwo}^{\minss}  }}

\newcommand{\thetaijkabs}[0]{\textcolor{\boundscolor}{\symVoltageAngle_{\indexGridNode \indexGridNodeTwo}^{\text{absmax}}  }}

\newcommand{\thetaijAAmin}[0]{\textcolor{\paramcolor}{\symVoltageAngle_{\indexGridNode\indexGridNodeTwo,\symPhaseA\symPhaseA }^{\minss}  }}
\newcommand{\thetaijAAmax}[0]{\textcolor{\paramcolor}{\symVoltageAngle_{\indexGridNode\indexGridNodeTwo,\symPhaseA\symPhaseA  }^{\maxss}  }}

\newcommand{\thetaijBBmin}[0]{\textcolor{\paramcolor}{\symVoltageAngle_{\indexGridNode\indexGridNodeTwo,\symPhaseB\symPhaseB }^{\minss}  }}
\newcommand{\thetaijBBmax}[0]{\textcolor{\paramcolor}{\symVoltageAngle_{\indexGridNode\indexGridNodeTwo,\symPhaseB\symPhaseB  }^{\maxss}  }}

\newcommand{\thetaijCCmin}[0]{\textcolor{\paramcolor}{\symVoltageAngle_{\indexGridNode\indexGridNodeTwo,\symPhaseC\symPhaseC }^{\minss}  }}
\newcommand{\thetaijCCmax}[0]{\textcolor{\paramcolor}{\symVoltageAngle_{\indexGridNode\indexGridNodeTwo,\symPhaseC\symPhaseC  }^{\maxss}  }}

\newcommand{\thetaijkmaxacc}[0]{\textcolor{\boundscolor}{\symVoltageAngle_{\indexGridNode \indexGridNodeTwo}^{'\maxss}  }}
\newcommand{\thetaijkminacc}[0]{\textcolor{\boundscolor}{\symVoltageAngle_{\indexGridNode \indexGridNodeTwo}^{'\minss}  }}

\newcommand{\thetaijkabsacc}[0]{\textcolor{\boundscolor}{\symVoltageAngle_{\indexGridNode \indexGridNodeTwo}^{'\text{absmax}}  }}

\newcommand{\phiijk}[0]{\textcolor{black}{\symPhaseDifference_{\indexGridNode\indexGridNodeTwo }^{}  }}
\newcommand{\phiikmax}[0]{\textcolor{\boundscolor}{\symPhaseDifference_{\indexGridNode }^{\maxss}  }}
\newcommand{\phiikmin}[0]{\textcolor{\boundscolor}{\symPhaseDifference_{\indexGridNode }^{\minss}  }}
\newcommand{\phiijkmin}[0]{\textcolor{\boundscolor}{\symPhaseDifference_{\indexGridNode   \indexGridNodeTwo }^{\minss}  }}
\newcommand{\phiijkmax}[0]{\textcolor{\boundscolor}{\symPhaseDifference_{\indexGridNode   \indexGridNodeTwo }^{\maxss}  }}

\newcommand{\phijik}[0]{\textcolor{black}{\symPhaseDifference_{\indexGridNodeTwo   \indexGridNode }^{}  }}
\newcommand{\phijikmin}[0]{\textcolor{\boundscolor}{\symPhaseDifference_{\indexGridNodeTwo   \indexGridNode    }^{\minss}  }}
\newcommand{\phijikmax}[0]{\textcolor{\boundscolor}{\symPhaseDifference_{\indexGridNodeTwo  \indexGridNode    }^{\maxss}  }}

\newcommand{\thetaiacck}[0]{\textcolor{black}{\symVoltageAngle_{\indexGridNode }^{'}  }}
\newcommand{\thetajacck}[0]{\textcolor{black}{\symVoltageAngle_{\indexGridNodeTwo }^{'}  }}
\newcommand{\thetaijacck}[0]{\textcolor{black}{\symVoltageAngle_{\indexGridNode  \indexGridNodeTwo }^{'}  }}

\newcommand{\thetaiacckl}[0]{\textcolor{black}{\symVoltageAngle_{\indexGridNode }^{\star}  }}
\newcommand{\thetajacckl}[0]{\textcolor{black}{\symVoltageAngle_{\indexGridNodeTwo }^{\star}  }}

\newcommand{\yijkshSDPseq}[0]{\textcolor{\complexparamcolor}{\mathbf{\symAdmittance}_{\indexGridLines \indexGridNode   \indexGridNodeTwo}^{\shuntss, \fortescuess}    }}
\newcommand{\yijksSDPseq}[0]{\textcolor{\complexparamcolor}{\mathbf{\symAdmittance}_{\indexGridLines \indexGridNode   \indexGridNodeTwo}^{\seriesss, \fortescuess}    }}

\newcommand{\yijkshSDPAdot}[0]{\textcolor{\complexparamcolor}{\mathbf{\symAdmittance}_{\indexGridLines \indexGridNode   \indexGridNodeTwo}^{\shuntss, \symPhaseA\!\cdot}    }}
\newcommand{\yijkshSDPBdot}[0]{\textcolor{\complexparamcolor}{\mathbf{\symAdmittance}_{\indexGridLines \indexGridNode   \indexGridNodeTwo}^{\shuntss, \symPhaseB\!\cdot}    }}
\newcommand{\yijkshSDPCdot}[0]{\textcolor{\complexparamcolor}{\mathbf{\symAdmittance}_{\indexGridLines \indexGridNode   \indexGridNodeTwo}^{\shuntss, \symPhaseC\!\cdot}    }}
\newcommand{\yijkshSDPNdot}[0]{\textcolor{\complexparamcolor}{\mathbf{\symAdmittance}_{\indexGridLines \indexGridNode   \indexGridNodeTwo}^{\shuntss, \symPhaseN\!\cdot}    }}

\newcommand{\yijkshSDPdotN}[0]{\textcolor{\complexparamcolor}{\mathbf{\symAdmittance}_{\indexGridLines \indexGridNode   \indexGridNodeTwo}^{\shuntss, \cdot \! \symPhaseN}    }}

\newcommand{\yijkshSDP}[0]{\textcolor{\complexparamcolor}{\mathbf{\symAdmittance}_{\indexGridLines \indexGridNode   \indexGridNodeTwo}^{\shuntss}    }}
\newcommand{\gijkshSDP}[0]{\textcolor{\paramcolor}{\mathbf{\symConductance}_{ \indexGridLines\indexGridNode   \indexGridNodeTwo}^{\shuntss}    }}
\newcommand{\bijkshSDP}[0]{\textcolor{\paramcolor}{\mathbf{\symSusceptance}_{ \indexGridLines\indexGridNode   \indexGridNodeTwo}^{\shuntss}    }}
\newcommand{\zijkshSDP}[0]{\textcolor{\complexparamcolor}{\mathbf{\symImpedance}_{\indexGridLines \indexGridNode   \indexGridNodeTwo}^{\shuntss}    }}

\newcommand{\yuunitSDP}[0]{\textcolor{\complexparamcolor}{\mathbf{\symAdmittance}_{    \indexUnit}  }}
\newcommand{\yuunitSDPDelta}[0]{\textcolor{\complexparamcolor}{\mathbf{\symAdmittance}_{    \indexUnit}^{\Delta}  }}

\newcommand{\clkshSDP}[0]{\textcolor{\paramcolor}{\mathbf{\symCapacitance}_{\indexGridLines    }^{\shuntss}    }}

\newcommand{\ybSDP}[0]{\textcolor{\complexparamcolor}{\mathbf{\symAdmittance}_{    \indexShunt}  }}
\newcommand{\gbSDP}[0]{\textcolor{\paramcolor}{\mathbf{\symConductance}_{    \indexShunt}  }}
\newcommand{\bbSDP}[0]{\textcolor{\paramcolor}{\mathbf{\symSusceptance}_{    \indexShunt}  }}

\newcommand{\yjikshSDP}[0]{\textcolor{\complexparamcolor}{\mathbf{\symAdmittance}_{ \indexGridLines\indexGridNodeTwo   \indexGridNode}^{\shuntss}    }}
\newcommand{\gjikshSDP}[0]{\textcolor{\paramcolor}{\mathbf{\symConductance}_{ \indexGridLines\indexGridNodeTwo   \indexGridNode}^{\shuntss}    }}
\newcommand{\bjikshSDP}[0]{\textcolor{\paramcolor}{\mathbf{\symSusceptance}_{ \indexGridLines\indexGridNodeTwo   \indexGridNode}^{\shuntss}    }}
\newcommand{\zjikshSDP}[0]{\textcolor{\complexparamcolor}{\mathbf{\symImpedance}_{\indexGridLines \indexGridNodeTwo   \indexGridNode}^{\shuntss}    }}

\newcommand{\yikSDP}[0]{\textcolor{\complexparamcolor}{\mathbf{\symAdmittance}_{ \indexGridNode   }^{}    }}
\newcommand{\gikSDP}[0]{\textcolor{\paramcolor}{\mathbf{\symConductance}_{ \indexGridNode   }^{}    }}
\newcommand{\bikSDP}[0]{\textcolor{\paramcolor}{\mathbf{\symSusceptance}_{ \indexGridNode   }^{}    }}

\newcommand{\zijksSDPkron}[0]{\textcolor{\complexparamcolor}{\mathbf{{\symImpedance}}_{ \indexGridLines}^{\seriesss, \text{Kron}}    }}
\newcommand{\zijksSDPphase}[0]{\textcolor{\complexparamcolor}{\mathbf{{\symImpedance}}_{ \indexGridLines}^{\seriesss, \text{phase}}    }}
\newcommand{\ylshSDPphase}[0]{\textcolor{\complexparamcolor}{\mathbf{{\symAdmittance}}_{ \indexGridLines}^{\shuntss, \text{phase}}    }}
\newcommand{\zijksSDPprim}[0]{\textcolor{\complexparamcolor}{\mathbf{\hat{\symImpedance}}_{ \indexGridLines}^{\seriesss, \text{circ}}    }}
\newcommand{\ylshSDPprim}[0]{\textcolor{\complexparamcolor}{\mathbf{\hat{\symAdmittance}}_{ \indexGridLines}^{\shuntss, \text{circ}}    }}

\newcommand{\zijksSDPvec}[0]{\textcolor{\complexparamcolor}{\mathbf{\hat{\symImpedance}}_{ \indexGridLines}^{\seriesss, \text{pn}}    }}
\newcommand{\zijksSDPvect}[0]{\textcolor{\complexparamcolor}{\mathbf{\hat{\symImpedance}}_{ \indexGridLines}^{\seriesss,  \text{np}}    }}
\newcommand{\zijksSDPhat}[0]{\textcolor{\complexparamcolor}{\mathbf{\hat{\symImpedance}}_{ \indexGridLines}^{\seriesss,  \text{pp}}    }}

\newcommand{\zijksSDPcond}[0]{\textcolor{\complexparamcolor}{\mathbf{\bar{\symImpedance}}_{ \indexGridLines}^{\seriesss, \text{cond}}    }}
\newcommand{\ylshSDPcond}[0]{\textcolor{\complexparamcolor}{\mathbf{\bar{\symAdmittance}}_{ \indexGridLines}^{\shuntss, \text{cond}}    }}
\newcommand{\zijksSDPseq}[0]{\textcolor{\complexparamcolor}{\mathbf{\symImpedance}_{ \indexGridLines}^{\seriesss, \fortescuess}    }}

\newcommand{\zijksSDPseqdiag}[0]{\textcolor{\complexparamcolor}{\mathbf{\symImpedance}_{ \indexGridLines}^{\seriesss, \fortescuess, \text{diag}}    }}

\newcommand{\zmat}[0]{\textcolor{\complexparamcolor}{\mathbf{\symImpedance}}}
\newcommand{\zmatP}[0]{\textcolor{\complexparamcolor}{\mathbf{\symImpedance}}^{\setPhases}}
\newcommand{\zmatF}[0]{\textcolor{\complexparamcolor}{\mathbf{\symImpedance}}^{\setFortescue}}

\newcommand{\thetaijSDPmax}[0]{\textcolor{\paramcolor}{\mathbf{\Theta}_{\indexGridNode\indexGridNodeTwo  }^{\maxss}  }}
\newcommand{\thetaijSDPmin}[0]{\textcolor{\paramcolor}{\mathbf{\Theta}_{\indexGridNode\indexGridNodeTwo  }^{\minss}  }}

\newcommand{\seq}[0]{\textcolor{\complexparamcolor}{\mathbf{A}}}
\newcommand{\seqinv}[0]{\textcolor{\complexparamcolor}{\mathbf{A}^{-1}}}
\newcommand{\seqH}[0]{\textcolor{\complexparamcolor}{\mathbf{A}^{\hermitiantranspose}}}
\newcommand{\seqconj}[0]{\textcolor{\complexparamcolor}{\mathbf{A}^{*}}}
\newcommand{\seqT}[0]{\textcolor{\complexparamcolor}{\mathbf{A}^{\transpose}}}

\newcommand{\seqvec}[0]{\textcolor{\complexparamcolor}{\mathbf{a}}}

\newcommand{\zslh}[0]{\textcolor{\complexparamcolor}{\mathbf{\symImpedance}_{ \indexGridLines, \symHarmonic}^{\seriesss}    }}
\newcommand{\rslh}[0]{\textcolor{\paramcolor}{\mathbf{\symResistance}_{  \indexGridLines, \symHarmonic   }^{\seriesss}    }}
\newcommand{\xslh}[0]{\textcolor{\paramcolor}{\mathbf{\symReactance}_{  \indexGridLines, \symHarmonic   }^{\seriesss}    }}

\newcommand{\yshlijh}[0]{\textcolor{\complexparamcolor}{\mathbf{\symAdmittance}_{ \indexGridLines\indexGridNode\indexGridNodeTwo, \symHarmonic}^{\shuntss}    }}
\newcommand{\gshlijh}[0]{\textcolor{\paramcolor}{\mathbf{\symConductance}_{  \indexGridLines\indexGridNode\indexGridNodeTwo, \symHarmonic   }^{\shuntss}    }}
\newcommand{\bshlijh}[0]{\textcolor{\paramcolor}{\mathbf{\symSusceptance}_{  \indexGridLines\indexGridNode\indexGridNodeTwo, \symHarmonic   }^{\shuntss}    }}

\newcommand{\yshljih}[0]{\textcolor{\complexparamcolor}{\mathbf{\symAdmittance}_{ \indexGridLines\indexGridNodeTwo\indexGridNode, \symHarmonic}^{\shuntss}    }}
\newcommand{\gshljih}[0]{\textcolor{\paramcolor}{\mathbf{\symConductance}_{  \indexGridLines\indexGridNodeTwo\indexGridNode, \symHarmonic   }^{\shuntss}    }}
\newcommand{\bshljih}[0]{\textcolor{\paramcolor}{\mathbf{\symSusceptance}_{  \indexGridLines\indexGridNodeTwo\indexGridNode, \symHarmonic   }^{\shuntss}    }}

\newcommand{\zijksSDP}[0]{\textcolor{\complexparamcolor}{\mathbf{\symImpedance}_{ \indexGridLines}^{\seriesss}    }}
\newcommand{\rijksSDP}[0]{\textcolor{\paramcolor}{\mathbf{\symResistance}_{  \indexGridLines   }^{\seriesss}    }}
\newcommand{\xijksSDP}[0]{\textcolor{\paramcolor}{\mathbf{\symReactance}_{  \indexGridLines   }^{\seriesss}    }}

\newcommand{\zjiksSDP}[0]{\textcolor{\complexparamcolor}{\mathbf{\symImpedance}_{ \indexGridLines   }^{\seriesss}    }}
\newcommand{\rjiksSDP}[0]{\textcolor{\paramcolor}{\mathbf{\symResistance}_{ \indexGridLines   }^{\seriesss}    }}
\newcommand{\xjiksSDP}[0]{\textcolor{\paramcolor}{\mathbf{\symReactance}_{ \indexGridLines   }^{\seriesss}    }}

\newcommand{\zijks}[0]{\textcolor{\complexparamcolor}{\symImpedance_{ \indexGridLines,\seriesss}    }}
\newcommand{\rijks}[0]{\textcolor{\paramcolor}{\symResistance_{  \indexGridLines,\seriesss}    }}
\newcommand{\xijks}[0]{\textcolor{\paramcolor}{\symReactance_{  \indexGridLines,\seriesss}    }}
\newcommand{\yijks}[0]{\textcolor{\complexparamcolor}{\symAdmittance_{  \indexGridLines}^{\seriesss}    }}
\newcommand{\bijks}[0]{\textcolor{\paramcolor}{\symSusceptance_{  \indexGridLines}^{\seriesss}    }}
\newcommand{\gijks}[0]{\textcolor{\paramcolor}{\symConductance_{  \indexGridLines}^{\seriesss}    }}

\newcommand{\zijksSDPH}[0]{\textcolor{\complexparamcolor}{(\mathbf{\symImpedance}_{ \indexGridLines}^{\seriesss})^{\hermitiantranspose}    }}
\newcommand{\rijksSDPH}[0]{\textcolor{\paramcolor}{(\mathbf{\symResistance}_{ \indexGridLines,\seriesss})^{\transpose}    }}
\newcommand{\xijksSDPH}[0]{\textcolor{\paramcolor}{(\mathbf{\symReactance}_{ \indexGridLines,\seriesss})^{\transpose}    }}

\newcommand{\yukSDP}[0]{\textcolor{\complexparamcolor}{\mathbf{\symAdmittance}_{ \indexUnit}    }}
\newcommand{\gukSDP}[0]{\textcolor{\paramcolor}{\mathbf{\symConductance}_{ \indexUnit}    }}
\newcommand{\bukSDP}[0]{\textcolor{\paramcolor}{\mathbf{\symSusceptance}_{ \indexUnit}    }}

\newcommand{\zukSDP}[0]{\textcolor{\complexparamcolor}{\mathbf{\symImpedance}_{ \indexUnit}    }}
\newcommand{\rukSDP}[0]{\textcolor{\paramcolor}{\mathbf{\symResistance}_{ \indexUnit}    }}
\newcommand{\xukSDP}[0]{\textcolor{\paramcolor}{\mathbf{\symReactance}_{ \indexUnit}    }}

\newcommand{\yijksSDP}[0]{\textcolor{\complexparamcolor}{\mathbf{\symAdmittance}_{ \indexGridLines}^{\seriesss}    }}
\newcommand{\gijksSDP}[0]{\textcolor{\paramcolor}{\mathbf{\symConductance}_{ \indexGridLines}^{\seriesss}    }}
\newcommand{\bijksSDP}[0]{\textcolor{\paramcolor}{\mathbf{\symSusceptance}_{ \indexGridLines}^{\seriesss}    }}

\newcommand{\yjiksSDP}[0]{\textcolor{\complexparamcolor}{\mathbf{\symAdmittance}_{  \indexGridLines\indexGridNodeTwo \indexGridNode }^{\seriesss}    }}
\newcommand{\gjiksSDP}[0]{\textcolor{\paramcolor}{\mathbf{\symConductance}_{ \indexGridLines \indexGridNodeTwo \indexGridNode }^{\seriesss}    }}
\newcommand{\bjiksSDP}[0]{\textcolor{\paramcolor}{\mathbf{\symSusceptance}_{  \indexGridLines\indexGridNodeTwo \indexGridNode }^{\seriesss}    }}

\newcommand{\zijksh}[0]{\textcolor{\complexparamcolor}{\symImpedance_{ \indexGridNode \indexGridNodeTwo,\shuntss}    }}
\newcommand{\rijksh}[0]{\textcolor{\paramcolor}{\symResistance_{ \indexGridNode \indexGridNodeTwo,\shuntss}    }}
\newcommand{\xijksh}[0]{\textcolor{\paramcolor}{\symReactance_{ \indexGridNode \indexGridNodeTwo,\shuntss}    }}
\newcommand{\yijksh}[0]{\textcolor{\complexparamcolor}{\symAdmittance_{ \indexGridNode \indexGridNodeTwo,\shuntss}    }}
\newcommand{\gijksh}[0]{\textcolor{\paramcolor}{\symConductance_{ \indexGridNode \indexGridNodeTwo,\shuntss}    }}
\newcommand{\bijksh}[0]{\textcolor{\paramcolor}{\symSusceptance_{ \indexGridNode \indexGridNodeTwo,\shuntss}    }}

\newcommand{\zjiksh}[0]{\textcolor{\complexparamcolor}{\symImpedance_{ \indexGridNodeTwo \indexGridNode ,\shuntss}    }}
\newcommand{\rjiksh}[0]{\textcolor{\paramcolor}{\symResistance_{\indexGridNodeTwo \indexGridNode ,\shuntss}    }}
\newcommand{\xjiksh}[0]{\textcolor{\paramcolor}{\symReactance_{ \indexGridNodeTwo\indexGridNode ,\shuntss}    }}
\newcommand{\yjiksh}[0]{\textcolor{\complexparamcolor}{\symAdmittance_{\indexGridNodeTwo \indexGridNode ,\shuntss}    }}
\newcommand{\gjiksh}[0]{\textcolor{\paramcolor}{\symConductance_{\indexGridNodeTwo \indexGridNode ,\shuntss}    }}
\newcommand{\bjiksh}[0]{\textcolor{\paramcolor}{\symSusceptance_{ \indexGridNodeTwo \indexGridNode ,\shuntss}    }}

\newcommand{\zlksh}[0]{\textcolor{\complexparamcolor}{\symImpedance_{  \indexGridLines,\shuntss}    }}
\newcommand{\rlksh}[0]{\textcolor{\paramcolor}{\symResistance_{  \indexGridLines,\shuntss}    }}
\newcommand{\xlksh}[0]{\textcolor{\paramcolor}{\symReactance_{  \indexGridLines,\shuntss}    }}
\newcommand{\ylksh}[0]{\textcolor{\complexparamcolor}{\symAdmittance_{  \indexGridLines,\shuntss}    }}
\newcommand{\glksh}[0]{\textcolor{\paramcolor}{\symConductance_{  \indexGridLines,\shuntss}    }}
\newcommand{\blksh}[0]{\textcolor{\paramcolor}{\symSusceptance_{  \indexGridLines,\shuntss}    }}

\newcommand{\zlkshspec}[0]{\textcolor{\complexparamcolor}{\dot{\symImpedance}_{  \indexGridLines,\shuntss}    }}
\newcommand{\rlkshspec}[0]{\textcolor{\paramcolor}{\dot{\symResistance}_{  \indexGridLines,\shuntss}    }}
\newcommand{\xlkshv}[0]{\textcolor{\paramcolor}{\dot{\symReactance}_{  \indexGridLines,\shuntss}    }}
\newcommand{\ylkshspec}[0]{\textcolor{\complexparamcolor}{\dot{\symAdmittance}_{  \indexGridLines,\shuntss}    }}
\newcommand{\glkshspec}[0]{\textcolor{\paramcolor}{\dot{\symConductance}_{  \indexGridLines,\shuntss}    }}
\newcommand{\blkshspec}[0]{\textcolor{\paramcolor}{\dot{\symSusceptance}_{  \indexGridLines,\shuntss}    }}

\newcommand{\zijksspec}[0]{\textcolor{\complexparamcolor}{\dot{\symImpedance}_{  \indexGridLines,\seriesss}    }}
\newcommand{\rijksspec}[0]{\textcolor{\paramcolor}{\dot{\symResistance}_{  \indexGridLines,\seriesss}    }}
\newcommand{\xijksspec}[0]{\textcolor{\paramcolor}{\dot{\symReactance}_{  \indexGridLines,\seriesss}    }}
\newcommand{\yijkshspec}[0]{\textcolor{\complexparamcolor}{\dot{\symAdmittance}_{  \indexGridNode \indexGridNodeTwo,\shuntss}    }}
\newcommand{\bijkshspec}[0]{\textcolor{\paramcolor}{\dot{\symSusceptance}_{  \indexGridNode \indexGridNodeTwo,\shuntss}    }}
\newcommand{\gijkshspec}[0]{\textcolor{\paramcolor}{\dot{\symConductance}_{  \indexGridNode \indexGridNodeTwo,\shuntss}    }}

\newcommand{\lijk}[0]{\textcolor{\paramcolor}{{\symLength}_{  \indexGridLines}    }}
\newcommand{\vprimvsecijk}[0]{\textcolor{black}{{\symRatio}_{ \indexGridNode \indexGridNodeTwo}    }}
\newcommand{\vprimvsecijkmin}[0]{\textcolor{\boundscolor}{{\symRatio}_{ \indexGridNode \indexGridNodeTwo}^{\minss}    }}
\newcommand{\vprimvsecijkmax}[0]{\textcolor{\boundscolor}{{\symRatio}_{ \indexGridNode \indexGridNodeTwo}^{\maxss}    }}

\newcommand{\vsecvprimijk}[0]{\textcolor{black}{{\symRatio}_{\indexGridNodeTwo \indexGridNode }    }}
\newcommand{\vsecvprimijkmin}[0]{\textcolor{\boundscolor}{{\symRatio}_{ \indexGridNodeTwo\indexGridNode }^{\minss}    }}
\newcommand{\vsecvprimijkmax}[0]{\textcolor{\boundscolor}{{\symRatio}_{\indexGridNodeTwo \indexGridNode }^{\maxss}    }}

\newcommand{\YPFconvexification}[0]{\textcolor{black}{{\symPenalty}^{\PFconvexss}   }}
\newcommand{\YGIconvexification}[0]{\textcolor{black}{{\symPenalty}^{\GIconvexss}   }}
\newcommand{\YPgridloss}[0]{\textcolor{black}{{\symPenalty}^{\symPower,\lossss}   }}
\newcommand{\YSgridloss}[0]{\textcolor{black}{{\symPenalty}^{\symApparentPower,\lossss}   }}
\newcommand{\Ytotcurrent}[0]{\textcolor{black}{{\symPenalty}^{\text{\symCurrent}}   }}

\newcommand{\Yobj}[0]{\textcolor{black}{{\symPenalty}^{\text{cost}}   }}
\newcommand{\Kobj}[0]{\textcolor{black}{{\symAnnualCost}^{\text{cost}}   }}

\newcommand{\NPFconvexification}[0]{\textcolor{\paramcolor}{{\symPenaltyWeight}_{\PFconvexss}   }}
\newcommand{\NGIconvexification}[0]{\textcolor{\paramcolor}{{\symPenaltyWeight}_{\GIconvexss}   }}

\newcommand{\wPFconvexification}[0]{\textcolor{\paramcolor}{{\symObjectiveWeight}_{\PFconvexss}   }}
\newcommand{\wPFconvexificationopt}[0]{\textcolor{\paramcolor}{{\symObjectiveWeight}_{\PFconvexss}^{\text{opt}}   }}

\newcommand{\zijkszero}[0]{\textcolor{\complexparamcolor}{\symImpedance_{ \indexGridLines , 0}^{\seriesss}    }}
\newcommand{\zijksone}[0]{\textcolor{\complexparamcolor}{\symImpedance_{ \indexGridLines , 1}^{\seriesss}    }}

\newcommand{\zijkssymmzerozero}[0]{\textcolor{\complexparamcolor}{\symImpedance_{ \indexGridLines , 0 0}^{\seriesss}    }}
\newcommand{\zijkssymmzeroone}[0]{\textcolor{\complexparamcolor}{\symImpedance_{ \indexGridLines , 0 1}^{\seriesss}    }}
\newcommand{\zijkssymmzerotwo}[0]{\textcolor{\complexparamcolor}{\symImpedance_{ \indexGridLines , 0 2}^{\seriesss}    }}

\newcommand{\zijkssymmonezero}[0]{\textcolor{\complexparamcolor}{\symImpedance_{ \indexGridLines , 1 0}^{\seriesss}    }}
\newcommand{\zijkssymmoneone}[0]{\textcolor{\complexparamcolor}{\symImpedance_{ \indexGridLines , 1 1}^{\seriesss}    }}
\newcommand{\zijkssymmonetwo}[0]{\textcolor{\complexparamcolor}{\symImpedance_{ \indexGridLines , 1 2}^{\seriesss}    }}

\newcommand{\zijkssymmtwozero}[0]{\textcolor{\complexparamcolor}{\symImpedance_{ \indexGridLines , 2 0}^{\seriesss}    }}
\newcommand{\zijkssymmtwoone}[0]{\textcolor{\complexparamcolor}{\symImpedance_{ \indexGridLines , 2 1}^{\seriesss}    }}
\newcommand{\zijkssymmtwotwo}[0]{\textcolor{\complexparamcolor}{\symImpedance_{ \indexGridLines , 2 2}^{\seriesss}    }}

\newcommand{\zjNG}[0]{\textcolor{\complexparamcolor}{\symImpedance_{ \indexGridNodeTwo , \symPhaseN \symPhaseG}   }}
\newcommand{\ziNG}[0]{\textcolor{\complexparamcolor}{\symImpedance_{ \indexGridNode , \symPhaseN \symPhaseG}   }}
\newcommand{\riNG}[0]{\textcolor{\paramcolor}{\symResistance_{ \indexGridNode , \symPhaseN \symPhaseG}   }}
\newcommand{\xiNG}[0]{\textcolor{\paramcolor}{\symReactance_{ \indexGridNode , \symPhaseN \symPhaseG}   }}

\newcommand{\zijksAA}[0]{\textcolor{\complexparamcolor}{\symImpedance_{ \indexGridLines , \symPhaseA \symPhaseA}^{\seriesss}    }}
\newcommand{\zijksAB}[0]{\textcolor{\complexparamcolor}{\symImpedance_{ \indexGridLines , \symPhaseA \symPhaseB}^{\seriesss}    }}
\newcommand{\zijksAC}[0]{\textcolor{\complexparamcolor}{\symImpedance_{ \indexGridLines , \symPhaseA \symPhaseC}^{\seriesss}    }}
\newcommand{\zijksBA}[0]{\textcolor{\complexparamcolor}{\symImpedance_{ \indexGridLines , \symPhaseB \symPhaseA}^{\seriesss}    }}
\newcommand{\zijksBB}[0]{\textcolor{\complexparamcolor}{\symImpedance_{ \indexGridLines , \symPhaseB \symPhaseB}^{\seriesss}    }}
\newcommand{\zijksBC}[0]{\textcolor{\complexparamcolor}{\symImpedance_{ \indexGridLines , \symPhaseB \symPhaseC}^{\seriesss}    }}
\newcommand{\zijksCA}[0]{\textcolor{\complexparamcolor}{\symImpedance_{ \indexGridLines , \symPhaseC \symPhaseA}^{\seriesss}    }}
\newcommand{\zijksCB}[0]{\textcolor{\complexparamcolor}{\symImpedance_{ \indexGridLines , \symPhaseC \symPhaseB}^{\seriesss}    }}
\newcommand{\zijksCC}[0]{\textcolor{\complexparamcolor}{\symImpedance_{ \indexGridLines , \symPhaseC \symPhaseC}^{\seriesss}    }}

 \newcommand{\rijksAA}[0]{\textcolor{\paramcolor}{\symResistance_{ \indexGridLines , \symPhaseA \symPhaseA}^{\seriesss}    }}
\newcommand{\rijksAB}[0]{\textcolor{\paramcolor}{\symResistance_{ \indexGridLines , \symPhaseA \symPhaseB}^{\seriesss}    }}
\newcommand{\rijksAC}[0]{\textcolor{\paramcolor}{\symResistance_{ \indexGridLines , \symPhaseA \symPhaseC}^{\seriesss}    }}
\newcommand{\rijksBA}[0]{\textcolor{\paramcolor}{\symResistance_{ \indexGridLines , \symPhaseB \symPhaseA}^{\seriesss}    }}
\newcommand{\rijksBB}[0]{\textcolor{\paramcolor}{\symResistance_{ \indexGridLines , \symPhaseB \symPhaseB}^{\seriesss}    }}
\newcommand{\rijksBC}[0]{\textcolor{\paramcolor}{\symResistance_{ \indexGridLines , \symPhaseB \symPhaseC}^{\seriesss}    }}
\newcommand{\rijksCA}[0]{\textcolor{\paramcolor}{\symResistance_{ \indexGridLines , \symPhaseC \symPhaseA}^{\seriesss}    }}
\newcommand{\rijksCB}[0]{\textcolor{\paramcolor}{\symResistance_{ \indexGridLines , \symPhaseC \symPhaseB}^{\seriesss}    }}
\newcommand{\rijksCC}[0]{\textcolor{\paramcolor}{\symResistance_{ \indexGridLines , \symPhaseC \symPhaseC}^{\seriesss}    }}
 \newcommand{\rijksAN}[0]{\textcolor{\paramcolor}{\symResistance_{ \indexGridLines , \symPhaseA \symPhaseN}^{\seriesss}    }}
\newcommand{\rijksBN}[0]{\textcolor{\paramcolor}{\symResistance_{ \indexGridLines , \symPhaseB \symPhaseN}^{\seriesss}    }}
\newcommand{\rijksCN}[0]{\textcolor{\paramcolor}{\symResistance_{ \indexGridLines , \symPhaseC \symPhaseN}^{\seriesss}    }}
\newcommand{\rijksNN}[0]{\textcolor{\paramcolor}{\symResistance_{ \indexGridLines , \symPhaseN \symPhaseN}^{\seriesss}    }}

 \newcommand{\xijksAA}[0]{\textcolor{\paramcolor}{\symReactance_{ \indexGridLines , \symPhaseA \symPhaseA}^{\seriesss}    }}
\newcommand{\xijksAB}[0]{\textcolor{\paramcolor}{\symReactance_{ \indexGridLines , \symPhaseA \symPhaseB}^{\seriesss}    }}
\newcommand{\xijksAC}[0]{\textcolor{\paramcolor}{\symReactance_{ \indexGridLines , \symPhaseA \symPhaseC}^{\seriesss}    }}
\newcommand{\xijksBA}[0]{\textcolor{\paramcolor}{\symReactance_{ \indexGridLines , \symPhaseB \symPhaseA}^{\seriesss}    }}
\newcommand{\xijksBB}[0]{\textcolor{\paramcolor}{\symReactance_{ \indexGridLines , \symPhaseB \symPhaseB}^{\seriesss}    }}
\newcommand{\xijksBC}[0]{\textcolor{\paramcolor}{\symReactance_{ \indexGridLines , \symPhaseB \symPhaseC}^{\seriesss}    }}
\newcommand{\xijksCA}[0]{\textcolor{\paramcolor}{\symReactance_{ \indexGridLines , \symPhaseC \symPhaseA}^{\seriesss}    }}
\newcommand{\xijksCB}[0]{\textcolor{\paramcolor}{\symReactance_{ \indexGridLines , \symPhaseC \symPhaseB}^{\seriesss}    }}
\newcommand{\xijksCC}[0]{\textcolor{\paramcolor}{\symReactance_{ \indexGridLines , \symPhaseC \symPhaseC}^{\seriesss}    }}
 \newcommand{\xijksAN}[0]{\textcolor{\paramcolor}{\symReactance_{ \indexGridLines , \symPhaseA \symPhaseN}^{\seriesss}    }}
\newcommand{\xijksBN}[0]{\textcolor{\paramcolor}{\symReactance_{ \indexGridLines , \symPhaseB \symPhaseN}^{\seriesss}    }}
\newcommand{\xijksCN}[0]{\textcolor{\paramcolor}{\symReactance_{ \indexGridLines , \symPhaseC \symPhaseN}^{\seriesss}    }}
\newcommand{\xijksNN}[0]{\textcolor{\paramcolor}{\symReactance_{ \indexGridLines , \symPhaseN \symPhaseN}^{\seriesss}    }}
 \newcommand{\xijksNA}[0]{\textcolor{\paramcolor}{\symReactance_{ \indexGridLines , \symPhaseN \symPhaseA}^{\seriesss}    }}
\newcommand{\xijksNB}[0]{\textcolor{\paramcolor}{\symReactance_{ \indexGridLines , \symPhaseN \symPhaseB}^{\seriesss}    }}
\newcommand{\xijksNC}[0]{\textcolor{\paramcolor}{\symReactance_{ \indexGridLines , \symPhaseN \symPhaseC}^{\seriesss}    }}

\newcommand{\zijksAAkron}[0]{\textcolor{\complexparamcolor}{\symImpedance_{ \indexGridLines , \symPhaseA \symPhaseA}^{\seriesss, \text{Kr}}    }}
\newcommand{\zijksABkron}[0]{\textcolor{\complexparamcolor}{\symImpedance_{ \indexGridLines , \symPhaseA \symPhaseB}^{\seriesss, \text{Kr}}    }}
\newcommand{\zijksACkron}[0]{\textcolor{\complexparamcolor}{\symImpedance_{ \indexGridLines , \symPhaseA \symPhaseC}^{\seriesss, \text{Kr}}    }}
\newcommand{\zijksBAkron}[0]{\textcolor{\complexparamcolor}{\symImpedance_{ \indexGridLines , \symPhaseB \symPhaseA}^{\seriesss, \text{Kr}}    }}
\newcommand{\zijksBBkron}[0]{\textcolor{\complexparamcolor}{\symImpedance_{ \indexGridLines , \symPhaseB \symPhaseB}^{\seriesss, \text{Kr}}    }}
\newcommand{\zijksBCkron}[0]{\textcolor{\complexparamcolor}{\symImpedance_{ \indexGridLines , \symPhaseB \symPhaseC}^{\seriesss, \text{Kr}}    }}
\newcommand{\zijksCAkron}[0]{\textcolor{\complexparamcolor}{\symImpedance_{ \indexGridLines , \symPhaseC \symPhaseA}^{\seriesss, \text{Kr}}    }}
\newcommand{\zijksCBkron}[0]{\textcolor{\complexparamcolor}{\symImpedance_{ \indexGridLines , \symPhaseC \symPhaseB}^{\seriesss, \text{Kr}}    }}
\newcommand{\zijksCCkron}[0]{\textcolor{\complexparamcolor}{\symImpedance_{ \indexGridLines , \symPhaseC \symPhaseC}^{\seriesss, \text{Kr}}    }}

\newcommand{\zijksAN}[0]{\textcolor{\complexparamcolor}{\symImpedance_{ \indexGridLines , \symPhaseA \symPhaseN}^{\seriesss}    }}
\newcommand{\zijksBN}[0]{\textcolor{\complexparamcolor}{\symImpedance_{ \indexGridLines , \symPhaseB \symPhaseN}^{\seriesss}    }}
\newcommand{\zijksCN}[0]{\textcolor{\complexparamcolor}{\symImpedance_{ \indexGridLines , \symPhaseC \symPhaseN}^{\seriesss}    }}
\newcommand{\zijksNA}[0]{\textcolor{\complexparamcolor}{\symImpedance_{ \indexGridLines , \symPhaseN \symPhaseA }^{\seriesss}    }}
\newcommand{\zijksNB}[0]{\textcolor{\complexparamcolor}{\symImpedance_{ \indexGridLines ,\symPhaseN \symPhaseB }^{\seriesss}    }}
\newcommand{\zijksNC}[0]{\textcolor{\complexparamcolor}{\symImpedance_{ \indexGridLines , \symPhaseN\symPhaseC }^{\seriesss}    }}

\newcommand{\zijksNN}[0]{\textcolor{\complexparamcolor}{\symImpedance_{ \indexGridLines , \symPhaseN \symPhaseN}^{\seriesss}    }}

\newcommand{\zijksAG}[0]{\textcolor{\complexparamcolor}{\symImpedance_{ \indexGridLines , \symPhaseA \symPhaseG}^{\seriesss}    }}
\newcommand{\zijksBG}[0]{\textcolor{\complexparamcolor}{\symImpedance_{ \indexGridLines , \symPhaseB \symPhaseG}^{\seriesss}    }}
\newcommand{\zijksCG}[0]{\textcolor{\complexparamcolor}{\symImpedance_{ \indexGridLines , \symPhaseC \symPhaseG}^{\seriesss}    }}
\newcommand{\zijksNG}[0]{\textcolor{\complexparamcolor}{\symImpedance_{ \indexGridLines , \symPhaseN \symPhaseG}^{\seriesss}    }}
\newcommand{\zijksGA}[0]{\textcolor{\complexparamcolor}{\symImpedance_{ \indexGridLines , \symPhaseG \symPhaseA }^{\seriesss}    }}
\newcommand{\zijksGB}[0]{\textcolor{\complexparamcolor}{\symImpedance_{ \indexGridLines ,\symPhaseG \symPhaseB }^{\seriesss}    }}
\newcommand{\zijksGC}[0]{\textcolor{\complexparamcolor}{\symImpedance_{ \indexGridLines , \symPhaseG\symPhaseC }^{\seriesss}    }}
\newcommand{\zijksGN}[0]{\textcolor{\complexparamcolor}{\symImpedance_{ \indexGridLines , \symPhaseG \symPhaseN}^{\seriesss}    }}

\newcommand{\zijksGG}[0]{\textcolor{\complexparamcolor}{\symImpedance_{ \indexGridLines , \symPhaseG \symPhaseG}^{\seriesss}    }}

\newcommand{\zijksNNhat}[0]{\textcolor{\complexparamcolor}{\hat{\symImpedance}_{ \indexGridLines , \symPhaseN \symPhaseN}^{\seriesss}    }}
\newcommand{\zijksANhat}[0]{\textcolor{\complexparamcolor}{\hat{\symImpedance}_{ \indexGridLines , \symPhaseA \symPhaseN}^{\seriesss}    }}
\newcommand{\zijksBNhat}[0]{\textcolor{\complexparamcolor}{\hat{\symImpedance}_{ \indexGridLines , \symPhaseB \symPhaseN}^{\seriesss}    }}
\newcommand{\zijksCNhat}[0]{\textcolor{\complexparamcolor}{\hat{\symImpedance}_{ \indexGridLines , \symPhaseC \symPhaseN}^{\seriesss}    }}

\newcommand{\zijksAAhat}[0]{\textcolor{\complexparamcolor}{\hat{\symImpedance}_{ \indexGridLines , \symPhaseA \symPhaseA}^{\seriesss}    }}
\newcommand{\zijksBAhat}[0]{\textcolor{\complexparamcolor}{\hat{\symImpedance}_{ \indexGridLines , \symPhaseB \symPhaseA}^{\seriesss}    }}
\newcommand{\zijksCAhat}[0]{\textcolor{\complexparamcolor}{\hat{\symImpedance}_{ \indexGridLines , \symPhaseC \symPhaseA}^{\seriesss}    }}

\newcommand{\zijksABhat}[0]{\textcolor{\complexparamcolor}{\hat{\symImpedance}_{ \indexGridLines , \symPhaseA \symPhaseB}^{\seriesss}    }}
\newcommand{\zijksBBhat}[0]{\textcolor{\complexparamcolor}{\hat{\symImpedance}_{ \indexGridLines , \symPhaseB \symPhaseB}^{\seriesss}    }}
\newcommand{\zijksCBhat}[0]{\textcolor{\complexparamcolor}{\hat{\symImpedance}_{ \indexGridLines , \symPhaseC \symPhaseB}^{\seriesss}    }}

\newcommand{\zijksAChat}[0]{\textcolor{\complexparamcolor}{\hat{\symImpedance}_{ \indexGridLines , \symPhaseA \symPhaseC}^{\seriesss}    }}
\newcommand{\zijksBChat}[0]{\textcolor{\complexparamcolor}{\hat{\symImpedance}_{ \indexGridLines , \symPhaseB \symPhaseC}^{\seriesss}    }}
\newcommand{\zijksCChat}[0]{\textcolor{\complexparamcolor}{\hat{\symImpedance}_{ \indexGridLines , \symPhaseC \symPhaseC}^{\seriesss}    }}

\newcommand{\zijksNAhat}[0]{\textcolor{\complexparamcolor}{\hat{\symImpedance}_{ \indexGridLines , \symPhaseN \symPhaseA }^{\seriesss}    }}
\newcommand{\zijksNBhat}[0]{\textcolor{\complexparamcolor}{\hat{\symImpedance}_{ \indexGridLines , \symPhaseN \symPhaseB }^{\seriesss}    }}
\newcommand{\zijksNChat}[0]{\textcolor{\complexparamcolor}{\hat{\symImpedance}_{ \indexGridLines , \symPhaseN \symPhaseC }^{\seriesss}    }}

\newcommand{\yijksAA}[0]{\textcolor{\complexparamcolor}{\symAdmittance_{ \indexGridLines , \symPhaseA \symPhaseA}^{\seriesss}    }}
\newcommand{\yijksAB}[0]{\textcolor{\complexparamcolor}{\symAdmittance_{ \indexGridLines , \symPhaseA \symPhaseB}^{\seriesss}    }}
\newcommand{\yijksAC}[0]{\textcolor{\complexparamcolor}{\symAdmittance_{ \indexGridLines , \symPhaseA \symPhaseC}^{\seriesss}    }}
\newcommand{\yijksBA}[0]{\textcolor{\complexparamcolor}{\symAdmittance_{ \indexGridLines , \symPhaseB \symPhaseA}^{\seriesss}    }}
\newcommand{\yijksBB}[0]{\textcolor{\complexparamcolor}{\symAdmittance_{ \indexGridLines , \symPhaseB \symPhaseB}^{\seriesss}    }}
\newcommand{\yijksBC}[0]{\textcolor{\complexparamcolor}{\symAdmittance_{ \indexGridLines , \symPhaseB \symPhaseC}^{\seriesss}    }}
\newcommand{\yijksCA}[0]{\textcolor{\complexparamcolor}{\symAdmittance_{ \indexGridLines , \symPhaseC \symPhaseA}^{\seriesss}    }}
\newcommand{\yijksCB}[0]{\textcolor{\complexparamcolor}{\symAdmittance_{ \indexGridLines , \symPhaseC \symPhaseB}^{\seriesss}    }}
\newcommand{\yijksCC}[0]{\textcolor{\complexparamcolor}{\symAdmittance_{ \indexGridLines , \symPhaseC \symPhaseC}^{\seriesss}    }}

\newcommand{\gijksAA}[0]{\textcolor{\paramcolor}{\symConductance_{ \indexGridLines , \symPhaseA \symPhaseA}^{\seriesss}    }}
\newcommand{\gijksAB}[0]{\textcolor{\paramcolor}{\symConductance_{ \indexGridLines , \symPhaseA \symPhaseB}^{\seriesss}    }}
\newcommand{\gijksAC}[0]{\textcolor{\paramcolor}{\symConductance_{ \indexGridLines , \symPhaseA \symPhaseC}^{\seriesss}    }}
\newcommand{\gijksBA}[0]{\textcolor{\paramcolor}{\symConductance_{ \indexGridLines , \symPhaseB \symPhaseA}^{\seriesss}    }}
\newcommand{\gijksBB}[0]{\textcolor{\paramcolor}{\symConductance_{ \indexGridLines , \symPhaseB \symPhaseB}^{\seriesss}    }}
\newcommand{\gijksBC}[0]{\textcolor{\paramcolor}{\symConductance_{ \indexGridLines , \symPhaseB \symPhaseC}^{\seriesss}    }}
\newcommand{\gijksCA}[0]{\textcolor{\paramcolor}{\symConductance_{ \indexGridLines , \symPhaseC \symPhaseA}^{\seriesss}    }}
\newcommand{\gijksCB}[0]{\textcolor{\paramcolor}{\symConductance_{ \indexGridLines , \symPhaseC \symPhaseB}^{\seriesss}    }}
\newcommand{\gijksCC}[0]{\textcolor{\paramcolor}{\symConductance_{ \indexGridLines , \symPhaseC \symPhaseC}^{\seriesss}    }}

\newcommand{\bijksAA}[0]{\textcolor{\paramcolor}{\symSusceptance_{ \indexGridLines , \symPhaseA \symPhaseA}^{\seriesss}    }}
\newcommand{\bijksAB}[0]{\textcolor{\paramcolor}{\symSusceptance_{ \indexGridLines , \symPhaseA \symPhaseB}^{\seriesss}    }}
\newcommand{\bijksAC}[0]{\textcolor{\paramcolor}{\symSusceptance_{ \indexGridLines , \symPhaseA \symPhaseC}^{\seriesss}    }}
\newcommand{\bijksBA}[0]{\textcolor{\paramcolor}{\symSusceptance_{ \indexGridLines , \symPhaseB \symPhaseA}^{\seriesss}    }}
\newcommand{\bijksBB}[0]{\textcolor{\paramcolor}{\symSusceptance_{ \indexGridLines , \symPhaseB \symPhaseB}^{\seriesss}    }}
\newcommand{\bijksBC}[0]{\textcolor{\paramcolor}{\symSusceptance_{ \indexGridLines , \symPhaseB \symPhaseC}^{\seriesss}    }}
\newcommand{\bijksCA}[0]{\textcolor{\paramcolor}{\symSusceptance_{ \indexGridLines , \symPhaseC \symPhaseA}^{\seriesss}    }}
\newcommand{\bijksCB}[0]{\textcolor{\paramcolor}{\symSusceptance_{ \indexGridLines , \symPhaseC \symPhaseB}^{\seriesss}    }}
\newcommand{\bijksCC}[0]{\textcolor{\paramcolor}{\symSusceptance_{ \indexGridLines , \symPhaseC \symPhaseC}^{\seriesss}    }}

\newcommand{\Ilijh}[0]{\textcolor{\complexcolor}{\mathbf{\symCurrent}_{\indexGridLines\indexGridNode \indexGridNodeTwo,\symHarmonic}^{}  }}
\newcommand{\Iljih}[0]{\textcolor{\complexcolor}{\mathbf{\symCurrent}_{\indexGridLines\indexGridNodeTwo\indexGridNode ,\symHarmonic}^{}  }}
\newcommand{\IlijAh}[0]{\textcolor{\complexcolor}{\symCurrent_{\indexGridLines\indexGridNode \indexGridNodeTwo, \symPhaseA,\symHarmonic}^{}  }}
\newcommand{\IlijBh}[0]{\textcolor{\complexcolor}{\symCurrent_{\indexGridLines\indexGridNode \indexGridNodeTwo, \symPhaseB,\symHarmonic}^{}  }}
\newcommand{\IlijCh}[0]{\textcolor{\complexcolor}{\symCurrent_{\indexGridLines\indexGridNode \indexGridNodeTwo, \symPhaseC,\symHarmonic}^{}  }}

\newcommand{\Slijh}[0]{\textcolor{\complexcolor}{\mathbf{\symApparentPower}_{\indexGridLines\indexGridNode \indexGridNodeTwo,\symHarmonic}^{}  }}
\newcommand{\Sljih}[0]{\textcolor{\complexcolor}{\mathbf{\symApparentPower}_{\indexGridLines \indexGridNodeTwo \indexGridNode,\symHarmonic}^{}  }}

\newcommand{\Sslijh}[0]{\textcolor{\complexcolor}{\mathbf{\symApparentPower}_{\indexGridLines\indexGridNode \indexGridNodeTwo,\symHarmonic}^{\seriesss}  }}
\newcommand{\Ssljih}[0]{\textcolor{\complexcolor}{\mathbf{\symApparentPower}_{\indexGridLines \indexGridNodeTwo \indexGridNode,\symHarmonic}^{\seriesss}  }}

\newcommand{\Islijh}[0]{\textcolor{\complexcolor}{\mathbf{\symCurrent}_{\indexGridLines\indexGridNode \indexGridNodeTwo,\symHarmonic}^{\seriesss}  }}
\newcommand{\Isljih}[0]{\textcolor{\complexcolor}{\mathbf{\symCurrent}_{\indexGridLines \indexGridNodeTwo \indexGridNode,\symHarmonic}^{\seriesss}  }}
\newcommand{\IslijAh}[0]{\textcolor{\complexcolor}{\symCurrent_{\indexGridLines\indexGridNode \indexGridNodeTwo, \symPhaseA,\symHarmonic}^{\seriesss}  }}
\newcommand{\IslijBh}[0]{\textcolor{\complexcolor}{\symCurrent_{\indexGridLines\indexGridNode \indexGridNodeTwo, \symPhaseB,\symHarmonic}^{\seriesss}  }}
\newcommand{\IslijCh}[0]{\textcolor{\complexcolor}{\symCurrent_{\indexGridLines\indexGridNode \indexGridNodeTwo, \symPhaseC,\symHarmonic}^{\seriesss}  }}

\newcommand{\Ishlijh}[0]{\textcolor{\complexcolor}{\mathbf{\symCurrent}_{\indexGridLines\indexGridNode \indexGridNodeTwo,\symHarmonic}^{\shuntss}  }}
\newcommand{\Ishljih}[0]{\textcolor{\complexcolor}{\mathbf{\symCurrent}_{\indexGridLines \indexGridNodeTwo\indexGridNode ,\symHarmonic}^{\shuntss}  }}
\newcommand{\IshlijAh}[0]{\textcolor{\complexcolor}{\symCurrent_{\indexGridLines\indexGridNode \indexGridNodeTwo, \symPhaseA,\symHarmonic}^{\shuntss}  }}
\newcommand{\IshlijBh}[0]{\textcolor{\complexcolor}{\symCurrent_{\indexGridLines\indexGridNode \indexGridNodeTwo, \symPhaseB,\symHarmonic}^{\shuntss}  }}
\newcommand{\IshlijCh}[0]{\textcolor{\complexcolor}{\symCurrent_{\indexGridLines\indexGridNode \indexGridNodeTwo, \symPhaseC,\symHarmonic}^{\shuntss}  }}

\newcommand{\IijkA}[0]{\textcolor{\complexcolor}{\symCurrent_{\indexGridLines\indexGridNode \indexGridNodeTwo, \symPhaseA}^{}  }}
\newcommand{\IijkB}[0]{\textcolor{\complexcolor}{\symCurrent_{\indexGridLines\indexGridNode \indexGridNodeTwo, \symPhaseB}^{}  }}
\newcommand{\IijkC}[0]{\textcolor{\complexcolor}{\symCurrent_{\indexGridLines\indexGridNode \indexGridNodeTwo, \symPhaseC}^{}  }}
\newcommand{\IijkN}[0]{\textcolor{\complexcolor}{\symCurrent_{\indexGridLines\indexGridNode \indexGridNodeTwo, \symPhaseN}^{}  }}
\newcommand{\IijkG}[0]{\textcolor{\complexcolor}{\symCurrent_{\indexGridLines\indexGridNode \indexGridNodeTwo, \symPhaseG}^{}  }}

\newcommand{\kifkG}[0]{\textcolor{\complexcolor}{\symCurrent_{k\indexGridNode f, \symPhaseG}^{}  }}

\newcommand{\IijkNreal}[0]{\textcolor{black}{\symCurrent_{\indexGridLines\indexGridNode \indexGridNodeTwo, \symPhaseN}^{\text{re}} }}
\newcommand{\IijkNimag}[0]{\textcolor{black}{\symCurrent_{\indexGridLines\indexGridNode \indexGridNodeTwo, \symPhaseN}^{\text{im}} }}

\newcommand{\IikNG}[0]{\textcolor{\complexcolor}{\symCurrent_{\indexGridNode , \symPhaseN\symPhaseG}}}
\newcommand{\IikNGreal}[0]{\textcolor{black}{\symCurrent_{\indexGridNode , \symPhaseN\symPhaseG}^{\text{re}} }}
\newcommand{\IikNGimag}[0]{\textcolor{black}{\symCurrent_{\indexGridNode , \symPhaseN\symPhaseG}^{\text{im}}}}

\newcommand{\Iijkzero}[0]{\textcolor{\complexcolor}{\symCurrent_{\indexGridLines\indexGridNode \indexGridNodeTwo, 0}^{}  }}
\newcommand{\Iijkone}[0]{\textcolor{\complexcolor}{\symCurrent_{\indexGridLines\indexGridNode \indexGridNodeTwo, 1}^{}  }}
\newcommand{\Iijktwo}[0]{\textcolor{\complexcolor}{\symCurrent_{\indexGridLines\indexGridNode \indexGridNodeTwo, 2}^{}  }}

\newcommand{\IjikA}[0]{\textcolor{\complexcolor}{\symCurrent_{\indexGridLines\indexGridNodeTwo\indexGridNode , \symPhaseA}^{}  }}
\newcommand{\IjikB}[0]{\textcolor{\complexcolor}{\symCurrent_{\indexGridLines\indexGridNodeTwo\indexGridNode  , \symPhaseB}^{}  }}
\newcommand{\IjikC}[0]{\textcolor{\complexcolor}{\symCurrent_{\indexGridLines\indexGridNodeTwo\indexGridNode , \symPhaseC}^{}  }}
\newcommand{\IjikN}[0]{\textcolor{\complexcolor}{\symCurrent_{\indexGridLines\indexGridNodeTwo\indexGridNode , \symPhaseN}^{}  }}
\newcommand{\IjikG}[0]{\textcolor{\complexcolor}{\symCurrent_{\indexGridLines\indexGridNodeTwo\indexGridNode , \symPhaseG}^{}  }}

\newcommand{\IijkrefA}[0]{\textcolor{\complexparamcolor}{\symCurrent_{\indexGridLines\indexGridNode \indexGridNodeTwo,  \symPhaseA}^{\refss}  }}
\newcommand{\IijkrefB}[0]{\textcolor{\complexparamcolor}{\symCurrent_{\indexGridLines\indexGridNode \indexGridNodeTwo, \symPhaseB}^{\refss}  }}
\newcommand{\IijkrefC}[0]{\textcolor{\complexparamcolor}{\symCurrent_{\indexGridLines\indexGridNode \indexGridNodeTwo, \symPhaseC}^{\refss}  }}
\newcommand{\IijksrefA}[0]{\textcolor{\complexparamcolor}{\symCurrent_{\indexGridLines\indexGridNode \indexGridNodeTwo, \seriesss, \symPhaseA}^{\refss}  }}
\newcommand{\IijksrefB}[0]{\textcolor{\complexparamcolor}{\symCurrent_{\indexGridLines\indexGridNode \indexGridNodeTwo, \seriesss, \symPhaseB}^{\refss}  }}
\newcommand{\IijksrefC}[0]{\textcolor{\complexparamcolor}{\symCurrent_{\indexGridLines\indexGridNode \indexGridNodeTwo, \seriesss, \symPhaseC}^{\refss}  }}

\newcommand{\IijksA}[0]{\textcolor{\complexcolor}{\symCurrent_{\indexGridLines\indexGridNode \indexGridNodeTwo, \symPhaseA}^{ \seriesss}  }}
\newcommand{\IijksB}[0]{\textcolor{\complexcolor}{\symCurrent_{\indexGridLines\indexGridNode \indexGridNodeTwo, \symPhaseB}^{\seriesss}  }}
\newcommand{\IijksC}[0]{\textcolor{\complexcolor}{\symCurrent_{\indexGridLines\indexGridNode \indexGridNodeTwo, \symPhaseC}^{\seriesss}  }}
\newcommand{\IijksN}[0]{\textcolor{\complexcolor}{\symCurrent_{\indexGridLines\indexGridNode \indexGridNodeTwo, \symPhaseN}^{\seriesss}  }}
\newcommand{\IijksG}[0]{\textcolor{\complexcolor}{\symCurrent_{\indexGridLines\indexGridNode \indexGridNodeTwo, \symPhaseG}^{\seriesss}  }}

\newcommand{\IjiksA}[0]{\textcolor{\complexcolor}{\symCurrent_{\indexGridLines\indexGridNodeTwo\indexGridNode , \symPhaseA}^{ \seriesss}  }}
\newcommand{\IjiksB}[0]{\textcolor{\complexcolor}{\symCurrent_{\indexGridLines\indexGridNodeTwo\indexGridNode , \symPhaseB}^{\seriesss}  }}
\newcommand{\IjiksC}[0]{\textcolor{\complexcolor}{\symCurrent_{\indexGridLines\indexGridNodeTwo\indexGridNode , \symPhaseC}^{\seriesss}  }}
\newcommand{\IjiksN}[0]{\textcolor{\complexcolor}{\symCurrent_{\indexGridLines\indexGridNodeTwo\indexGridNode , \symPhaseN}^{\seriesss}  }}

\newcommand{\IijkshA}[0]{\textcolor{\complexcolor}{\symCurrent_{\indexGridLines\indexGridNode \indexGridNodeTwo, \symPhaseA}^{ \shuntss}  }}
\newcommand{\IijkshB}[0]{\textcolor{\complexcolor}{\symCurrent_{\indexGridLines\indexGridNode \indexGridNodeTwo, \symPhaseB}^{\shuntss}  }}
\newcommand{\IijkshC}[0]{\textcolor{\complexcolor}{\symCurrent_{\indexGridLines\indexGridNode \indexGridNodeTwo, \symPhaseC}^{\shuntss}  }}
\newcommand{\IijkshN}[0]{\textcolor{\complexcolor}{\symCurrent_{\indexGridLines\indexGridNode \indexGridNodeTwo, \symPhaseN}^{\shuntss}  }}
\newcommand{\IijkshG}[0]{\textcolor{\complexcolor}{\symCurrent_{\indexGridLines\indexGridNode \indexGridNodeTwo, \symPhaseG}^{\shuntss}  }}

\newcommand{\IjikshA}[0]{\textcolor{\complexcolor}{\symCurrent_{\indexGridLines\indexGridNodeTwo\indexGridNode , \symPhaseA}^{ \shuntss}  }}
\newcommand{\IjikshB}[0]{\textcolor{\complexcolor}{\symCurrent_{\indexGridLines\indexGridNodeTwo\indexGridNode , \symPhaseB}^{\shuntss}  }}
\newcommand{\IjikshC}[0]{\textcolor{\complexcolor}{\symCurrent_{\indexGridLines\indexGridNodeTwo\indexGridNode , \symPhaseC}^{\shuntss}  }}\newcommand{\IjikshN}[0]{\textcolor{\complexcolor}{\symCurrent_{\indexGridLines \indexGridNodeTwo \indexGridNode , \symPhaseN}^{\shuntss}  }}

\newcommand{\IikA}[0]{\textcolor{\complexcolor}{\symCurrent_{\indexGridNode , \symPhaseA}^{}  }}
\newcommand{\IikB}[0]{\textcolor{\complexcolor}{\symCurrent_{\indexGridNode , \symPhaseB}^{}  }}
\newcommand{\IikC}[0]{\textcolor{\complexcolor}{\symCurrent_{\indexGridNode , \symPhaseC}^{}  }}

\newcommand{\yijkshAA}[0]{\textcolor{\complexparamcolor}{\symAdmittance_{\indexGridLines \indexGridNode \indexGridNodeTwo, \symPhaseA \symPhaseA}^{\shuntss}    }}
\newcommand{\yijkshAB}[0]{\textcolor{\complexparamcolor}{\symAdmittance_{\indexGridLines \indexGridNode \indexGridNodeTwo, \symPhaseA \symPhaseB}^{\shuntss}    }}
\newcommand{\yijkshAC}[0]{\textcolor{\complexparamcolor}{\symAdmittance_{\indexGridLines \indexGridNode \indexGridNodeTwo, \symPhaseA \symPhaseC}^{\shuntss}    }}
\newcommand{\yijkshAN}[0]{\textcolor{\complexparamcolor}{\symAdmittance_{\indexGridLines \indexGridNode \indexGridNodeTwo, \symPhaseA \symPhaseN}^{\shuntss}    }}
\newcommand{\yijkshAG}[0]{\textcolor{\complexparamcolor}{\symAdmittance_{\indexGridLines \indexGridNode \indexGridNodeTwo, \symPhaseA \symPhaseG}^{\shuntss}    }}

\newcommand{\yijkshBA}[0]{\textcolor{\complexparamcolor}{\symAdmittance_{\indexGridLines \indexGridNode \indexGridNodeTwo, \symPhaseB \symPhaseA}^{\shuntss}    }}
\newcommand{\yijkshBB}[0]{\textcolor{\complexparamcolor}{\symAdmittance_{\indexGridLines \indexGridNode \indexGridNodeTwo, \symPhaseB \symPhaseB}^{\shuntss}    }}
\newcommand{\yijkshBC}[0]{\textcolor{\complexparamcolor}{\symAdmittance_{\indexGridLines \indexGridNode \indexGridNodeTwo, \symPhaseB \symPhaseC}^{\shuntss}    }}
\newcommand{\yijkshBN}[0]{\textcolor{\complexparamcolor}{\symAdmittance_{\indexGridLines \indexGridNode \indexGridNodeTwo, \symPhaseB \symPhaseN}^{\shuntss}    }}
\newcommand{\yijkshBG}[0]{\textcolor{\complexparamcolor}{\symAdmittance_{\indexGridLines \indexGridNode \indexGridNodeTwo, \symPhaseB \symPhaseG}^{\shuntss}    }}

\newcommand{\yijkshCA}[0]{\textcolor{\complexparamcolor}{\symAdmittance_{\indexGridLines \indexGridNode \indexGridNodeTwo, \symPhaseC \symPhaseA}^{\shuntss}    }}
\newcommand{\yijkshCB}[0]{\textcolor{\complexparamcolor}{\symAdmittance_{\indexGridLines \indexGridNode \indexGridNodeTwo, \symPhaseC \symPhaseB}^{\shuntss}    }}
\newcommand{\yijkshCC}[0]{\textcolor{\complexparamcolor}{\symAdmittance_{\indexGridLines \indexGridNode \indexGridNodeTwo, \symPhaseC \symPhaseC}^{\shuntss}    }}\newcommand{\yijkshCN}[0]{\textcolor{\complexparamcolor}{\symAdmittance_{\indexGridLines \indexGridNode \indexGridNodeTwo, \symPhaseC \symPhaseN}^{\shuntss}    }}\newcommand{\yijkshCG}[0]{\textcolor{\complexparamcolor}{\symAdmittance_{\indexGridLines \indexGridNode \indexGridNodeTwo, \symPhaseC \symPhaseG}^{\shuntss}    }}

\newcommand{\yijkshNA}[0]{\textcolor{\complexparamcolor}{\symAdmittance_{\indexGridLines \indexGridNode \indexGridNodeTwo, \symPhaseN \symPhaseA}^{\shuntss}    }}
\newcommand{\yijkshNB}[0]{\textcolor{\complexparamcolor}{\symAdmittance_{\indexGridLines \indexGridNode \indexGridNodeTwo, \symPhaseN \symPhaseB}^{\shuntss}    }}
\newcommand{\yijkshNC}[0]{\textcolor{\complexparamcolor}{\symAdmittance_{\indexGridLines \indexGridNode \indexGridNodeTwo, \symPhaseN \symPhaseC}^{\shuntss}    }}\newcommand{\yijkshNN}[0]{\textcolor{\complexparamcolor}{\symAdmittance_{\indexGridLines \indexGridNode \indexGridNodeTwo, \symPhaseN \symPhaseN}^{\shuntss}    }}\newcommand{\yijkshNG}[0]{\textcolor{\complexparamcolor}{\symAdmittance_{\indexGridLines \indexGridNode \indexGridNodeTwo, \symPhaseN \symPhaseG}^{\shuntss}    }}

\newcommand{\yijkshGA}[0]{\textcolor{\complexparamcolor}{\symAdmittance_{\indexGridLines \indexGridNode \indexGridNodeTwo, \symPhaseG \symPhaseA}^{\shuntss}    }}
\newcommand{\yijkshGB}[0]{\textcolor{\complexparamcolor}{\symAdmittance_{\indexGridLines \indexGridNode \indexGridNodeTwo, \symPhaseG \symPhaseB}^{\shuntss}    }}
\newcommand{\yijkshGC}[0]{\textcolor{\complexparamcolor}{\symAdmittance_{\indexGridLines \indexGridNode \indexGridNodeTwo, \symPhaseG \symPhaseC}^{\shuntss}    }}\newcommand{\yijkshGN}[0]{\textcolor{\complexparamcolor}{\symAdmittance_{\indexGridLines \indexGridNode \indexGridNodeTwo, \symPhaseG \symPhaseN}^{\shuntss}    }}\newcommand{\yijkshGG}[0]{\textcolor{\complexparamcolor}{\symAdmittance_{\indexGridLines \indexGridNode \indexGridNodeTwo, \symPhaseG \symPhaseG}^{\shuntss}    }}

\newcommand{\gijkshAA}[0]{\textcolor{\paramcolor}{\symConductance_{ \indexGridLines\indexGridNode \indexGridNodeTwo, \symPhaseA \symPhaseA}^{\shuntss}    }}
\newcommand{\gijkshAB}[0]{\textcolor{\paramcolor}{\symConductance_{\indexGridLines \indexGridNode \indexGridNodeTwo, \symPhaseA \symPhaseB}^{\shuntss}    }}
\newcommand{\gijkshAC}[0]{\textcolor{\paramcolor}{\symConductance_{\indexGridLines \indexGridNode \indexGridNodeTwo, \symPhaseA \symPhaseC}^{\shuntss}    }}

\newcommand{\gijkshBA}[0]{\textcolor{\paramcolor}{\symConductance_{ \indexGridLines\indexGridNode \indexGridNodeTwo, \symPhaseB \symPhaseA}^{\shuntss}    }}
\newcommand{\gijkshBB}[0]{\textcolor{\paramcolor}{\symConductance_{\indexGridLines \indexGridNode \indexGridNodeTwo, \symPhaseB \symPhaseB}^{\shuntss}    }}
\newcommand{\gijkshBC}[0]{\textcolor{\paramcolor}{\symConductance_{\indexGridLines \indexGridNode \indexGridNodeTwo, \symPhaseB \symPhaseC}^{\shuntss}    }}

\newcommand{\gijkshCA}[0]{\textcolor{\paramcolor}{\symConductance_{ \indexGridLines\indexGridNode \indexGridNodeTwo, \symPhaseC \symPhaseA}^{\shuntss}    }}
\newcommand{\gijkshCB}[0]{\textcolor{\paramcolor}{\symConductance_{ \indexGridLines\indexGridNode \indexGridNodeTwo, \symPhaseC \symPhaseB}^{\shuntss}    }}
\newcommand{\gijkshCC}[0]{\textcolor{\paramcolor}{\symConductance_{\indexGridLines \indexGridNode \indexGridNodeTwo, \symPhaseC \symPhaseC}^{\shuntss}    }}

\newcommand{\bijkshAA}[0]{\textcolor{\paramcolor}{\symSusceptance_{ \indexGridLines\indexGridNode \indexGridNodeTwo, \symPhaseA \symPhaseA}^{\shuntss}    }}
\newcommand{\bijkshAB}[0]{\textcolor{\paramcolor}{\symSusceptance_{ \indexGridLines\indexGridNode \indexGridNodeTwo, \symPhaseA \symPhaseB}^{\shuntss}    }}
\newcommand{\bijkshAC}[0]{\textcolor{\paramcolor}{\symSusceptance_{ \indexGridLines\indexGridNode \indexGridNodeTwo, \symPhaseA \symPhaseC}^{\shuntss}    }}

\newcommand{\bijkshBA}[0]{\textcolor{\paramcolor}{\symSusceptance_{ \indexGridLines\indexGridNode \indexGridNodeTwo, \symPhaseB \symPhaseA}^{\shuntss}    }}
\newcommand{\bijkshBB}[0]{\textcolor{\paramcolor}{\symSusceptance_{ \indexGridLines\indexGridNode \indexGridNodeTwo, \symPhaseB \symPhaseB}^{\shuntss}    }}
\newcommand{\bijkshBC}[0]{\textcolor{\paramcolor}{\symSusceptance_{ \indexGridLines\indexGridNode \indexGridNodeTwo, \symPhaseB \symPhaseC}^{\shuntss}    }}

\newcommand{\bijkshCA}[0]{\textcolor{\paramcolor}{\symSusceptance_{ \indexGridLines\indexGridNode \indexGridNodeTwo, \symPhaseC \symPhaseA}^{\shuntss}    }}
\newcommand{\bijkshCB}[0]{\textcolor{\paramcolor}{\symSusceptance_{\indexGridLines \indexGridNode \indexGridNodeTwo, \symPhaseC \symPhaseB}^{\shuntss}    }}
\newcommand{\bijkshCC}[0]{\textcolor{\paramcolor}{\symSusceptance_{ \indexGridLines\indexGridNode \indexGridNodeTwo, \symPhaseC \symPhaseC}^{\shuntss}    }}

\newcommand{\yjikshAA}[0]{\textcolor{\complexparamcolor}{\symAdmittance_{\indexGridLines\indexGridNodeTwo \indexGridNode , \symPhaseA \symPhaseA}^{\shuntss}    }}
\newcommand{\yjikshAB}[0]{\textcolor{\complexparamcolor}{\symAdmittance_{\indexGridLines\indexGridNodeTwo \indexGridNode , \symPhaseA \symPhaseB}^{\shuntss}    }}
\newcommand{\yjikshAC}[0]{\textcolor{\complexparamcolor}{\symAdmittance_{\indexGridLines \indexGridNodeTwo\indexGridNode , \symPhaseA \symPhaseC}^{\shuntss}    }}
\newcommand{\yjikshAN}[0]{\textcolor{\complexparamcolor}{\symAdmittance_{\indexGridLines \indexGridNodeTwo\indexGridNode , \symPhaseA \symPhaseN}^{\shuntss}    }}
\newcommand{\yjikshAG}[0]{\textcolor{\complexparamcolor}{\symAdmittance_{\indexGridLines \indexGridNodeTwo\indexGridNode , \symPhaseA \symPhaseG}^{\shuntss}    }}

\newcommand{\yjikshBA}[0]{\textcolor{\complexparamcolor}{\symAdmittance_{\indexGridLines\indexGridNodeTwo \indexGridNode , \symPhaseB \symPhaseA}^{\shuntss}    }}
\newcommand{\yjikshBB}[0]{\textcolor{\complexparamcolor}{\symAdmittance_{\indexGridLines\indexGridNodeTwo \indexGridNode , \symPhaseB \symPhaseB}^{\shuntss}    }}
\newcommand{\yjikshBC}[0]{\textcolor{\complexparamcolor}{\symAdmittance_{\indexGridLines \indexGridNodeTwo\indexGridNode , \symPhaseB \symPhaseC}^{\shuntss}    }}
\newcommand{\yjikshBN}[0]{\textcolor{\complexparamcolor}{\symAdmittance_{\indexGridLines \indexGridNodeTwo\indexGridNode , \symPhaseB \symPhaseN}^{\shuntss}    }}
\newcommand{\yjikshBG}[0]{\textcolor{\complexparamcolor}{\symAdmittance_{\indexGridLines \indexGridNodeTwo\indexGridNode , \symPhaseB \symPhaseG}^{\shuntss}    }}

\newcommand{\yjikshCA}[0]{\textcolor{\complexparamcolor}{\symAdmittance_{\indexGridLines\indexGridNodeTwo \indexGridNode , \symPhaseC \symPhaseA}^{\shuntss}    }}
\newcommand{\yjikshCB}[0]{\textcolor{\complexparamcolor}{\symAdmittance_{\indexGridLines\indexGridNodeTwo \indexGridNode , \symPhaseC \symPhaseB}^{\shuntss}    }}
\newcommand{\yjikshCC}[0]{\textcolor{\complexparamcolor}{\symAdmittance_{\indexGridLines \indexGridNodeTwo\indexGridNode , \symPhaseC \symPhaseC}^{\shuntss}    }}
\newcommand{\yjikshCN}[0]{\textcolor{\complexparamcolor}{\symAdmittance_{\indexGridLines \indexGridNodeTwo\indexGridNode , \symPhaseC \symPhaseN}^{\shuntss}    }}
\newcommand{\yjikshCG}[0]{\textcolor{\complexparamcolor}{\symAdmittance_{\indexGridLines \indexGridNodeTwo\indexGridNode , \symPhaseC \symPhaseG}^{\shuntss}    }}

\newcommand{\yjikshNA}[0]{\textcolor{\complexparamcolor}{\symAdmittance_{\indexGridLines\indexGridNodeTwo \indexGridNode , \symPhaseN \symPhaseA}^{\shuntss}    }}
\newcommand{\yjikshNB}[0]{\textcolor{\complexparamcolor}{\symAdmittance_{\indexGridLines\indexGridNodeTwo \indexGridNode , \symPhaseN \symPhaseB}^{\shuntss}    }}
\newcommand{\yjikshNC}[0]{\textcolor{\complexparamcolor}{\symAdmittance_{\indexGridLines \indexGridNodeTwo\indexGridNode , \symPhaseN \symPhaseC}^{\shuntss}    }}
\newcommand{\yjikshNN}[0]{\textcolor{\complexparamcolor}{\symAdmittance_{\indexGridLines \indexGridNodeTwo\indexGridNode , \symPhaseN \symPhaseN}^{\shuntss}    }}
\newcommand{\yjikshNG}[0]{\textcolor{\complexparamcolor}{\symAdmittance_{\indexGridLines \indexGridNodeTwo\indexGridNode , \symPhaseN \symPhaseG}^{\shuntss}    }}

\newcommand{\yjikshGA}[0]{\textcolor{\complexparamcolor}{\symAdmittance_{\indexGridLines\indexGridNodeTwo \indexGridNode , \symPhaseG \symPhaseA}^{\shuntss}    }}
\newcommand{\yjikshGB}[0]{\textcolor{\complexparamcolor}{\symAdmittance_{\indexGridLines\indexGridNodeTwo \indexGridNode , \symPhaseG \symPhaseB}^{\shuntss}    }}
\newcommand{\yjikshGC}[0]{\textcolor{\complexparamcolor}{\symAdmittance_{\indexGridLines \indexGridNodeTwo\indexGridNode , \symPhaseG \symPhaseC}^{\shuntss}    }}
\newcommand{\yjikshGN}[0]{\textcolor{\complexparamcolor}{\symAdmittance_{\indexGridLines \indexGridNodeTwo\indexGridNode , \symPhaseG \symPhaseN}^{\shuntss}    }}
\newcommand{\yjikshGG}[0]{\textcolor{\complexparamcolor}{\symAdmittance_{\indexGridLines \indexGridNodeTwo\indexGridNode , \symPhaseG \symPhaseG}^{\shuntss}    }}

\newcommand{\gijkshpp}[0]{\textcolor{\paramcolor}{\symConductance_{\indexGridLines \indexGridNode \indexGridNodeTwo, \indexPhases\indexPhases}^{\shuntss}    }}
\newcommand{\gijkshph}[0]{\textcolor{\paramcolor}{\symConductance_{ \indexGridLines\indexGridNode \indexGridNodeTwo, \indexPhases\indexPhasesTwo}^{\shuntss}    }}

\newcommand{\gijkspp}[0]{\textcolor{\paramcolor}{\symConductance_{ \indexGridLines  , \indexPhases\indexPhases}^{\seriesss}    }}
\newcommand{\gijksph}[0]{\textcolor{\paramcolor}{\symConductance_{ \indexGridLines  , \indexPhases\indexPhasesTwo}^{\seriesss}    }}

\newcommand{\bijkshpp}[0]{\textcolor{\paramcolor}{\symSusceptance_{\indexGridLines \indexGridNode \indexGridNodeTwo, \indexPhases\indexPhases}^{\shuntss}    }}
\newcommand{\bijkspp}[0]{\textcolor{\paramcolor}{\symSusceptance_{ \indexGridLines, \indexPhases\indexPhases}^{\seriesss}    }}
\newcommand{\bijksph}[0]{\textcolor{\paramcolor}{\symSusceptance_{ \indexGridLines, \indexPhases\indexPhasesTwo}^{\seriesss}    }}

\newcommand{\bijkshph}[0]{\textcolor{\paramcolor}{\symSusceptance_{ \indexGridNode \indexGridNodeTwo, \indexPhases\indexPhasesTwo}^{\shuntss}    }}

\newcommand{\zbental}[0]{{z}}
\newcommand{\ybental}[0]{{y}}
\newcommand{\xbental}[0]{{x}}
\newcommand{\xonebental}[0]{\textcolor{black}{\xbental_{1} }}
\newcommand{\xtwobental}[0]{\textcolor{black}{\xbental_{2} }}
\newcommand{\xthreebental}[0]{\textcolor{black}{\xbental_{3} }}
\newcommand{\xfourbental}[0]{\textcolor{black}{\xbental_{4} }}
\newcommand{\xfivebental}[0]{\textcolor{black}{\xbental_{5} }}
\newcommand{\xsixbental}[0]{\textcolor{black}{\xbental_{6} }}
\newcommand{\xsevenbental}[0]{\textcolor{black}{\xbental_{7} }}
\newcommand{\xibental}[0]{\textcolor{black}{\xbental_{\counterbental} }}
\newcommand{\upperbental}[0]{{U}}
\newcommand{\lowerbental}[0]{{L}}
\newcommand{\xilowerbental}[0]{\textcolor{\paramcolor}{\xibental^{ \lowerbental } }}
\newcommand{\xiupperbental}[0]{\textcolor{\paramcolor}{\xibental^{ \upperbental } }}
\newcommand{\xlowerbental}[1]{\textcolor{\boundscolor}{\xbental_{#1}^{ \lowerbental } }}
\newcommand{\xupperbental}[1]{\textcolor{\boundscolor}{\xbental_{#1}^{ \upperbental } }}

\newcommand{\sigmabental}[1]{\textcolor{black}{\sigma_{(#1)}}}
\newcommand{\etabental}[1]{\textcolor{black}{\eta_{(#1)}}}
\newcommand{\counterbental}[0]{\textcolor{\paramcolor}{i}}
\newcommand{\countbental}[0]{\textcolor{\paramcolor}{\nu}}
\newcommand{\errorbental}[0]{\textcolor{\paramcolor}{\epsilon}}
\newcommand{\countpoly}[0]{\textcolor{\paramcolor}{\mu}}

\newcommand{\circumscribedopoly}[0]{\textcolor{\paramcolor}{\symSetting_{\circumss}}}
\newcommand{\rotatedopoly}[0]{\textcolor{\paramcolor}{\symSetting_{\rotatedss}}}
\newcommand{\genericbinary}[1]{\textcolor{black}{b_{#1}}}


\newcommand{\pipediameterijk}[0]{ \textcolor{\paramcolor}{  \symFluidPipeDiameter_{\indexGridLines}     }}
\newcommand{\pipeareaijk}[0]{ \textcolor{\paramcolor}{  \symFluidPipeArea_{\indexGridLines}     }}

\newcommand{\fluidpressure}[1]{\textcolor{black}{\symPressure_{#1}}}
\newcommand{\fluidpressureik}[0]{\textcolor{black}{\symPressure_{\indexGridNode }    }}
\newcommand{\fluidpressurejk}[0]{\textcolor{black}{\symPressure_{\indexGridNodeTwo }    }}

\newcommand{\fluidpressureikmin}[0]{\textcolor{\paramcolor}{\symPressure_{\indexGridNode }^{\minss}    }}
\newcommand{\fluidpressureikmax}[0]{\textcolor{\paramcolor}{\symPressure_{\indexGridNode }^{\maxss}    }}
\newcommand{\fluidpressurejkmin}[0]{\textcolor{\paramcolor}{\symPressure_{\indexGridNodeTwo }^{\minss}    }}
\newcommand{\fluidpressurejkmax}[0]{\textcolor{\paramcolor}{\symPressure_{\indexGridNodeTwo }^{\maxss}    }}

\newcommand{\fluidpressuresqik}[0]{\textcolor{black}{\symPressureSquared_{\indexGridNode }    }}
\newcommand{\fluidpressuresqjk}[0]{\textcolor{black}{\symPressureSquared_{\indexGridNodeTwo }    }}

\newcommand{\fluidpressuresqikmin}[0]{\textcolor{\paramcolor}{\symPressureSquared_{\indexGridNode }^{\minss}    }}
\newcommand{\fluidpressuresqjkmin}[0]{\textcolor{\paramcolor}{\symPressureSquared_{\indexGridNodeTwo }^{\minss}    }}
\newcommand{\fluidpressuresqikmax}[0]{\textcolor{\paramcolor}{\symPressureSquared_{\indexGridNode }^{\maxss}    }}
\newcommand{\fluidpressuresqjkmax}[0]{\textcolor{\paramcolor}{\symPressureSquared_{\indexGridNodeTwo }^{\maxss}    }}

\newcommand{\fluidflow}[1]{\textcolor{black}{\symVolumeFlow_{#1}}}
\newcommand{\fluidflowijk}[0]{\textcolor{black}{\fluidflow{\indexGridNode \indexGridNodeTwo}   }}

\newcommand{\fluidmassflow}[1]{\textcolor{black}{\symMassFlow_{#1}}}
\newcommand{\fluidmassflowijk}[0]{\textcolor{black}{\fluidmassflow{\indexGridNode \indexGridNodeTwo} }}

\newcommand{\fluidmassflowdirectionone}[0]{\textcolor{black}{y_{ij}^{+}}}
\newcommand{\fluidmassflowdirectiontwo}[0]{\textcolor{black}{y_{ij}^{-}}}

\newcommand{\fluiddensity}[1]{{\symDensity_{#1}  }}
\newcommand{\fluiddensityijk}[0]{\textcolor{\paramcolor}{\fluiddensity{\indexGridLines}   }}

\newcommand{\fluidviscosity}[1]{{\symViscosity_{#1}  }}
\newcommand{\fluidviscosityijk}[0]{\textcolor{\paramcolor}{\fluidviscosity{\indexGridLines}   }}

\newcommand{\heatflow}[1]{\textcolor{black}{\symHeatFlow_{#1}}}
\newcommand{\heatflowijk}[0]{\textcolor{black}{\symHeatFlow_{\indexGridNode \indexGridNodeTwo} }}

\newcommand{\heatresistance}[1]{\textcolor{\paramcolor}{\symHeatResistance_{#1}}}
\newcommand{\heatresistanceijk}[0]{ \textcolor{\paramcolor}{ \heatresistance{\indexGridLines}  }}   

\newcommand{\heattempijk}[0]{\textcolor{black}{\symTemperature_{\indexGridNode \indexGridNodeTwo}   }}
\newcommand{\heattempik}[0]{\textcolor{black}{\symTemperature_{\indexGridNode }   }}
\newcommand{\heattempjk}[0]{\textcolor{black}{\symTemperature_{ \indexGridNodeTwo}   }}

\newcommand{\fluidhagenresistance}[0]{ \textcolor{\paramcolor}{ w^{\hagenposeuilless}_{\indexGridLines}  }    }
\newcommand{\fluiddarcyresistance}[0]{ \textcolor{\paramcolor}{ w^{\darcyweisbachss}_{\indexGridLines}   }    }

\newcommand{\heatcapacity}[0]{ \textcolor{\paramcolor}{ \symHeatCapacity^{\text{\symPressure}} }}   


\newcommand{\symFlowVariable}[0]{   F }
\newcommand{\symFlowVariableSquare}[0]{   f }
\newcommand{\symPotentialVariable}[0]{ U }
\newcommand{\symPotentialAngleVariable}[0]{ \theta }
\newcommand{\symPotentialVariableSquare}[0]{ u }
\newcommand{\setFlowVariable}[0]{   \mathcal{F} }
\newcommand{\setPotentialVariable}[0]{ \mathcal{U} }
\newcommand{\indexFlowVariable}[0]{  f }
\newcommand{\indexPotentialVariable}[0]{ u }
\newcommand{\indexDomains}[0]{ d }
\newcommand{\setDomains}[0]{\mathcal{D}  }
\newcommand{\setDomainsExt}[0]{\mathcal{D}_0  }
\newcommand{\potentiali}[0]{ \symPotentialVariable_{\indexGridNode} }
\newcommand{\potentialz}[0]{ \symPotentialVariable_{\indexZIP} }
\newcommand{\potentialzlin}[0]{ \symPotentialVariable_{\indexZIP}^{\linss} }
\newcommand{\potentialilin}[0]{ \symPotentialVariable_{\indexGridNode}^{\linss} }
\newcommand{\potentialj}[0]{ \symPotentialVariable_{\indexGridNodeTwo} }
\newcommand{\potentialanglei}[0]{ \symPotentialAngleVariable_{\indexGridNode} }
\newcommand{\potentialanglej}[0]{ \symPotentialAngleVariable_{\indexGridNodeTwo} }
\newcommand{\potentialref}[0]{ \textcolor{\paramcolor}{ \symPotentialVariable^{\refss} }}
\newcommand{\potentialiref}[0]{ \textcolor{\paramcolor}{ \symPotentialVariable_{\indexGridNode}^{\refss} }}
\newcommand{\potentialzref}[0]{ \textcolor{\paramcolor}{ \symPotentialVariable_{\indexZIP}^{\refss} }}
\newcommand{\flowij}[0]{ \symFlowVariable_{\indexGridNode \indexGridNodeTwo} }
\newcommand{\flowi}[0]{ \symFlowVariable_{\indexGridNode } }
\newcommand{\flowz}[0]{ \symFlowVariable_{\indexZIP } }
\newcommand{\flowcomplexij}[0]{ \textcolor{\complexcolor}{\symFlowVariable_{\indexGridNode \indexGridNodeTwo} }}
\newcommand{\flowcomplexji}[0]{ \textcolor{\complexcolor}{\symFlowVariable_{\indexGridNodeTwo\indexGridNode } }}
\newcommand{\flowangleij}[0]{ {}{\symFlowVariable_{\indexGridNode \indexGridNodeTwo}^{\imagangle} }}
\newcommand{\flowji}[0]{ \symFlowVariable_{\indexGridNodeTwo \indexGridNode } }
\newcommand{\flowresistancetij}[0]{\textcolor{\paramcolor}{ r_{\indexGridLines } }}
\newcommand{\flowimpedancetij}[0]{\textcolor{\complexcolor}{ z_{\indexGridLines } }}
\newcommand{\flowreactancetij}[0]{\textcolor{\paramcolor}{ x_{\indexGridLines } }}
\newcommand{\potentialsqi}[0]{ \symPotentialVariableSquare_{\indexGridNode} }
\newcommand{\potentialsqz}[0]{ \symPotentialVariableSquare_{\indexZIP} }
\newcommand{\potentialsqj}[0]{ \symPotentialVariableSquare_{\indexGridNodeTwo} }
\newcommand{\flowsqij}[0]{ \symFlowVariableSquare_{\indexGridNode \indexGridNodeTwo} }
\newcommand{\potentialimin}[0]{ \textcolor{\paramcolor}{\symPotentialVariable_{\indexGridNode}^{\minss} }}
\newcommand{\potentialimax}[0]{ \textcolor{\paramcolor}{\symPotentialVariable_{\indexGridNode}^{\maxss} }}
\newcommand{\potentialirated}[0]{ \textcolor{\sizingcolor}{\symPotentialVariable_{\indexGridNode}^{\ratedss} }}
\newcommand{\potentialjmin}[0]{ \textcolor{\paramcolor}{\symPotentialVariable_{\indexGridNodeTwo}^{\minss} }}
\newcommand{\potentialjmax}[0]{ \textcolor{\paramcolor}{\symPotentialVariable_{\indexGridNodeTwo}^{\maxss} }}
\newcommand{\potentialjrated}[0]{ \textcolor{\sizingcolor}{\symPotentialVariable_{\indexGridNodeTwo}^{\ratedss} }}
\newcommand{\potentialijrated}[0]{ \textcolor{\sizingcolor}{\symPotentialVariable_{\indexGridNode\indexGridNodeTwo}^{\ratedss} }}
\newcommand{\potentialjirated}[0]{ \textcolor{\sizingcolor}{\symPotentialVariable_{\indexGridNodeTwo \indexGridNode}^{\ratedss} }}
\newcommand{\flowijmin}[0]{ \textcolor{\paramcolor}{\symFlowVariable_{\indexGridNode \indexGridNodeTwo}^{\minss} }}
\newcommand{\flowijmax}[0]{\textcolor{\paramcolor}{ \symFlowVariable_{\indexGridNode \indexGridNodeTwo}^{\maxss} }}
\newcommand{\flowijrated}[0]{\textcolor{\sizingcolor}{ \symFlowVariable_{\indexGridNode \indexGridNodeTwo}^{\ratedss} }}
\newcommand{\flowjirated}[0]{\textcolor{\sizingcolor}{ \symFlowVariable_{ \indexGridNodeTwo \indexGridNode }^{\ratedss} }}

\newcommand{\powerbalancefactorij}[0]{\textcolor{\paramcolor}{c_{\indexGridLines } }}

\newcommand{\enthalpycombustion}{\textcolor{\paramcolor}{\Delta \symEnthalpy^{\circ}_{\combustionss}}}
\newcommand{\enthalpyformation}{\textcolor{\paramcolor}{\Delta \symEnthalpy^{\circ}_{\formationss}}}

\newcommand{\slackconvexificationij}[0]{{\symSlack_{\indexGridNode \indexGridNodeTwo }^{\PFconvexss} }}
\newcommand{\slackconvexification}[0]{{\symSlack^{\PFconvexss} }}
\newcommand{\slackmaxconvexification}[0]{{\symSlack_{  }^{\maxss} }}
\newcommand{\slackminconvexification}[0]{{\symSlack_{  }^{\minss} }}
\newcommand{\slackabsconvexification}[0]{{\symSlack_{  }^{\absss} }}
\newcommand{\slackcomplexmagnitudeconvexification}[0]{{\symSlack_{  }^{\absss} }}
\newcommand{\slackbental}[0]{{\symSlack_{  }^{\bentallinearization} }}


\newcommand{\mapping}[0]{   \to }
\newcommand{\injection}[0]{   \rightarrowtail }
\newcommand{\surjection}[0]{   \twoheadrightarrow }
\newcommand{\bijection}[0]{   \leftrightarrow }

\newcommand{\symbollistlabelwith}{$\PPNeleclossk$}
\newcommand{\imagnumber}[0]{\textcolor{\paramcolor}{j}}
\newcommand{\imagangle}[0]{\textcolor{\paramcolor}{\angle}}

\title{Computational Analysis of Impedance Transformations for Four-Wire Power Networks with Sparse Neutral Grounding}

\author{Frederik Geth}
\email{frederik.geth@gmail.com}
\orcid{0000-0002-9345-2959}
\affiliation{%
  \institution{CSIRO Energy}
  \city{Newcastle}
  \state{NSW}
  \country{Australia}
}

\author{Rahmat Heidari}
\email{rahmat.heidarihaei@csiro.au}
\orcid{0000-0003-4835-2614} 
\affiliation{%
  \institution{CSIRO Energy}
  \city{Newcastle}
  \country{Australia}
}

\author{Arpan Koirala}
\email{arpan.koirala@kuleuven.be}
\orcid{0000-0003-4826-7137} 
\affiliation{%
  \institution{KU Leuven / EnergyVille}
  \city{Leuven / Genk}
  \country{Belgium}
}

\renewcommand{\shortauthors}{Geth et al.}

\begin{abstract}
In low-voltage distribution networks, the integration of novel energy technologies can be accelerated through advanced optimization-based analytics such as network state estimation and network-constrained dispatch engines for distributed energy resources.
The scalability of distribution network optimization models is challenging due to phase unbalance and neutral voltage rise effects necessitating the use of 4 times as many voltage variables per bus than in transmission systems.
This paper proposes a novel technique to limit this to a factor 3, exploiting common physical features of low-voltage networks specifically, where neutral grounding is  sparse, as it is in many parts of the world.
We validate the proposed approach in \textsc{OpenDSS}, by translating a number of published test cases to the reduced form, and observe that the proposed ``phase-to-neutral'' transformation is highly accurate for the common single-grounded low-voltage network configuration, and provides a high-quality approximation for  other  configurations. 
We finally provide numerical results for unbalanced power flow optimization problems using \textsc{PowerModelsDistribution.jl}, to illustrate the computational speed benefits of a factor of about 1.42. 
\end{abstract}




\maketitle

\section{Introduction}
\subsection{Background on Public Four-wire Networks}

Low-voltage (LV) networks are dealing with the integration of new technologies, including solar PV, home batteries and electric vehicles.
Solving physics-based mathematical optimization problems in such networks is an interesting technology development opportunity\,\cite{Claeys2021b}, with problems such as unbalanced network state estimation, distribution locational marginal pricing, network-constrained battery dispatch and EV charging coordination receiving significant attention by researchers. 

Around the globe, many LV networks have a dual-voltage three-phase + neutral configuration (e.g. 400/230 V or 208/120 V), with grounding of the neutral only at the MV/LV transformer.
This leads to situations where there can be a nontrivial amount of neutral voltage, a.k.a. neutral voltage \emph{shift}, downstream of the transformer.
Neutral voltage shift tends to negatively impact the voltage quality of the network, as it exacerbates \emph{phase-to-neutral} voltage magnitude issues. 
Accurately capturing neutral voltage shift is therefore important, as it may limit the uptake of distributed energy resources. 



Where the neutral is grounded,  neutral current flows via the earth through a grounding rod\footnote{Sometimes first through an impedance.}. 
The regularity of grounding depends on the network design principles, with``multi-grounded'' and ``single-grounded'' being the two main  approaches. 
In three-phase dual-voltage residential networks, grounding is probably near-universal on the wye-side of transformers, where the neutral point would otherwise be floating, which itself would not enable stable phase-to-neutral voltages.
Nevertheless, neutral grounding can be performed additionally throughout the network when network codes allow for it -- or demand it.

Furthermore, in the real world, cables and overhead lines induce current into the ground due to capacitive coupling with earth (nF/km) \cite{Kersting2001b}\cite{Olivier2018}.
However, at 50 or 60 Hz the shunt admittance in LV networks is typically very small, and  the voltage is low as well, so the shunt current contribution is limited in normal operation.
Note that it is not expected that electrical devices owned by end-users inject a significant amount of current directly into the ground - that is typically considered a fault situation and is protected against by residual current devices a.k.a. ground fault circuit interrupters.  
Nevertheless, at medium voltage or higher, the line shunt currents may become significant in the overall current balance of the network at low load.

Industrial engineering software such as \textsc{Digsilent PowerFactory} and open-source power system research platforms such as \textsc{OpenDSS} \cite{Dugan2011} include scalable algorithms for steady-state \emph{simulation} of voltage unbalance, neutral voltage rise in single and multi-grounded systems with capacitive coupling to ground \cite{Ciric2003}.


\subsection{Research Question and Contributions}

In the context of unbalanced power flow \emph{optimization} problems, scalability remains crucial but challenging and commercial solutions are lacking  \cite{Claeys2020}.
Capturing neutral voltage rise usually means representing four complex-valued voltage variables per bus, one for each terminal $t \in \{a,b,c,n\}$.
When one assumes the neutral is pervasively grounded (e.g. justifiable in a multi-grounded philosophy), Kron's reduction of the neutral \cite{Dorfler2011a_edit} introduces minimal error, and enables elimination of the neutral voltage variables, which implies that the optimization problem may be solved faster (due to 25\% fewer variables). 
However, this transformation can not capture neutral voltage rise accurately, and therefore it is not an attractive option for  networks with a sparsely grounded neutral (e.g. single-grounded neutral) with significant load unbalance.

The central research question addressed in this paper is: `Is there an alternative to Kron's reduction of the neutral that can capture neutral voltage rise, but equally enables a reduction in problem size from 4 to 3 complex voltage variables per bus?'
We therefore:
\begin{itemize}
    \item derive a novel ``phase-to-neutral'' transformation  from first principles;
    \item postulate conditions under which it is accurate;
    \item analyze some of its properties;
    \item run numerical experiments, using \textsc{OpenDSS} \cite{Dugan2011}, to validate the approximation quality across a selection of LV and MV networks;
    \item demonstrate the scalability advantages in a power flow optimization context, using \textsc{PowerModelsDistribution} \cite{Fobes2020}. 
\end{itemize}
 
   \section{Phase-to-neutral Transformation}
 We first introduce our notation and derive a self-contained mathematical description of four-wire power flow with neutral grounding.
 The remainder of the section develops a linear transformation $\mathcal{T}: \setComplex^{4\times 4} \rightarrow \setComplex^{3\times 3}$, which is a generalization of ideas presented first in 
 Koirala et al. \cite{Koirala2019} when applied to series impedance matrices.
 
\subsection{Symbols}
Table \ref{table_symbols} summarizes the \emph{vector} and \emph{matrix} symbols used in this article. 
A blue color indicates a complex variable, a brown color indicates a complex parameter, a red color indicates a real parameter.
A bold typeface indicates a vector or matrix variables, a normal typeface is a scalar variable. 
\begin{table}[tph]
\setlength{\tabcolsep}{1pt}
  \centering 
   \caption{Vector and matrix symbols} \label{table_symbols}
    \begin{tabular}{l l  l}
\hline
$\VikSDP$ &(V) & Phase-to-ground voltage vector at bus $i$\\
$\VikSDPwye$& (V) & Phase-to-neutral voltage vector at bus $i$\\
$\IijkSDP$ &(A) & Total current vector in line $l$ in the direction $i \rightarrow j$\\
$\IijskSDP$ &(A) & Series current vector in line $l$ in the direction $i \rightarrow j$\\
$\IijshkSDP$& (A) & Shunt current vector in line $l$ in the direction $i \rightarrow j$\\
$\IloadSDP$& (A) & Load $d$ current vector\\
$\IgenSDP$& (A) & Generator $g$ current vector\\
$\IgroundSDP$& (A) & Bus $i$ grounding current vector\\
$\zijksSDP$& (A) & Series impedance matrix of line $l$\\
$\textcolor{\paramcolor}{\mathbf{T}}$& (-) & Phase-to-neutral transformation matrix \\
\hline
\end{tabular}
\end{table}

   \subsection{Preliminaries} 
\begin{figure}[h]
  \centering
    \includegraphics[scale=0.6]{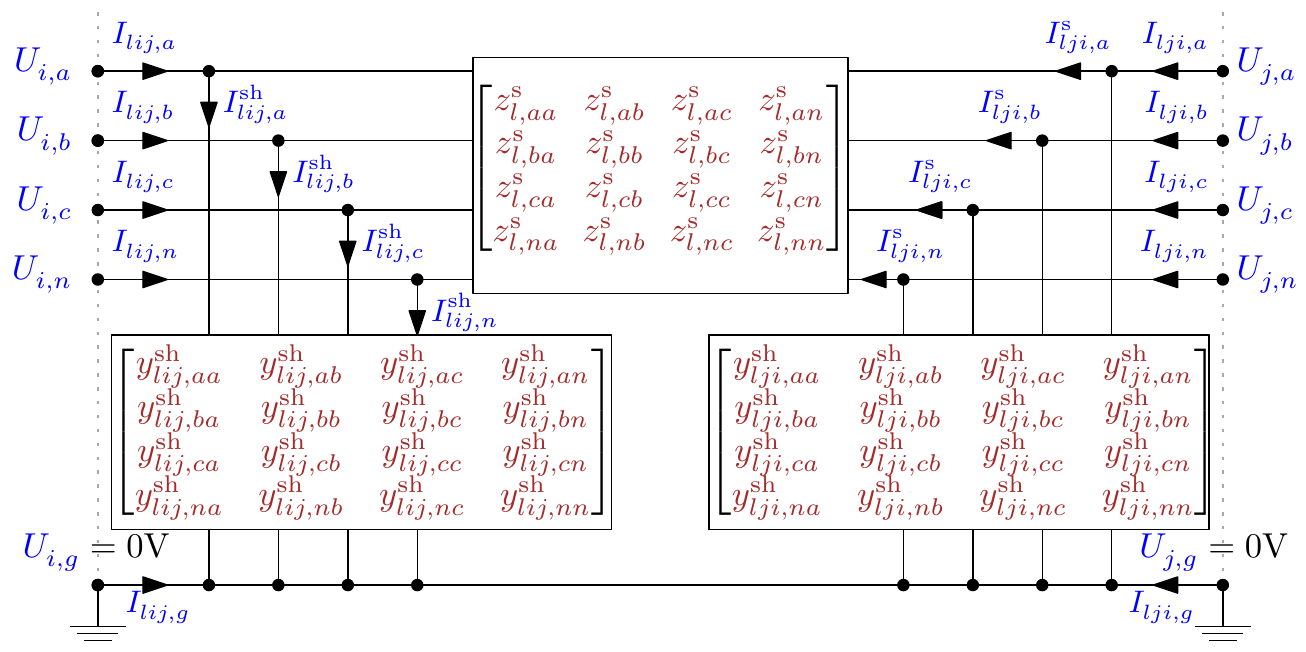}
  \caption{Four-wire line model}  \label{fig_four_wire}
\end{figure}

 Fig. \ref{fig_four_wire} defines the \emph{scalar} symbols for the total, series and shunt current variables, phase-to-ground voltages, impedances for the (asymmetric) $\Pi$-section with phases $\setPhases = \{a,b,c \}$ and neutral $n$, for a line $l$ between buses $i$ and $j$.

The \textbf{s}eries impedance matrix of line $l$ is,
        \begin{equation}
\zijksSDP
=
\begin{bmatrix} 
\zijksAA   &\zijksAB  &\zijksAC &\zijksAN  \\
\zijksBA   &\zijksBB  &\zijksBC &\zijksBN  \\
\zijksCA   &\zijksCB  &\zijksCC &\zijksCN  \\
\zijksNA   &\zijksNB  &\zijksNC &\zijksNN  \\
 \end{bmatrix} . 
  \end{equation}  
In the real world, both the real and the imaginary part of the series impedance are symmetric and positive semidefinite, however we don't need to assume this in the remainder of this work, and therefore define separate parameters for mutual induction, e.g. $a\rightarrow b$ and $b\rightarrow a$. 
The diagonal elements are called the self-impedances and the off-diagonal the mutual impedances. 
  For convenient comparison of impedance transformations later, we define the \emph{do-nothing} transformation $\mathcal{O}(\zijksSDP) = \zijksSDP$.
  
We stack the scalar current variables (Fig. \ref{fig_four_wire}) in vectors to  obtain the vector expression of conservation of current, linking the total current $\IijkSDP $, series current $\IijskSDP$ and shunt current $\IijshkSDP$ in line $l$ in the direction $i\rightarrow j$,
           \begin{equation}
\IijkSDP =  \begin{bmatrix} 
\IijkA    \\
\IijkB     \\
\IijkC  \\
\IijkN  
 \end{bmatrix} = 
 \underbrace{\begin{bmatrix} 
\IijksA    \\
\IijksB     \\
\IijksC  \\
\IijksN  
 \end{bmatrix}}_{\IijskSDP}
 +
  \underbrace{\begin{bmatrix} 
\IijkshA    \\
\IijkshB     \\
\IijkshC  \\
\IijkshN  
 \end{bmatrix}}_{\IijshkSDP} = \IijskSDP + \IijshkSDP .
   \end{equation}   
   
The ground current for line $l$ in the direction  $i\rightarrow j$ is,
      \begin{equation}
\IijkG = -(\IijkA + \IijkB + \IijkC + \IijkN ) = -\textcolor{\paramcolor}{\mathbf{1}}^{\text{T}} \IijkSDP \label{eq_ground_current}.
\end{equation}
   
Next, we stack the depicted scalar  voltage-to-ground\footnote{i.e. w.r.t. $\VikG = 0$ V} variables at bus $i$  (Fig. \ref{fig_four_wire}) in the vector $\VikSDP$ and derive the phase-to-neutral voltages $\VikSDPwye$,
\begin{equation}
\VikSDP = \begin{bmatrix} 
\VikA    \\
\VikB     \\
\VikC  \\
\VikN  
 \end{bmatrix} ,
    \VikSDPwye = 
    \begin{bmatrix} 
\VikAN    \\
\VikBN    \\
\VikCN 
\end{bmatrix} =
  \underbrace{ \begin{bmatrix}
   1 & 0 & 0 &-1\\
    0 & 1 & 0& -1\\
   0 & 0 & 1 & -1\\
   \end{bmatrix}}_{\textcolor{\paramcolor}{\mathbf{T}}}
   \begin{bmatrix} 
\VikA    \\
\VikB     \\
\VikC  \\
\VikN  
 \end{bmatrix}  =
 \textcolor{\paramcolor}{\mathbf{T}} \VikSDP \label{eq_voltage_transf}.
   \end{equation}   
   Note that this linear transformation $ \textcolor{\paramcolor}{\mathbf{T}}$ \emph{does} lose information, i.e. there is no way to recover $\VikN  $ from the variable $\VikSDPwye$.


 

   Ohm's law in matrix form between adjacent buses $i$ and $j$ connected through line $l$ is,
           \begin{equation}
\VjkSDP  =   \VikSDP - \zijksSDP \IijskSDP. \label{eq_ohms_4w}
\end{equation}  
Kirchhoff's current law is applied at all buses $i \in \mathcal{I}$,
           \begin{equation}
\forall i \in \mathcal{I}: \underbrace{\sum_{lij} \IijkSDP}_{\text{branches}} + \underbrace{\sum_{d} \IloadSDP}_{\text{loads}}- \underbrace{\sum_{g} \IgenSDP}_{\text{generators}} +  \underbrace{ \IgroundSDP}_{\text{bus grounding}}= 0. \label{eq_4w_kcl}
   \end{equation}   
At a subset of  buses $i \in \mathcal{I}^{\text{g}} \subset \mathcal{I}$, grounding of the neutral is performed, in which case the neutral current entry in $\IgroundSDP$ can be nonzero, however the phases are not grounded:
           \begin{equation}
 \IgroundSDP =
 \begin{bmatrix}
 0&0&0& I_{i,n}^{\text{g}}
 \end{bmatrix}^{\text{T}}.
   \end{equation}   
  At the other buses (the ungrounded ones), $ \IgroundSDP$ is the zero-vector.
Loads with constant power set points $\SloadarefMAT$ are defined, 
   \begin{equation}
\SloadarefMAT =
\begin{bmatrix} 
\Sloadaref    \\
\Sloadbref    \\
\Sloadcref 
\end{bmatrix} 
=
\diag\left(\VikSDPwye (\IloadSDP)^{\hermitiantranspose}  \right) =
\begin{bmatrix} 
\VikAN    \\
\VikBN    \\
\VikCN 
\end{bmatrix} 
\circ
\begin{bmatrix}
\Iloada \\
\Iloadb\\
\Iloadc \\
\end{bmatrix}^{*}, \label{eq_load_model}
 \end{equation}   
 where  $\circ$ is the element-wise multiplication, superscript * is complex conjugate, and superscript H is conjugate transpose.

  \subsection{Novel Impedance Transformation} \label{sec_transform_derivation}
  We now assume we  inject current into the ground only when grounding the neutral, which occurs only at a subset of buses $i \in \mathcal{I}^{\text{g}} \subset \mathcal{I}$:
  \begin{enumerate}
  \item the neutral grounding is perfect where it occurs, i.e. fixes the voltage at the ground voltage of 0 V;
  \item we don't inject any current into the \emph{ground} through capacitive effects, i.e. the line shunt admittances $ \yijkshSDP$ are small and therefore shunt currents are $\IijshkSDP \approx \IjishkSDP \approx 0$ and therefore we can subsitute series current for total current $ \IijkSDP \approx  \IijskSDP$.
      \item  power consumption or generation devices connected to the network inject a negligible amount of current into the ground.
  \end{enumerate}

Under these assumptions, the total current through the ground equals zero everywhere - except for at buses with neutral grounding. Consequently,
\begin{equation}
\IijkG = 0 = -(\IijkA + \IijkB + \IijkC + \IijkN ),
\end{equation}
and therefore we can express the neutral current as a linear combination of the phase currents,
\begin{equation}
\IijkN =  -(\IijkA + \IijkB + \IijkC  ). \label{eq_neutral_current_expr}
\end{equation}
We can write all current variables as a linear transformation of the phase current vector,
\begin{equation}
\IijkSDP =  \begin{bmatrix} 
\IijkA    \\
\IijkB     \\
\IijkC  \\
\IijkN  
 \end{bmatrix}
 =  \underbrace{ \begin{bmatrix}
   1 & 0 & 0 \\
    0 & 1 & 0\\
   0 & 0 & 1 \\
   -1 & -1 &- 1 \\
   \end{bmatrix}}_{\textcolor{\paramcolor}{\mathbf{T}^{\text{T}}}}
   \underbrace{   \begin{bmatrix} 
\IijkA    \\
\IijkB     \\
\IijkC  
 \end{bmatrix} 
 }_{\IijkSDP[\setPhases]}
 = 
 \textcolor{\paramcolor}{\mathbf{T}^{\text{T}}} \IijkSDP[\setPhases] \label{eq_current_transf}.
   \end{equation}   
   Note that this transformation does not lose any information.
We can now transform Ohm's law \eqref{eq_ohms_4w}, multiplying it on the left with $\textcolor{\paramcolor}{\mathbf{T}}$,
        \begin{equation}
\textcolor{\paramcolor}{\mathbf{T}} \VjkSDP  =  \textcolor{\paramcolor}{\mathbf{T}} \VikSDP - \textcolor{\paramcolor}{\mathbf{T}} \zijksSDP \IijkSDP,
\end{equation}   
and perform the substitutions  \eqref{eq_voltage_transf}, \eqref{eq_current_transf} .
 In phase-to-neutral voltages ($3 \times 1 $ vectors), the novel expression for Ohm's law is,
        \begin{equation}
\VjkSDPwye = \VikSDPwye - \underbrace{ \textcolor{\paramcolor}{\mathbf{T}} \zijksSDP  \textcolor{\paramcolor}{\mathbf{T}^{\text{T}}} }_{\mathcal{T}(\zijksSDP)} \IijkSDP[\setPhases] .
\label{eq_transformed_ohms}
\end{equation}   
The explicit form of the  transformed  matrix $\mathcal{T}(\zijksSDP)= \textcolor{\paramcolor}{\mathbf{T}} \zijksSDP  \textcolor{\paramcolor}{\mathbf{T}^{\text{T}}}$ is, 
\resizebox{0.48\textwidth}{!}{
$
\begin{bmatrix} 
\zijksAA - \zijksNA - \zijksAN + \zijksNN  &\zijksAB  - \zijksNB - \zijksAN + \zijksNN  & \zijksAC  - \zijksNC - \zijksAN + \zijksNN  \\
\zijksBA - \zijksNA - \zijksBN + \zijksNN  &\zijksBB - \zijksNB - \zijksBN  + \zijksNN   & \zijksBC - \zijksNC - \zijksBN + \zijksNN \\
\zijksCA - \zijksNA - \zijksCN + \zijksNN  & \zijksCB - \zijksNB - \zijksCN+  \zijksNN  & \zijksCC - \zijksNC - \zijksCN  + \zijksNN  \\
 \end{bmatrix} .
 $
  }
\begin{equation}
\phantom{i} \label{eq_pn_impedance_explicit}
\end{equation}


 Koirala et al.  propose a related transformation $\mathcal{U}(\zijksSDP)$ \cite{Koirala2019}.
We can now cast that approach as a special case of \eqref{eq_pn_impedance_explicit}, using the further assumption that all mutual impedances are negligible, i.e. $\zijksAB \approx \zijksAC \approx \zijksAN \approx \zijksBC \approx \zijksBN \approx \zijksCN \approx \zijksBA \approx \zijksCA \approx \zijksNA \approx \zijksCB \approx \zijksNB \approx \zijksNC \approx0$ , and obtain their result,
        \begin{equation}
\mathcal{U}(\zijksSDP) =
\begin{bmatrix} 
\zijksAA + \zijksNN  &\zijksNN  &\zijksNN  \\
\zijksNN  &\zijksBB  + \zijksNN   &\zijksNN \\
\zijksNN  & \zijksNN  & \zijksCC  + \zijksNN  \\
 \end{bmatrix} . 
  \end{equation}  

  \subsection{Neutral Voltage Recovery Algorithm}
\label{sec_transf_properties}
Comparing \eqref{eq_ohms_4w} and \eqref{eq_transformed_ohms} we observe we have lost the expression for the neutral voltage rise, i.e. 
     \begin{equation}
\VjkN = \VikN - (\zijksNA - \zijksNN ) \IijkA - (\zijksNB - \zijksNN ) \IijkB - (\zijksNC - \zijksNN )\IijkC. \label{eq_neutral_rise_relaxation}
\end{equation}   
Without this constraint, the transformation is effectively  a relaxation\footnote{also an approximation when the assumptions in \S \ref{sec_transform_derivation} do not hold.}. 
   Although this tells us that we may not be able to recover a unique solution in general if we only use the transformed Ohm's law expression \eqref{eq_transformed_ohms}, it will be shown  all hope is not lost.
  The neutral voltage recovery  is performed through  Algorithm \ref{algo_neutral_recovery}.
  The algorithm requires the network graph to not have any isolated subgraphs.
It starts by assigning $0$ V to the neutral voltage variable $\VikN$ for all grounded buses $i \in \mathcal{I}^{g}$ (step 1), and adds them to a set collecting the buses with known neutral voltage $\mathcal{I}^{\text{visited}}$ (step 2). 
From this set of buses with known voltage magnitude, we search for  adjacent buses where the neutral voltage is not yet known, using the topology information (step 3-4). 
We can now solve \eqref{eq_neutral_rise_relaxation} for the receiving end voltage $\VjkN$ given the sending end voltage $\VikN$ and the currents $\IijkA,\IijkB,\IijkC$ (step 5).
We now add the receiving end  bus $j$ to set $\mathcal{I}^{\text{visited}}$ (step 6), and try to find a new adjacent bus (step 3).
Note that the algorithm would also work when the voltage of the grounded buses is assigned a nonzero value.
  \begin{algorithm}[tb]
 \caption{Recover neutral voltage  $\VikN$ from  $\IijkSDP[\setPhases]$ and $\zijksSDP$} \label{algo_neutral_recovery}
 \begin{algorithmic}[1]
 \renewcommand{\algorithmicrequire}{\textbf{Input:}}
 \renewcommand{\algorithmicensure}{\textbf{Output:}}
 \REQUIRE $\mathcal{I}^{\text{g}},\mathcal{I}, \mathcal{T}, \forall lij \in \mathcal{T}:\IijkA, \IijkB , \IijkC , \zijksNA, \zijksNB, \zijksNC, \zijksNN$
 \ENSURE  $\forall i \in \mathcal{I} :\VikN$
  \STATE $\forall i \in \mathcal{I}^{\text{g}}: \VikN \gets 0$ V
  \STATE $\mathcal{I}^{\text{visited}} \gets \mathcal{I}^{\text{g}}$
  \WHILE {$ \mathcal{I}^{\text{visited}} \neq  \mathcal{I}  $} 
    \STATE   $lij \gets \exists\, lij  \in \mathcal{T}: i \in \mathcal{I}^{\text{visited}}, j \notin \mathcal{I}^{\text{visited}}$
\STATE  $\VjkN \gets \VikN - (\zijksNA - \zijksNN ) \IijkA - (\zijksNB - \zijksNN ) \IijkB - (\zijksNC - \zijksNN )\IijkC
 $
  \STATE $\mathcal{I}^{\text{visited}}\gets \mathcal{I}^{\text{visited}} \cup \{j \}$
  \ENDWHILE
 \end{algorithmic} 
 \end{algorithm}

\begin{figure}[tbh]
  \centering
    \includegraphics[scale=0.75]{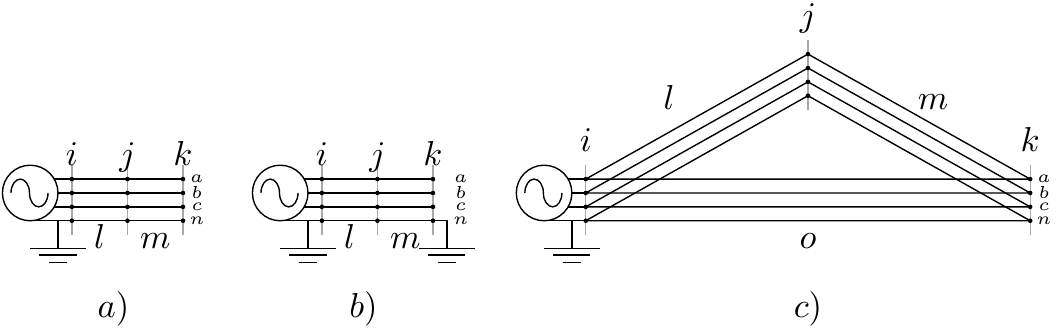}
  \caption{Examples of network configuration variants: a) single-grounded radial, b) multi-grounded radial, c) single-grounded meshed.}  \label{fig_network_variants}
\end{figure}

We note that the solution this algorithm provides is, in general, not unique, e.g. for multi-grounded topologies. 
We illustrate 3 configurations in Fig.\,\ref{fig_network_variants}.
For instance, in a 3-bus meshed system with 4-wire branches, and neutral grounding only at bus $i$, it is possible to recover the neutral voltage at bus $j$ in two different ways: $i \rightarrow j$ or $i \rightarrow k \rightarrow j$ (Fig.\,\ref{fig_network_variants}c). 
A similar problem occurs where the neutral is grounded at two buses, but not on in-between ones. For instance, a radial system $i-j-k$ where the neutral is grounded on bus $i$ and $k$ but not on  $j$  (Fig.\,\ref{fig_network_variants}b).
In all cases the algorithm terminates.

In networks where each continuous neutral section is grounded only once (e.g. Fig.\,\ref{fig_network_variants}a, Fig.\,\ref{fig_network_variants}c), the algorithm's solution is not only unique, but also identical to that of the original 4x4 Ohm's law  \eqref{eq_ohms_4w} problem!
The reason is that KCL is satisfied everywhere, and only the single grounding allows the neutral current to flow. 
The neutral current therefore must be uniquely determined through the inductive coupling with the phases, and the current balance w.r.t loads/generators.
When the network is grounded on multiple locations, the current will divide across the different paths to ground, and therefore the constraint \eqref{eq_neutral_rise_relaxation} is needed to guarantee the correct solution.

\subsection{Discussion}
Note that the load behavior \eqref{eq_load_model} is still valid, and generalizes to ZIP and exponential load models.
KCL \eqref{eq_4w_kcl} only requires the phase-currents now. 
We refer to \cite{Geth2020b} for a comprehensive description of the feasible set of 3-wire (O)PF. 

The proposed phase-to-neutral impedance transformation is exact when the following conditions hold simultaneously:
\begin{enumerate}
    \item networks with branches without shunt admittance between the conductors and earth;
    \item networks where none of the connected devices put current into the ground;
    \item networks where neutral grounding is only performed once for each continuous neutral section.
\end{enumerate}
In any other circumstance, it is either a true approximation, a relaxation, or both.
The three conditions are realistic assumptions in single-grounded LV networks. 
However, at higher voltage levels, and/or in multi-grounded systems, the assumptions are harder to justify. 
Therefore in the next section we will perform numerical experiments to understand the implications.

In the context of optimization problems, by performing the transformation, we have  eliminated the bus neutral voltage and line neutral current variables, thereby reducing the problem size significantly. 
In post-processing, we run Algorithm 1 to recover the neutral voltage values.
We can now directly impose magnitude bounds on the phase-to-neutral voltages, however, we lose the ability to impose neutral-to-ground voltage limits\footnote{that would require adding neutral voltage variables and \eqref{eq_neutral_rise_relaxation}.}. 
Neutral current limits are still enforcable using \eqref{eq_neutral_current_expr},
\begin{equation}
|\IijkN| \rightarrow  |\IijkA + \IijkB + \IijkC| \leq \textcolor{\paramcolor}{I_{l,n}^{\text{max}}}.
\end{equation}

Note that we do \emph{not} suggest this transformation technique should be used in \emph{simulation} contexts without justification, despite high accuracy in contexts that may be easy to identify a-priori. 
As Kersting  writes in his seminal  book \cite{Kersting2001b} ``\textit{When the [simulation] analysis is being done using a computer, the approach to take is to go ahead and model the shunt admittance for both overhead and underground lines. Why make a simplifying assumption when it is not necessary?}''. 
The scalability of power flow solvers using for example the fixed-point iteration variant of the current injection method is vastly better than the scalability of nonlinear unbalanced \emph{optimal} power flow (UBOPF). Therefore, the impedance transformation approach can be justified more easily for UBOPF. 

\subsection{Contrast: Kron's Reduction of the Neutral}
To enable the computational comparison between different impedance transformations in the upcoming section, we now briefly summarize Kron's reduction of the neutral.
Kron's reduction of the neutral  is defined for a partition of the series impedance matrix based on the neutral self-impedance, 
        \begin{equation}
\zijksSDP
=
\left[
\begin{array}{c|c}
\textcolor{\complexparamcolor}{\mathbf{z}^{\text{PP}}_l} &  \textcolor{\complexparamcolor}{\mathbf{z}^{\text{PN}}_l}\\
\hline
\textcolor{\complexparamcolor}{\mathbf{z}^{\text{NP}}_l} &  \zijksNN
\end{array}
\right]
 \in \setComplex^{4\times 4}, \textcolor{\complexparamcolor}{\mathbf{z}^{\text{PP}}_l} \in \setComplex^{3\times 3},  \textcolor{\complexparamcolor}{\mathbf{z}^{\text{PN}}_l}, (\textcolor{\complexparamcolor}{\mathbf{z}^{\text{NP}}_l})^{\text{T}}  \in \setComplex^{3\times 1}.
  \end{equation}  
We also partition Ohm's law \eqref{eq_ohms_4w} accordingly,
  \begin{equation}
\left[
\begin{array}{c} 
\VjkSDP[\setPhases] \\ \hline
\VjkN  
 \end{array} 
 \right] = 
 \left[
\begin{array}{c} 
\VikSDP[\setPhases] \\ \hline
\VikN  
 \end{array} 
 \right] -
 \left[
\begin{array}{c|c}
\textcolor{\complexparamcolor}{\mathbf{z}^{\text{PP}}_l} &  \textcolor{\complexparamcolor}{\mathbf{z}^{\text{PN}}_l}\\
\hline
\textcolor{\complexparamcolor}{\mathbf{z}^{\text{NP}}_l} &  \zijksNN
\end{array}
\right]
\left[
\begin{array}{c} 
\IijskSDP[\setPhases] \\ \hline
\IijksN  
 \end{array} 
 \right] \label{eq_ohms_krons}.
  \end{equation}
  Under the assumption $\VikN \approx 0 \approx \VjkN$ we derive an expression for the neutral current from the lower entries of \eqref{eq_ohms_krons}:
            \begin{equation}
\IijksN =  -\frac{1}{\zijksNN}  \textcolor{\complexparamcolor}{\mathbf{z}^{\text{NP}}_l} \IijskSDP[\setPhases] .
  \end{equation} 
We substitute this in the top entries of \eqref{eq_ohms_krons} and obtain,
\begin{equation}
    \VjkSDP[\setPhases] = \VikSDP[\setPhases] - 
    \left(
    \textcolor{\complexparamcolor}{\mathbf{z}^{\text{PP}}_l} -\frac{1}{\zijksNN} \textcolor{\complexparamcolor}{\mathbf{z}^{\text{PN}}_l}  \textcolor{\complexparamcolor}{\mathbf{z}^{\text{NP}}_l}
    \right)
    \IijskSDP[\setPhases] .
\end{equation}
Kron's reduction of the neutral is now defined as the  transformation $\mathcal{K}$,
          \begin{equation}
 \mathcal{K}(\zijksSDP) = \textcolor{\complexparamcolor}{\mathbf{z}^{\text{PP}}_l} -\frac{1}{\zijksNN} \textcolor{\complexparamcolor}{\mathbf{z}^{\text{PN}}_l}  \textcolor{\complexparamcolor}{\mathbf{z}^{\text{NP}}_l} .
  \end{equation}  
  Note that this transformation is nonlinear in the entries of $\zijksSDP$; for further details we refer the reader to  \cite{Dorfler2011a_edit}.  
 For completeness, we provide the corresponding transformation for the shunt admittance in the $\Pi$-section,
\begin{equation}
  \left[
  \begin{array}{c} 
\IijshkSDP[\setPhases] \\ \hline
\IijkshN  
 \end{array} 
 \right] =
 \left[
\begin{array}{c|c}
\textcolor{\complexparamcolor}{\mathbf{y}^{\text{PP}}_l} &  \textcolor{\complexparamcolor}{\mathbf{y}^{\text{PN}}_l}\\
\hline
\textcolor{\complexparamcolor}{\mathbf{y}^{\text{NP}}_l} &  \yijkshNN
\end{array}
\right]
\left[
\begin{array}{c} 
\VikSDP[\setPhases] \\ \hline
\VikN \approx 0 
 \end{array} 
 \right],
\end{equation}
 which tells us the neutral shunt current is $\IijkshN  = \textcolor{\complexparamcolor}{\mathbf{y}^{\text{NP}}_l} \VikSDP[\setPhases] $.

Whereas Kron's reduction of the neutral assumes the neutral is perfectly grounded everywhere ($\VikN = 0\, \forall i$),  the phase-to-neutral transform allows the neutral voltage to be \emph{nonzero} everywhere, except for a limited set of buses. Both transformations generate $3\times 3$ impedance matrices, with voltage variables after the transformation to be interpreted as phase-to-neutral voltages.




\section{Numerical experiments}

Table \ref{tab:transformations} summarizes the formulations discussed in the previous sections.
These transformations will  be compared numerically in this section. 
\begin{table}[htb]
{
  \centering
  \caption{Summary of Impedance Transformations}
    \begin{tabular}{l l l l}
     & Transform &  Result & Assumptions   \\
    \hline
$\mathcal{O}$ & Original (do-nothing) & $4\times 4$ & -     \\
$\mathcal{K}$ & Kron's reduction & $3\times 3$&  neutral voltage 0      \\
$\mathcal{T}$ & Phase-to-neutral (\textit{proposed}) & $3\times 3$&  ground current 0   \\
$\mathcal{U}$ & Modified phase-to-neutral & $3\times 3$&  ground current 0,      \\
 &   &  &   no mutual impedance     \\
\hline
    \end{tabular}%
  \label{tab:transformations}%
  }
\end{table}%

\subsection{Data Sets}

We have adapted a previously-developed library of 128 single-grounded networks \cite{Claeysthesis2021}, based on originals from the Low-Voltage Network Solutions project \cite{LVNSclosedown14}. 
\begin{itemize}
    \item explicit-neutral versions of `ENWL' networks \cite{LVNSclosedown14} with impedance data replaced with data developed by Urquhart and Thomson  \cite{Urquhart2015} ($\mathcal{O}$, $4 \times 4$);
    \item the Kron-reduced  equivalents ($\mathcal{K}, 3 \times 3$, previously explored in \cite{Claeysthesis2021});
    \item novel  phase-to-neutral transformed versions ($\mathcal{T}, 3 \times 3$);
    \item novel  modified phase-to-neutral transformed versions ($\mathcal{U}, 3 \times 3$).
\end{itemize}
Note that these networks start at the 11 kV bus, and then include a delta-grounded-wye transformer, providing the sink for the neutral current. Each of these test cases  have only a single neutral grounding, and are lacking shunt admittance data, and therefore are expected to be exact under the phase-to-neutral transformation. Note that all this networks exclusively contain cables, not overhead lines.


Moreover, we include a simple 4-wire overhead line MV 2-bus single-grounded validation system from Kersting and Philips \cite{Kersting1994}.
This example includes shunt admittance data, and different levels of load unbalance, and therefore serves as a good sensitivity study. 
Finally, we also develop a similar 2-bus single-grounded system for an Australian LV overhead line.

\subsection{Software}

We use \textsc{OpenDSS} \cite{Dugan2011} as a well-tested unbalanced power flow solver that supports explicit neutral wires. 
We  apply the proposed transformation to the \textsc{OpenDSS} data and write it to a new  file, and compare the voltages and currents w.r.t. the originals.
The  \textsc{Julia} toolbox \textsc{PowerModelsDistribution} (\textsc{PMD}) \cite{Fobes2020}  is used to build the UBOPF problems.
As \textsc{PMD} has a built-in \textsc{OpenDSS} parser, we also generate the transformed `DSS' files for the different transformations (Table~\ref{tab:transformations}), allowing us to validate w.r.t OpenDSS as well as to solve UBOPF case studies.

\textsc{PMD}  contains high-quality implementations of the 3-wire  (detailed in \cite{Geth2020b}) and  4-wire UBOPF models (detailed in \cite{Claeys2020}). 
These implementations are consistent with seminal work \cite{Araujo2013}. 
We specifically use the rectangular current-voltage formulation, which results in nonconvex quadratic programming models, that are solved through \textsc{JuMP} \cite{Dunning2015}   (for automatic differentiation) and \textsc{Ipopt}  \cite{Wachter2006} (nonlinear programming solver).
In our experience, voltage variable initialization is crucial for reliability and convergence speed, however it is not as crucial for the other variables (current, power).
Therefore in our experiments, the voltage variables are initialized at a balanced 1-pu phasor, with the neutral voltage set to 0 (where applicable). 
Note that when traversing transformers, those initialized voltage phasors are rotated according to the transformer vector group.

We configure both OpenDSS and Ipopt to a numerical tolerance of 1E-8.

\subsection{Validation and Approximation Error}
In this section, we compare power flow results between the original 4-wire data and the phase-to neutral transform.
First, we compare cases without shunts to ground, to serve as validation. 
Secondly, we compare 4-wire cases with nonzero shunt admittance, to understand the impact of neglecting the shunts when doing the phase-to-neutral transformation.

\subsubsection{Line Shunts Negligible}

\begin{figure}[tbh]
  \centering
    \includegraphics[scale=0.4]{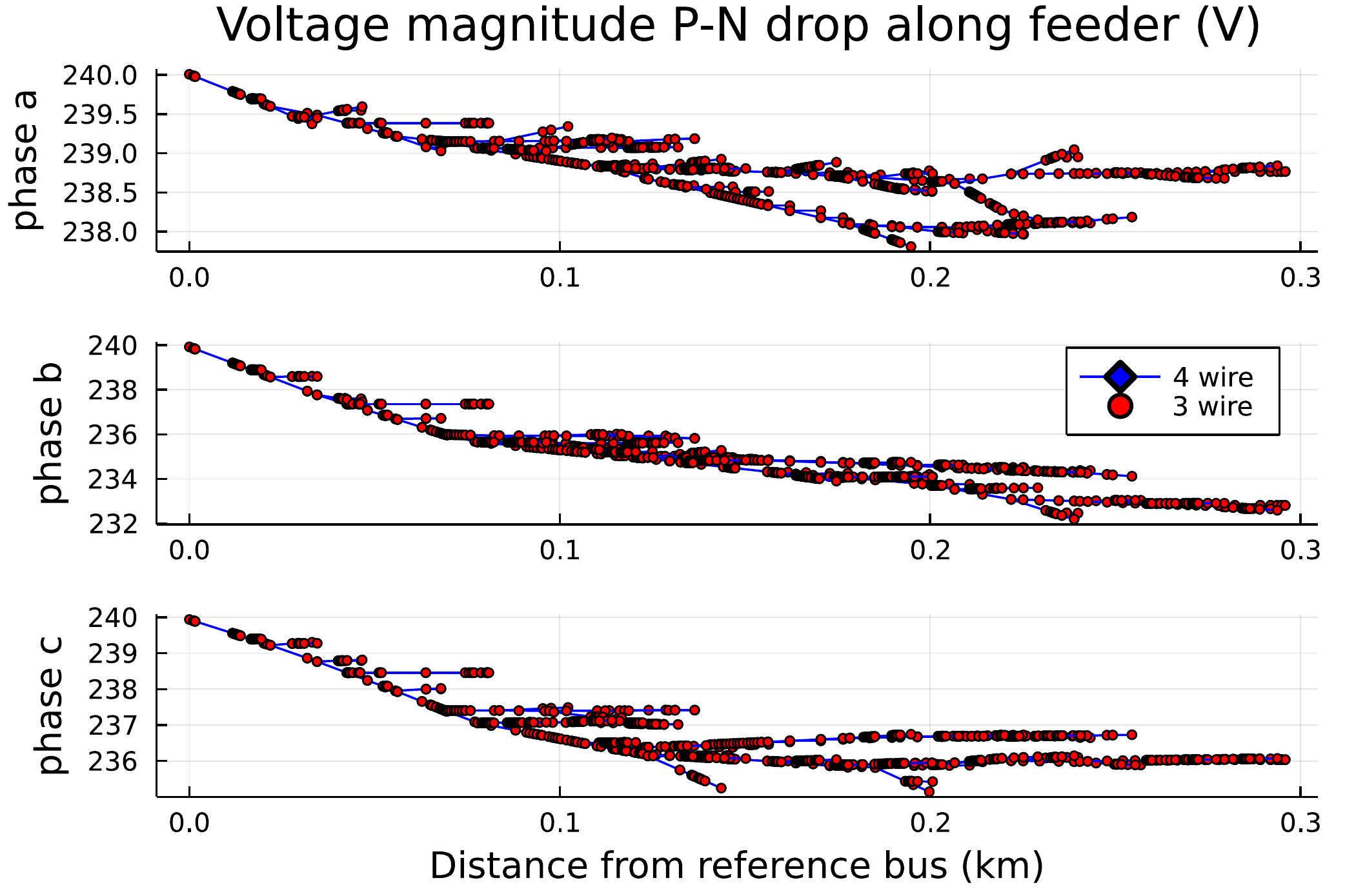}
  \caption{Comparison of original and phase-to-neutral impedance models on voltage magnitudes along the feeder 1 of network 1  (906 buses. The 4-wire $\mathcal{O}$ and 3-wire $\mathcal{T}$ voltages overlap in all three phases across the network).}  
  \label{fig_voltage_compare_Sander}
\end{figure}

Fig.~\ref{fig_voltage_compare_Sander} shows phase to neutral voltage magnitude drop along a feeder (specifically network 1 feeder 1) of one of the ENWL test cases  for both 4-wire $\mathcal{O}$ and 3-wire $\mathcal{T}$. We see that all values match. 
We note that the voltages are identical to numerical tolerance for the phase-to-neutral transformation across all 128 networks. 

\begin{figure}[tbh]
  \centering
    \includegraphics[scale=0.4]{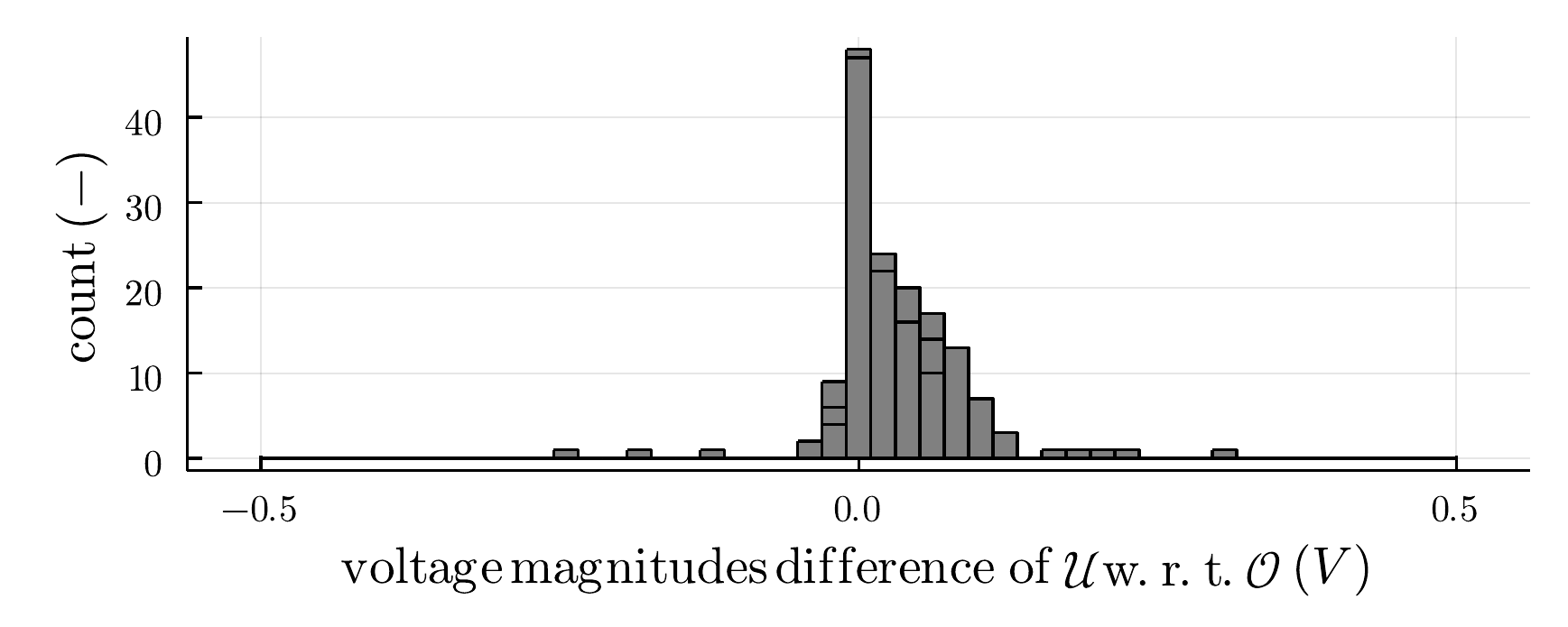}
  \caption{Histogram depicting voltage magnitudes difference of the modified phase-to-neutral transformation $\mathcal{U}$ across the library of 128 networks.}  
  \label{fig_gap_3wMpn_accuracy}
\end{figure}

\begin{figure}[tbh]
  \centering
    \includegraphics[scale=0.4]{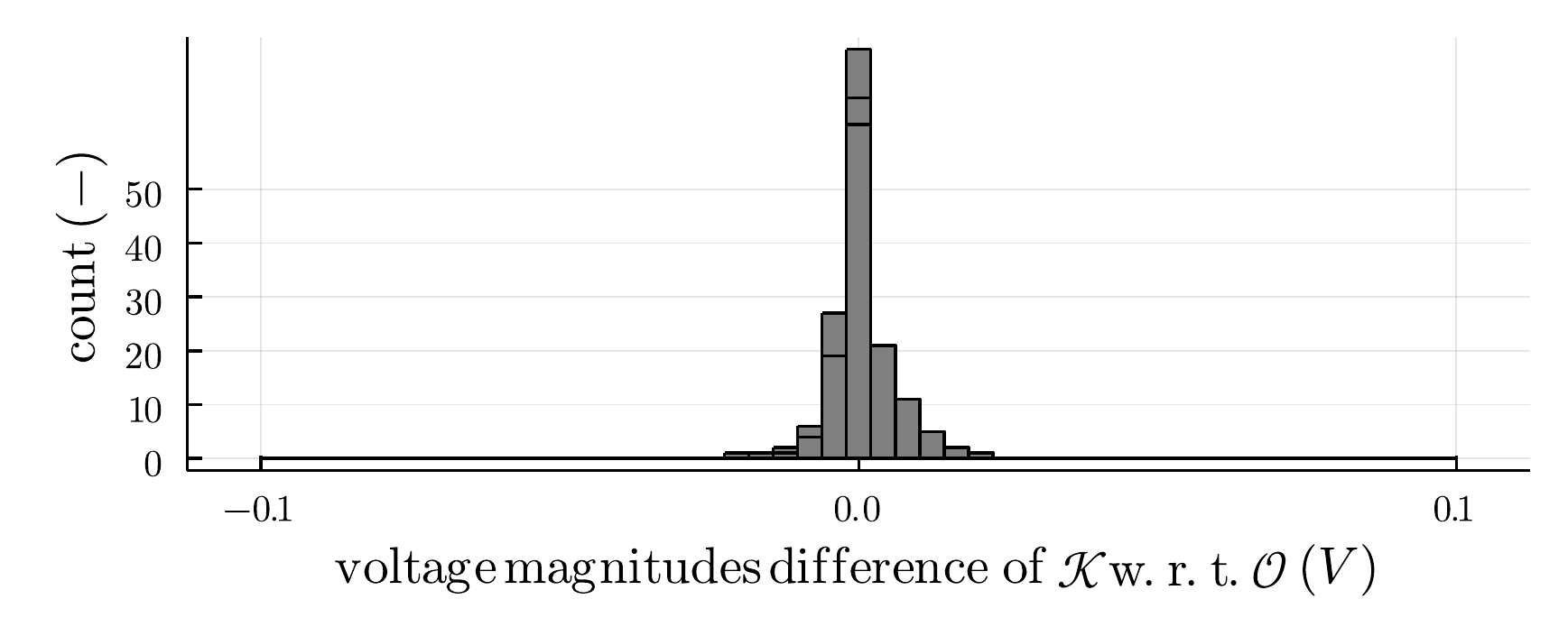}
  \caption{Histogram depicting voltage magnitudes difference of the Kron's reduction $\mathcal{K}$ across the library of 128 networks.}  
  \label{fig_gap_3wkr_accuracy}
\end{figure}

Figures~\ref{fig_gap_3wMpn_accuracy}--\ref{fig_gap_3wkr_accuracy} depict the voltage magnitude difference for the modified phase-to-neutral transformation $\mathcal{U}$ and Kron's reduction $\mathcal{K}$. 
It is noted that  $\mathcal{U}$ is the least accurate and that $\mathcal{K}$ has small, but non-zero difference.


\subsubsection{Line Shunts Nonnegligible}
The goal here is to understand the transformation where it is a true approximation -  i.e. cases \emph{with} shunt admittance.
We compare two 2-bus single-grounded overhead-line cases: a MV one previously published by Kersting, and a novel one, based on a common Australian LV overhead line configuration. 
The source bus has a fixed, balanced voltage phasor with grounded neutral, the second bus has a constant power load with ungrounded neutral. 
We analyse three scenarios for the load: balanced, unbalanced and very unbalanced. 
Note that, due to the detailed impedance models used, even with a balanced load, the  current/voltage phasors aren't balanced, as the line impedance isn't balanced. 
Table~\ref{table_kersting_basecase_4wire} reports the power flow results for the simple 4-wire 12.47 $kV$ system described in Kersting \cite{Kersting1994} with nonnegligible line shunts. 
Comparison of the phase-to-neutral voltage magnitudes show that the $\mathcal{T}$ has a small difference with respect to the original 4-wire line model, despite  the nonnegligible line shunt.

\begin{table}[tph]
\setlength{\tabcolsep}{3pt}
  \centering 
   \caption{Power flow results Kersting: voltage magnitudes and comparisons. Balanced, unbalanced, and very unbalanced load setpoints are respectively: [3, 3, 3] MVA, [4, 3, 2] MVA, and [4, 0.3, 2] MVA, PF 0.9 lagging.} \label{table_kersting_basecase_4wire}
    \begin{tabular}{l l  r r r r}
Model	&	& $|\VikAN|$  & $|\VikBN|$ & $|\VikCN|$  \\ \hline
$\mathcal{O}$  
& bal. & 0.93241815 & 0.96099983 & 0.94272734 & (pu) \\
& unb. & 0.85158645 & 1.00804238 & 0.98165957 & (pu) \\
&very unb. & 0.87569142 & 1.09761676 & 0.92216782 & (pu)  \\
\hline
$\mathcal{T}$ 
& bal. & 0.93241001 & 0.96099242 & 0.94272011 & (pu) \\
& unb. & 0.85157804 & 1.00803523 & 0.98165219 & (pu) \\
&very unb. & 0.87568319 & 1.09760949 & 0.92216047 & (pu)  \\
\hline
$\mathcal{U}$ 
&bal. & 0.87975849 & 0.87975849 & 0.87975849 & (pu) \\
& unb. & 0.73327937 & 0.9535492 & 0.9783037 & (pu) \\
&very unb. & 0.56748513 & 1.20372151 & 0.93725998 & (pu)  \\
\hline
$\mathcal{K}$ 
&bal. & 0.93270236 & 0.95880124 & 0.94469055 & (pu)  \\
& unb. & 0.88275064 & 0.98653945 & 0.97183049 & (pu) \\
&very unb. & 0.90210055 & 1.05734142 & 0.93567207 & (pu)  \\
\hline \hline
$\mathcal{T}$ - $\mathcal{O}$ 
&bal. & -8E-06 & -7E-06 & -7E-06 & (pu) \\
& unb. & -8E-06 & -7E-06 & -7E-06 & (pu) \\
&very unb. & -8E-06 & -7E-06 & -7E-06 & (pu)  \\
\hline
$\mathcal{U}$ - $\mathcal{O}$ 
&bal. & -0.05266 & -0.081241 & -0.062969 & (pu) \\
& unb. & -0.118307 & -0.054493 & -0.003356 & (pu) \\
&very unb. & -0.308206 & 0.106105 & 0.015092 & (pu)  \\
\hline
$\mathcal{K}$ - $\mathcal{O}$ 
&bal. & 0.000284 & -0.002199 & 0.001963 & (pu)       \\
& unb. & 0.031164 & -0.021503 & -0.009829 & (pu) \\
&very unb. & 0.026409 & -0.040275 & 0.013504 & (pu) \\
\hline
\end{tabular}
\end{table}

A similar analysis has been carried out for an Australian distribution two-bus system with an overhead line (modeled from Carson's equations), and the results are reported in Table~\ref{table_AusOverhead_basecase_4wire}.
\begin{table}[tph]
\setlength{\tabcolsep}{3pt}
  \centering 
   \caption{Power flow results Australian overhead line: voltage magnitudes and comparisons. Balanced, unbalanced, and very unbalanced load setpoints are respectively: [30, 30, 30] kVA, [40, 30, 20] kVA, and [40, 3, 20] kVA, PF 0.9 lagging.} \label{table_AusOverhead_basecase_4wire}
    \begin{tabular}{l l  r r r r}
Model	&	& $|\VikAN|$  & $|\VikBN|$ & $|\VikCN|$  \\ \hline
$\mathcal{O}$ 
& bal. & 0.94692051 & 0.97494785 & 0.96777093 & (pu) \\
&unb. & 0.89880375 & 1.00170312 & 0.99102676 & (pu) \\
&very unb. & 0.91174448 & 1.06698015 & 0.95470596 & (pu) \\
\hline
$\mathcal{T}$
&bal.       & 0.94691995 & 0.97494315 & 0.96777538 & (pu) \\
& unb.      & 0.89881128 & 1.0016897 & 0.99103186 & (pu) \\
&very unb.  & 0.91175119 & 1.06697602 & 0.95470524 & (pu) \\
\hline
$\mathcal{U}$
&bal.       & 0.92496348 & 0.92496348 & 0.92496348 & (pu) \\
& unb.      & 0.79688733 & 0.99131649 & 0.98129524 & (pu) \\
&very unb.  & 0.82994038 & 1.1092935 & 0.90006353 & (pu) \\
\hline
$\mathcal{K}$
&bal.       & 0.95221817 & 0.97128645 & 0.96623487 & (pu) \\
& unb.      & 0.92136173 & 0.98786491 & 0.98295279 & (pu) \\
&very unb.  & 0.93358972 & 1.03739209 & 0.95948488 & (pu) \\
\hline \hline
$\mathcal{T}$ - $\mathcal{O}$
&bal.       & -1.00E-06 & -5.00E-06 & 4.00E-06 & (pu) \\
& unb.      & 8.00E-06 & -1.30E-05 & 5.00E-06 & (pu) \\
&very unb.  & 7.00E-06 & -4.00E-06 & -1.00E-06 & (pu) \\
\hline
$\mathcal{U}$ - $\mathcal{O}$
&bal.       & -0.021957 & -0.049984 & -0.042807 & (pu) \\
& unb.      & -0.101916 & -0.010387 & -0.009732 & (pu) \\
&very unb.  & -0.081804 & 0.042313 & -0.054642 & (pu) \\
\hline
$\mathcal{K}$ - $\mathcal{O}$
&bal.       & 0.005298 & -0.003661 & -0.001536 & (pu) \\
& unb.      & 0.022558 & -0.013838 & -0.008074 & (pu) \\
&very unb.  & 0.021845 & -0.029588 & 0.004779 & (pu) \\
\hline
\end{tabular}
\end{table}
We conclude that for both cases, the error due to neglecting the shunts and then doing the phase-to-neutral transformation is small.
Transformation $\mathcal{U}$ should be avoided in MV, as the mutual impedance values are significant.
The phase-to-neutral transformation is the best $3\times 3$ representation overall - for these cases with a single-grounded neutral. If both buses would have a grounded neutral, Kron's reduction would be the best. 

\subsection{Scalability in UBOPF}
OpenDSS files describe power flow simulation case studies, not UBOPF ones. To use these network case studies as UBOPF examples we therefore now add bounds and degrees of freedom. 
\subsubsection{Adding Voltage Bounds and Generators}
The voltage magnitude bounds are set to 
\begin{equation}
\begin{bmatrix}
0.9 pu\\
0.9 pu\\
0.9 pu\\
\end{bmatrix}
=
\begin{bmatrix}
\VikANmin   \\
\VikBNmin    \\
\VikCNmin   \\
\end{bmatrix}
\leq 
\begin{bmatrix} 
    |\VikAN|    \\
    |\VikBN|    \\
    |\VikCN| 
\end{bmatrix} \leq 
\begin{bmatrix}
\VikANmax   \\
\VikBNmax    \\
\VikCNmax   \\
\end{bmatrix}
=
\begin{bmatrix}
    1.1 pu\\
    1.1 pu\\
    1.1 pu\\
\end{bmatrix}
.
\end{equation}

Note that in the four-wire case, the phase-to-ground voltages are still bounded. 
However, for the purposes of formal formulation of an optimisation problem, these bounds are weakly expressed at values that cannot become binding before the phase-to-neutral limits, even at maximal neutral voltage shift (i.e. >1.1pu).
We fix the voltage phasor at the reference bus, 
\begin{equation}
 \VikSDP
 =
 \begin{bmatrix}
\VikANref   \\
\VikBNref    \\
\VikCNref   \\
\VikNref   \\
\end{bmatrix}
 =
 \begin{bmatrix}
1 \imagangle 0^\circ pu  \\
1 \imagangle -120^\circ   pu \\
1 \imagangle 120^\circ   pu \\
0  pu
\end{bmatrix}.
\end{equation}

We read in the 128 networks, add a slack generator to represent the external system (capable of absorbing power if needed), and then add a distributed generator (DG) to the bus of each fourth load\footnote{Load IDs are numbers in the original data.}. 
The objective is to minimize the total generation cost.
Each DG has the same cost, but the slack generator is twice as expensive. 
Therefore, we expect the DGs to be dispatched until the voltage limits become binding.
Note that due to losses in the power flow, and up to one generator per bus, the solution is not degenerate despite identical costs, as the influence of the generators on the network losses will be different.
Note that for a few test cases, \textsc{Ipopt} did not converge for the modified phase-to-neutral transformation (the others transformations did not pose problems), therefore we removed those from the UBOPF testing, ending up with a final set of 111 networks.

\subsubsection{Computation Speed Results}

\begin{figure}[tbh]
  \centering
    \includegraphics[scale=0.4]{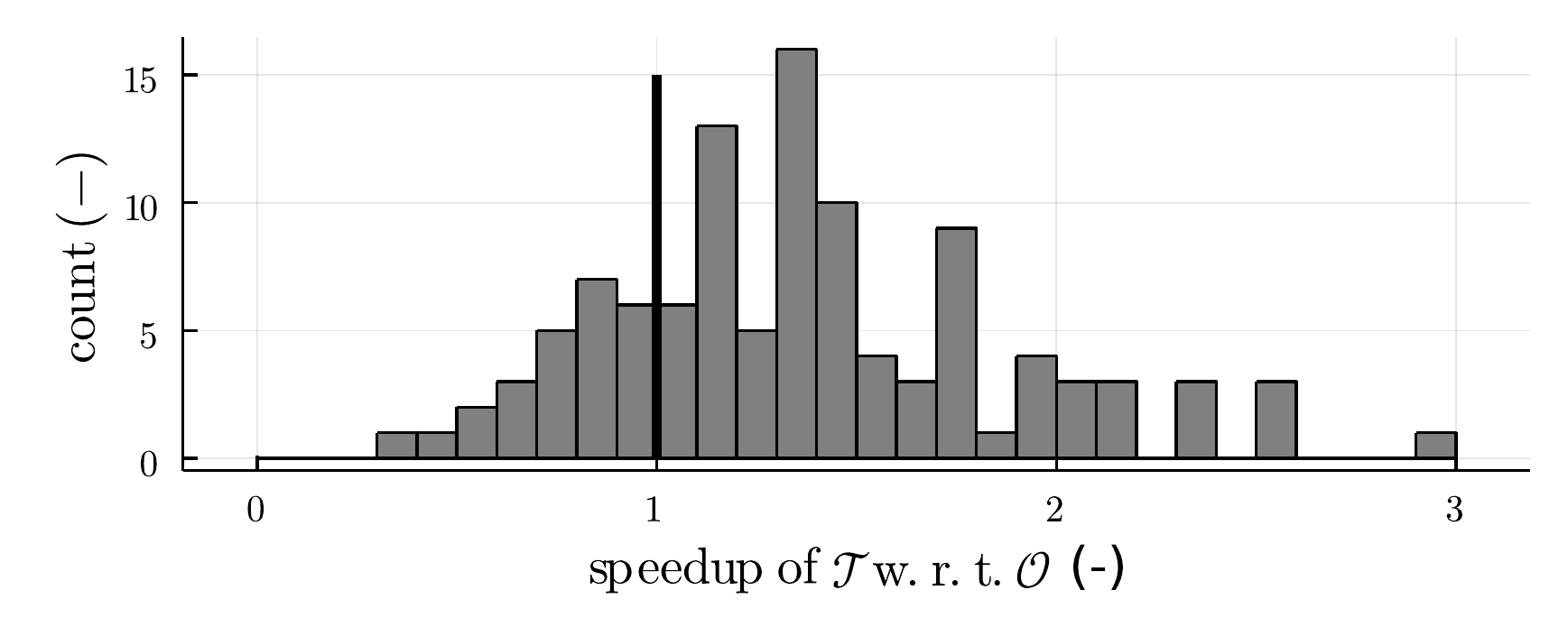}
  \caption{Histogram depicting the computational speedup of the phase-to-neutral transformation $\mathcal{T}$ across  111 networks.}  \label{fig_speedup_3wpn}
\end{figure}

Fig.\,\ref{fig_speedup_3wpn}  depicts the ratio of computation time of the 4-wire case studies before and after the phase-to-neutral transformation. 
Most cases have a speedup above 1, therefore we conclude an overall improvement. The mean speedup is 1.42 with 25 cases solving slower than the 4-wire original.

\begin{figure}[tbh]
  \centering
    \includegraphics[scale=0.4]{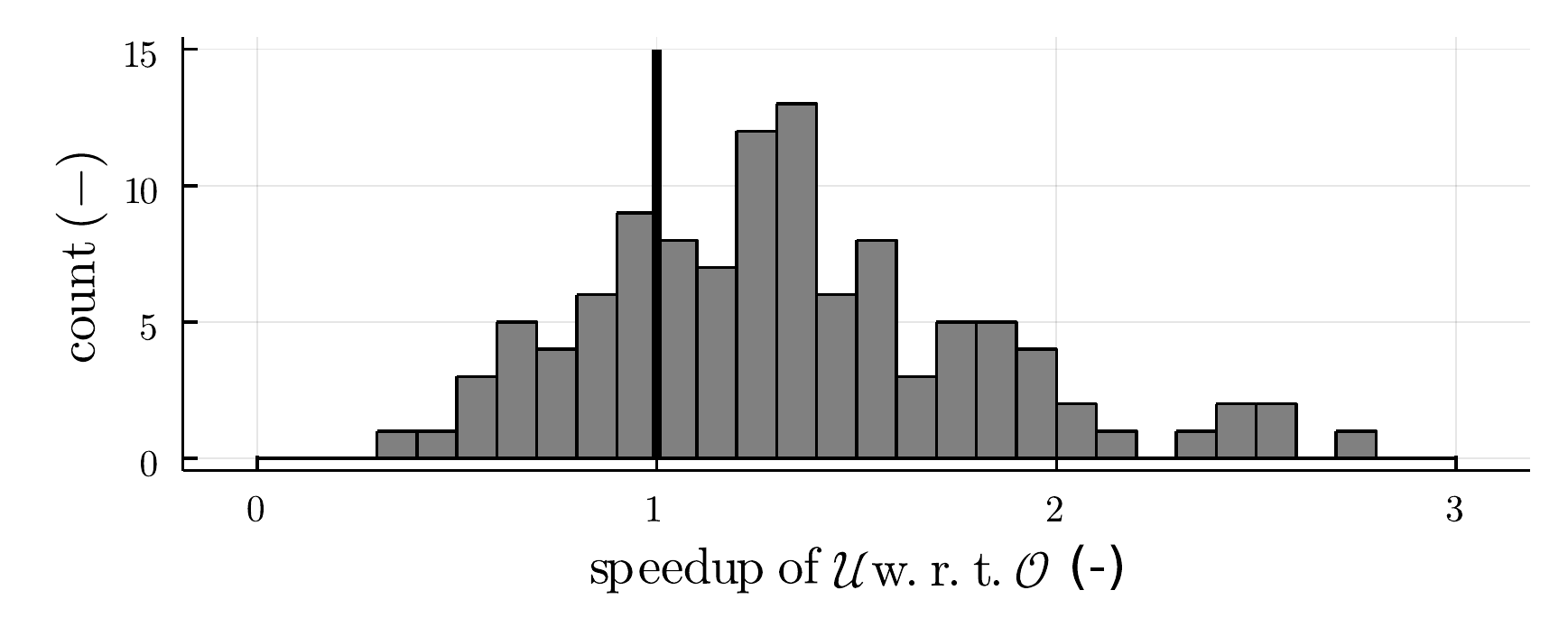}
  \caption{Histogram depicting the computational speedup of the modified phase-to-neutral transformation $\mathcal{U}$ across 111 networks.}  \label{fig_speedup_3wMpn}
\end{figure}

Fig.\,\ref{fig_speedup_3wMpn}  depicts the ratio of computation time of the 4-wire case studies before and after the modified phase-to-neutral transformation. The mean speedup is 1.37 with 29 cases solving slower than the 4-wire original.

\begin{figure}[tbh]
  \centering
    \includegraphics[scale=0.4]{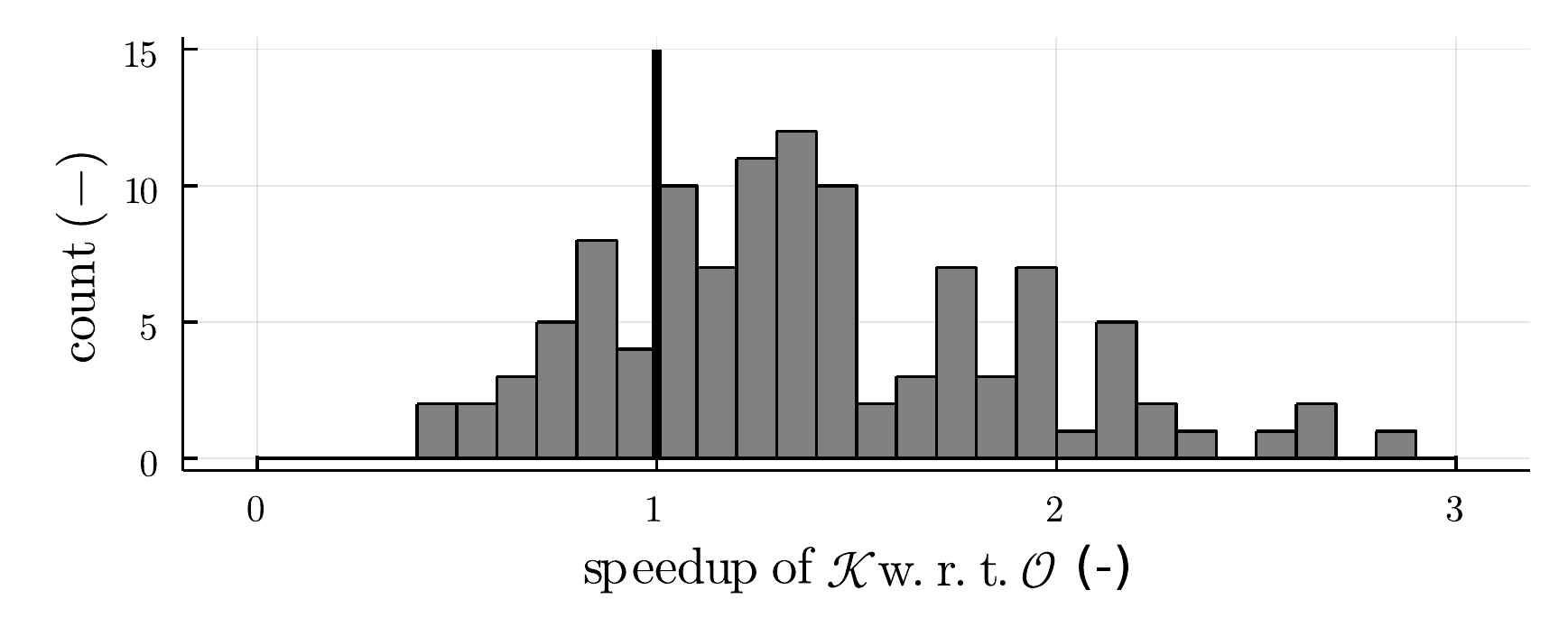}
  \caption{Histogram depicting the computational speedup of Kron's reduction $\mathcal{K}$ across 111 networks.}  \label{fig_speedup_3wKr}
\end{figure}

Fig.\ref{fig_speedup_3wKr}  depicts the ratio of computation time of the 4-wire case studies before and after Kron's reduction. The mean speedup is 1.43 with 24 cases solving slower than the 4-wire original.

\begin{figure}[tbh]
  \centering
    \includegraphics[scale=0.4]{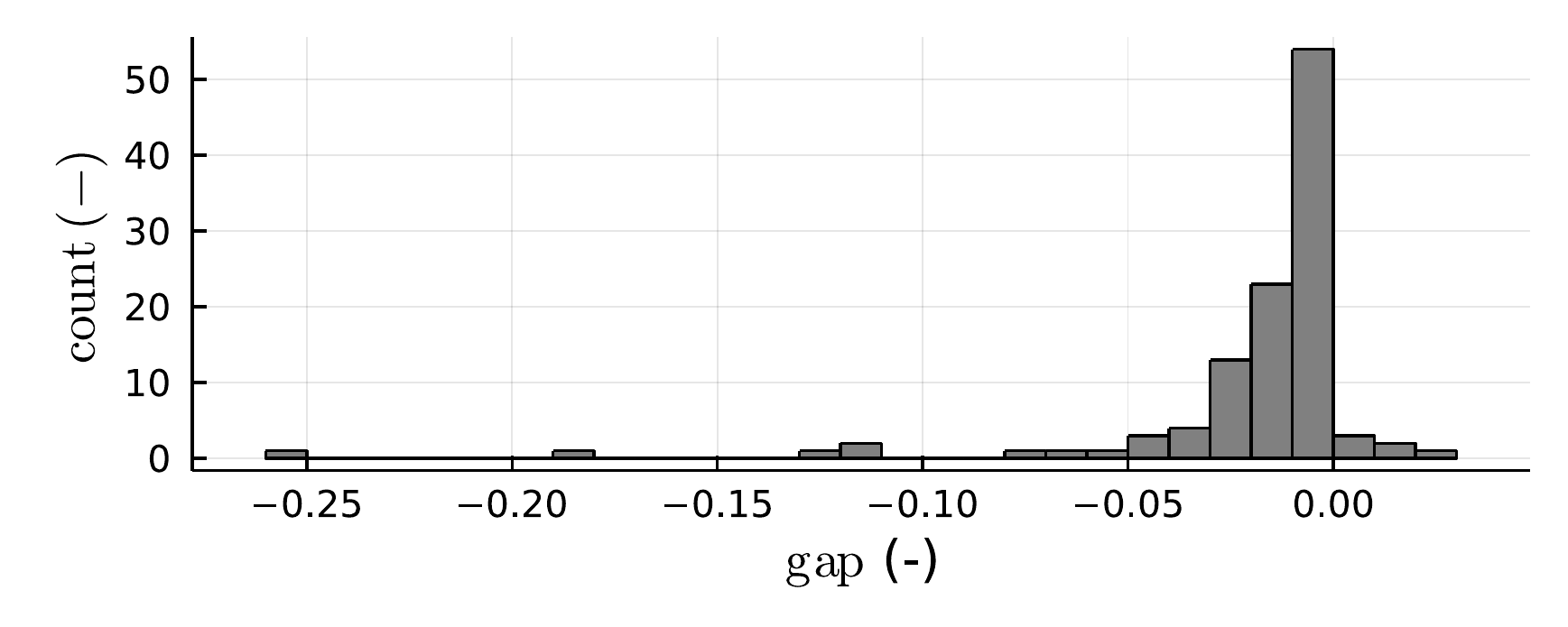}
  \caption{Histogram depicting optimality gap of the modified phase-to-neutral transformation $\mathcal{U}$ across the 111 networks.}  \label{fig_gap_3wMpn}
\end{figure}

\begin{figure}[tbh]
  \centering
    \includegraphics[scale=0.4]{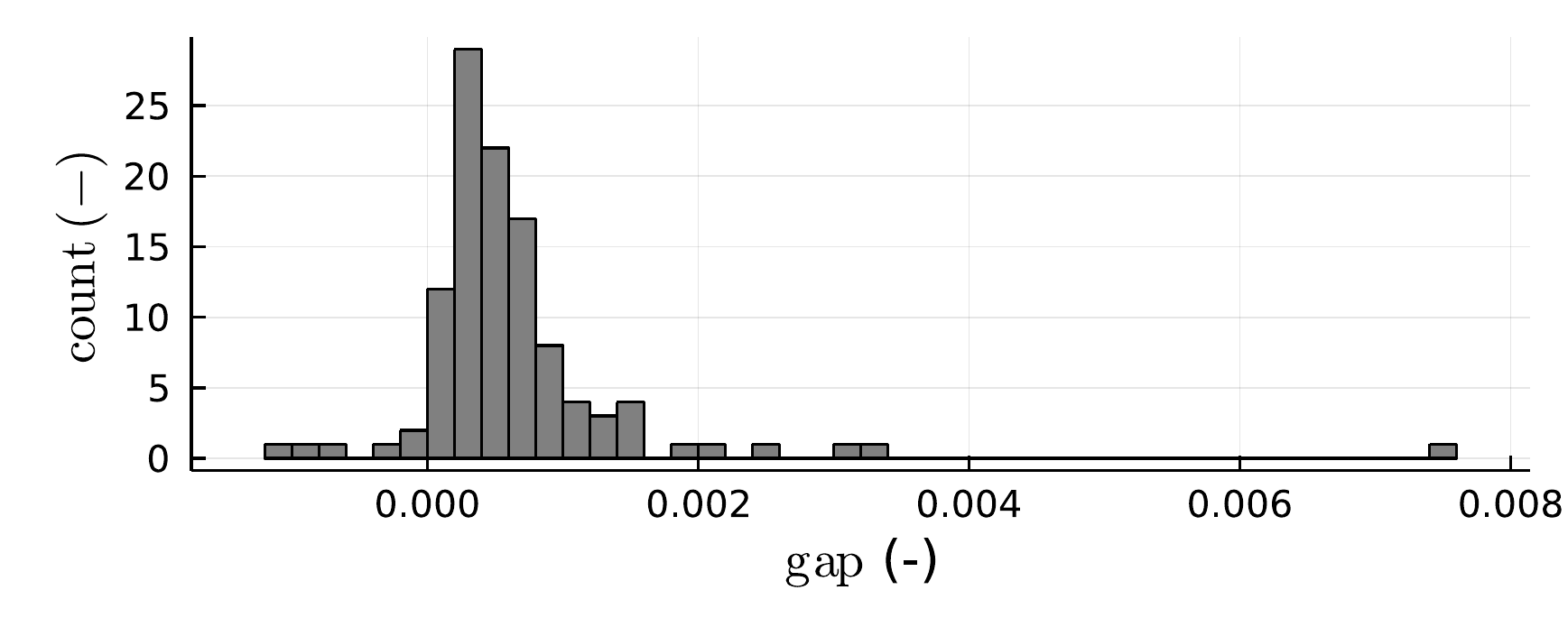}
  \caption{Histogram depicting optimality gap of Kron's reduction $\mathcal{K}$ across the 111 networks.}  \label{fig_gap_3wKr}
\end{figure}

The modified phase-to-neutral transformation and Kron's reduction result in different dispatch decisions and therefore optimal objective values. We quantify the gap for the former in Fig. \ref{fig_gap_3wMpn} and the latter in Fig. \ref{fig_gap_3wKr}. We note that the gap with Kron's reduction is typically small, however, it is nonzero. Therefore, the phase-to-neutral transformation should be considered preferentially in single-grounded networks.


We note that the scalability improvements of the proposed transformation are nice-to-have, as part of set of techniques available to practitioners looking to stretch performance.
Nevertheless, it is not game-changing, as we are still bound by the limitations of interior-point methods for nonlinear programming.

\subsubsection{Discussion}
Note that in the four-wire context phase-to-neutral voltage bounds require  an explicit phase-to-neutral variable and constraint for each phase \eqref{eq_voltage_transf}, 
or the following quadratic expression in phase-to-ground variables:\\
\resizebox{0.48\textwidth}{!}{
$\begin{bmatrix}
\VikANmin   \\
\VikBNmin    \\
\VikCNmin   \\
\end{bmatrix}  
\leq 
\begin{bmatrix}
     (\VikAreal)^2 + (\VikAimag)^2 + (\VikNreal)^2 + (\VikNimag)^2 - 2 \VikAreal \VikNreal +2 \VikAimag \VikNimag  \\
     (\VikBreal)^2 + (\VikBimag)^2 + (\VikNreal)^2 + (\VikNimag)^2 - 2 \VikBreal \VikNreal +2 \VikBimag \VikNimag  \\
     (\VikCreal)^2 + (\VikCimag)^2 + (\VikNreal)^2 + (\VikNimag)^2 - 2 \VikCreal \VikNreal +2 \VikCimag \VikNimag  \\
\end{bmatrix}
\leq 
\begin{bmatrix}
\VikANmax   \\
\VikBNmax    \\
\VikCNmax   \\
\end{bmatrix} .
$
}
\begin{equation}
\phantom{i} \label{eq_voltage_pn_quad_bounds}
\end{equation}
The phase-to-neutral transform $\mathcal{T}$ has results in more simple expressions, as it is defined in phase-to-neutral variables,
\begin{equation}
\begin{bmatrix}
\VikANmin   \\
\VikBNmin    \\
\VikCNmin   \\
\end{bmatrix}
\leq 
\begin{bmatrix}
  (\VikANreal)^2 +   (\VikANimag)^2 \\
  (\VikBNreal)^2 +   (\VikBNimag)^2 \\
  (\VikCNreal)^2 +   (\VikCNimag)^2 \\
\end{bmatrix}
\leq 
\begin{bmatrix}
\VikANmax   \\
\VikBNmax    \\
\VikCNmax   \\
\end{bmatrix}.
\end{equation}
We believe the improvement in computation speed can be partially explained by the simpler expressions for the voltage magnitude bounds. 
Nevertheless, the availability of diverse benchmarks for UBOPF is still problematic, so the computational speedup may not be as significant with different network configurations.

\section{Conclusions}

Similar to Kron's reduction, the proposed phase-to-neutral transformation reduces the size of the line matrices from $4 \times 4$ to $3 \times 3$, and can be solved with 3 voltage variables per bus instead of 4.
The novel transformation has  properties which are distinct from Kron's reduction though:
\begin{itemize}
    \item it can  \emph{facilitate} some, but only a \emph{few}, buses with neutral grounded, versus \emph{requiring} the neutral at \emph{all} buses to be grounded;
    \item  establishes neutral voltage shift, instead of approximating it as 0.
\end{itemize}

When the shunt current induced in the ground is negligible throughout the network except for at one bus where grounding enables current to flow into the ground, this work shows that you can uniquely recover the neutral voltage rise without needing explicit neutral voltage variables.
In MV networks, or in sparsely-grounded networks, the transformation may still provide a high-quality approximation that offers computational benefits. 
We believe the proposed transformation is fundamentally superior to Kron's reduction for single-grounded MV and LV networks, due to inherently limited branch shunt current contributions to the power flow at those voltage levels and due to the ability to capture neutral voltage shift accurately. 
Furthermore, the computation speed of the proposed transformation is slightly better than that of Kron's reduction and the original 4-wire description.

Throughout the development of this paper, the authors noted a lack of network data sets that are 1) explicitly  four-wire, 2) with unbalanced impedance matrices, 3) with line shunt data, 4) with complete neutral grounding descriptions. 
Often the data sets 1) are missing shunt admittance values, 2) assume balanced series impedances (e.g. because defined in terms of positive and zero sequence impedance only), 3) are only available in Kron-reduced form. 


Future work therefore includes
1) the development of unbalanced OPF data model and data sets, specifically with  multiple neutral grounding setups;
2) broader comparison of impedance transformations and coordinate spaces (e.g. sequence components).
Finally, we conjecture an application where the proposed transformation may improve algorithmic scalability  more impressively.
For four-wire UBOPF in the power-voltage variable space, the neutral voltage does not have a magnitude lower bound, which   causes convergence issues in nonlinear programming (due to $S=UI^*$ where $U\approx 0$) \cite{Claeys2020}. 
With the proposed transformation, one does not need explicit neutral voltage variables, thereby likely enabling more reliable solving of the  power-voltage equations.

\begin{acks}
The authors are grateful to Sander Claeys for doing a lot of the ground work preparing the test cases. 
We also appreciate Andrew Urquhart and Murray Thomson's help in the generation of the impedance data  \cite{Urquhart2015} for the  test case library.
Finally, the authors are grateful to all the contributors to the open source packages that enabled this work, including but not limited to \textsc{JuMP} \cite{Dunning2015}, \textsc{Ipopt}  \cite{Wachter2006}, \textsc{PowerModelsDistribution}\cite{Fobes2020}.
\end{acks}

\bibliographystyle{ACM-Reference-Format}

\appendix

\section{Appendix: Reduced Phases}
In this section we illustrate how the phase-to-neutral transformation works in a single-phase branch, as these are relatively common in practice.
We perform the transformation for a two-wire branch with phase $a$ and neutral $n$.
Ohm's law  now uses a $2 \times 2$ matrix,
  \begin{equation}
   \begin{bmatrix} 
\VikA    \\
\VikN  
 \end{bmatrix}  
 = 
    \begin{bmatrix} 
\VjkA    \\
\VjkN  
 \end{bmatrix}  
 - 
 \underbrace{\begin{bmatrix}
 \zijksAA & \zijksAN \\
 \zijksNA & \zijksNN \\
 \end{bmatrix}}_{\zijksSDP}
 \begin{bmatrix} 
\IijkA    \\
\IijkN  
 \end{bmatrix}.
   \end{equation}   
   The phase-to-neutral voltage is,
   \begin{equation}
\VikAN     =
  \underbrace{ \begin{bmatrix}
   1  &-1\\
   \end{bmatrix}}_{\textcolor{\paramcolor}{\mathbf{T}}}
   \begin{bmatrix} 
\VikA    \\
\VikN  
 \end{bmatrix}  .
   \end{equation}   
The neutral current is the negation of the phase current, 
\begin{equation}
\begin{bmatrix} 
\IijkA    \\
\IijkN  
 \end{bmatrix}
 =   \begin{bmatrix}
   1   \\
   -1  \\
   \end{bmatrix}
\IijkA    
 = 
 \textcolor{\paramcolor}{\mathbf{T}^{\text{T}}} \IijkA
   \end{equation}  
The phase-to-neutral transformed Ohm's law epxression is,
        \begin{equation}
\VikAN   
= 
\VjkAN   
- 
\underbrace{
\textcolor{\paramcolor}{\mathbf{T}} \zijksSDP  \textcolor{\paramcolor}{\mathbf{T}^{\text{T}}} }_{\mathcal{T}(\zijksSDP)} \IijkA ,
\end{equation}   
finally delivering us the single-phase impedance transformation,
   \begin{equation}
    \mathcal{T}(\zijksSDP)
    =
 \begin{bmatrix}
 \zijksAA - \zijksNA - \zijksAN + \zijksNN 
 \end{bmatrix}.
   \end{equation}   
  Note that this is just the top diagonal element of \eqref{eq_pn_impedance_explicit}.
The two-phase+neutral transformation is generated with, 
   \begin{equation}
\textcolor{\paramcolor}{\mathbf{T}}
=
\begin{bmatrix}
     1 & 0& -1\\
    0 & 1 & -1\\
   \end{bmatrix}.
   \end{equation}  
but we leave the derivation of the explicit form for the interested reader.

\end{document}